\documentclass[a4j,10pt]{article}
\usepackage{amsmath,amscd,amssymb}

\newcommand{\nbiga}{\mathcal{A}}

\newcommand{\nbigc}{\mathcal{C}}
\newcommand{\nbigd}{\mathcal{D}}
\newcommand{\nbige}{\mathcal{E}}
\newcommand{\nbigf}{\mathcal{F}}

\newcommand{\nbigh}{\mathcal{H}}
\newcommand{\nbigi}{\mathcal{I}}

\newcommand{\nbigk}{\mathcal{K}}
\newcommand{\nbigl}{\mathcal{L}}

\newcommand{\nbign}{\mathcal{N}}
\newcommand{\nbigo}{\mathcal{O}}

\newcommand{\nbigs}{\mathcal{S}}
\newcommand{\nbigt}{\mathcal{T}}
\newcommand{\nbigu}{\mathcal{U}}
\newcommand{\nbigv}{\mathcal{V}}
\newcommand{\nbigw}{\mathcal{W}}
\newcommand{\nbigx}{\mathcal{X}}

\newcommand{\proj}{\Bbb{P}}
\newcommand{\seisuu}{\Bbb{Z}}
\newcommand{\rnum}{{\bf Q}}

\newcommand{\cnum}{{\bf C}}
\newcommand{\real}{{\bf R}}
\newcommand{\hyperh}{\Bbb{H}}

\newcommand{\DD}{\Bbb{D}}

\newcommand{\gbigl}{\frak L}

\newcommand{\gbigs}{\frak S}

\newcommand{\gminiq}{\frak q}

\newcommand{\vece}{{\boldsymbol e}}

\newcommand{\vecv}{{\boldsymbol v}}
\newcommand{\vecu}{{\boldsymbol u}}
\newcommand{\vecw}{{\boldsymbol w}}

\newcommand{\vecl}{{\boldsymbol l}}

\newcommand{\vecalpha}{{\boldsymbol \alpha}}
\newcommand{\veca}{{\boldsymbol a}}
\newcommand{\vecb}{{\boldsymbol b}}

\newcommand{\vecdelta}{{\boldsymbol \delta}}
\newcommand{\vecs}{{\boldsymbol s}}
\newcommand{\vect}{{\boldsymbol t}}
\newcommand{\vecc}{{\boldsymbol c}}
\newcommand{\vecd}{{\boldsymbol d}}
\newcommand{\vech}{{\boldsymbol h}}
\newcommand{\veck}{{\boldsymbol k}}

\newcommand{\rarr}{\rightarrow}
\newcommand{\lrarr}{\longrightarrow}

\newcommand{\pf}{{\bf Proof}\hspace{.1in}}
\newcommand{\qed}{\mbox{\rule{1.2mm}{3mm}}}

\def\Cok{\mathop{\rm Cok}\nolimits}
\def\cok{\mathop{\rm Cok}\nolimits}
\def\Image{\mathop{\rm Im}\nolimits}

\def\rank{\mathop{\rm rank}\nolimits}

\def\Ker{\mathop{\rm Ker}\nolimits}
\def\ker{\mathop{\rm Ker}\nolimits}

\def\Sym{\mathop{\rm Sym}\nolimits}
\def\sym{\mathop{\rm sym}\nolimits}

\def\Res{\mathop{\rm Res}\nolimits}
\def\ord{\mathop{\rm ord}\nolimits}

\def\Ric{\mathop{\rm Ric}\nolimits}
\def\tr{\mathop{\rm tr}\nolimits}
\def\vol{\mathop{\rm dvol}\nolimits}

\newcommand{\del}{\partial}
\newcommand{\delbar}{\overline{\del}}
\newcommand{\bardel}{\overline{\del}}

\newcommand{\nbar}{\underline{n}}
\newcommand{\dbar}{\underline{d}}
\newcommand{\jbar}{\underline{j}}
\newcommand{\mbar}{\underline{m}}
\newcommand{\kbar}{\underline{k}}
\newcommand{\jitibar}{\underline{j+1}}
\newcommand{\jminusitibar}{\underline{j-1}}
\newcommand{\mitibar}{\underline{m+1}}

\newcommand{\mzeroijbar}{\underline{m_0(i,j)}}

\newcommand{\lbar}{\underline{l}}
\newcommand{\itibar}{\underline{1}}
\newcommand{\nibar}{\underline{2}}
\newcommand{\shikaku}{\sharp}
\newcommand{\sankaku}{\triangle}
\newcommand{\twoprime}{\prime\prime}

\newcommand{\harmonicbundle}{(E,\delbar_E,\theta,h)}
\newcommand{\harmonicbundledual}{(E^{\lor},\delbar_{E^{\lor}},\theta^{\lor},h^{\lor})}
\newcommand{\conjugate}{\clubsuit}
\newcommand{\barz}{\bar{z}}
\newcommand{\zbar}{\barz}
\newcommand{\zetabar}{\bar{\zeta}}
\newcommand{\baralpha}{\bar{\alpha}}
\newcommand{\alphabar}{\baralpha}
\newcommand{\barlambda}{\bar{\lambda}}
\newcommand{\lambdabar}{\barlambda}
\newcommand{\fbar}{\bar{f}}
\newcommand{\varphibar}{\bar{\varphi}}
\newcommand{\etabar}{\bar{\eta}}

\newcommand{\tbar}{\bar{t}}
\newcommand{\xbar}{\bar{x}}
\newcommand{\sbar}{\bar{s}}

\newcommand{\modelbundle}[2]{E({#1},{#2})}
\newcommand{\modeldeform}{{\mathcal Mod}}
\newcommand{\transmatrix}{B}

\newcommand{\volume}{\Omega}

\newcommand{\dotvecw}{\dot{\vecw}}

\newcommand{\dotF}{\dot{F}}
\newcommand{\dotvece}{\dot{\vece}}

\newcommand{\dote}{\dot{e}}

\newcommand{\dotPhi}{\dot{\Phi}}
\newcommand{\doth}{\dot{h}}
\newcommand{\dottheta}{\dot{\theta}}
\newcommand{\dotn}{\dot{n}}
\newcommand{\dotdelbar}{\dot{\delbar}}
\newcommand{\dotvecv}{\dot{\vecv}}
\newcommand{\dotv}{\dot{v}}

\newcommand{\soeji}{I}

\newcommand{\Poin}{{\bf p}}
\newcommand{\poin}{\Poin}
\newcommand{\prolong}[1]{{}^{\diamond}{#1}}

\newcommand{\resddlambda}{\Res(\DD^{\lambda})}
\newcommand{\residuei}{\Res_{\nbigd_i}(\DD)}
\newcommand{\residueidagger}{\Res_{\nbigd_i^{\dagger}}(\DD^{\dagger})}
\newcommand{\residueiflat}{\Res_{\nbigd_i}(\DD^f)}
\newcommand{\residueiflatdagger}{\Res_{\nbigd^{\dagger}_i}(\DD^{\dagger\,f})}

\newcommand{\weightfilt}{W}

\newcommand{\laplacian}{\Delta''}

\newcommand{\doublelangle}{\langle\langle}
\newcommand{\doublerangle}{\rangle\rangle}

\newcommand{\graded}{{\mathcal Gr}}

\newcommand{\deformoverorigin}{\prolong{\nbige}_{|\cnum_{\lambda}\times O}}
\newcommand{\deformoverorigindagger}{\prolong{\nbige^{\dagger}}_{|\cnum_{\mu}\times O}}

\newcommand{\tiisai}{\nbigu}

\newcommand{\blowup}{\widetilde}

\newcommand{\projection}{\gminiq}
\newcommand{\directsummand}{\nbigu}

\newcommand{\superinfty}{^{(\infty)}}

\newcommand{\partialkoszul}{\Pi(N_1,\ldots,N_n)}
\newcommand{\partialkoszulinfty}{\Pi(N_1\superinfty,\ldots,N_n\superinfty)}

\newcommand{\gradedrestrictedtoni}{\graded_{b|\nbigd_{\nibar}}^{(1)}}

\newcommand{\bdmath}{\begin{displaymath}}
\newcommand{\edmath}{\end{displaymath}}

\newcommand{\beqn}{\begin{equation}}
\newcommand{\eeqn}{\end{equation}}

\newcommand{\beqnarray}{\begin{eqnarray}}
\newcommand{\eeqnarray}{\end{eqnarray}}

\newcommand{\bitemize}{\begin{itemize}}
\newcommand{\eitemize}{\end{itemize}}

\newcommand{\benumerate}{\begin{enumerate}}
\newcommand{\eenumerate}{\end{enumerate}}

\newcommand{\bdescriprion}{\begin{description}}
\newcommand{\edescriprion}{\end{description}}

\newtheorem{thm}{Theorem}[section]
\newtheorem{cor}{Corollary}[section]

\newtheorem{rem}{Remark}[section]
\newtheorem{lem}{Lemma}[section]
\newtheorem{prop}{Proposition}[section]
\newtheorem{df}{Definition}[section]

\newtheorem{condition}{Condition}[section]
\newtheorem{assumption}{Assumption}[section]

\def\longuparrow{\vcenter{%
   \lineskip0pt\lineskiplimit0pt\baselineskip0pt\ialign{%
   \hfil{##}\hfil\crcr\hbox to 0pt{\hss$\uparrow$\hss}\cr%
   \hbox to 0pt{\hss\vrule width.4pt depth 0pt height 1em\hss}\cr}}}
\def\longdownarrow{\vcenter{%
   \lineskip0pt\lineskiplimit0pt\baselineskip0pt\ialign{\hfil{##}\hfil%
   \crcr\hbox to 0pt{\hss\vrule width.4pt depth 0pt height 1em\hss}\cr%
   \hbox to 0pt{\hss$\downarrow$\hss}\cr}}}

\begin{document}

\title{Asymptotic behaviour of tame nilpotent
harmonic bundles\\
with trivial parabolic structure
}
\author{Takuro Mochizuki}
\date{}

\maketitle

\abstract{
Let $E$ be a holomorphic vector bundle.
Let $\theta$ be a Higgs field,
that is a holomorphic section of
$End(E)\otimes\Omega^{1,0}_X$
satisfying $\theta^2=0$.
Let $h$ be a pluriharmonic metric of
the Higgs bundle $(E,\theta)$.
The tuple $(E,\theta,h)$ is called a harmonic bundle.

Let $X$ be a complex manifold,
and $D$ be a normal crossing divisor of $X$.
In this paper, we study
the harmonic bundle $(E,\theta,h)$ over $X-D$.
We regard $D$ as the singularity of $(E,\theta,h)$,
and we are particularly interested 
in the asymptotic behaviour of the harmonic bundle around $D$.
We will see that it is similar to
the asymptotic behaviour of complex variation
of polarized Hodge structures,
when the harmonic bundle is tame and nilpotent
with the trivial parabolic structure.
For example,
we prove constantness of general monodromy weight filtrations,
compatibility of the filtrations,
norm estimates, and the purity theorem.

For that purpose,
we will obtain a limiting mixed twistor structure
from a tame nilpotent harmonic bundle with trivial parabolic structure,
on a punctured disc.
It is a partial solution of a conjecture of Simpson.

\vspace{.1in}


}


\section{Introduction}

\subsection{Harmonic bundles}

Let $X$ be a complex manifold.
Let $(E,\delbar_E)$ be a holomorphic bundle.
Let $\theta$ be a Higgs field of $E$,
namely, it is a holomorphic section of 
$End(E)\otimes\Omega^{1,0}$ satisfying $\theta\wedge\theta=0$.
Let $h$ be a hermitian metric of $E$.
Let $\del_E$ denote the $(1,0)$-part of the metric connection
of $(E,\delbar_E,h)$.
We also have the adjoint $\theta^{\dagger}$
of $\theta$ with respect to the metric $h$.
Then we obtain the following connection:
\begin{equation}\label{eq;12.15.51}
 \DD^1:=\delbar_E+\del_E+\theta+\theta^{\dagger}:
C^{\infty}(X,E)\lrarr C^{\infty}(X,E\otimes\Omega^1_X).
\end{equation}
More generally, we obtain the following $\lambda$-connection
for any $\lambda\in\cnum$:
\[
 \DD^{\lambda}:=\delbar_E+\theta+\lambda\cdot(\del_E+\theta^{\dagger}):
C^{\infty}(X,E)\lrarr C^{\infty}(X,E\otimes\Omega^1_X).
\]
\begin{df}
The metric $h$ is called pluriharmonic, if
the connection $\DD^1$ is flat,
that is, $\DD^1\circ\DD^1=0$.
The tuple $(E,\delbar_E,\theta,h)$
is called a harmonic bundle.
\hfill\qed
\end{df}
Note that the condition is equivalent to
``$\DD^{\lambda}\circ \DD^{\lambda}=0$ for all of $\lambda$''.
We have $\DD^0=\delbar_E+\theta$,
and the condition $\DD^0\circ \DD^0=0$
is equivalent to the condition that $\theta$ is a Higgs field.

\begin{rem}
Probably, such object should be called
a `pluriharmonic bundle'.
However we use `harmonic bundle' for simplicity.
\hfill\qed
\end{rem}

Let $D$ be a normal crossing divisor of a complex manifold $X$.
In this paper, our main interest is a harmonic bundle over $X-D$,
and we investigate the asymptotic behaviour of the harmonic bundle
around $D$.
We will impose the following conditions
(see the subsection \ref{subsection;12.13.1} for more detail):
\begin{condition}\mbox{{}}
Let $P$ be any point of $X$, and $(\nbigu,\varphi)$
be an admissible coordinate around $P$
(Definition {\rm\ref{df;12.9.1}}).
On $\nbigu$, we have the description:
\[
 \theta=
 \sum_{j=1}^l f_j\cdot\frac{dz_j}{z_j}
+\sum_{j=l+1}^{n}g_j\cdot dz_j.
\]
\begin{description}
\item[(Tameness)]
 Let $t$ be a formal variable.
 We have the polynomials
 $\det(t-f_j)$ and $\det(t-g_j)$ of $t$,
 whose coefficients are holomorphic functions
 defined over $\nbigu-\bigcup_{j=1}^l D_{i_j}$.
 When the functions are extended to the holomorphic functions
 over $\nbigu$,
 the harmonic bundle is called tame at $P$.
\item[(Nilpotentness)]
Assume that the harmonic bundle is tame at $P$.
When $\det(t-f_j)|_{\nbigu\cap D_{i_j}}=t^r$,
then the harmonic bundle is called nilpotent at $P$.

When $(E,\delbar_E,h,\theta)$ is a tame nilpotent at any point $P\in X$,
then it is called a tame nilpotent harmonic bundle.
\item[(Trivial parabolic structure)]
We say that the parabolic structure of $(E,\delbar_E,\theta,h)$ is trivial,
if the parabolic structure of 
the prolongment of the restriction $(E,\delbar_E,\theta,h)_{|C}$ is trivial
for any holomorphic curve $C$ transversal with $D$.
(See Condition {\rm\ref{condition;11.28.1}}
and Definition {\rm\ref{df;12.9.2}}).
\end{description}
\end{condition}

In the words of the flat bundle $(E,\DD^1)$,
the combination of the nilpotentness condition
and the triviality of the parabolic structures
are described as follows:
\begin{condition}\mbox{{}}
\begin{enumerate}
\item
 The monodromies around the components of $D$ 
 are unipotent.
\item
 Let $s$ be a multiple-flat section.
 Let $(\nbigu,\varphi)$ be an admissible coordinate
 around $P$.
 Then we have
 $C_1\cdot\prod_{i=1}^l|z_i|^{\epsilon}
 \leq |s|_h
 \leq C_2\cdot \prod_{i=1}^l |z_i|^{-\epsilon}$.
(Precisely, we need only the estimate on curves.)
\hfill\qed
\end{enumerate}
\end{condition}

Recall that the harmonic bundle can be regarded
as a generalization of
the complex variation of polarized Hodge structures
(CVHS).
On the CVHS, 
the highly developed theories for the asymptotic behaviour 
are well known
due to Cattani-Kaplan-Schmid and Kashiwara-Kawai.
Briefly and imprecisely speaking,
their results say that
we have some nice relations between the monodromies,
and that the monodromy weight filtrations describe
the asymptotic behaviour.
Although their results indicate the direction of our study,
it seems difficult to directly apply their method in our case.

When the base manifold $X$ is one dimensional,
such behaviour was deeply studied by Simpson.
Moreover, he proposed the `mixed twistor structure',
which is quite important for the study 
in the case $X$ is higher dimensional.
In fact,
most of the essential ideas contained in this paper are due to Simpson
(see \cite{s1}, \cite{s2}, \cite{s3} and \cite{s4}):
We will heavily
owe to many results and methods
that he developed
in \cite{s1} and \cite{s2}.
We will often use them without mention his name.
The papers are fundamental for our study of harmonic bundles.
The mixed twistor structure was introduced in \cite{s3}.
The mixed twistor structure
permits us to obtain some compatibilities
on the relations between the monodromy weight filtrations
at the intersection points of divisors.
(Such compatibilities are well known for the complex variation of the 
polarized Hodge structure.)
It seems difficult for the author to obtain such results,
if we use only some rather classical elliptic analytic argument
without mixed twistor structure.

\subsection{Main results}

Recall that the harmonic bundle can be regarded
as a generalization of the complex variation of Hodge structure.
Briefly speaking,
our final but unreached purpose in this study
is to see the following:
\begin{quote}
The asymptotic behaviour of a tame harmonic bundle
around the singularity is similar
to the behaviour of CVHS around the singularity,
in some sense.
\end{quote}
As is already noted, we will mainly investigate
the tame harmonic bundles under the assumptions
of nilpotentness and the triviality of parabolic structures.
We explain what results are obtained in the case $\dim X=2$.

\subsubsection{Flat connection}

Since we are interested in the asymptotic behaviour
around the singularity of harmonic bundles,
we can assume that $X=\Delta^2=\{(z_1,z_2)\,|\,|z_i|<1\}$
and $D=D_1\cup D_2$.
Here we put $D_i:=\{z_i=0\}$.
Let $\harmonicbundle$ be a tame nilpotent harmonic bundle
with trivial parabolic structure over $X-D$.
Let $P$ be a point of $X-D$.
We have the loop $\gamma_i:[0,1]\lrarr X-D$ defined as follows:
\[
 z_j(\gamma_i(t))=
 \left\{
 \begin{array}{ll}
 z_i(P)\cdot \exp(2\pi\sqrt{-1}t) &(j=i)\\
 z_j(P) &(j\neq i).
 \end{array}
 \right.
\]
We put $V=E_{|P}$. We have the monodromy
$M(\gamma_i)\in End(E_{|P})$
with respect to the flat connection $\DD^1$
given in (\ref{eq;12.15.51}).
Due to our assumption, it is unipotent.
Thus we have the logarithm $N_i=\log M(\gamma_i)$.
We put $N(\veca)=\sum_{i=1}^2 a_i\cdot N_i$
for $\veca=(a_1,a_2)\in\real^2_{>0}$.
Let $W(\veca)$ denote the weight filtration of $N(\veca)$.
We have the constantness of the filtration on the positive cone.
Namely, the following holds.
\begin{thm}
We have $W(\veca_1)=W(\veca_2)$ for any $\veca_i\in\real_{>0}^2$.
\hfill\qed
\end{thm}

We put $N(\itibar)=N_1$ and $N(\nibar)=N_1+N_2$.
We denote the weight filtration of $N(\jbar)$ by $W(\jbar)$.
Let $\graded^{(1)}$ denote the associated graded space
of $W(\itibar)$.
Then we have the induced filtration $W^{(1)}(\nibar)$ on $\graded^{(1)}$.
On the other hand,
we have the induced action $N^{(1)}(\nibar)$ on $\graded^{(1)}$.
Let $W(N^{(1)}(\nibar))$ denote the weight filtration
of $N^{(1)}(\nibar)$.
Namely, we have two filtrations $W^{(1)}(\nibar)$
and $W(N^{(1)})(\nibar)$
on the graded vector space
$\graded^{(1)}=\bigoplus_h\graded^{(1)}_h$.

\begin{thm}
We have
$W^{(1)}(\nibar)_{a+h}\cap \graded^{(1)}_{a}
=W(N^{(1)}(\nibar))_{h}\cap \graded^{(1)}_a$.
\hfill\qed
\end{thm}

We put $\hyperh:=\{\zeta=x+\sqrt{-1}y\in\cnum\,|\,y>0\}$.
Then we have the universal covering
$\pi:\hyperh^2\lrarr X-D$ defined by
$\pi(\zeta_1,\zeta_2)
=\bigl(
   \exp(2\pi\sqrt{-1}\zeta_1),
   \exp(2\pi\sqrt{-1}\zeta_2)
 \bigr)$.
We put as follows:
\[
 \blowup{Z}(id,2,C,A):=
 \bigl\{
 (x_1+\sqrt{-1}y_1,x_2+\sqrt{-1}y_2)
 \in\hyperh^2
\,\big|\,
 |x_1|<A,\,\,
 |x_2|<A,\,\,
 C\cdot y_{1}>y_2,
 \bigr\}.
\]
Let $\blowup{P}$ be a point of $\blowup{Z}(id,2,C,A)$.
We have the pull back $\pi^{\ast}(E,\DD^1,h)$.
The fibers $\pi^{\ast}(E)_{|\blowup{P}}$
and $E_{|P}$ are naturally identified.

Let $u$ be a non-zero element of $E_{|P}$.
Then we have the numbers
$h_j=\deg^{W(\jbar)}(u)$ for $j=1,2$.
Let $f$ be a flat section of $(\pi^{\ast}E,\DD^1)$
such that $f_{|\blowup{P}}=u$.
\begin{thm}
There exist positive numbers $C_1$ and $C_2$
such that the following inequality holds
on $\blowup{Z}(id,2,C,A)$:
\[
 0<C_1\leq
 |f|_h^2\cdot y_1^{h_1}\cdot y_2^{h_2-h_1}
 \leq C_2.
\]
\hfill\qed
\end{thm}

From the tuple $(V,N_1,N_2)$,
we obtain the following complex:
\[
 \Pi(N_1,N_2):
 V\stackrel{d}{\lrarr}
 \Image(N_1)\oplus \Image(N_2)
 \stackrel{d}{\lrarr} \Image(N_1N_2).
\]
It is easy to see $H^0(\Pi(N_1,N_2))=\ker(N_1)\cap\ker(N_2)$
and $H^2(\Pi(N_1,N_2))=0$.
The filtration $W(\nibar)$ on $V$ induces the filtration $W$
of the complex as follows:
\[
 W_k(V)=W(\nibar)_k,
\quad
 W_k(\Image(N_i))=N_i(W(\nibar)_k),
\quad
 W_k(\Image(N_1N_2))=N_1N_2(W(\nibar)_k).
\]
It induces the filtration $W$
on the cohomology group $H^{\ast}(\Pi(N_1,N_2))$.
\begin{thm}[Purity theorem]
Assume that $\harmonicbundle$ has a real structure.
We have $W_kH^k(\Pi(N_1,N_2))=H^k(\Pi(N_1,N_2))$.
In other words,
the naturally defined morphisms
$\ker(d)\cap W_k(\Pi(N_1,N_2)^k)\lrarr H^k(\Pi(N_1,N_2))$
are surjective.
\hfill\qed
\end{thm}

\begin{rem}
Due to Cattani-Kaplan-Schmid and Kashiwara-Kawai,
the intersection cohomology and the $L^2$-cohomology
of CVHS are isomorphic.
The purity theorem was crucially used for the proof of such coincidence.
Although we do not discuss the relation between
the intersection cohomology and the $L^2$-cohomology
of harmonic bundles in this paper,
it seems appropriate to include the purity theorem here.
The author intends to study the relations between $L^2$-cohomology and
the intersection cohomology in another paper.
\hfill\qed
\end{rem}

\subsubsection{Holomorphic bundles}

In the previous subsubsection, we state the results
in terms of the flat bundles and the monodromies.
It is reworded in terms of the holomorphic flat bundles
and the residues.
The $(0,1)$-part $\delbar+\lambda\cdot\theta^{\dagger}$
of the connection $\DD^1$ gives a holomorphic structure
$d^{\prime\prime\,1}$ to the $C^{\infty}$-vector bundle $E$.
We denote the holomorphic bundle $(E,d^{\twoprime\,1})$
by $\nbige^1$.
Let $\prolong{\nbige^1}$
denote the prolongment of $\nbige^1$ by an increasing order.
(See subsection  \ref{subsubsection;12.9.60}.)
We can show that $\prolong{\nbige^1}$ is locally free
and same as the canonical extension,
i.e.,
$\DD^1$ is of log type on $\prolong{\nbige^1}$.

We have the residues
$N_i:=\Res_{D_i}(\DD^1)\in\Gamma(D_i,End(\prolong{\nbige^1}_{|D_i}))$
for $i=1,2$.
We put $D_{\itibar}=D_1$ and $D_{\nibar}=D_1\cap D_2$.
For $\veca\in\real^2_{>0}$,
we put $N(\veca):=\sum_{i=1}^2 a_i\cdot N_{i\,|\,D_{\nibar}}$.
Let $W(\veca)$ denote the weight filtration of $N(\veca)$.
\begin{thm}[Theorem \ref{thm;12.7.200}]
We have $W(\veca_1)=W(\veca_2)$
for any $\veca_i\in\real_{>0}^2$.
\hfill\qed
\end{thm}

On $D_{\itibar}$, we put $N(\itibar)=N_1$.
Let $W(\itibar)$ denote the weight filtration of $N(\itibar)$,
which is a filtration of $\prolong{\nbige^1}_{D_{\itibar}}$
by vector subbundles.
We obtain the graded vector bundle
$\graded^{(1)}$ on $D_{\itibar}$.

On $D_{\nibar}$, we put $N(\nibar)=\sum_{i=1}^2 N_{i\,|\,D_i}$.
Let $W(\nibar)$ denotes the weight filtration of $N(\nibar)$,
which is a filtration of $\prolong{\nbige^1}_{D_{\nibar}}$.

Let consider $\graded^{(1)}_{|D_{\itibar}}$.
The filtration $W(\nibar)$ of $\prolong{\nbige^1}_{|D_{\nibar}}$
induces the induced filtration $W^{(1)}(\nibar)$
of $\graded^{(1)}_{|D_{\itibar}}$.
On the other hand, $N(\nibar)$ induces the endomorphism $N^{(1)}(\nibar)$
of $\graded^{(1)}_{|D_{\nibar}}$.
Then we obtain the weight filtration
$W(N^{(1)}(\nibar))$.

\begin{thm}[Theorem \ref{thm;12.13.50}] \label{thm;12.9.15}
We have
$W^{(1)}(\nibar)_{h+a}\cap \graded^{(1)}_a
=W(N^{(1)}(\nibar))_h\cap \graded^{(1)}_a$
for any integers $a$ and $h$.
\hfill\qed
\end{thm}
We will see more strong compatibility
in Theorem \ref{thm;12.10.250}.
Clearly, we can replace the role of $N_1$ and $N_2$
(Theorem \ref{thm;12.15.110}).

Let take a holomorphic frame $\vecv=(v_1,\ldots,v_r)$ of $\prolong{\nbige^1}$
over $X$, compatible
with the sequence of the filtrations $(W(\itibar),W(\nibar))$.
Namely it satisfies the following:
\begin{itemize}
\item
 $\vecv_{|D_{\jbar}}$
 is compatible with the filtration $W(\jbar)$ over $D_{\jbar}$.
 In particular,
 we obtain the induced frame $\vecv^{(1)}$ over $D_{\itibar}$.
\item
 $\vecv^{(1)}_{|D_{\nibar}}$ is compatible with $W^{(1)}(\nibar)$.
\item
 We have $\deg^{W(\nibar)}(v_i)=\deg^{W^{(1)}(\nibar)}(v^{(1)}_i)$.
\end{itemize}
We put $2\cdot k_1(v_i):=\deg^{W(\itibar)}(v_i)$
and $2\cdot k_2(v_i):=\deg^{W(\nibar)}(v_i)-\deg^{W(\itibar)}(v_i)$.
We obtain the $C^{\infty}$-frame $\vecv'=(v_1',\ldots v_r')$
of $\nbige^1$ over $X-D$ defined as follows:
\[
 v_i':=v_i\cdot
 (-\log|z_1|)^{-k_1(v_i)}\cdot (-\log|z_2|)^{-k_2(v_i)}.
\]
For a positive number $C$, we put as follows:
\[
 Z(id,2,C):=\{(z_1,z_2)\in X-D\,|\,|z_1|^C<|z_2|\}.
\]
\begin{thm}[Theorem \ref{thm;12.7.110}]
On $Z(id,2,C)$, the $C^{\infty}$-frame $\vecv'$ is adapted.
Namely the hermitian matrix-valued functions
$H(h,\vecv')$ and the inverse $H(h,\vecv')^{-1}$
are bounded over $Z(id,2,C)$.
\hfill\qed
\end{thm}
Clearly we can replace the roles of $1$ and $2$
(Theorem \ref{thm;12.15.111}).

We put $\nbigv=\prolong{\nbige^1}_{|D_{\nibar}}$
and $\nbign_i:=N_{i\,|\,D_{\nibar}}$.
\begin{thm} [Theorem \ref{thm;12.7.40}]
Assume that $\harmonicbundle$ has a real structure.
The purity theorem for $(\nbigv,\nbign_1,\nbign_2)$ holds.
\hfill\qed
\end{thm}

Let $\vecv$ be a holomorphic frame of $\prolong{\nbige^1}$,
which is compatible with $W(\nibar)$ on $D_{\nibar}$.
We put $2\cdot k(v_i):=\deg^{W(\nibar)}(v_i)$.
Let consider the morphism $\psi_{m,\nibar}:X-D\lrarr X-D$
defined by $\psi_{m,\nibar}(z_1,z_2)=(z_1^m,z_2^m)$.
We have the pull back $\psi_{m,\nibar}^{\ast}(\nbige^1,\DD^1,h)$.
We have the frame $\vecv^{(m)}$ of $\psi_{m,\nibar}^{\ast}\nbige^1$
defined as follows:
\[
 v^{(m)}_i=\psi_{m,\nibar}^{\ast}(v_i)\cdot m^{-k(v_i)}.
\]
Let $F=\bigoplus_{i=1}^r\nbigo_{\Delta^{\ast}}\cdot u_i$
denote the trivial holomorphic bundle
with the frame $\vecu=(u_i)$.
We denote the holomorphic structure of $F$
by $d''_F$.
Due to the frames $\vecv^{(m)}$ and $\vecu$,
we obtain the isomorphism
$\Phi_m:\psi_{m,\nibar}^{\ast}\nbige^1\lrarr F$.
Then we obtain the sequences of
the metrics $\{h^{(m)}\}$,
the connections $\{\DD^{1\,(m)}\}$,
the (non-holomorphic) Higgs fields $\{\theta^{(m)}\}$
and the conjugates $\{\theta^{(m)\,\dagger}\}$.
We also obtain the sequences of
the holomorphic structures
$\delbar_F^{(m)}:=d''_F-\theta^{(m)\,\dagger}$.

\begin{thm}[Theorem \ref{thm;12.15.120}]
\mbox{{}}
\begin{itemize}
\item
We can pick a subsequence $\{m_i\}$ of $\{m\}$
such that the corresponding sequences
$\{h^{(m_i)}\}$, $\{\DD^{1\,(m_i)}\}$,
$\{\theta^{(m_i)}\}$, $\{\theta^{(m_i)\,\dagger}\}$,
$\{\delbar_F^{(m_i)}\}$ converge
in the $L_l^p$-sense for any $l$ and for any sufficiently large $p$.
The limits are denoted by
$h\superinfty$, $\DD\superinfty$,
$\theta\superinfty$, $\theta^{\dagger\,(\infty)}$,
and $\delbar_F\superinfty$.
\item
The tuple $(F,\delbar_F\superinfty,\theta\superinfty,h\superinfty)$
is a CVHS.
\hfill\qed
\end{itemize}
\end{thm}

\subsection{Mixed twistor structure}

In the previous subsection, we stated the results
for the holomorphic flat bundle $(\nbige^1,\DD^1)$.
In fact, we will consider the $\lambda$-connections
$(\nbige^{\lambda},\DD^{\lambda})$,
and the conjugates $(\nbige^{\dagger\,\mu},\DD^{\dagger\,\mu})$
for all $\lambda$ and $\mu$.
We will show similar results, at once.
It is the merit to consider all of $\lambda$ and $\mu$
that we obtain a limiting mixed twistor structure,
which is a partial solution of a conjecture by Simpson in \cite{s3}.
(See Definition \ref{df;12.9.11} for the definition of mixed twistor.)

Let $\harmonicbundle$ be a tame nilpotent harmonic bundle
with trivial parabolic structure over $\Delta^{\ast\,2}$.
Let $P$ be a point of $\Delta^{\ast\,2}$,
and $O$ be the origin of $\Delta^2$.
We obtain the vector bundle $S(O,P)$ over $\proj^1$,
and the morphisms
$N_i^{\sankaku}:S(O,P)\lrarr S(O,P)\otimes\nbigo_{\proj^1}(2)$
for $i=1,2$.
(See the subsubsection \ref{subsubsection;12.9.10} for the construction.)
We put $N^{\sankaku}(\nibar):=\sum_{i=1}^2N^{\sankaku}_i$.
We denote the weight filtration of $N^{\sankaku}(\nibar)$
by $W(\nibar)$.

\begin{thm}[A limiting mixed twistor theorem, Theorem
 \ref{thm;12.15.121}]
\mbox{{}} \label{thm;12.9.20}
\begin{itemize}
\item
For any neighborhood $U$ of $O$ in $\Delta^2$,
we can take an appropriate point $P\in U\cap\Delta^{\ast\,2}$
such that the filtered vector bundle $(S(O,P),W(\nibar))$
is a mixed twistor.
\item
The morphisms $N_i^{\sankaku}$ is a morphism of mixed twistors.
\hfill\qed
\end{itemize}
\end{thm}

The mixed twistor structure is essentially used in the proof of
the compatibility of the filtrations,
that is, Theorem \ref{thm;12.9.15}.
Briefly speaking,
the proof of Theorem \ref{thm;12.9.15} is divided
into three steps.
\begin{description}
\item[Step 1.]
 Let $b$ denote the bottom number of the filtration $W(\itibar)_b$.
 We see that we only have to prove the coincidence
 in the bottom part,
 i.e.,
 $W^{(1)}(\nibar)_{h+b}\cap \graded^{(1)}_b
=W(N^{(1)}(\nibar))_h\cap \graded^{(1)}_b$.
\item[Step 2.]
 We see the implication
 $W^{(1)}(\nibar)_{h+b}\cap \graded^{(1)}_b
\supset W(N^{(1)}(\nibar))_h\cap \graded^{(1)}_b$.
\item[Step 3.]
 We see the coincidence
 $W^{(1)}(\nibar)_{h+b}\cap \graded^{(1)}_b
=W(N^{(1)}(\nibar))_h\cap \graded^{(1)}_b$.
\end{description}

The step 1 is rather elementary.
We will only use the linear algebra.
In the step 2, we use a comparison method and norm estimate
in one dimensional case. However we need only rather classical analysis.
We use the mixed twistor in the step 3.
We need the following:
\begin{itemize}
\item
The filtration $W(\nibar)$ of $S(O,P)$ is a mixed twistor structure.
This is a consequence of a limiting mixed twistor theorem.
\item
We have the vector bundle:
\[
 \graded^{W(N^{(1)}(\nibar))}_{b,h}:=
 \frac{\graded^{(1)}_b\cap W(N^{(1)}(\nibar))_h}
 {\graded^{(1)}_{b}\cap W(N^{(1)}(\nibar))_{h-1}}.
\]
The equalities
$c_1(\graded^{W(N^{(1)}(\nibar))}_{h,b})
=(h+b)\cdot
 \rank(\graded^{W(N^{(1)}(\nibar))}_{h,b})$ hold
for any $h$.
Here $c_1(\nbigf)$ denotes the first Chern class
of a coherent sheaf $\nbigf$ on $\proj^1$.
\end{itemize}
From these two facts and the implication
$W^{(1)}(\nibar)_{h+b}\cap \graded^{(1)}_b
\subset W(N^{(1)}(\nibar))_h\cap \graded^{(1)}_b$,
we obtain the coincidence,
due to some general lemma for mixed twistors.

However, the consideration of 
all $(\nbige^{\lambda},\DD^{\lambda})$
and $(\nbige^{\dagger\,\mu},\DD^{\dagger\,\mu})$
raises some difficulties.
It is a principle that
the arguments over $\cnum_{\lambda}^{\ast}=\cnum_{\mu}^{\ast}$
are not different from the arguments for $\lambda=1$.
On the other hand, we need some additional argument
around $\lambda=0$.

One big difference is the existence of normalizing frames.
Let $\vecv$ be a holomorphic frame of $\nbige^{\lambda}$ over $X-D$.
We obtain the $\lambda$-connection form
$\nbiga=\sum A_j\cdot dz_j/z_j\in\Gamma(X-D,M(r)\otimes\Omega^{1,0}_{X-D})$,
determined by the relation $\DD\vecv=\vecv\cdot \nbiga$.
If $A_j$ are constant, 
the frame $\vecv$ is called a normalizing frame.
When $\lambda\neq 0$,
we can always take a normalizing frame.
On the contrary, we do not have a normalizing frame in general,
in the case $\lambda=0$.

We will see that $\prolong{\nbige^{\lambda}}$ is locally free.
Simpson has already shown
that $\prolong{\nbige}^{\lambda}$ is locally free,
if the base manifold is one dimensional.
By using the fact, it is easy to see
that the normalizing frame of the prolongment $\nbige^{\lambda}$
gives, in fact, the frame of $\prolong{\nbige}^{\lambda}$.
However,
we have to prove some extension property
of holomorphic sections on hyperplanes,
if $\lambda=0$, 

One more point which we should care the conjugacy classes of
the residues $N_i$.
When the base manifold is one dimensional,
Simpson showed that the conjugacy classes
of $N_{|(\lambda,O)}$ are independent of $\lambda$.
Thus it is easy to see that 
the conjugacy classes of $N_{1\,|\,(\lambda,Q)}$ 
are independent of $\lambda$,
when $Q$ is contained in $D_{\itibar}-D_{\nibar}$.
On the other hand, it is easy to see
the conjugacy classes of $N_{1\,|\,(\lambda,Q)}$
are independent of $Q$ due to the existence of a normalizing frame,
when we fix $\lambda\neq 0$.
As a result, we can immediately obtain that
the conjugacy classes of $N_{1\,|\,(\lambda,Q)}$ are independent
of $(\lambda,Q)\in
  \cnum_{\lambda}\times D_{\itibar}
   -\{0\}\times D_{\nibar}$.
However we need some argument to see
that the degeneration of the conjugacy classes
does not occur at $\{0\}\times D_{\nibar}$.
Interestingly, we can show that the conjugacy classes
of $N_{1\,|\,(\lambda,Q)}$ are independent of $\lambda$
for any $Q\in D_{\itibar}$,
by using a limiting mixed twistor theorem.

\subsection{Outline of the paper}

\subsubsection{Section \ref{section;12.15.1}}
In the subsections \ref{subsection;12.15.2}--\ref{subsection;12.15.3},
we make a preparation on  the commuting tuples of nilpotent maps
and the weight filtrations,
on vector spaces or vector bundles.
The most part is well known and standard.
What we would like to see is summarized 
in Corollary \ref{cor;12.15.100},
Proposition \ref{prop;12.10.80}
and Corollary \ref{cor;12.15.101}.
In the subsection \ref{subsection;12.15.4},
we recall mixed twistor structure
and give lemmas which will be used later.
The author feels that the mixed twistor is useful
when we would like to obtain a lower bound of the degree
with respect to a filtration.

\subsubsection{Section \ref{section;12.15.70}}
In the subsection \ref{subsection;12.15.5},
we recall the definition of harmonic bundles.
Following Simpson (\cite{s3} and \cite{s4}),
we consider the deformed holomorphic bundle
and the conjugate.
In the subsection \ref{subsection;11.27.27},
we recall some easy examples of harmonic bundles
over the punctured disc $\Delta^{\ast}$.
In particular,
the model bundles $Mod(l,a,C)$ will be used later
as a convenient tool.
These examples can be closely investigated by direct calculations.
Some results,
for example, the corollaries
\ref{cor;12.15.110},
\ref{cor;11.26.5},
\ref{cor;11.27.41},
\ref{cor;12.2.5} and \ref{cor;11.27.42}
will be used.
In the subsection \ref{subsection;12.4.10},
we discuss the `convergency' of a sequence of the harmonic bundles.
We will see that we can pick a nice `convergent' subsequence,
as is naturally expected from the ellipticity of harmonic bundles.

\subsubsection{Section \ref{section;12.15.71}}
In the subsection \ref{subsubsection;12.9.60},
we prepare some words used in this paper.
In the subsection \ref{subsection;12.13.1},
we recall definitions of tameness and nilpotentness.
We will recall an estimate of the norms of Higgs fields
with respect to the metric.
In the subsection \ref{subsection;12.15.10},
we will recall some results of Simpson in one dimensional case.
In the subsection \ref{subsection;12.15.11},
we will recall the definition of triviality of parabolic structures.
We will also see that
the tame nilpotent rank one harmonic bundle
with trivial parabolic structure is smooth.

{\em
After the subsection {\rm\ref{subsection;12.15.11}},
the harmonic bundles will be always assumed to be
tame, nilpotent and has trivial parabolic structure.}

In the remaining of the section 4,
we will see that the prolongment $\prolong{\nbige}$ of
the deformed holomorphic bundle by an increasing order
are locally free.
We will use some ideas of Cornalba-Griffiths (\cite{cg}).
In the subsection \ref{subsection;12.15.15},
we recall something from their paper.
In the subsubsections
\ref{subsubsection;12.8.6}--\ref{subsubsection;12.15.17},
a normalizing frame gives a frame over
$\nbigx^{\shikaku}$,
and we show that $\prolong{\nbige}$ is locally free
if $\prolong{\nbige^0}$ is locally free.
In particular, it implies that
the $\prolong{\nbige}$ is locally free when
the base manifold is one dimensional.
In the subsubsection \ref{subsubsection;12.15.18},
we state the family versions of the results in
the subsection \ref{subsection;12.15.10}.
In the subsection \ref{subsection;12.15.19},
we will show the extendability
of the sections on hyperplane in the case $\lambda=0$.
It immediately implies that $\prolong{\nbige^0}$ is locally free,
and thus that the $\prolong{\nbige}$ is locally free.
In the subsection \ref{subsection;12.15.20},
we see the functoriality of the prolongment.

\subsubsection{Section \ref{section;12.15.75}}
In the subsections \ref{subsection;12.15.21}--\ref{subsection;12.15.22},
we recall the construction of the vector bundle $S(Q,P)$
of Simpson.
In the subsection \ref{subsection;12.15.23},
a limiting mixed twistor theorem is given and proved,
in the one dimensional case.
In the subsubsection \ref{subsubsection;12.15.24},
a refinement for higher dimensional case is given.
Once we obtain a mixed twistor structure,
briefly speaking,
we know that `a degeneration at $\lambda=0$ does not occur'.
Some easy and useful consequences of such type
are given in the subsubsections
\ref{subsubsection;12.15.25} and \ref{subsubsection;12.10.40}.
In the subsubsection \ref{subsubsection;12.15.26},
we give a weak constantness of the filtrations
as an easy consequence of a limiting mixed twistor structure,
although it can be shown without mixed twistor structure.

\subsubsection{Section \ref{section;12.15.76}}
In the subsection \ref{subsection;12.4.1},
we explain our method of comparison to obtain some estimate for metrics.
The method will be used in the beginning of all of the latter sections.
Briefly speaking, the method reduces the estimate
over the region to the estimate over the boundary.
Since the dimension of the boundary is lower,
we can use an estimate for lower dimensional case.
(In the subsection \ref{subsection;12.4.1}, we only need
 a rough estimate on the boundary.)
In the subsection \ref{subsection;12.4.70},
we consider the morphisms
$\psi_{m,\itibar}:\Delta^n\lrarr \Delta^n$
defined by $(z_1,\ldots,z_n)\longmapsto (z_1^m,z_2,\ldots,z_n)$.
For a harmonic bundle $\harmonicbundle$
over $\Delta^{\ast\,l}\times \Delta^{n-l}$,
we obtain the sequence
$\bigl\{\psi_{m,\itibar}^{\ast}\harmonicbundle\bigr\}$
of harmonic bundles.
We can apply the result of the subsection \ref{subsection;12.4.10},
due to the rough norm estimate obtained
in the subsection \ref{subsection;12.4.70}.
In the subsection \ref{subsection;12.15.30},
we see some orthogonality in the limit.
In particular, we see that the limiting harmonic bundle
is a CVHS in one dimensional case.
In the subsection \ref{subsection;12.15.31},
we investigate the first Chern class of the vector bundle
$\graded^{W(N^{\sankaku(1)}(\nibar)) }_{h_1,h_2}$ over $\proj^1$,
obtained in the subsubsection \ref{subsubsection;12.10.40}.

\subsubsection{Section \ref{section;12.15.32}}
In the section \ref{section;12.15.32},
we will prove the constantness of the filtrations
on the positive cones for the tuple of residues of harmonic bundles.
As a preliminary,
we give a norm estimate in some special case
in the subsection \ref{subsection;12.5.100}.
Then we consider the morphisms
$\psi_{m,\nbar}:\Delta^n\lrarr\Delta^n$
defined by $(z_1,\ldots,z_n)\longmapsto (z_1^m,\ldots,z_n^m)$
in the subsection \ref{subsection;12.15.33}.
By investigating the limiting harmonic bundle,
we obtain the constantness of the filtrations
on the positive cones in some special case.
Although the theorem is stated in the subsection
\ref{subsection;12.15.34},
the main part of the proof is done
in the subsection \ref{subsection;12.15.33}.

\subsubsection{Section \ref{section;12.15.35}}

In the section \ref{section;12.15.35},
we see the strongly sequential compatibility of the residues.
In the subsubsection \ref{subsubsection;12.15.37},
we see some compatibility in the bottom part
in the two dimensional case.
By using a method of comparison,
we obtain some implication of the filtrations.
Then we obtain the coincidence
by using the result in the subsubsections \ref{subsubsection;12.10.40}
and Lemma  \ref{lem;1.8.10} for mixed twistor structures.
Once we know such compatibility in two dimensional case,
a similar compatibility in the higher dimensional case is easy to
obtain, which is shown in the subsubsection \ref{subsubsection;12.15.40}.
Then we obtain the theorems in the subsection \ref{subsection;12.15.41}.

\subsubsection{Section \ref{section;12.15.80}}

In the subsection \ref{subsection;12.15.45},
we obtain a norm estimate.
As a preliminary,
we consider the pull back of harmonic bundle $\harmonicbundle$
via the `blowup' $\phi_N:\blowup{X}\lrarr X$,
and we obtain the norm estimate for
$\phi_N^{\ast}\harmonicbundle$
in the subsubsection \ref{subsubsection;12.7.100}.
The method is same as that in \ref{subsection;12.4.1}.
Since we have already shown the strongly sequential compatibility of
the residues in the subsection \ref{subsection;12.15.41},
we can obtain a stronger estimate.
By translating such a result,
we obtain the theorem in the subsubsection
\ref{subsubsection;12.15.46}.
In the subsubsection \ref{subsubsection;12.14.1},
we consider the pull backs of the harmonic bundles
via the morphism $\psi_{m,\nbar}$,
as in the subsection \ref{subsection;12.15.33}.
Then we obtain a limiting harmonic bundle.
We see that it is, in fact, a CVHS.
As an application of limiting CVHS,
we see the purity theorem
in the subsubsection \ref{subsubsection;12.15.50}.

\subsection{Some remarks}

Unfortunately, this paper looks rather long.
However, the reader will know that
much of the part is elementary and
not new for both of the reader and the author.
Many of the definitions, the lemmas and the propositions are
more or less standard, familiar and obvious.
They are included to clarify what the author would like to say.

This paper is the revision of \cite{mochi}.
The main difference is as follows:
\begin{itemize}
\item
 In \cite{mochi}, the dimensions of the base manifolds
 are assumed to be less than two.
 In this paper, we discuss higher dimensional case.
 The much part of the preliminary for the filtrations
 are added for that purpose.
\item
 The explanation for the norm estimate in one dimensional case
 is added. In \cite{s2}, such estimates are stated for
 the cases $\lambda=0$ and $\lambda=1$.
 Clearly his argument works for the other $\lambda$.
 We only indicate how to change.
\item
 The explanation for the prolongation of the deformed holomorphic bundle
 is added. 
\item
 The mixed twistor structure is used more efficiently.
 As a result, some arguments for the filtrations on the 
 divisors are simplified.
\item
 The author hopes that the discussion and the explanation
 in this paper are more clear
 than the previous version \cite{mochi}.
\end{itemize}

The author's original motivation of the study is
to generalize the Kobayashi-Hitchin correspondence
(see for example, \cite{b}, \cite{don}, \cite{s1}, \cite{s2} and \cite{uy}).
Namely he would like to clarify the relation
of stable Higgs bundles and harmonic bundles.
For that purpose
it is clearly important to characterize the residues
of the Higgs fields.
Unfortunately, the understanding seems insufficient
to solve such problem in this stage,
for the author.
Probably, one direction of the study
is a more precise comparison of a limiting CVHS
and the original harmonic bundle.

\subsection{Acknowledgement}

The author sincerely thanks the referee
for his effort to read the paper
and his useful suggestions.
He kindly informed the author
on problems mentioned in the paper of Jost-Li-Zuo \cite{jlz}.
It seems quite related to the original purpose of the author.
Surely it will be one of our next subjects in this area.

The author is grateful to the colleagues in Osaka City University.
In particular, he thanks Mikiya Masuda for his encouragement.

The author thanks the financial supports by Japan Society for
the Promotion of Science and the Sumitomo Foundation.
The paper was written during his stay at the Institute
for Advanced Study. The author is sincerely grateful
to their excellent hospitality.
He also acknowledges National Scientific Foundation
for a grant DMS 9729992,
although any opinions, findings and conclusions or recommendations
expressed in this material are those of the author.

\subsection{Some sets}

We will use the following notation:

\begin{tabular}{llll}
$\seisuu$: & the set of the integers, &
$\seisuu_{>0}$: &the set of the positive integers,\\
$\rnum$: &the set of the rational numbers,&
$\rnum_{>0}$: &the set of the positive rational numbers,\\
$\real$: & the set of the real numbers, &
$\real_{>0}$: &the set of the positive real numbers,\\
$\cnum$: &the set of the complex numbers, &
$\nbar$: & the set $\{1,2,\ldots,n\}$,\\
$\gbigs_l$: & the $l$-th symmetric group.
\end{tabular}

We put as follows:
\[
 \Delta(C):=\{z\in\cnum\,|\,|z|<C\},\quad
 \Delta^{\ast}(C):=\{z\in\cnum\,|\,0<|z|<C\},\quad
 \cnum^{\ast}=\{z\in \cnum\,|\,z\neq 0\}.
\]
When $C=1$, we often omit to denote $C$,
i.e.,
$\Delta=\Delta(1)$ and $\Delta^{\ast}=\Delta^{\ast}(1)$.
If we emphasize the variable,
we describe as $\Delta_z$, $\Delta_i$.
For example,
$\Delta_z\times \Delta_w=\{(z,w)\in\Delta\times \Delta\}$,
and $\Delta_1\times \Delta_2=\{(z_1,z_2)\in\Delta\times \Delta\}$.
We often use the notation $\cnum_{\lambda}$ and $\cnum_{\mu}$.

Unfortunately, the notation $\Delta$ is also used to denote the
Laplacian. The author hopes that there will be no confusion.

The set of $r\times r$-matrices with $\cnum$-coefficient
is denoted by $M(r)$,
and the set of $r\times r$-hermitian matrices by $\nbigh(r)$.

In general, $q_i:X^n\lrarr X$ denotes the projection onto
the $i$-th component,
and $\pi_i:X^n\lrarr X^{n-1}$
denotes the projection omitting the $i$-th component.
However we will often use $\pi$ to denote some other projections.

Let $l$ be a positive integer.
We have the decomposition of $\real_{\geq 0}^l$ 
into $\coprod_{I\subset\lbar}\real_{>0}^{I}$, defined as follows:
\[
 \real_{>0}^I:=
 \{(a_1,\ldots,a_l)\in\real_{\geq\,0}^l\,|\,
 a_i>0 \Longleftrightarrow i\in I\}.
\]

Let $l$ be a positive integer, $C$ be a positive real number,
and $\sigma$ be an element of $\gbigs_l$.
Then we put as follows:
\[
 Z(\sigma,l,C):=
 \bigl\{(z_1,\ldots,z_n)\in\Delta^{\ast\,l}\times \Delta^{n-l}\,\big|\,
 |z_{\sigma(i-1)}|^C<|z_{\sigma(i)}|\,\,
 i=2,\ldots,l\bigr\}.
\]

\section{Preliminary for filtrations}
\label{section;12.15.1}

\subsection{Vector space and filtrations}
\label{subsection;12.15.2}

\subsubsection{Base and metric}

Let $V$ be an $n$-dimensional vector space over $\cnum$.
To describe a base of $V$,
or more generally, to describe a tuple of elements of $V$,
we use a notation $\vecv=(v_1,\ldots,v_n)$.
Let $\vecv$ and $\vecw$ be two bases of $V$.
We obtain the matrix $A=(A_{i\,j})$
determined by the following formula:
\[
 v_j=\sum_i A_{i\,j}\cdot w_i.
\]
In that case, it is described as $\vecv=\vecw\cdot A$.

Let $h$ be a hermitian metric of $V$.
Then we have the hermitian matrix $H(h,\vecv)=(H_{i\,j})$
determined as follows:
\[
 H_{i\,j}:=h(v_i,v_j).
\]
The  $H(h,\vecv)\in\nbigh(n)$ is called the hermitian matrix
of the metric $h$ with respect to $\vecv$.

\subsubsection{Compatibility with direct sum}

Let $V$ be a finite dimensional vector space
with a direct sum decomposition $V=\bigoplus_i V_i$.

Let $v$ be a non-zero element of $V$.
It is called compatible with the decomposition
if there exists an $i$ such that $v\in V_i$.
The number $i$ is called the degree of $v$.

Let $\vecv=(v_1,\ldots,v_n)$ be a base of $V$.
It is compatible with the decomposition
if each $v_i$ is compatible with the decomposition.

Let $W$ be an increasing filtration of $V$.
It is called compatible with the decomposition
if $W_j=\bigoplus_i (W_j\cap V_i)$
for any $j$.
We denote the induced filtration of $V_i$
by $W\cap V_i$.

Let $f$ be an endomorphism of $V$.
It is called compatible with the decomposition if $f(V_i)\subset V_i$.

For a tuple $(u_1,\ldots,u_l)$ of elements of $V$,
$\langle u_1,\ldots,u_l\rangle$ denotes the vector subspace
generated by $(u_1,\ldots,u_l)$.

\subsubsection{Filtration}

In this paper, we mainly use increasing filtrations.
Thus `filtration' means an `increasing filtration'
if we do not notice.

Let $V$ be a vector space with a filtration $W$,
the associated graded vector space is 
denoted by $Gr^W=\bigoplus_i Gr^W_i$,
where $Gr^{W}_i:=W_i/W_{i-1}$.

For a filtration $W$,
we have the number $b(W)$ determined as follows:
\[
 b(W):=\min\{\,h\,|\,Gr^W_h\neq 0\}.
\]
The number $b(W)$ is called the bottom number of $W$.

For a non-zero element $v\in V$,
the number $\deg^W(v)$ is defined as follows:
\[
 \deg^W(v):=\min\{h\,|\,v\in W_h\}.
\]
The number $\deg^W(v)$ is called the degree of $v$
with respect to the filtration $W$.
We have the induced element $v^{(1)}$ of
$Gr^W_{\deg^W(v)}(V)$.

Let $\vecv=(v_1,\ldots,v_n)$ be a base of $V$.
We say that $\vecv$ is compatible with the filtration $W$,
if the following is satisfied:
\begin{quote}
For any $i$, we have a subset $I\subset \{1,\ldots,n\}$
such that $\{v_j\,|\,j\in I\}$ gives a base of $W_i$.
\end{quote}
In that case,
the induced elements $\vecv^{(1)}=\{v_1^{(1)},\ldots,v^{(1)}_n\}$
gives a base of $Gr^W(V)$ compatible with the natural decomposition.

An endomorphism $f$ of $V$ is called compatible 
with the filtration $W$ if $f(W_i)\subset W_i$.

\subsubsection{The induced filtrations}

Let $W(\itibar)$ and $W(\nibar)$ be a filtrations on $V$.
We have the graded associated space of $W(\itibar)$:
\[
 Gr^{(1)}:=\bigoplus_{a}Gr^{(1)}_a,
\quad
 Gr^{(1)}_a:=Gr^{W(\itibar)}_a.
\]
We have 
the induced filtration $W^{(1)}(\nibar)$
by $W(\nibar)$ on $Gr^{(1)}$, which is defined as follows:
\[
 W^{(1)}(\nibar)_l
:=\bigoplus_a W^{(1)}(\nibar)_l\cap Gr_a^{(1)},
\quad
 W^{(1)}(\nibar)_l\cap Gr_a^{(1)}
:=\frac{W(\nibar)_{l+a}\cap W(\itibar)_a }{W(\itibar)_{a-1}}
  \subset Gr_a^{(1)}.
\]

Let $(W(\itibar),W(\nibar),\ldots, W(\nbar))$ be a filtrations on $V$.
We have the induced filtrations
$(W^{(1)}(\nibar),\ldots,W^{(1)}(\nbar))$
on $Gr^{(1)}$.
Inductively, we obtain the filtrations $W^{(m)}(\jbar)$
on $Gr^{(m)}$ for $1\leq m<j\leq n$ as follows:
\begin{enumerate}
\item
On $Gr^{(m)}$, 
we have the filtrations $W^{(m)}(\jbar)$
for $(j=m+1,\ldots,n)$.
\item
Then we put $Gr^{(m+1)}:=Gr^{W^{(m)}(\mitibar)}$.
We have the filtrations $W^{(m+1)}(\jbar)$ $(j=m+2,\ldots,n)$
induced by $W^{(m)}(\jbar)$.
\end{enumerate}

\subsubsection{Compatible sequence of filtrations}

Let $(W(\itibar),W(\nibar),\ldots, W(\nbar))$ be a filtrations on $V$.
Let $\vech=(h_1,\ldots,h_n)$  be a tuple of integers.
We have the following morphism:
\[
 \pi_{\vech}:
 \bigcap_{j=1}^nW(\jbar)_{h_j}
 \lrarr Gr_{h_1}^{(1)}.
\]
The image of $\pi_{\vech}$ is always contained
in $Gr_{h_1}^{(1)}\cap \bigcap_{j=2}^n W^{(1)}(\jbar)_{h_j}$.

\begin{df}
A sequence of the filtration $(W(\itibar),\ldots,W(\nbar))$
is called compatible,
if the following holds, inductively:
\begin{enumerate}
\item
$(W^{(1)}(\nibar),\ldots,W^{(1)}(\nbar))$ is a compatible sequence.
\item
For any $\vech\in \seisuu^n$,
the image of $\pi_{\vech}$ is same as
$Gr_{h_1}^{(1)}\cap \bigcap_{j=2}^n W^{(1)}(\jbar)_{h_j}$.
\hfill\qed
\end{enumerate}
\end{df}

\begin{rem}
When $n\leq 2$, the condition is trivial.
\end{rem}

\begin{df}
Let $(W(\itibar),\ldots,W(\nbar))$ is a compatible sequence of filtrations.
A non-zero element $f\in V$ is called compatible with
the sequence $(W(\itibar),\ldots,W(\nbar))$ if the following holds,
inductively:
\begin{enumerate}
\item The induced element $f^{(1)}\in Gr^{(1)}$ is compatible
with the sequence $(W^{(1)}(\nibar),\ldots,W^{(1)}(\nbar))$.
\item
For any $j\geq 2$,
we have $\deg^{W(\jbar)}(f)=\deg^{W^{(1)}(\jbar)}(f^{(1)})$.
\hfill\qed
\end{enumerate}
\end{df}

\begin{df}
Let $(W(\itibar),\ldots,W(\nbar))$ is a compatible sequence of filtrations.
A base $\vecv=(v_i)$ of $V$ is called compatible,
if the following holds, inductively:
\begin{itemize}
\item
 For each $i$, the element $v_i$ is compatible
 with the sequence $(W(\itibar),\ldots,W(\nbar))$.
\item
 $\vecv$ is a compatible base with $W(\itibar)$.
\item
 The induced base $\vecv^{(1)}$ of $Gr^{(1)}$ 
 is compatible with the sequence
 $(W^{(1)}(\nibar),\ldots,W^{(1)}(\nbar))$.
\hfill\qed
\end{itemize}
\end{df}

Let consider a decomposition of $V$
into $\bigoplus_{\vech\in\seisuu^{n}}U_{\vech}$.
\begin{df}
Let $(W(\itibar),\ldots,W(\nbar))$ is a compatible sequence of
 filtrations.
A decomposition $V=\bigoplus_{\vech\in\seisuu^{n}}U_{\vech}$
is called a splitting compatible with
the sequence $(W(\itibar),\ldots,W(\nbar))$,
if the following holds for any $\vech\in\seisuu^n$:
\[
 \bigcap_{j=1}^nW(\itibar)_{q_j(\vech)}
=\bigoplus_{\veck\in \nbigt(\vech)}
 U_{\veck},
\quad
 \nbigt(\vech):=
 \{\veck\in\seisuu^n\,|\,q_j(\veck)\leq q_j(\vech)\}.
\]
Here $q_j$ denotes the projection $\seisuu^n\lrarr \seisuu$
onto the $j$-th component.
\hfill\qed
\end{df}

The dimension of $U_{\vech}$ is denoted by $d(\vech)$.

\begin{lem}
Let $(W(\itibar),\ldots,W(\nbar))$ is a compatible sequence
of filtrations.
\begin{enumerate}
\item
Let $\vecv$ be a base of $V$ compatible with
the decomposition $V=\bigoplus_{\vech\in\seisuu^n}U_{\vech}$.
Assume that the decomposition is compatible with
$(W(\itibar),\ldots,W(\nbar))$.
Then the base $\vecv$ is compatible with
$(W(\itibar),\ldots,W(\nbar))$
\item
Let $\vecv$ be a base of $V$ compatible with
the sequence $(W(\itibar),\ldots,W(\nbar))$.
We put as follows:
\[
 U_{\vech}:=
 \langle v_i\,|\, \deg^{W(\jbar)}(v_i)=q_j(\vech)
 \rangle.
\]
Then the decomposition $V=\bigoplus_{\vech\in\seisuu^n}U_{\vech}$
is compatible with the sequence $(W(\itibar),\ldots,W(\nbar))$.
\end{enumerate}
\end{lem}
\pf
Let see the first claim.
We use an induction on $n$.
Let $\vecv$ be a base compatible with the decomposition.
Clearly it is compatible with the filtration
$W(\itibar)$.

Let $\pi_{h}'$ denote the composite of the following morphisms.
Here we put $h_1=h$ for simplicity of notation:
\[
 \bigoplus_{
 \substack{\veck\in \seisuu^n\\ q_1(\veck)=h
 }}U_{\veck}
\stackrel{\subset}{\lrarr}
 \bigcap_{j=1}^nW(\jbar)_{h_j}
\stackrel{\pi_{h}}{\lrarr}
 Gr^{(1)}_{h}.
\]
Since the decomposition is compatible with the sequence of the filtrations,
$\pi_h'$ is isomorphic.
Thus we obtain the decomposition of $Gr^{(1)}_{h}$ as follows:
\[
 Gr^{(1)}_{h}=
 \bigoplus_{
 \substack{\veck\in \seisuu^n\\ q_1(\veck)=h
 }}\pi_h'(U_{\veck}).
\]
Due to the following isomorphism,
the decomposition is compatible with the sequence
$(W^{(1)}(\nibar),\ldots,W^{(1)}(\nbar))$:
\[
  \bigoplus_{
 \substack{\veck'\in \nbigt(\vecl)
 }}\pi_h'(U_{(h,\veck')})
\simeq
 \frac{W(\itibar)_h\cap \bigcap_{j=1}^{n-1} W(\jitibar)_{l_{j}}}
 {W(\itibar)_{h-1}\cap \bigcap_{j=1}^{n-1} W(\jitibar)_{l_{j}}}
\simeq
 Gr^{(1)}_{h}\cap
 \bigcap_{j=1}^{n-1}
 W^{(1)}(\jitibar)_{l_j},
\]
Here we put
$\nbigt(\vecl)=
 \{
 \veck'\in \seisuu^{n-1},
 q_j(\veck)\leq l_j,\,\,j=1,\ldots n-1
\}$,
and $(h,\veck')$ denotes $(h,k_1',\ldots,k_{n-1}')$
for $\veck'=(k_1',\ldots,k_{n-1}')$.
By our assumption of the induction,
the induced base $\vecv^{(1)}$ is compatible
with the sequence $(W^{(1)}(\nibar),\ldots,W^{(1)}(\nbar))$.
We also have $\deg^{W(\jbar)}(v_i)=\deg^{W^{(1)}(\jbar)}(v_i)$
for $j\geq 2$.
Thus we obtain the first claim.

The second claim can be shown similarly.
\hfill\qed

\begin{lem}
Let $(W(\itibar),\ldots,W(\nbar))$ be a compatible sequence
of the filtrations.
There exists a decomposition
$V=\bigoplus_{\vech\in\seisuu^n}U_{\vech}$
compatible with the sequence $(W(\itibar),\ldots,W(\nbar))$.
\end{lem}
\pf
We use an induction on $n$.
Let consider the induced filtrations 
$Gr^{(1)}_h\cap W^{(1)}(\jbar)$ on $Gr^{(1)}_h$,
which is easily checked to be compatible.
By our assumption of the induction,
we can take a compatible decomposition:
\[
 Gr^{(1)}_h=\bigoplus_{\veck\in\seisuu^{n-1}}
 U'_{h,\veck},
\quad\quad
 Gr_{h}^{(1)}\cap
 \bigcap_{j=1}^{n-1}
  W^{(1)}(\jitibar)_{l_j}
=\bigoplus_{\veck\in \nbigt(\vecl)}
 U'_{h,\veck}.
\]
For an element $\vech=(h_1,\ldots,h_n)\in\seisuu^n$,
we put $\vech'=(h_2,\ldots,h_n)$.
Since $(W(\itibar),\ldots,W(\nbar))$ is compatible,
$U'_{h_1,\vech'}$ is contained in the following morphism:
\[
 \pi_{\vech}:
 W(\itibar)_{h_1}\cap
 \bigcap_{j=2}^nW(\jbar)_{h_j}
\lrarr
 Gr^{(1)}_{h_1}\cap
 \bigcap_{j=2}^nW^{(1)}(\jbar)_{h_j}.
\]
We pick the subspace $U_{\vech}$ of
$ W(\itibar)_{h_1}\cap \bigcap_{j=2}^nW(\jbar)_{h_j}$,
which is isomorphic to $U'_{h_1,\vech'}$
via the morphism $\pi_{\vech}$.
Then we obtain the decomposition
$V=\bigoplus_{\vech\in\seisuu^n}U_{\vech}$.

By our construction,
we have the following:
\[
 \left[
  W(\itibar)_{h_1}\cap
 \bigcap_{j=2}^nW(\jbar)_{h_j}
 \right]
= \left[
 W(\itibar)_{h_1-1}\cap
 \bigcap_{j=2}^nW(\jbar)_{h_j}
 \right]
\oplus
 \left[
 \bigoplus_{\veck\in\nbigt(\vech')}
 U_{(h_1,\veck)}\right].
\]
Thus the induction can proceed.
\hfill\qed

\subsection{A commuting tuple of nilpotent maps}

\subsubsection{Tensor product, symmetric products and exterior products}
\label{subsubsection;12.4.145}

Let $V$ be a finite dimensional vector space over $\cnum$
and $N$ be an endomorphism of $V$.
Then we have the tensor product $V^{\otimes\,h}$,
the symmetric products $\Sym^h(V)$ and
the exterior products $\bigwedge^h(V)$.
We also have the endomorphism of $V^{\otimes\,h}$:
\[
 \tilde{N}:=
 \sum \overbrace{1\otimes\cdots \otimes 1}^{j-1}\otimes N
 \otimes
 \overbrace{1\otimes\cdots\otimes 1}^{h-j}.
\]
We often denote $\tilde{N}$ by $N^{\otimes\,n}$.
We will not use the endomorphism
$\overbrace{N\otimes\cdots\otimes N}^n$.
Thus the author hopes that any confusion does not occur.
The morphism $N^{\otimes\,n}$ preserves
the subspaces $\Sym^h(V)$ and $\bigwedge^h(V)$.
Thus it induces the endomorphisms
of $\Sym^h(V)$ and $\bigwedge^h(V)$.
We denote them by $N^{\sym\,h}$ and $N^{\wedge\,h}$
respectively.

\subsubsection{The weight filtration of nilpotent maps}

Let $V$ be a finite dimensional vector space over $k$
and $N$ be a nilpotent map of $V$.
Recall that $N$ induces the weight filtration 
$W(N)$ of $V$,
which is characterized by the following properties:
\begin{itemize}
\item
 $N\cdot W_{l}(N)\subset W_{l-2}(N)$.
\item
 The induced morphism $N^k:Gr^{W(N)}_k\lrarr Gr^{W(N)}_{-k}$
 is isomorphic for any $k\geq 0$.
\end{itemize}
We have the following obvious lemma.
\begin{lem}
Let $\vecv$ be a base of $V$ compatible with the filtration $W(N)$.
Then we have the following equality:
\[
 \sum_i \deg^{W(N)}(v_i)=0.
\]
\hfill\qed
\end{lem}
For any $l\geq 0$, we put
$ P_lGr_l^{W(N)}:=
 \Ker\bigl(
 N^{l+1}:Gr_l^{W(N)}\lrarr Gr_{-l-2}^{W(N)}
 \bigr)$.
When $l-a=2m\geq 0$ for some non-negative integer $m$,
we put $P_lGr_a^{W(N)}:=\Image(N_1^m:P_lGr_l^{W(N)}\lrarr Gr_a)$.
Then we obtain the decomposition
$Gr_a:=\bigoplus_{0\leq h}P_{a+2h}Gr_a$,
which is called the primitive decomposition.

\subsubsection{Splittings of the weight filtrations}
\label{subsubsection;12.9.50}

Let $sl_2$ is a Lie subalgebra of the $(2\times 2)$-matrix algebra
$M(2)$ with the following base:
\[
 N_0=
 \left(
 \begin{array}{cc}
 0 & 1 \\ 0 & 0
 \end{array}
 \right),
\quad
 N_1=
  \left(
 \begin{array}{cc}
 0 & 0 \\ 1 & 0
 \end{array}
 \right),
\quad
C=\left(
 \begin{array}{cc}
 1 & 0 \\ 0 & -1
 \end{array}
 \right).
\]
Assume that we have a homomorphism of Lie algebras
$\eta:sl_2\lrarr End(V)$.
Then we obtain the weight filtration 
$W(\eta(N_1))$,
and the decomposition of $V$ into the eigenspaces of $\eta(C)$:
\[
 V=\bigoplus_{\alpha}V_{\alpha}.
\]
Here $\alpha$ runs through the eigenvalues of $\eta(C)$.
Then we have the following:
\[
 W_l(\eta(N_1))=\bigoplus_{h\leq l}V_{h}.
\]

For any non-negative integer $n$,
we have the naturally induced representation
$\eta^{\otimes\,n}:sl_2\lrarr End(V^{\otimes\,n})$.
For an $n$-tuple $\vecalpha=(\alpha_1,\ldots,\alpha_n)$
of the eigenvalues of $\eta(C)$,
we put $V_{\vecalpha}:=V_{\alpha_1}\otimes\cdots\otimes V_{\alpha_n}$.
We have the decomposition of $V^{\otimes\,n}$ as follows:
\[
 V^{\otimes\,n}=
 \bigoplus_{\vecalpha}V_{\vecalpha}.
\]
Here $\vecalpha$ runs through the set of $n$-tuples 
of the eigenvalues of $\eta(C)$.
It is clear that $V_{\vecalpha}$ is contained
in the eigenspace of $\eta^{\otimes\,n}(C)$ with
the eigenvalue $\rho(\vecalpha)=\sum_{i=1}^n\alpha_i$.
Thus we have the eigen decomposition of $V^{\otimes\,n}$:
\begin{equation} \label{eq;12.2.25}
 V^{\otimes\,n}=
 \bigoplus_{\alpha}
 \Bigl(
 \bigoplus_{\substack{\vecalpha\\ \rho(\vecalpha)=\alpha}}
 V_{\vecalpha}\Bigr).
\end{equation}
Here $\alpha$ runs through the eigenvalues of $\eta^{\otimes\,n}(C)$.
Then we obtain the following:
\begin{equation} \label{eq;12.2.27}
 W_l(\eta(N_1)^{\otimes\,n})=
 \bigoplus_{\substack{\vecalpha\\\rho(\vecalpha)\leq l}}
 V_{\vecalpha}
=\Bigl\{
 \sum x_1\otimes\cdots\otimes x_n
 \,\Big|\,\sum \deg^{W(\eta(N_1))}(x_i)\leq l
 \Bigr\}.
\end{equation}

Let consider the case of the symmetric product
and the exterior product.
We put as follows for a set $I$:
\[
 \nbigs(I,n):=
 \Bigl\{
 f:I\lrarr \seisuu_{\geq\,0}
 \,\Big|\,
 \sum_{i\in I}f(i)=n
 \Bigr\}.
\]
For $f\in\nbigs(I,n)$,
we put $\rho(f):=\sum_{i\in I}i\cdot f(i)$.

Let $I$ be the set of the eigenvalues of $\eta(C)$.
Then we have the following decomposition:
\begin{equation} \label{eq;12.2.26}
 \Sym^n(V)=\bigoplus_{f\in\nbigs(I,n)}
 \bigotimes_{\alpha\in I}\Sym^{f(\alpha)} V_{\alpha},
\quad
 \bigwedge^n(V)=\bigoplus_{f\in\nbigs(I,n)}
 \bigotimes_{\alpha\in I}\bigwedge^{f(\alpha)}V_{\alpha}.
\end{equation}
By considering the eigenvalues
of $\eta(C)^{\sym\,n}$ and $\eta(C)^{\wedge\,n}$,
we obtain the following:
\begin{equation} \label{eq;12.2.28}
\begin{array}{l}
{\displaystyle
 W_l\bigl(\eta(N_1)^{\sym\,n}\bigr)=
 \bigoplus_{\substack{f\in\nbigs(I,n),\\ \rho(f)\leq l}}
 \bigotimes_{\alpha\in I}\Sym^{f(\alpha)} V_{\alpha}=
\Bigl\{
\sum x_1\cdot \cdots\cdot x_n\,\Big|\,
 \sum\deg^{W(\eta(N_1))}(x_i)\leq l
\Bigr\}
} \\
{\displaystyle
 W_l\bigl(\eta(N_1)^{\wedge\,n}\bigr)=
 \bigoplus_{\substack{f\in\nbigs(I,n),\\ \rho(f)\leq l}}
 \bigotimes_{\alpha\in I}\bigwedge^{f(\alpha)}V_{\alpha}
=\Bigl\{
 \sum x_1\wedge\cdots\wedge x_n\,\Big|\,
 \sum\deg^{W(\eta(N_1))}(x_i)\leq l
 \Bigr\}.
}
\end{array}
\end{equation}

Let $N$ be a nilpotent map on a finite dimensional vector space $V$.
We can pick a representation $\eta:sl_2\lrarr End(V)$
such that $\eta(N_1)=N$.
Thus we obtain the following:
\[
 \begin{array}{l}
 W_l(N^{\otimes\,n})=
 \Bigl\{
 \sum x_1\otimes\cdots\otimes x_n\,\Big|\,
 \sum \deg^{W(N)}x_i\leq l
 \Bigr\},\\
 \mbox{{}}\\
 W_l(N^{\sym\,n})=
 \Bigl\{
 \sum x_1\cdot\cdots\cdot x_n\,\Big|\,
 \sum \deg^{W(N)}x_i\leq l
 \Bigr\},\\
 \mbox{{}}\\
  W_l(N^{\wedge\,n})=
 \Bigl\{
 \sum x_1\wedge\cdots\wedge x_n\,\Big|\,
 \sum \deg^{W(N)}x_i\leq l
 \Bigr\}.\\
 \end{array}
\]

Let assume that we have a splitting of the weight filtration $W(N)$,
i.e.,
we have a decomposition $V=\bigoplus_h U_h$
such that $W(N)_l=\bigoplus_{h\leq l}U_h$.
Then we have the decomposition of
the products $V^{\otimes\,n}$, $\Sym^n(V)$
and $\bigwedge^n V$ by the same formula as those
(\ref{eq;12.2.25}) and (\ref{eq;12.2.26}),
although the meaning is slightly different.
They give the splitting of the filtrations
$W(N^{\otimes\,n})$,
$W(N^{\sym\,n})$ and $W(N^{\wedge\,n})$
by the same formula as those (\ref{eq;12.2.27})
and (\ref{eq;12.2.28}).

\subsubsection{Compatibility of a commuting tuple of nilpotent maps}
Let $V$ be a finite dimensional vector space
with a decomposition $V=\bigoplus V_i$.
Let $N$ be a nilpotent endomorphism of $V$.
Then it induces the weight filtration,
which we denote by $W(N)$.
Recall that $\nbar$ denote the set $\{1,\ldots,n\}$.

\begin{df}
Let $N_1,\ldots,N_n$ be a tuple of nilpotent maps of $V$.
It is called a commuting tuple,
if $N_i$ and $N_j$ are commutative for any $i,j\in\nbar$.
\hfill\qed
\end{df}

\begin{df}
Let $(N_1,\ldots,N_n)$ be a commuting tuple of
nilpotent maps.
We say that the constantness of the induced filtration
 on the positive cones holds,
if the following holds:
\begin{quote}
For any subset $I\subset \nbar$,
the filtration
$W(\sum_{i\in I} a_i N_i)$ is independent of
$(a_i\,|\,i\in I)\in \real_{>0}^I$.
\end{quote}
When the constantness of the filtrations on the positive cone
holds,
we denote the filtration $W\bigl(\sum_{i\in I}a_iN_i\bigr)$ $(a_i>0)$
by $W(I)$.
\hfill\qed
\end{df}

Assume that the constantness of the filtrations
on the positive cones holds.
We put $N(\jbar)=\sum_{i\leq j}N_i$.
Let $W(\jbar)$ denote
the weight filtration of $N(\jbar)$.
We denote the graded vector space associated
with $W(\itibar)$ by $Gr^{(1)}$.
We have the projection $\pi_{h_1}:W(\itibar)_{h_1}\lrarr Gr^{(1)}_{h_1}$.
Let $\vech=(h_1,\ldots,h_n)$ denote an $n$-tuple of integers.
Then we have the following morphism:
\[
 \pi_{\vech}:
 W(\itibar)_{h_1}\cap
 W(\nibar)_{h_2}\cap\cdots\cap W(\nbar)_{h_n}
\lrarr 
 Gr^{(1)}_{h_1}.
\]
On the other hand,
$N(\jbar)$ induces the morphism $N^{(1)}(\jbar)$
on $Gr^{(1)}$.
Let $W(N^{(1)}(\jbar))$ denote the weight filtration
of $N^{(1)}(\jbar)$.
Then we have the subspace
$Gr^{(1)}_{h_1}\cap\bigcap_{j=2}^n W(N^{(1)}(\jbar))_{h_j-h_1}$
of $Gr^{(1)}_{h_1}$.

\begin{df} \label{df;12.3.1}
Let $(N_1,N_2,\ldots,N_l)$ be a commuting tuple of nilpotent maps.
It is called sequentially compatible,
if the following holds inductively:
\begin{itemize}
\item
The constantness of the filtrations on the positive cones holds.
\item 
 We have the induced tuple
 $(N_2^{(1)},\ldots,N_l^{(1)})$
 of the commuting nilpotent maps
 on $Gr^{W(N_1)}$.
 It is sequentially compatible.
\item
 For any $\vech=(h_1,\ldots,h_n)$, we have
 $\Image(\pi_{\vech})
=Gr^{(1)}_{h_1}\cap\bigcap_{j=2}^n W(N^{(1)}(\jbar))_{h_j-h_1}$.
\hfill\qed
\end{itemize}
\end{df}

\begin{rem} \label{rem;12.10.200}
The third condition in Definition {\rm\ref{df;12.3.1}}
can be reworded as follows:
\begin{itemize}
\item
Let $W^{(1)}(\jbar)$ denote the filtration of $Gr^{(1)}$
induced by $W(\jbar)$.
Then we have the following:
\[
 W^{(1)}(\jbar)_{l+a}\cap Gr^{(1)}_a
=W(N^{(1)}(\jbar))_l\cap Gr^{(1)}_a.
\]
\item
We have
$\Image(\pi_{\vech})=
 Gr^{(1)}_{h_1}\cap \bigcap_{j=2}^n W^{(1)}(\jbar)_{h_j}$
for any $\vech=(h_1,\ldots,h_n)\in\seisuu^n$.
\hfill\qed
\end{itemize}
\end{rem}

\begin{rem}
On $Gr^{(1)}$,
we have $N^{(1)}(\jbar)=\sum_{i\leq j}N^{(1)}_i=\sum_{i=2}^jN^{(1)}_i$.
\hfill\qed
\end{rem}

\begin{lem}
Assume that $(N_1,\ldots,N_n)$ is compatible.
Then $(W(\itibar),\ldots,W(\nbar))$ is a compatible sequence
of filtrations.
\end{lem}
\pf
We only have to note Remark \ref{rem;12.10.200}.
\hfill\qed

\begin{df}
Assume that $(N_1,\ldots,N_n)$ is sequentially compatible.
A base $\vecv$ is called compatible with the sequence
 $(N_1,\ldots,N_n)$,
if it is compatible with the sequence
$(W(\itibar),\ldots,W(\nbar))$.
\hfill\qed
\end{df}

\vspace{.1in}

For any $n$-tuple $\vech=(h_1,\ldots,h_n)$,
we have the following morphism:
\[
 P\pi_{\vech}:
\left[
 \Ker(N(\itibar)^{h_1+1})
 \cap \bigcap_{j=2}^n
 W(\jbar)_{h_j}\right]
\lrarr
 P_{h_1}Gr_{h_1}^{(1)}.
\]
\begin{lem}
When $(N_1,\ldots,N_n)$ is a sequentially compatible,
we have the following implication:
\[
\Image(P\pi_{\vech})\subset
 \left[
 P_{h_1}Gr_{h_1}^{(1)}\cap \bigcap_{j=2}^n W^{(1)}(\jbar)_{h_j}
 \right].
\]
\end{lem}
\pf
We have the following implication:
\[
 \Image(P\pi_{\vech})
 \subset
 \left[
 P_{h_1}Gr_{h_1}^{(1)}
\cap
 \pi_{\vech}\Bigl(
 \bigcap_{j=1}^n W(\jbar)_{h_j}
 \Bigr)
 \right]
\subset
 \left[
  P_{h_1}Gr_{h_1}^{(1)}
\cap
\bigcap_{j=2}^n W^{(1)}(\jbar)_{h_j}\right].
\]
Thus we are done.
\hfill\qed

\begin{df}
A tuple $(N_1,\ldots,N_n)$ is called strongly sequentially compatible
if the following holds:
\begin{itemize}
\item
  $(N_1,\ldots,N_n)$ is sequentially compatible.
\item
 For any tuple $\vech=(h_1,\ldots,h_n)$,
 we have the following:
 \[
  \Image(P\pi_{\vech})=
 \left[
 P_{h_1}Gr_{h_1}^{(1)}\cap \bigcap_{j=2}^n W^{(1)}(\jbar)_{h_j}
 \right].
 \]
\end{itemize}
\end{df}

\begin{df}
A commuting tuple $(N_1,\ldots,N_n)$
of nilpotent maps are called of Hodge type,
if the following holds:
\begin{itemize}
\item
 For any permutation $\sigma$ of $\nbar$,
 the tuple $(N_{\sigma(1)},\ldots,N_{\sigma(n)})$ 
 is strongly sequentially compatible.
\hfill\qed
\end{itemize}
\end{df}


\subsubsection{Sequential compatibility in the level $h$}

\begin{df}\label{df;12.2.30}
Let $(N_1,N_2,\ldots,N_n)$ be a commuting tuple of nilpotent maps.
It is called sequentially compatible in the level $h$,
if the following holds:
\begin{enumerate}
\item \label{number;12.2.31}
The constantness of the filtrations on the positive cones holds.
\item  \label{number;12.2.32}
 The induced tuple
 $(N_2^{(1)},\ldots,N_n^{(1)})$
 on $Gr^{W(N_1)}$
 is sequentially compatible.
\item \label{number;12.2.33}
 For any $\vech=(h_1,\ldots,h_n)$ such that $h_1\leq h$,
 we have
 $\Image(\pi_{\vech})
=Gr^{(1)}_{h_1}\cap\bigcap_{j=2}^n W(N^{(1)}(\jbar))_{h_j-h_1}$.
\end{enumerate}
In particular, when $h$ is the bottom number 
of the filtration $W(\itibar)$,
we say that $(N_1,\ldots,N_n)$ are sequentially compatible
in the bottom part.
\hfill\qed
\end{df}

\begin{lem}
Let $(N_1,\ldots,N_n)$ be a commuting tuple of nilpotent maps
satisfying the conditions {\rm\ref{number;12.2.31}} and
{\rm\ref{number;12.2.32}}
in Definition {\rm\ref{df;12.2.30}}.
When we check whether $(N_1,\ldots,N_n)$
is sequentially compatible in the bottom part,
we only have to check the following instead of the condition
{\rm\ref{number;12.2.33}} in Definition {\rm\ref{df;12.2.30}}.
\begin{description}
\item [4]
 For any $h$,
 we have
 $W(\itibar)_b\cap W(\jbar)_{h}=Gr_b^{(1)}\cap W(N^{(1)}(\jbar))_{h-b}$,
where $b$ denotes the bottom number of $W(\itibar)$.
\end{description}
\end{lem}
\pf
It is clear that the condition \ref{number;12.2.33}
implies the condition 4.
In the bottom part,
we have the equality:
\[
 W(\itibar)_b\cap \bigcap_{j=2}^{n}W(\jbar)_{h_j}
=\bigcap_{j=2}^n
 \bigl(
  W(\itibar)_b
 \cap W(\jbar)_{h_j}
 \bigr).
\]
Thus the condition 4 implies
the condition \ref{number;12.2.33}.
\hfill\qed

\begin{df}
Let $(N_1,\ldots,N_n)$ be a commuting tuple of nilpotent maps.
We say that it is universally sequentially compatible
in the bottom part,
if $(N_1^{\wedge\,m},\ldots,N_n^{\wedge\,m})$ is sequentially compatible
in the bottom part
for any non-negative integer $m$.
\hfill\qed
\end{df}

\subsubsection{Splitting of the sequentially compatible nilpotent maps}

\begin{lem}
Let $(N_1,\ldots,N_n)$ be a commuting tuple,
which is sequentially compatible in the level $h$.
Then there exists the decomposition of $W_h(\itibar)$:
\[
 W_h(\itibar)=
 \bigoplus_{\veck\in\seisuu^n} U_{\veck}.
\]
The decomposition satisfies the following
for any $\vech=(h_1,\ldots,h_n)$ such that $h_1\leq h$:
\begin{equation} \label{eq;12.2.35}
 \bigcap_{j=1}^n W(\jbar)_{h_j}=
 \bigoplus_{
 \veck\in \tiisai(\vech)
 } U_{\veck},
\quad
 \tiisai(\vech):=
 \{\veck\in\seisuu^n\,|\,\rho_j(\veck)\leq h_j\,\,j=1,\ldots,n\}.
\end{equation}
Here we put $\rho_{j}(\veck):=\sum_{i\leq j}k_i$
for $\veck=(k_1,\ldots,k_n)$.
\end{lem}
\pf
We have the filtrations $W(\itibar)_h\cap W(\jbar)$
of $W(\itibar)_h$ for any $j$.
The sequence
$\bigl(W(\itibar)_h\cap W(\itibar),W(\itibar)_h\cap W(\nibar),
\ldots, W(\itibar)_h\cap W(\nbar)\bigr)$
is compatible.
Thus we obtain the compatible decomposition:
\[
 W(\itibar)_h=
 \bigoplus_{\veck\in\seisuu^n}
 \overline{U}_{\veck}.
\]
For any $\vech\in \seisuu^n$,
we put
$\mu(\vech):=(\rho_1(\vech),\rho_2(\vech),\ldots,\rho_n(\vech))$.
Then we put as follows:
\[
 U_{\vech}=
 \overline{U}_{\mu(\vech)}.
\]
Then the decomposition $V=\bigoplus_{\vech\in\seisuu}U_{\vech}$
has the desired property.
\hfill\qed

\vspace{.1in}

Let $(N_1,\ldots,N_n)$ be a sequentially compatible
in the level $h$.
We put as follows:
\begin{equation} \label{eq;12.9.51}
 R=\dim W(\itibar)_h,
\quad
b=\sum_{a\leq h} a\cdot \dim Gr_a^{(1)}.
\end{equation}
We put $\nbigv:=\bigwedge^R V$
and $\nbign_i:=N_i^{\bigwedge\,R}$,
and $\nbign(\jbar)=\sum_{i\leq j}\nbign_i$.

\begin{lem}
Let consider the weight filtration
$\nbigw(\jbar):=\nbigw(\nbign(\jbar))$.
\begin{itemize}
\item
The bottom number of the filtration $\nbigw(\itibar)$
is $b$.
\item
We have the natural isomorphism
$\nbigw(\itibar)_b\simeq \det(W(\itibar)_h)$.
\item
Let $e$ be a non-zero element of $\nbigw(\itibar)_b$.
Then we have $\deg^{\nbigw(\mbar)}(e)=b$
for any $m$.
\end{itemize}
\end{lem}
\pf
The first two claims are clear, and we do not need
the sequentially compatibility in the level $h$.
(See the subsubsection \ref{subsubsection;12.9.50}.)
The sequence $(W(\itibar),W(\mbar))$ is compatible.
Let take a frame $\vecv$ which is 
compatible with the sequence $(W(\itibar),W(\mbar))$.
Then we have the following equality:
\[
 \deg^{W(\mbar)}(v_i)
=\deg^{W^{(1)}(\mbar)}(v_i)
=\deg^{W(\itibar)}(v_i)+\deg^{W(N^{(1)}(\mbar))}(v_i^{(1)}).
\]
We also have the following:
\[
 \sum_{\deg^{W(\itibar)}(v_i)=a}\deg^{W(N^{(1)}(\mbar))}(v_i^{(1)})=0.
\]
Then we obtain the equality
$\sum_i\deg^{W(\mbar)}(v_i)
=\sum_{a\leq h} a\cdot \dim Gr^{(1)}_a=b$
\hfill\qed

\subsubsection{Splitting of strongly sequentially compatible nilpotent maps}

Let $q_1:\seisuu^n\lrarr \seisuu$ denote the projection
onto the first component.
\begin{lem} \label{lem;12.3.10}
Let $(N_1,\ldots,N_n)$ be commuting tuple of nilpotent maps,
which is strongly sequentially compatible.
Then we have the decomposition:
\[
 V=\bigoplus_{k\geq 0}\bigoplus_{\veck\in\seisuu^n}P_{k}U_{\veck}.
\]
It satisfies the following:
\begin{enumerate}
\item
It gives a splitting of the filtrations $W(\jbar)$,
that is we have the following:
\[
 \bigcap_{j=1}^nW(\jbar)_{h_j}=
 \bigoplus_{k\geq \,0}
 \bigoplus _{\veck\in \nbigu(\vech)}
 P_kU_{\veck}.
\]
\item
 $P_kU_{\veck}=0$ unless $|q_1(\veck)|\leq k$ and $k-q_1(\veck)$ is even.
\item
 When $-k< q_1(\veck)\leq k$,
 we have
 $N_1\bigl(
  P_kU_{\veck}\bigr)
 =P_kU_{\veck-2\vecdelta_1}$.
 Here we put $\veck-2\vecdelta_1=(k_1-2,k_2,\ldots,k_n)$
 for $\veck=(k_1,\ldots,k_n)$.
\item
 When $q_1(\veck)=-k$,
 $N_1\bigl(  P_kU_{\veck}\bigr)=0$.
\end{enumerate}
\end{lem}
\pf
We have the sequentially compatible tuple
$(N_2^{(1)},\ldots,N_n^{(1)})$ on $P_{h}Gr^{(1)}_{h}$.
Thus we have the decomposition of $P_{h}Gr^{(1)}_{h}$:
\[
 P_{h}Gr^{(1)}_{h}=
 \bigoplus_{\veck\in\seisuu^{n-1}}
 P_hU_{h,\veck},
\quad
 P_{h}Gr^{(1)}_{h}\cap
 \bigcap_{j=1}^{n-1}
 W(N^{(1)}(\jitibar))_{h'_j}=
 \bigoplus_{\veck\in\nbigu(\vech')}
 P_hU_{h,\veck}.
\]
For a tuple $\vech=(h_1,\ldots,h_n)$,
we put $\chi(\vech):=(h_1,h_2-h_1,h_3-h_2,\ldots,h_n-h_{n-1})$,
and $\chi'(\vech)=(h_2-h_1,h_3-h_2,\ldots,h_n-h_{n-1})$.
Note that $\chi(\vech)\in\nbigu(\vech)$.
When we put $\vech'=(h_2-h_1,\ldots,h_{n}-h_1)\in\seisuu^{n-1}$,
we have $\chi'(\vech)\in\nbigu(\vech')$.

Since $(N_1,\ldots,N_n)$ is strongly sequentially compatible,
the space $P_{h_1}U_{h_1,\chi'(\vech)}$ is contained in the image
of $P\pi_{\vech}$:
\[
 P\pi_{\vech}\left[
 \Ker(N(1)^{h_1+1})\cap
 \bigcap_{j=2}^n W(\jbar)_{h_j}\right]
=\left[
 P_{h_1}Gr^{(1)}_{h_1}\cap
 \bigcap_{j=2}^n W(N^{(1)}(\jbar))_{h_j-h_1}
 \right]
 \supset P_{h_1}U_{h_1,\chi'(\vech)}.
\]
Thus we can take a subspace $P_{h_1}U_{\chi(\vech)}$
of $\Ker(N(1)^{h_1+1})\cap \bigcap_{j=2}^nW(\jbar)_{h_j}$
such that $P_{h_1}U_{\chi(\vech)}$ is isomorphic to
$P_{h_1}U_{h_1,\chi'(\vech)}$
via the morphism $P\pi_{\vech}$.

For an integer $m$ such that $0\leq m\leq h_1$,
we put as follows:
\[
 P_{h_1}U_{\vech-2m\vecdelta_1}:=
 N_1^m\Bigl(
 P_{h_1}U_{\vech}
 \Bigr).
\]
Here we put $\vech-2m\cdot\vecdelta_1=(h_1-2m,h_2,\ldots,h_n)$
for $\vech=(h_1,h_2,\ldots,h_n)$.
By our choice,
we have $N_1^{h_1+1}\bigl(P_{h_1}U_{\vech} \bigr)=0$.
Then we obtain the desired decomposition.
\hfill\qed

For any tuple $\vech\in\seisuu^n$ and $k\geq 0$,
we have the number $d(k,\vech):=\dim P_kU_{\vech}$.
Clearly we have $d(k,\vech)=d(k-2,\vech)$
if $-k<q_1(\vech)\leq k$.

\begin{cor}\label{cor;12.10.201}
Assume that $(N_1,\ldots,N_n)$ is strongly sequentially compatible.
Then we can take a base $\vecv$ of $V$
satisfying the following:
\begin{enumerate}
\item
$\vecv=
 \Bigl(v_{k,\vech,\eta}\,
 \Big|\,k\geq 0,\vech\in\seisuu^n,
 \eta=1,\ldots,d(k,\vech) \Bigr)$.
\item
 We have
 $N_1(v_{k,\vech,\eta})=v_{k,\vech-2\vecdelta_1,\eta}$
 when $-k<q_1(\vech)\leq k$.
\item
 $N_1(v_{k,\vech,\eta})=0$
 if $q_1(\vech)=-k$.
\item
 $\deg^{W(\jbar)}(v_{k,\vech,\eta})=\rho_j(\vech)$.
\item
$\vecv$ is compatible with the sequence
$\bigl(W(\itibar),\ldots,W(\nbar)\bigr)$.
\end{enumerate}
\end{cor}
\pf
We take a base 
$\Bigl(v_{h_1,\vech,\eta}\,\Big|\,\eta=1,\ldots,d(h_1,\vech)\Bigr)$
of $P_{h_1}U_{\vech}$
in the case $\vech=(h_1,\ldots,h_n)$.
We put
$v_{h_1,\vech-2m\vecdelta_1,\eta}:=N_1^m(v_{h_1,\vech,\eta})$.
Then we obtain the frame desired.
\hfill\qed

\begin{df}
A frame $\vecv$ satisfying the condition in
Corollary {\rm\ref{cor;12.10.201}} is called
strongly compatible with
$\bigl(N(\itibar),N(\nibar),\ldots,N(\nbar)\bigr)$.
\end{df}

\subsubsection{A reduction of the sequential compatibility}

\begin{prop} \label{prop;12.10.80}
Let $(N_1,\ldots,N_n)$ be a commuting tuple.
Assume the following:
\begin{itemize}
\item
$(N_1,N_2,\ldots,N_n)$
is universally sequentially compatible in the bottom parts.
\item
$(N_1+N_2,N_3,\ldots,N_n)$
is sequentially compatible.
\end{itemize}
Then $(N_1,N_2,\ldots,N_n)$ is sequentially compatible.
\end{prop}
\pf
We only have to show that $(N_1,\ldots,N_n)$
is sequentially compatible in the level $h$
for any $h$.
We use an induction on $h$.
We assume that we have already known that
$(N_1,\ldots,N_n)$ is sequentially compatible
in the level $h-1$.

We put as follows:
\[
 R_1:=\dim W(\itibar)_{h-1}+1,
\quad
 b_1:=\sum_{a\leq h-1} a\cdot \dim Gr_a^{(1)}+h.
\]
We put $\nbigv_1:=\bigwedge^{R_1} V$,
$\nbign_i:=N_i^{\wedge\,R_1}$
and $\nbign(\jbar)=\sum_{i\leq j}\nbign_i$.
We obtain the weight filtration $\nbigw(\jbar)$ of $\nbign(\jbar)$.
We denote the associated graded vector space
of $\nbigw(\itibar)$ by $\graded^{(1)}$.
Then the bottom number of the filtration $\nbigw(\itibar)$
is $b_1$ above.
The tuple $(\nbign_1,\ldots,\nbign_n)$ is sequentially compatible
in the bottom part.
Thus we have the following equality for any $(h_2,\ldots,h_n)$:
\[
 \nbigw(\itibar)_{b_1}\cap
 \bigcap_{j=2}^n
 \nbigw(\jbar)_{h_j}
=\graded^{(1)}_{b_1}
\cap \bigcap_{j=2}^n
 \nbigw(\nbign^{(1)}(\jbar))_{h_j-b_1}.
\]

Due to our choice of $R_1$,
we have the natural isomorphism:
$\graded_{b_1}^{(1)}
=\nbigw(\itibar)_{b_1}\simeq \det(W(\itibar)_{h-1})\otimes Gr_h^{(1)}$.
Under the isomorphism,
the morphisms $\nbign^{(1)}_i$ and $N^{(1)}_i$ also correspond.
Thus we obtain the following equality
for any $h$ and $(h_2,\ldots,h_n)$
under the isomorphism:
\[
 \graded_{b_1}^{(1)}
\cap
 \bigcap_{j=2}^n
 \nbigw(\nbign^{(1)}(\jbar))_{h_j-h}
=\det(W(\itibar)_{h-1})\otimes 
 \left[
 Gr_h^{(1)}
 \cap
 \bigcap_{j=2}^n
 W(N^{(1)}(\jbar))_{h_j-h}
 \right].
\]
Thus we only have to check the following coincidence.
(Here we put $h_1=h$ for simplicity of notation):
\[
 \bigcap_{j=1}^n
 \nbigw(\jbar)_{h_j-h+b_1}
=\det(W(\itibar)_{h-1})\otimes 
\left[
 \Image\Bigl(
 \pi_{\vech}:
 \bigcap_{j=1}^n
 W(\jbar)_{h_j}
\lrarr Gr_{h}^{(1)}
 \Bigr)
 \right].
\]

First we see that
$\det(W(\itibar)_{h-1})\otimes\Image(\pi_{\vech})$ is contained
in $\bigcap_{j=1}^n\nbigw(\jbar)_{h_j-h_1+b}$.
Let take a non-zero element $e$ of $\det(W(\itibar)_{h-1})$.
We know that the degree of $e$
in $\bigwedge^{R}V$ with respect to
the filtration $W(N(\jbar)^{\wedge\,R})$ is 
$b=b_1-h_1$ for any $j$.
Here $R$ and $b$ are given as follows:
(See (\ref{eq;12.9.51}). We use $h-1$ instead of $h$.)
\[
 R:=\dim W(\itibar)_{h-1},
\quad
 b:=\sum_{a\leq h-1} a\cdot \dim Gr_a^{(1)}=b_1-h_1.
\]
Thus we obtain the following inequality:
\[
 \deg^{\nbigw(\jbar)}(e\wedge y)
\leq
 b_1-h_1+h_j.
\]
It implies that $e\wedge y\in
\bigcap_{j=1}^n\nbigw(\jbar)_{b_1-h_1+h_j}$.

Let see the implication
$\bigcap_{j=1}^n\nbigw(\jbar)_{b_1-h_1+h_j}
\subset
 \det(W(\itibar)_{h-1})\otimes \Image(\pi_{\vech})$.
Any element of $\bigcap_{j=1}^n\nbigw(\jbar)_{b_1-h_1+h_j} $
is described as follows:
\[
 e\wedge y,\quad
 y\in W(\itibar)_{h_1}.
\]
Let consider the splitting
$V=\bigoplus_{\veck\in\seisuu^{n-1}}U_{\veck}$
compatible with the sequentially compatible tuple
$(N_1+N_2,N_3,\ldots,N_1)$.
Then we have the following decomposition:
\[
 y=\sum_{\veck\in\seisuu^{n-1}} y_{\veck},
\quad
 y_{\veck}\in U_{\veck}.
\]
Due to the condition $\deg^{\nbigw(\nbar)}(e\wedge y)\leq h_n-h_1+b$,
we have the following vanishing for any $l>h_{n}$:
\[
 e\wedge \sum_{
 \substack{\veck\in\seisuu^{n-1},\\ \rho_n(\veck)=l}
 } y_{\veck}=0.
\]
It implies that
$\sum_{\rho_n(\veck)=l}y_{\veck}\in W(\itibar)_{h-1}$.
We put as follows:
\[
 y'=y-\sum_{l>h_n}\sum_{
 \substack{\veck\in\seisuu^{n-1}\\
 \rho_n(\veck)=l}}y_{\veck}.
\]
Then we know the following:
\[
 e\wedge y'=e\wedge y,
\quad
  y'\in W(\itibar)_{h}.
\]
Thus we can assume that $y_{\veck}=0$
if $\rho_h(\veck)>h_n$, from the beginning.

By an inductive argument, we can assume that
$y_{\veck}=0$ if there exist $2\leq j\leq n$
such that $\rho_j(\veck)>h_j$.
In that case,
$y$ is contained in $\bigcap_{j=1}^nW(\jbar)_{h_j}$.
It implies the implication desired.
Thus we are done.
\hfill\qed

\subsubsection{A lemma of Cattani-Kaplan}

Let $N_1,\ldots,N_n$ be a commuting tuple of nilpotent maps on $V$.
Let $t_1\ldots,t_n$ be formal variables and we put
$N(\vect):=\sum t_i\cdot N_i$.
Let $K$ denote the rational function field
$\cnum(t_1,\ldots,t_n)$ with variables $t_1,\ldots,t_n$.
Then $N(\vect)$ gives a nilpotent map over $V\otimes_{\cnum} K$.
The weight filtration induced by $N(\vect)$
is denoted by $W(\vect)$.

Let $\veca=(a_1,\ldots,a_n)$ be an element of $\cnum^n$.
Then we have $N(\veca)=\sum a_i\cdot N_i$.
We denote the weight filtration of $N(\veca)$
by $W(\veca)$.

\begin{df} \label{df;12.3.40}
When we have $\dim_{\cnum}(W(\veca)_l)=\dim_{K}(W(\vect)_l)$
for any $l$,
we say that $\veca$ is general,
or that $N(\veca)$ is general.
\hfill\qed
\end{df}

Since $N_i$ are commuting,
we always have $N_i\cdot W(\veca)_l\subset W_l(\veca)$.
\begin{lem}[Cattani-Kaplan, \cite{ck}] \label{lem;12.3.41}
When $N(\veca)$ is general,
we have $N_i\cdot W(\veca)_l\subset W(\veca)_{l-1}$.
\hfill\qed
\end{lem}

\subsubsection{A lemma for the conjugacy classes of the nilpotent maps}

We recall some general result on the conjugacy classes
of nilpotent maps.
We will use the result later without mention.
Let $R$ be a discrete valuation ring.
Let $K$ and $k$ denote the quotient field and the residue field
of $R$ respectively.
Let $V$ be a free module over a discrete valuation ring $R$.
Let $N$ be a nilpotent maps of $V$
defined over $R$.
We put $V_K:=V\otimes_R K$ and $V_k:=V\otimes_Rk$.
We have the induced nilpotent maps
$N_K\in End(V_K)$ and $N_k\in End(V_k)$.
They induce the weight filtrations
$W_K$ and $W_k$ of $V_K$ and $V_k$ respectively.

\begin{lem}\label{lem;2.6.10}
We put $l_0:=\min\{l\,|\,\dim W_{K,l}\neq \dim W_{k,l}\}$.
Then we have the following inequality:
\[
 \dim_K W_{K,l_0}>\dim_k W_{k,l_0}.
\]
\end{lem}
\pf
First observe the following:
Let $b(K)$ and $b(k)$ be the bottom numbers
of the filtrations $W_K$ and $W_k$.
If $N_K^l=0$, then $N^l=0$
and thus $N_k^l=0$.
It implies that $b(K)\leq b(k)$.
If $b(K)=b(k)=b$,
then we have the following inequality:
\[
  \dim W_{K,b(K)}=\dim\Image(N_K^{b})\geq \dim\Image(N_k^b)=\dim W_{k,b(k)}.
\]

We put $D=\sum_{l<l_0}\dim W_{K,l}$.
By considering the exterior product $\bigwedge^{D+1}V$,
we can reduce the problem to the comparison
of the dimension of the bottom parts.
Thus we are done.
\hfill\qed

\subsection{Vector bundles and filtrations}

\subsubsection{Words and Notation}
Let $X$ be a complex manifold
and $E$ be a $C^{\infty}$-vector bundle over $\cnum$.
The space of $C^{\infty}$-sections
are denoted by $C^{\infty}(X,E)$.
When $E$ is a holomorphic bundle,
the space of holomorphic sections
is denoted by $\Gamma(X,E)$.
For frame $\vecv$ and $\vecw$,
we have the transformation matrices $B$
determined by $\vecv=\vecw\cdot B$.

Let $h$ be a hermitian metric of $E$,
and $\vecv$ be a frame of $E$ of rank $r$.
Then we obtain the $\nbigh(r)$-valued function
$H(h,\vecv)$ and $H(h,\vecv)^{-1}$.

\begin{df}
A frame $\vecv$ is called adapted,
if $H(h,\vecv)$ and $H(h,\vecv)^{-1}$ are bounded
over $X$.
\hfill\qed
\end{df}

Let $Y$ be a subset of $X$.
The restriction of $E$ to $Y$ is denoted by $E_{|Y}$.

Assume that a decomposition of $E_{|Y}$ 
into a direct sum of vector bundles $\bigoplus E_i$ is given.
The restriction of a $C^{\infty}$-section $f$ of $E$
to $Y$ is denoted by $f_{|Y}$.
It is called compatible with the decomposition,
if there is an $i$ such that 
$f_{|Y}$ is a section of $E_{i}$.
A frame $\vecv=(v_1,\ldots,v_n)$ of $E$ is compatible
with the decomposition,
if each $v_i$ is compatible with the decomposition.

\subsubsection{Filtrations}

A filtration $W$ of $E_{|Y}$ by vector bundle
is defined to be a finite increasing sequence
of vector subbundles:
\[
W_a\subset W_{a+1}\subset \cdots \subset W_{a+h}\subset E_{|Y}.
\]
The associated graded vector bundle on $Y$
is denoted by $Gr^W(E_{|Y})$.

If $E$ is a holomorphic vector bundle,
then a filtration of $E$ by subsheaves
is defined to be a similar finite increasing sequence
of subsheaves.

When a decomposition of $E_{|Y}$ is given,
a filtration $W$ of $E_{|Y}$ is called compatible with the decomposition
if $W_l=\bigoplus_i E_i\cap W_l$.

\begin{df}
Let $W$ be a filtration of $E_{|Y}$.
A $C^{\infty}$-section $f$ of $E$
is called compatible with the filtration $W$,
if the numbers $\deg^{W_{|P}}(f(P))$ are independent of $P\in Y$.
In that case, we put $\deg^W(f):=\deg^{W_{|P}}(f(P))$
for some $P\in Y$.
\end{df}

\begin{df}
Let $\vecv=(v_1,\ldots,v_n)$ be a $C^{\infty}$-frame
of $E$.
It is called compatible with the filtration $W$ of $E_{|Y}$,
if the following conditions are satisfied:
\begin{enumerate}
\item
Each $v_i$ is compatible with the filtration $W$.
\item
For any point $P\in Y$,
the frame $\vecv_{|P}$ is compatible
with the filtration $W_{|P}$.
\end{enumerate}
\end{df}

The induced sections $\vecv^{(1)}=(v^{(1)}_1,\ldots,v^{(1)}_n)$
on $Y$
gives a frame of $Gr^{W}(E_{|Y})$
compatible with the decomposition.

We put $X=\Delta^n$, $D_i=\{z_i=0\}$ and $D=\bigcup_{i=1}^l D_i$
for $l\leq n$.
We put $D_I=\bigcup_{i\in I}D_i$ for $I\subset \lbar$.
Let $\sigma$ be an element of the $l$-th symmetric group $\gbigs_l$.
Then we obtain the sequence of the subset:
\[
 I_1\subset I_2\subset\cdots\subset I_l,
\quad
I_j=\bigl\{\sigma(i)\,\big|\,i\leq j\bigr\}.
\]
Then we obtain the sequence
$D_{I_1}\supset D_{I_2}\supset\cdots\supset D_{I_l}$.

\begin{df}\label{df;12.10.206}
Let $E$ be a vector bundle over $X$.
A sequence $\bigl(W(I_1),W(I_2),\ldots,W(I_l)\bigr)$
of the filtration of $E$ is the following data:
\begin{enumerate}
\item
 For any point $Q\in D_{I_m}$,
 the sequence $\bigl(W(I_1)_{|Q},W(I_2)_{|Q},\ldots, W(I_m)_{|Q}\bigr)$
 is a compatible sequence of filtrations.
\item
 $W(I_j)$ is a filtration of $E_{|D_{I_j}}$ 
 by vector subbundles.
\item \label{number;12.10.205}
 Let $\vech=(h_1,\ldots,h_m)$ denote a tuple of integers.
 We have a vector spaces
 $\bigcap_{j=1}^mW(I_j)_{h_j\,|\,Q}$
 for any $Q\in D_{I_m}$.
 Then they form a vector subbundle of
 $E_{|D_{I_m}}$.
 Namely we have a vector subbundle
 $\bigcap_{j=1}^mW(I_j)_{h_j\,|\,D_{I_m}}$
 of $E_{|D_{I_m}}$.
\hfill\qed
\end{enumerate}
\end{df}

Let $Q\in D_{I_m}$.
Then for any tuple $\vech\in\seisuu^m$,
we have the number $d(\vech)=\dim U_{\vech}$
for a splitting $E_{Q}=\bigoplus_{\veck\in\seisuu}U_{\veck}$
compatible with the filtrations
$\bigl(W(I_1)_{|D_{I_m}},\ldots,W(I_m)_{|D_{I_m}}\bigr)$.

\begin{lem}
For a tuple $\vech\in \seisuu^m$,
the number $d(\veck)$ is independent of 
a choice of $Q\in D_{I_m}$.
\end{lem}
\pf
As is easily seen,
the number $d(\veck)$ is determined by
the following numbers:
\[
 \dim \left[
 \bigcap_{j=1}^m W(I_j)_{h_j\,|\,Q}
 \right],
\quad
 (h_1,\ldots,h_m)\in\seisuu^m.
\]
By our assumption, the numbers are independent of a choice of $Q$.
\hfill\qed

On $D(I_1)$,
we obtain the graded vector space $\graded^{(1)}$
associated with the filtration $W(I_1)$.
For any point $Q\in D_{\mbar}$,
we have the induced filtration
$W^{(1)}(\mbar)_{|Q}$ of $\graded^{(1)}_{|Q}$.
\begin{lem} \label{lem;12.13.2}
$\bigl\{W^{(1)}(I_m)_{|Q}\,\big|\,Q\in D_{I_m}\bigr\}$
gives a filtration of $E_{|D_{\mbar}}$
by vector subbundles.
\end{lem}
\pf
We only have to check that the dimension
of $W^{(1)}(I_m)_{|Q}$ is independent of a choice of $Q\in D_{I_m}$.
It follows from the fact that
the numbers $d(\vech)$ are independent of a choice of $Q\in D_{I_m}$.
\hfill\qed

For any $\vecl\in\seisuu^{m-1}$,
we have the vector subspaces
$\bigl\{\bigcap_{j=1}^{m-1}W^{(1)}(I_{j+1})_{l_j\,|\,Q}\,|\,Q\in D_{I_m}
\bigr\}$.
\begin{lem} \label{lem;12.13.3}
For each $\vecl\in\seisuu^{m-1}$,
$\bigl\{
 \bigcap_{j=1}^{m-1}W^{(1)}(I_{j+1})_{l_j\,|\,Q}\,|\,Q\in D_{I_m}
\bigr\}$ forms a vector subbundle over $D_{I_m}$.
\end{lem}
\pf
Again we only have to see the independence of the dimension
of the vector spaces $\bigcap_{j=1}^{m-1}W^{(1)}(I_{j+1})_{l_j\,|\,Q}$.
It follows from the fact that
the numbers $d(\vech)$ are independent of a choice of $Q\in D_{I_m}$.
\hfill\qed

We have the vector bundle
$\graded^{(1)}$ on $D_{I_1}$,
and the filtrations
$W^{(1)}(I_j)$ on $D_{I_j}$ for $j\geq 1$.
We have already seen the following proposition
in the lemmas \ref{lem;12.13.2} and \ref{lem;12.13.3}.
\begin{prop}
$\bigl(W^{(1)}(I_2),\ldots,W^{(1)}(I_l) \bigr)$
is a compatible sequence of filtrations.
\hfill\qed
\end{prop}

\begin{df}
Let $f$ be a section of $E$ over $X$.
We say that $f$ is compatible with the sequence of the filtrations
$\bigl(W(I_1),W(I_2),\ldots,W(I_l)\bigr)$, if the following is satisfied:
\begin{itemize}
\item
$f$ is compatible with $W(I_j)$ for any $j$.
\item
Let $P$ be a point of $D_{I_j}$.
Then $f_{|P}$ is compatible with the filtrations
$\bigl(W(I_1)_{|P},\ldots, W(I_j)_{|P}\bigr)$.
\hfill\qed
\end{itemize}
\end{df}

\begin{df}
Let $\vecv$ be a frame of $E$.
We say $\vecv$ is compatible with the sequence
$\bigl(W(I_1),\ldots,W(I_l)\bigr)$ if the following holds:
\begin{itemize}
\item
 $\vecv$ is compatible with the filtration $W(I_j)$
 for any $j$.
\item
 Let $P$ be a point of $I_j$.
 Then $\vecv_{|P}$ is compatible with the sequence
 $\bigl(W(I_1)_{|P},\ldots,W(I_j)_{|P}\bigr)$.
\hfill\qed
\end{itemize}
\end{df}

\subsubsection{The existence of compatible splitting}

Let $\bigl(W(I_1),\ldots, W(I_l)\bigr)$ be a compatible sequence
of filtrations.
For simplicity of notation,
we assume that $I_j=\jbar=\{1,\ldots,j\}$.

\begin{lem}
For any $1\leq m\leq l$,
there are decompositions of $E_{|D_{\mbar}}$:
\[
 E_{|D_{\mbar}}=
 \bigoplus_{\vech\in\seisuu^m}
 \nbigk_{\vech}.
\]
They satisfy the following:
\begin{enumerate}
\item
 For any $\vech\in\seisuu^m$,
 we have
 $\bigcap_{j=1}^m W(\jbar)_{h_j}
=\bigoplus_{\veck\in\nbigt(\vech)}\nbigk_{\veck}$
on $D_{\mbar}$.
Here $\nbigt(\vech)$ denotes the set of $\veck\in\seisuu^m$
satisfying $q_j(\veck)\leq h_j$ for any $1\leq j\leq m$.
\item
 We have
 $\nbigk_{\vech\,|\,D_{\mitibar}}=
 \bigoplus _{k}\nbigk_{(\vech,k)}$.
Here $(\vech,k)=(h_1,\ldots,h_m,k)$
for $\vech=(h_1,\ldots,h_m)$.
\end{enumerate}
\end{lem}
\pf
We use an induction on $l$.
We have the vector bundle
$\graded^{(1)}_h$
on $D_{\itibar}$,
and the filtrations $\graded^{(1)}_h\cap W^{(1)}(\jbar)$
on $D_{\jbar}$.
Since the sequence of the filtration is compatible,
we have the compatible splitting:
\[
 \graded^{(1)}_{h\,|\,D_{\mbar}}
=\bigoplus_{\veck'\in\seisuu^{m-1}}
 \nbigk_{h,\veck'}.
\]

We construct $\nbigk_{\vech}$ on $D_{\mbar}$
by using an descending induction on $m$.
Assume that we have already constructed $\nbigk_{\vech}$
on $D_{\mitibar}$,
and consider the decomposition on $D_{\mbar}$.

For a tuple $\vech=(h_1,\ldots,h_m)\in\seisuu^m$,
we put $\vech'=(h_2,\ldots,h_m)$.
Then $\nbigk_{h_1,\vech'}$ is contained in
the image of the following morphism on $D_{\mbar}$, by our assumption:
\[
 \pi_{\vech}:
 \bigcap_{j=1}^mW(\jbar)_{h_j}
\lrarr
 \graded^{(1)}_{h_1}\cap
 \bigcap_{j=2}^mW^{(1)}(\jbar)_{h_j}
\supset \nbigk_{h_1,\vech'}.
\]
On $D_{\mitibar}$,
we already have $\bigoplus_k\nbigk_{(\vech,k)}$.
By extending it,
we can take a subbundle $\nbigk_{\vech}$
of $\bigcap_{j=1}^mW(\jbar)_{h_j}$
on $D_{\mbar}$,
satisfying the following:
\begin{itemize}
\item
$\nbigk_{\vech}$ is isomorphic
to $\nbigk_{h_1,\vech'}$ via the morphism $\pi_{\vech}$.
\item
We have
$\nbigk_{\vech|D_{\mitibar}}= \bigoplus_k\nbigk_{(\vech,k)}$.
\end{itemize}
Thus the induction can proceed.
\hfill\qed

\begin{df}
Such tuple
$\{\nbigk_{\vech}\,|\,\vech\in\seisuu^{m},\,\,m=1,\ldots,l\}$
is called a compatible splitting of the sequence
$(W(\itibar),\ldots,W(\lbar))$.
\hfill\qed
\end{df}

\begin{lem}
Let $(W(\itibar),\ldots,W(\lbar))$ be a compatible sequence
of filtrations of $E$.
There exists a frame $\vecv$ of $E$
compatible with the sequence $(W(\itibar),\ldots,W(\lbar))$.
\end{lem}
\pf
We take a compatible splitting
$\{\nbigk_{\vech}\,|\,\vech\in\seisuu^{m},\,\,m=1,\ldots,l\}$
of the sequence $(W(\itibar),\ldots,W(\lbar))$.
We can take a frame $\vecv$
compatible with the splitting 
$\{\nbigk_{\vech}\,|\,\vech\in\seisuu^{m},\,\,m=1,\ldots,l\}$.
Thus we are done.
\hfill\qed
\subsection{Commuting tuple for a vector bundle}
\label{subsection;12.15.3}

\subsubsection{Constantness of the filtrations on the positive cones}

Let $E$ be a holomorphic vector bundle over $\Delta^n$.
Let $l$ be a natural number such that $l\leq n$.
We put $D_i:=\{(z_1,\ldots,z_n)\in\Delta^n\,|\,z_i=0\}$,
and $D_{\mbar}=\bigcap_{i=1}^mD_i$.
Let $N_i$ be an element of
$\Gamma(D_i,End(E)_{|D_i})$
for $i=1,\ldots,l$.
For $m\leq l$,
we have the nilpotent maps
$N_1,\ldots,N_m$ of $End(E_{|D_{\mbar}})$
on $D_{\mbar}$.
Then we put as follows, for any $\veca\in\real_{\geq 0}^m$:
\[
 N(\veca):=
 \sum_{j=1}^m a_j\cdot N_{j\,|\,D_{\mbar}}.
\]
\begin{df}
We say that the constantness of the the filtrations
on the positive cones for $(N_1,\ldots,N_l)$ holds,
if the following holds:
\begin{itemize}
\item
For any $m\leq l$,
for any $Q\in D_{\mbar}$
and for any $I\subset\mbar$,
the filtration $W(\veca)_{|Q}$ is independent
of a choice of $\veca\in\real_{>0}^I$.
\item
$\bigl\{W(\veca)_{|Q}\,\big|\,Q\in D_{\mbar}\bigr\}$
forms the vector bundle on $D_{\mbar}$.
\hfill\qed
\end{itemize}
\end{df}

\subsubsection{Sequential compatibility}

Let $E$ be a holomorphic vector bundle over $\Delta^n$.
We put $D_i:=\{(z_1,\ldots,z_n)\in\Delta^n\,|\,z_i=0\}$,
and $D_{\mbar}=\bigcap_{i=1}^mD_i$.
Let $l$ be a number less than $n$.
Let $N_i$ be an element of
$\Gamma(D_i,End(E)_{|D_i})$
for $i=1,\ldots,l$.
On $D_{\mbar}$,
we have the nilpotent maps
$N_{1\,|\,D_{\mbar}},\ldots,N_{m\,|\,D_{\mbar}}$
of $End(E_{|D_{\mbar}})$.

\begin{df} \label{df;12.10.208}
A commuting tuple $(N_1,\ldots,N_l)$ is called
sequentially compatible,
if the following holds:
\begin{enumerate}
\item
Let $P$ be a point of $D_{\jbar}$.
Then $(N_{1\,|\,P},\ldots,N_{j\,|\,P})$ is sequentially compatible.
\item
We put $N(\jbar)=\sum_{i\leq j}N_{i\,|\,D_{\jbar}}$.
Then the conjugacy classes of $N(\jbar)_{|Q}$
are independent of $Q\in D_{\jbar}$.
\item \label{number;12.10.207}
Let $W(\jbar)$ denote the weight filtration of $N(\jbar)$.
Then $(W(\itibar),\ldots,W(\lbar))$ is a compatible sequence
of filtrations.
\hfill\qed
\end{enumerate}
\end{df}

\begin{rem}
When we check whether a commuting tuple $(N_1,\ldots,N_n)$
is sequentially compatible,
we only have to check
the condition {\rm\ref{number;12.10.205}}
in Definition {\rm\ref{df;12.10.206}}
instead of the condition {\rm\ref{number;12.10.207}}
in Definition {\rm\ref{df;12.10.208}}.
\hfill\qed
\end{rem}

\begin{df}
Let $(N_1,\ldots,N_l)$ be a sequentially compatible commuting tuple.
A frame $\vecv$ is called compatible with $(N_1,\ldots,N_l)$,
if $\vecv$ is compatible with
the sequence $\bigl(W(\itibar),\ldots,W(\lbar)\bigr)$.
\hfill\qed
\end{df}

\begin{lem}
There are decompositions of $E_{|D_{\mbar}}$ for $1\leq m\leq l$:
\[
 E_{|D_{\mbar}}=
 \bigoplus_{\vech\in\seisuu^m}
 \nbigk_{\vech}.
\]
They satisfy the following:
\begin{enumerate}
\item
 For any $\vech\in\seisuu^m$,
 we have
 $\bigcap_{j=1}^m W(\jbar)_{h_j}
=\bigoplus_{\veck\in\nbigu(\vech)}\nbigk_{\veck}$
on $D_{\mbar}$.
Here $\nbigu(\vech)$ denotes the set of $\veck\in\seisuu^n$
satisfying $\rho_j(\veck)=\sum_{i\leq j}q_i(\veck)\leq h_j$
for any $1\leq j\leq m$.
\item
 We have
 $\nbigk_{\vech\,|\,D_{\mitibar}}=
 \bigoplus _{k}\nbigk_{(\vech,k)}$.
Here $(\vech,k)=(h_1,\ldots,h_m,k)$
for $\vech=(h_1,\ldots,h_m)$.
\end{enumerate}
\end{lem}
\pf
Since $(W(\itibar),\ldots,W(\lbar))$ is compatible,
we have the compatible splitting
$\{\overline{\nbigk}_{\vech}\,|\,
 \vech\in\seisuu^{m},\,\,m=1,\ldots,l\}$.
For any $\vech=(h_1,\ldots,h_m)\in\seisuu^m$,
we put $\mu(\vech):=(\rho_1(\vech),\rho_2(\vech),\ldots,\rho_n(\vech))$.
Then we put as follows:
\[
 \nbigk_{\vech}=\overline{\nbigk}_{\mu(\vech)}.
\]
Then the tuple $\{\nbigk_{\vech}\,|\,
 \vech\in\seisuu^{m}\,\,m=1,\ldots,l\}$
has the desired property.
\hfill\qed

\begin{df}
Such tuple
$\{\nbigk_{\vech}\,|\,\vech\in\seisuu^{m},\,\,m=1,\ldots,l\}$
is called a compatible splitting of the commuting tuple $(N_1,\ldots,N_l)$.
\hfill\qed
\end{df}

\begin{cor}\label{cor;12.15.100}
Let $E$ and $N_1,\ldots, N_l$ be as above.
We can take a holomorphic frame $\vecv$ of $E$,
which is compatible with $(N_1,\ldots,N_l)$
\hfill\qed
\end{cor}

\subsubsection{Strongly sequential compatibility}

Let $E$ be a holomorphic vector bundle over $\Delta^n$.
We put $D_i:=\{(z_1,\ldots,z_n)\in\Delta^n\,|\,z_i=0\}$,
and $D_{\mbar}=\bigcap_{i=1}^mD_i$.
Let $l$ be a number less than $n$.
Let $N_i$ be an element of
$\Gamma(D_i,End(E)_{|D_i})$
for $i=1,\ldots,l$.
For any $m\leq l$,
we have the nilpotent maps
$N_{1\,|\,D_{\mbar}},\ldots,N_{m\,|\,D_{\mbar}}$
of $E_{|D_{\mbar}}$.

\begin{df}
We say that the tuple $(N_1,\ldots,N_l)$ is strongly sequentially
 compatible, if the following holds:
\begin{enumerate}
\item
All the assumptions in Definition {\rm\ref{df;12.10.208}}
are satisfied.
\item
Moreover
$(N_{1\,|\,P},\ldots, N_{m\,|\,P})$ is assumed to be
strongly sequentially compatible for each $P\in D_{\mbar}$.
\hfill\qed
\end{enumerate}
\end{df}

\begin{lem} \label{lem;12.3.31}
There are decompositions of $E_{|D_{\mbar}}$ for $1\leq m\leq l$:
\[
 E_{|D_{\mbar}}=
 \bigoplus_{k\geq 0}
 \bigoplus_{\vech\in\seisuu^m}
 P_k\nbigk_{\vech}.
\]
They satisfy the following:
\begin{enumerate}
\item
 For any $\vech\in\seisuu^m$ and $k\geq 0$,
 we have
 $\bigcap_{j=1}^m W(\jbar)_{h_j}
=\bigoplus_{k\geq 0}\bigoplus_{\veck\in\nbigu(\vech)}P_k\nbigk_{\veck}$
on $D_{\mbar}$.
Here $\nbigu(\vech)$ denotes the set of $\veck\in\seisuu^n$
satisfying $\rho_j(\veck)\leq h_j$ for any $1\leq j\leq m$.
\item
 We have
 $P_k\nbigk_{\vech\,|\,D_{\mitibar}}=
 \bigoplus _{a}P_k\nbigk_{(\vech,a)}$.
Here $(\vech,a)=(h_1,\ldots,h_m,a)$
for $\vech=(h_1,\ldots,h_m)$.
\item
 $P_k\nbigk_{\vech}=0$ unless $|q_1(\vech)|\leq k$ and $k-q_1(\vech)$
is even.
\item
 When $-k< q_1(\vech)\leq k$,
 we have
 $N_1\bigl(
  P_k\nbigk_{\vech}\bigr)
 =P_k\nbigk_{\vech-2\vecdelta_1}$.
 Here we put $\vech-2\vecdelta_1=(h_1-2,h_2,\ldots,h_n)$.
\item
 When $h_1=-k$,
 $N_1\bigl(  P_k\nbigk_{\vech}\bigr)=0$.
\end{enumerate}
\end{lem}
\pf
On $D_{\itibar}$,
we have the graded vector space
$P_hGr_h^{(1)}$.
On $D_{\jbar}$ for $2\leq j\leq l$,
we have the nilpotent morphisms $N_{j}^{(1)}$
of $P_hGr_h^{(1)}$.
Since $(N_2^{(1)},\ldots,N_l^{(2)})$ is sequentially compatible,
we obtain the compatible decompositions
$P_hGr^{(1)}_{h\,|\,D_{\mbar}}
=\bigoplus_{\veck'\in\seisuu^{m}}
 \nbigk_{h,\veck'}$ for any $2\leq m\leq l$.

Then we construct $P_k\nbigk_{\vech}$ on $D_{\mbar}$
by using an induction on $m$.
Assume that we already have the decomposition
on $D_{\mitibar}$.
Let $\vech$ be a tuple $(h_1,\ldots,h_m)$.
We construct $P_{h_1}\nbigk_{\vech}$ on $D_{\mbar}$
in the following.

For a tuple $\vech=(h_1,\ldots,h_m)$,
we put $\chi(\vech):=(h_1,h_2-h_1,h_3-h_2,\ldots,h_m-h_{m-1})$,
and $\chi'(\vech)=(h_2-h_1,h_3-h_2,\ldots,h_m-h_{m-1})$.
Note that $\chi(\vech)\in\nbigu(\vech)$.
If $\vech'=(h_2-h_1,\ldots,h_{m}-h_1)\in\seisuu^{m-1}$,
then we have $\chi'(\vech)\in\nbigu(\vech')$.

By our assumption,
$\nbigk_{h_1,\chi'(\vech)}$ is contained in the image
of $P\pi_{\vech}$:
\[
 P\pi_{\vech}\left[
 \Ker(N(\itibar)^{h_1+1})\cap
 \bigcap_{j=2}^m W(\jbar)_{h_j}\right]
=\left[
 P_{h_1}Gr^{(1)}_{h_1}\cap\bigcap_{j=2}^n W(N^{(1)}(\jbar))_{h_j-h_1}
 \right]
 \supset P_{h_1}\nbigk_{h_1,\chi'(\vech)}.
\]
On $D_{\mitibar}$,
we have $\bigoplus_k\nbigk_{(\chi(\vech),k)}$.
By extending it,
we can take a subbundle $P_{h_1}\nbigk_{\chi(\vech)}$
of $\Ker(N(\itibar)^{h_1+1})\bigcap_{j=2}^mW(\jbar)_{h_j}$
satisfying the following:
\begin{itemize}
\item
$P_{h_1}\nbigk_{\chi(\vech)}$ is isomorphic
to $\nbigk_{h_1,\chi'(\vech)}$ via the morphism $P\pi_{\vech}$.
\item
We have
$P_{h_1}\nbigk_{\chi(\vech)|D_{\mitibar}}
=\bigoplus_aP_{h_1}\nbigk_{(\chi(\vech),a)}$.
\end{itemize}
By an inductive argument, 
we obtain $P_{h}\nbigk_{\vech}$ on $D_{\mbar}$
for any $m$ and $\vech\in\seisuu^m$ such that $h_1=h$.

For an integer $m$ such that $0\leq m\leq h_1$,
we put as follows:
\[
 P_{h_1}\nbigk_{\vech-2m\vecdelta_1}:=
 N_1^m\Bigl(
 P_{h_1}\nbigk_{\vech}
 \Bigr).
\]
Here we put $\vech-2m\vecdelta_1=(h_1-2m,h_2,\ldots,h_n)$.
By our choice,
we have $N_1^{h_1+1}\bigl(P_{h_1}\nbigk_{\vech} \bigr)=0$.
Then we obtain the desired decomposition.
\hfill\qed

\begin{df}
Such tuple
$\{P_k\nbigk_{\vech}\,|\,k\geq 0,\,\,\vech\in\seisuu^{m},\,m=1,\ldots,l\}$
is called the strongly compatible splitting
of $(N_1,\ldots,N_l)$.
\hfill\qed
\end{df}

For any tuple $\vech\in\seisuu^l$ and $k\geq 0$,
we have the number $d(k,\vech):=\rank P_k\nbigk_{\vech}$.
Clearly we have $d(k,\vech)=d(k-2,\vech)$
if $-k<q_1(\vech)\leq k$.

\begin{cor} \label{cor;12.15.101}
Let $E$ and $(N_1,\ldots,N_l)$ be as above.
Then we can take a holomorphic frame $\vecv$ of $E$,
around the origin $O$ of $\Delta^n$,
satisfying the following:
\begin{enumerate}
\item
$\vecv=
 \Bigl(v_{k,\vech,\eta}\,
 \Big|\,k\geq 0,\vech\in\seisuu^l,
 \eta=1,\ldots,d(k,\vech) \Bigr)$.
\item
 We have
 $N_1(v_{k,\vech,\eta})=v_{k,\vech-2\vecdelta_1,\eta}$
 when $-k<h_1\leq k$, on $D_{\itibar}$.
\item
 $N_1(v_{k,\vech,\eta})=0$
 if $h_1=-k$, on $D_{\itibar}$.
\item
 $\deg^{W(\jbar)}(v_{k,\vech,\eta})=
 \sum_{i\leq j}h_i$.
\item
$\vecv$ is compatible with the sequence
$\bigl(W(\itibar),\ldots,W(\lbar)\bigr)$.
\hfill\qed
\end{enumerate}
Such frame is called strongly compatible with
$(N_1,\ldots,N_l)$.
\end{cor}

\begin{df}
Such frame is called a frame strongly compatible
with $(N_1,\ldots,N_l)$.
\hfill\qed
\end{df}

\label{subsubsection;12.7.50}

\begin{df} \label{df;12.10.90}
A tuple $(N_1,\ldots,N_l)$ is called of Hodge,
if $(N_{\sigma(1)},\ldots,N_{\sigma(l)})$
is strongly sequentially compatible
for any $\sigma\in\gbigs_l$.
\hfill\qed
\end{df}

\subsection{Mixed twistor structure}
\label{subsection;12.15.4}

\subsubsection{Definition}
The harmonic metric can be regarded as a generalization
of the variation of polarized Hodge structure.
To regard the harmonic metric as a variation of {\em some structure},
Simpson introduced the twistor structure.
\begin{df}[Simpson \cite{s3}] \label{df;12.9.11}
The pure twistor structure and the mixed twistor structure
are defined as follows:
\begin{enumerate}
\item
A holomorphic vector bundle on the projective line $\proj^1$
is called a pure twistor structure of weight $i$
if it is holomorphically isomorphic to
a direct some of the line bundle $\nbigo_{\proj^1}(i)$.
\item
A holomorphic vector bundle $V$ with an ascending filtration $W$
by holomorphic vector subbundles
is called a mixed twistor structure
if $Gr_i^W$ is a pure twistor structure of weight $i$.
\end{enumerate}
\end{df}
He also introduced the variation of the twistor structure
and observed that a harmonic bundle can be regarded as 
the variation of the pure twistor structure.
(See \cite{s3} for more detail.)

In the next subsubsections,
we explain a method to use the mixed twistor structure
in this paper.

\subsubsection{Lower bound of the degree}
We will use the mixed twistor structure
to obtain a lower bound of the degree.
We consider the following situation.
Let $(V,W)$ be a mixed twistor structure.
Let $L$ be a holomorphic vector subbundle of $V$
with the filtration $W_L$.
If we have the upper bound of the degree
of any non-zero element of $W_{L,l}(L)$,
i.e., $\deg^W(s)\leq l+a$ for any
non-zero $s\in W_{L,l}(L)_{\lambda},(\lambda\in\proj^1)$,
then we have the inclusion $W_{L,l}(L)\subset W_{l+a}$.
We have the following easy lemma.

\begin{lem}\label{lem;1.8.10}
Let $(V,W)$ be a mixed twistor structure.
Let $a$ be an integer.
Let $L$ be a holomorphic vector subbundle of $V$
with the filtration $W_L$
such that $W_{L,l}(L)\subset W_{l+a}$
and that the first Chern classes
$c_1(Gr^{W_{L}}_l(L)))$ are $(l+a)\cdot\rank (Gr^{W_{L}}_l(L))$
for any $l$.

Then the induced morphisms
$Gr^{W_{L}}_l(L)\lrarr Gr^{W}_{l+a}$ is an injection of
vector bundles,
that is,
$Gr_l^{W_L}$ naturally gives a subbundle of $Gr^W_{l+a}$.
Moreover $Gr^{W_{L}}_l(L)$ is a pure twistor of weight $l+a$.

In particular, the degree of any non-zero element
$s\in W_{L,l}(L_{\lambda})-W_{L,l-1}(L_{\lambda})$
is exactly $l+a$.
\end{lem}
\pf
Let $b$ be the bottom number of the filtration $W_L$.
First we consider the bottom part $W_{L,b}$.
It is well known that
a holomorphic vector bundle $W_{L,b}$ over $\proj^1$
is holomorphically
a direct sum $\bigoplus_i\nbigo(i)^{\oplus n_i}$.
Assume that
$W_{L,b}$ is not a pure twistor structure of weight $b+a$.
Since we have the equality $\sum n_i i=(b+a)\cdot \rank(W_{L,b})$,
there is an $i>b+a$ such that $n_i\neq 0$.
It is easy to see that $Hom(\nbigo(i),W_{b+a})=0$
for any $i>b+a$.
Hence we have no injective morphism $W_{L,b}\rarr W_{b+a}$,
which contradicts the fact that $L$ is a subbundle of $W_{b+a}$.
Thus $W_{L,b}$ is a pure twistor structure of weight $b+a$.
Moreover the composition $W_{L,b}\lrarr W_{b+a}\lrarr Gr_{b+a}$ is
injective of the vector bundles.

We use an induction.
We assume that we have proved the claims of the theorem
for $W_{L,l}$ for any $l<l_0$
and then we prove that the claim for $W_{L,l_0}$ holds.
We have the morphism
$\phi_{l_0}:Gr_{l_0}^{W_L}\lrarr Gr^W_{l_0+a}$.
Assume that $\phi_{l_0}$ is not injective.
Note that $c_1(\ker(\phi_{l_0}))$
is larger than $(l_0+a)\cdot\rank \ker(\phi_{l_0})$.
We put $K_{l_0-1}:= \pi_{l_0-1}^{-1}(\ker\phi_{l_0})$,
where $\pi_{l_0-1}$ denotes the projection
$W_{L,l_0}\lrarr W_{L,l_0}/W_{L,l_0-1}$.
Then we have the morphism
$\phi_{l_0-1}:K_{l_0-1}/W_{L,l_0-2}\lrarr Gr_{l_0+a-1}$.
We have the nontrivial kernel $\ker(\phi_{l_0-1})$
such that 
$c_1(\ker(\phi_{l_0-1}))>(l_0-1+a)\cdot \rank(\ker(\phi_{l_0-1}))$.
We put $K_{l_0-2}:=\pi_{l_0-2}^{-1}(\ker\phi_{l_0-1})$,
where $\pi_{l_0-2}$ denotes the projection
$W_{L,l_0}\lrarr W_{L,l_0}/W_{l_0-2}$.
In general,
we denote the projection $W_{L,l_0}\lrarr W_{l_0}/ W_{L,i}$
by $\pi_i$.
Inductively,
we can construct the vector subbundles $K_i$
of $W_{L,l_0}$ as follows:
\begin{quote}
Assume that we have $W_{L,i}\subset K_{i}\subset W_{i+a}$.
Then we have the morphism 
$\phi_i:K_i/W_{L,i-1}\lrarr Gr_{i+a}$.
We put $K_{i-1}:=\pi_{i-1}^{-1}(\phi_i)$.
\end{quote}
Then we can check that
$c_1(K_i)>(i+a)\rank K_i$.
For sufficiently small $i$,
we obtain the inequality
$0=c_1(K_i)>(i+a)\cdot 0=0$.
Thus we arrive at the contradiction
if we assume that $\phi_{l_0}$ is not injective.
Thus we obtain the injectivity
of the morphism $\phi_{l_0}:Gr_{L,l_0}\lrarr Gr_{l_0+a}$
as a vector bundles,
namely $\phi_{l_0|P}$ is a injection for each point $P\in\proj^1$.
Since $Gr_{l_0+a}$ is pure twistor of weight $l_0+a$,
we have the inequality
$c_1(Gr_{L,l_0})\leq (l_0+a)\rank Gr_{l_0+a}$,
which is in fact equality by assumption.
Thus $Gr_{L,l_0}$ is a pure twistor of weight $l_0+a$.
\hfill\qed

\begin{rem}
It is remarkable that
we obtain the lower bound of degree
from some topological information, that is, Chern class.
\end{rem}

\subsubsection{Morphism of mixed twistors}

\begin{df}
Let $(V^{(i)},W^{(i)})$ $(i=1,2)$ be mixed twistors.
A morphism of mixed twistors are
the morphism of locally free coherent sheaves $V^{(1)}\lrarr V^{(2)}$
preserving the filtrations.
\hfill\qed
\end{df}

Let $f$ be a morphism of locally free coherent sheaves 
$V^{(1)}\lrarr V^{(2)}$.
We have the morphism $f_{|P}:V^{(1)}_{|P}\lrarr V^{(2)}_{|P}$
for any $P\in\proj^1$.
Then the rank of $f_{|P}$ are not constant, in general.
However, when $f$ is a morphism of mixed twistors,
then the rank of $f_{|P}$ is constant,
as we will see in the following lemma.

\begin{lem} \label{lem;12.7.30}
Let $(V^{(i)},W^{(i)})$ $(i=1,2)$ be mixed twistors,
and $f$ be a morphism of mixed twistors.
\begin{enumerate}
\item
The rank of $f_{|P}$ is constant,
and thus the $\Cok(f)$ is locally free.
Hence $\{\Ker(f)_{|P}\,|\,P\in\proj^1\}$ and
$\{\Image(f)\,|\,P\in\proj^1\}$ form
subbundles of $V^{(1)}$ and $V^{(2)}$ respectively.
\item
We put $W_l(\ker(f))=W_l^{(1)}\cap \ker(f)$.
The filtration $W_{\cdot}(\ker(f))$ induces the mixed twistor
structure to $\ker(f)$.
\item
We put $W_l(\cok(f))=\pi(W_l^{(2)})$,
where $\pi$ denotes the projection $V^{(2)}\lrarr \cok(f)$.
Then the filtration $W_{\cdot}(\cok(f))$ gives the mixed twistor structure
to $\cok(f)$.
\item \label{number;11.26.10}
We have the equality $f(W_l^{(1)})=\Image(f)\cap W_l^{(2)}$.
We put $W_l(\Image(f))=f(W_l^{(2)})$,
and then the filtration $W_{\cdot}(\Image(f))$
gives the mixed twistor structure
to $\Image(f)$.
\end{enumerate}
\end{lem}
\pf
Note that $W_l^{(i)}$ $(i=1,2)$ are naturally mixed twistor structures,
and that we have the morphism $f_l:W_l^{(1)}\lrarr W_l^{(2)}$
of the mixed twistors.
Thus we can use the induction on $l$.

Assume that $l$ is the bottom number $b_1$ of $W^{(1)}$.
We have the morphism $f_{b_1}:W^{(1)}_{b_1}\lrarr W^{(2)}_{b_1}$
and $W^{(1)}_{b_1}$ is isomorphic to
a direct sum of $\nbigo(b_1)$.
We denote the projection $W^{(2)}_{b_1}\lrarr Gr^{(2)}_{b_1}$
by $\pi^{(2)}_{b_1}$.
Then we obtain the following morphisms
\[
 \begin{CD}
 W^{(1)}_{b_1}@>{f_{b_1}}>>
 W^{(2)}_{b_1}@>{\pi^{(2)}_{b_1}}>>
 Gr^{(2)}_{b_1}.
 \end{CD}
\]
We have the composite $Gr_{b_1}(f):=\pi^{(2)}_{b_1}\circ f_{b_1}$.
Since $W^{(1)}_{b_1}$ and $Gr^{(2)}_{b_1}$ are pure twistor
of weight $b_1$,
it is easy to see the following:
\begin{itemize}
\item
The ranks of $Gr_{b_1}(f)_{|P}$
are independent of $P\in\proj^1$.
\item
The kernel, the image and the cokernel of $Gr_{b_1}(f)$
are pure twistors of weight $b_1$.
The kernel is a subbundle of $W^{(1)}_{b_1}$.
\end{itemize}
We have the naturally defined morphism
$\Ker(Gr_{b_1}(f))\lrarr W_{b_1-1}^{(2)}$.
Then it is easy to see that the morphism is in fact $0$.
Thus $\Ker(Gr_{b_1}(f))$ and $\ker(f_{b_1})$ are same.
We also obtain the following exact sequence:
\[
 0\lrarr W_{b_1-1}^{(2)}\lrarr \Cok(f_{b_1})\lrarr
 \Cok(Gr(f_{b_1}))\lrarr 0.
\]
Thus the $\Cok(f_{b_1})$ is locally free,
and the ranks of $f_{b_1\,|\,P}$ are independent of $P\in\proj^1$.
We also know that
the image $\Image(f_{b_1})$ is a subbundle of $W_{b_1}$,
and $\Image(f_{b_1})$ is a pure twistor of weight $b_1$.
Thus we can show the claim \ref{number;11.26.10}.
In all, we obtain the claims in the case that $l$ is the bottom number
of $W^{(1)}$.

\vspace{.1in}
We assume that the claims hold for
$f_{l-1}:W^{(1)}_{l-1}\lrarr W^{(2)}_{l-1}$,
and we will prove that the claims hold for
$f_l:W_l^{(1)}\lrarr W_l^{(2)}$.
Since $f_l$ preserves the filtration,
we have the natural morphism
$Gr_l(f):Gr_l^{(1)}\lrarr Gr_l^{(2)}$.
Because $Gr_l^{(i)}$ are pure twistors of weight $l$,
it is easy to see the following:
\begin{itemize}
\item The ranks of $Gr_l(f)_{|P}$ are independent of $P\in\proj^1$.
\item
The kernel, the image and the cokernel
of $Gr_l(f)$ is pure twistors of weight $l$.
The kernel is a subbundle of $Gr^{(1)}_l$,
and the image is a subbundle of $Gr^{(2)}_l$.
\end{itemize}
We have the natural morphism
$\phi:\ker(f_l)\lrarr \ker(Gr_l(f))$.
By an easy diagramm chasing,
we obtain the injection $\cok(\phi)\lrarr \Cok(f_{l-1})$
of coherent sheaves.
By our assumption,
$\Cok(f_{l-1})$ is mixed twistor,
such that $Gr_l(\Cok(f_{l-1}))=0$.
On the other hand, $\cok(\phi)$ is a quotient of
a pure twistor $\ker(Gr_l(f))$ of weight $l$.
Thus the morphism $\cok(\phi)\lrarr \cok(f_{l-1})$
must be $0$,
in other words,
$\phi$ must be surjective.
Thus we obtain the exact sequence:
\[
 0\lrarr\ker(f_{l-1})\lrarr \ker(f_l)\lrarr \ker(Gr_l(f))\lrarr 0.
\]
It implies the assertions for $\ker(f_l)$.
We also obtain the exact sequences:
\[
 \begin{array}{l}
 0\lrarr \Image(f_{l-1})\lrarr \Image(f_l)
\lrarr \Image(Gr_l(f))\lrarr 0.\\
\mbox{{}}\\
 0\lrarr \Cok(f_{l-1})\lrarr \Cok(f_l)
\lrarr \Cok(Gr_l(f))\lrarr 0.
 \end{array}
\]
It implies the assertions for $\Image(f_l)$ and $\cok(f_l)$.
Thus we are done.
\hfill\qed

\begin{rem}
Note that the property $4$ is remarkable.
If we consider a morphism of filtered vector spaces
$(V^{(1)},W^{(1)})\lrarr (V^{(2)},W^{(2)})$,
the properties $1$, $2$ and $3$ obviously holds.
However $4$ does not hold, in general.
\end{rem}

Let $(V,W)$ be a mixed twistor,
and $a$ be an integer.
Then we have the naturally defined mixed twistor
of $V\otimes\nbigo(a)$,
as follows:
\[
 W_l(V\otimes\nbigo(a))=W_{l-a}(V)\otimes\nbigo(a).
\]
It is easy to check that $(V\otimes\nbigo(a),W)$
gives a mixed twistor.

Let $f:V\lrarr V\otimes\nbigo(a)$ be a morphism
of locally free coherent sheaves.
Then $f_{|P}:V\lrarr V$ is determined
as an element of $\proj(M(r)^{\lor})$.
Here $r$ denotes a rank of $V$ and $M(r)^{\lor}$
denotes the dual space of the vector space of $r$-matrices.
When $f$ is nilpotent,
then it induces the filtration $W(f_{|P})$
on the fiber $V_{|P}$ for any point $P\in \proj^1$.

\begin{lem}
Let $(V,W)$ be a mixed twistor
and $f:(V,W)\lrarr (V\otimes\nbigo(2),W)$
be a nilpotent morphism of mixed twistor.
Then the conjugacy classes of the endomorphisms
$f_{|P}$ are independent of $P\in\proj^1$.
\end{lem}
\pf
Let $\eta$ be a generic point of $\proj^1$.
Then we have the filtration
$W(f_{|\eta})$ of $V_{|\eta}$ induced by the nilpotent maps $f_{|\eta}$.
Let $b$ be the bottom number of the filtration $W(f_{|\eta})$.
We use the induction on $b$.

When $b=0$,
then $f_{|\eta}=0$. Thus $f=0$. Thus we have nothing to prove.
Assume that we have proved the claim in the case $b>b_0$,
and we will prove the claim in the case $b=b_0$.
Note that the bottom part of the filtration $W(f_{|\eta})$
is same as $\Image(f_{|\eta}^{-b_0})$
in $V_{|\eta}$,
and $W_{-b_0-1}(f_{|\eta})$ is same as 
$\ker(f_{|\eta}^{-b_0})$.
We also have $f^{-b_0+1}=0$.

We have the morphism
$f^{-b_0}:V\lrarr V\otimes\nbigo(-2b_0)$,
which is a morphism of mixed twistors.
Then $\ker(f^{-b_0})$ and $\Image(f^{-b_0})$ are
subbundles of $V$ and $V\otimes\nbigo(-2b_0)$,
and they have the naturally induced mixed twistors.
We put as follows:
\[
 V'=\frac{\ker(f^{-b_0})}{\Image(f^{-b_0})\otimes\nbigo(2b_0)}.
\]
Then it has the naturally induced mixed twistor structure.

We have the naturally induced morphism
$\tilde{f}:V'\lrarr V'\otimes\nbigo(2)$.
We can apply the claim for $b<b_0$ to $\tilde{f}$,
and thus the conjugacy class of $\tilde{f}_{|P}$
are independent of $P\in\proj^1$.
Then we obtain the independence of the conjugacy class
of $\tilde{f}_{|P}$.
\hfill\qed

\subsubsection{Sub mixed twistors}

Let $(V,W)$ be mixed a twistor.
Let $V_1$ be a subbundle of $V$.
Then we obtain the filtration of $V_1$
by the coherent subsheaves $W_h\cap V_1$.

\begin{df}
We say that $V_1$ is a sub mixed twistor of $(V,W)$
if the filtration $\{W_h\cap V_1\,|\,h\in\seisuu\}$
gives a mixed twistor structure.
\hfill\qed
\end{df}

We have already shown the following:
\begin{lem}
Let $f:(V_1,W_1)\lrarr (V_2,W_2)$ be a morphism of mixed twistors.
Then the kernel and image are sub mixed twistors
of $(V_1,W_1)$ and $(V_2,W_2)$ respectively.
\hfill\qed
\end{lem}

\begin{lem}
Let $(V,W)$ be a mixed twistor,
and $V_i$ $(i=1,2)$ be sub mixed twistors.
Then $V_1+V_2$ and $V_1\cap V_2$ are also sub mixed twistors.
\end{lem}
\pf
We can regard $V_1\cap V_2$ is the kernel of
the morphism of mixed twistors $V_1\oplus V_2\lrarr V$.
We can regard $V_1+V_2$ is the image of the morphism.
\hfill\qed

\vspace{.1in}
Let $N:(V,W)\lrarr (V,W)\otimes\nbigo(2)$
be a morphism of mixed twistors.
Since the conjugacy classes of $N$ are independent
of $\lambda\in\proj^1$,
we obtain the weight filtration $W(N)$ of $N$
by vector subbundles.
The following lemma is easy to see.
\begin{lem}
For any $h$,
the vector subbundle $W(N)_h$ is a sub mixed twistor.
\hfill\qed
\end{lem}

Let $N_i:(V,W)\lrarr (V,W)\otimes\nbigo(2)$
be a morphism of mixed twistors for $i=1,\ldots,n$.
\begin{lem}
For any tuple $\vech=(h_1,\ldots,h_n)\in\seisuu^n$,
$\bigcap_{j=1}^n W(N_i)_{h_j}$ is sub mixed twistor 
of $(V,W)$
\hfill\qed
\end{lem}

\subsubsection{Commuting tuple of nilpotent maps}

Let $V$ be a vector bundle over $\proj^1$.
Let $(N_1,\ldots,N_n)$ be a commuting tuple of
nilpotent morphisms $V\lrarr V\otimes\nbigo_{\proj^1}(2)$.

\begin{prop} \label{prop;12.10.70}
Assume the following:
\begin{enumerate}
\item \label{number;12.3.22}
 The weight filtration $W(\nbar)$ of $N(\nbar)$ 
 is a mixed twistor.
\item \label{number;12.3.23}
 $N(\jbar):(V,W(\nbar))\lrarr (V,W(\nbar))\otimes\nbigo(2)$
 $(1\leq j\leq n)$ gives a morphism of mixed twistor structure.
\item \label{number;12.3.21}
 $(N_{1\,|\,\lambda},\ldots,N_{n-1\,|\,\lambda},N_{n\,|\,\lambda})$ 
 is sequentially compatible at any point $\lambda\in\proj^1$.
\item \label{number;12.3.20}
 $(N_{1\,|\,\lambda},\ldots,N_{n-1\,|\,\lambda})$
 is strongly sequentially compatible 
at any point $\lambda\in\proj^1$.
\end{enumerate}
Then $(N_{1\,|\,\lambda},\ldots,N_{n-1\,|\,\lambda}, N_{n\,|\,\lambda})$
is strongly sequentially compatible at any $\lambda\in\proj^1$.
\end{prop}
\pf
First we note that
we have
$W^{(1)}(\jbar)_{a+l}\cap Gr^{(1)}_a 
=W(N^{(1)}(\jbar))_{l}\cap Gr^{(1)}_a$,
by the assumption \ref{number;12.3.21}.
By the assumption \ref{number;12.3.20},
$P\pi_{\vech}$ induces the following isomorphism
for any $\vech\in\seisuu^{n-1}$:
\begin{equation} \label{eq;12.3.35}
  P\pi_{\vech}:
 \Ker(N(\itibar)^{h_1+1})\cap
 \bigcap_{j=2}^{n-1} W(\jbar)_{h_j}
\lrarr
 P_{h_1}Gr^{(1)}_{h_1}
 \cap
 \bigcap_{j=2}^{n-1}W^{(1)}(\jbar)_{h_j}
\end{equation}

Due to the conditions \ref{number;12.3.22} and \ref{number;12.3.23},
the filtrations $W(\nbar)$ and $W^{(1)}(\nbar)$
induces the mixed twistor structures
on the both sides of (\ref{eq;12.3.35}).
Due to the condition \ref{number;12.3.21},
The morphism $P\pi_{\vech}$ preserves the mixed twistor structures.
Thus we can conclude that
the morphism $P\pi_{\vech}$ gives an isomorphism
of the filtered vector bundles over $\proj^1$.
It implies that 
$(N_{1\,|\,\lambda},\ldots,N_{n-1\,|\,\lambda}, N_{n\,|\,\lambda})$
is strongly sequentially compatible at any $\lambda\in\proj^1$.
\hfill\qed

\section{Preliminary for harmonic bundles}
\label{section;12.15.70}

\subsection{Harmonic bundles and deformed holomorphic bundles}
\label{subsection;12.15.5}

\subsubsection{harmonic bundles}

Let $X$ be a complex manifold.
Let $(E,\delbar_E)$ be a holomorphic bundle.
Here $E$ denotes a $C^{\infty}$-vector bundle
and $\delbar_E$ denotes an operator
$\delbar_E:C^{\infty}(X,E)\lrarr C^{\infty}(X,E\otimes\Omega^{0,1}_X)$,
such that $(\delbar_E)^2=0$ and that
$\delbar_E(fv)=\delbar(f)\cdot v+f\!\cdot\!\delbar_E(v)$
for any $f\in C^{\infty}(X)$ and $v\in C^{\infty}(X,E)$.
Let $h$ be a hermitian metric of $E$.
We denote the inner product of $h$ by $(\cdot,\cdot)_h$.
We often omit $h$ if there is no confusion.
For a holomorphic vector bundle $(E,\delbar_E)$
with a hermitian metric,
we obtain 
$\del_E:C^{\infty}(X,E)\lrarr C^{\infty}(X,E\otimes\Omega^{1,0})$
satisfying
$\delbar(f,g)_h=(\delbar_E(f),g)_h+(f,\del(g))_h$.
We denote the curvature of the unitary connection
$\del_E+\delbar_E$
by $R(\del_E+\delbar_E)$.

Let $\theta$ be a section of $C^{\infty}(X,End(E)\otimes \Omega^{1,0})$.
It is called a (holomorphic) Higgs field
if $\delbar_E\theta=0$ and $\theta\wedge\theta=0$.
The tuple $(E,\delbar_E,h)$ is called a Higgs bundle.

We have the adjoint of $\theta$ with respect to $h$,
which we denote by $\theta^{\dagger}$,
namely
$(\theta\!\cdot\! f,g)_h=(f,\theta^{\dagger}\!\cdot\! g)_h$.
Then
$\theta^{\dagger}$ is an element of
$C^{\infty}(X,End(E)\otimes\Omega^{0,1})$
satisfying $\del_E(\theta^{\dagger})=0$ and
$\theta^{\dagger}\wedge\theta^{\dagger}=0$.

From a Higgs bundle $(E,\delbar_E,\theta)$
with a hermitian metric $h$,
we obtain the following connection:
\[
 \DD^1:=\delbar_E+\del_E+\theta+\theta^{\dagger}:
 C^{\infty}(X,E)\lrarr C^{\infty}(X,E\otimes\Omega^1_X).
\]
\begin{df}
A tuple $(E,\delbar_E,h,\theta)$ is called a harmonic bundle,
if $\DD^1$ is flat,
namely $\DD^1\circ \DD^1=0$.
\hfill\qed
\end{df}

\begin{rem}
Probably, such object should be called
a pluriharmonic bundle.
But we use `harmonic bundle' for simplicity.
\hfill\qed
\end{rem}

The condition $\DD^1\circ\DD^1=0$ is equivalent to
the following:
\[
 (\delbar_E+\theta^{\dagger})^2=
 (\del_E+\theta)^2=
 R(\del_E+\delbar_E)
+\theta\wedge\theta^{\dagger}
+\theta^{\dagger}\wedge\theta=0.
\]

\begin{lem}[Simpson, \cite{s4}]
Let $(E,\delbar_E,h,\theta)$
be a harmonic bundle.
Then we have $\del_E(\theta)=\delbar_E(\theta^{\dagger})=0$.
\end{lem}
\pf
We know that $\del_E^2=\theta^2=(\del_E+\theta)^2=0$.
It implies that $\del_E(\theta)=0$.
Similarly we obtain the equality
$\delbar_E(\theta^{\dagger})=0$.
\hfill\qed

\subsubsection{The deformed holomorphic bundle}

Let $(E,\delbar_E,h,\theta)$ be a harmonic bundle
over $X$.
We denote $\cnum_{\lambda}\times X$ 
by $\nbigx$.
We denote the projection $\cnum_{\lambda}\times X\lrarr X$
by $p_{\lambda}$.
We have the $C^{\infty}$-bundle $p_{\lambda}^{-1}(E)$
over $\nbigx$.
We have the operator
$d'':
 C^{\infty}(\nbigx,p_{\lambda}^{-1}(E))
\lrarr 
 C^{\infty}(\nbigx,p_{\lambda}^{-1}(E)\otimes \Omega_{\nbigx}^{0,1})$.
\[
 d'':=\delbar_E+\lambda\cdot\theta^{\dagger}+\delbar_{\lambda}.
\]

\begin{lem}[Simpson, \cite{s4}]
The operator $d''$ gives a holomorphic structure
of $p_{\lambda}^{-1}(E)$.
\end{lem}
\pf
We only have to see that $d''\circ d''=0$,
which follows from
$\delbar_E(\theta^{\dagger})=0$.
\hfill\qed

The holomorphic bundle $(p^{-1}_{\lambda}(E),d'')$
is denoted by $\nbige$,
which we call the deformed holomorphic bundle.
We have the pull back of the hermitian metric $h$.
The metric connection is given by the following:
\[
 d_{\lambda}+\delbar_E+\del_E
+\lambda\cdot\theta^{\dagger}
-\overline{\lambda}\cdot \theta.
\]
The curvature of the metric is as follows:
\begin{equation} \label{eq;11.27.25}
 -d\bar{\lambda}\cdot\theta+d\lambda\cdot\theta^{\dagger}+
R(\bar{\del}_E+\del_E)
-|\lambda|^2[\theta,\theta^{\dagger}].
\end{equation}

We put as follows:
\[
 \nbigx^{\lambda}:=\{\lambda\}\times X,
\quad
 \nbigx^{\shikaku}:=\cnum_{\lambda}^{\ast}\times X.
\]
The restrictions
$(\nbige,d'')_{|\nbigx^{\lambda}}$
and $(\nbige,d'')_{|\nbigx^{\shikaku}}$
are denoted by
$(\nbige^{\lambda},d^{\twoprime\,\lambda})$
and
$(\nbige^{\shikaku},d^{\twoprime\,\shikaku})$.

\subsubsection{The $\lambda$-connection}

Let $(E,\delbar_E)$ be a holomorphic bundle.
Let $\lambda$ be a complex number.
In general,
an operator
$\nabla^{\lambda}:
 C^{\infty}(X,E)\lrarr C^{\infty}(X,E\otimes\Omega^1)$
is called a $\lambda$-connection
if the following holds for any $f\in C^{\infty}(X)$
and $v\in C^{\infty}(X,E)$:
\[
 \nabla^{\lambda}(f\cdot v)=
 (\lambda\del(f)+\bardel(f))\cdot v
+f\cdot\nabla^{\lambda}(v).
\]
It is called holomorphic if the $(0,1)$-part is same as
$\delbar_E$.
It is called flat if $\nabla^{\lambda}\circ\nabla^{\lambda}=0$.

It is easy to see that 
a flat holomorphic $0$-connection is equivalent
to a pair of holomorphic structure and
a holomorphic Higgs field.
And a flat holomorphic $1$-connection
is equivalent to an ordinary holomorphic flat connection.

Let $(E,\delbar_E,\theta,h)$
be a harmonic bundle over $X$.
Then we have the operator
$\DD^{\lambda}:C^{\infty}(X,E)\lrarr C^{\infty}(X,E\otimes\Omega^{1}_X)$
defined as follows:
\[
 \DD^{\lambda}=\delbar_E+\theta+\lambda(\del_E+\theta^{\dagger}).
\]
Recall that we have the holomorphic bundle
$\nbige^{\lambda}$ on $\nbigx^{\lambda}$,
whose holomorphic structure is given by
$\delbar_E+\lambda\theta^{\dagger}$.

\begin{lem}
The operator $\DD^{\lambda}$
is a flat holomorphic $\lambda$-connection of
$\nbige^{\lambda}$.
\end{lem}
\pf
It is clear from the definition
that $\DD^{\lambda}$ gives a holomorphic $\lambda$-connection
of $\nbige^{\lambda}$.
We have the following equality:
\[
 \DD^{\lambda}\circ\DD^{\lambda}=
 \bigl(\delbar_E+\lambda\cdot\theta^{\dagger}\bigr)^2
+\bigl(\lambda\del_E+\theta\bigr)^2
+\lambda\cdot
 \Bigl(
 R(\del_E+\delbar_E)+[\theta,\theta^{\dagger}]
 \Bigr).
\]
Since we have $\del(\theta)=\delbar_E(\theta^{\dagger})=0$,
we can obtain the desired flatness.
\hfill\qed

The $\DD^{\lambda}$ is called the $\lambda$-connection
associated with the harmonic bundle $\harmonicbundle$.
We have the operator
$\DD:
 C^{\infty}(\nbigx,\nbige)
 \lrarr
 C^{\infty}(\nbigx,\nbige\otimes p_{\lambda}^{\ast}\Omega^1_X)$
defined by 
$\DD=\delbar_E+\theta+\lambda(\del_E+\theta^{\dagger})$.
The operator $\DD$ is also called the $\lambda$-connection
associated with $\harmonicbundle$.

Note that $\DD$ and $\delbar_{\lambda}$ are commutative.
Thus $\DD(v)$ is holomorphic if $v$ is holomorphic section
of $\nbige$.
Let $\vecv$ be a holomorphic frame of
$\nbige$ on an open subset $U$ of $\nbigx$.
Then the $\lambda$-connection form
$\nbiga=(\nbiga_{i\,j})$
is defined by the following relation:
\[
 \DD v_j=
 \sum_{i}\nbiga_{i\,j}\cdot v_i.
\]
Obviously $\nbiga_{i\,j}$ are holomorphic sections
of $p_{\lambda}^{\ast}\Omega_X^{1,0}$,
i.e.,
$\nbiga$ is an element of
$\Gamma(U,M(r)\otimes p^{-1}_{\lambda}\Omega_X^{1,0})$,
where $r$ is a rank of $E$.
We describe as $\DD\vecv=\vecv\cdot\nbiga$.

\subsubsection{The associated flat connections}

For any $\lambda\neq 0$,
we have the holomorphic flat connection
of $\nbigx^{\lambda}$:
\[
 \DD^{\lambda,f}:=
 \delbar_E+\lambda\cdot\theta^{\dagger}
+\del_E+\lambda^{-1}\theta.
\]
Again the flatness follows from the equalities
$\del_E(\theta)=\delbar_E(\theta^{\dagger})=0$.

We have the operator
$ \DD^{f}:=\delbar_E+\del_E+\lambda\theta^{\dagger}+\lambda^{-1}\theta:
 C^{\infty}(\nbigx^{\shikaku},\nbige^{\shikaku})\lrarr
 C^{\infty}(\nbigx^{\shikaku},\nbige^{\shikaku}
                  \otimes p_{\lambda}^{\ast}\Omega^1_X)$.
We call it the associated family of the flat connections.

Let $\vecv$ be a holomorphic frame of $\nbige^{\shikaku}$
on some open subset $U$ of $\nbigx^{\shikaku}$.
Then we obtain the holomorphic section
$\nbiga^{f}=(\nbiga^{f}_{i\,j})$
of $\Gamma(U,M(r)\otimes p_{\lambda}^{\ast}\Omega_X^{1,0})$
defined as follows:
\[
 \DD^{f}v_j=\sum_i \nbiga^f_{i\,j}\cdot v_i,
\quad
\mbox{i.e.,}
\quad
 \DD^f\vecv=\vecv\cdot\nbiga^{f}.
\]

\begin{lem}
We have the relation
$\nbiga^f=\lambda^{-1}\cdot\nbiga$.
\end{lem}
\pf
The $(0,1)$-parts of $\DD^f$ and $\DD$ are same,
we have the following relation between
the $(1,0)$-parts of $\DD^f$ and $\DD$:
\[
 \DD^{f\,(1,0)}
=\del_E+\lambda^{-1}\theta
=\lambda^{-1}(\lambda\cdot \del_E+\theta)
=\lambda^{-1}\cdot \DD^{(1,0)}.
\]
Thus we are done.
\hfill\qed

\subsubsection{Conjugate}

We denote the conjugate of $X$ by $X^{\dagger}$.
Namely $X^{\dagger}$ denotes the complex manifold
whose underlying $C^{\infty}$-manifold is same as $X$
and whose holomorphic structure is given by $\del$.
If $(E,\delbar_E)$ is a holomorphic bundle over $X$,
then $(E,\del_E)$ is a holomorphic bundle over $X^{\dagger}$.
If $\theta$ is a holomorphic Higgs field of $(E,\delbar_E)$
over $X$,
then $\theta^{\dagger}$ is a holomorphic Higgs field
of $(E,\del_E)$ over $X^{\dagger}$.

Let $\harmonicbundle$ be a harmonic bundle over $X$.
Then the conjugate $(E,\del_E,h,\theta^{\dagger})$ is
a harmonic bundle over $X^{\dagger}$.
We put $\nbigx^{\dagger}:=\cnum_{\mu}\times X^{\dagger}$.
We denote the projection $\nbigx^{\dagger}\lrarr X^{\dagger}$
by $p_{\mu}^{\dagger}$.
From a harmonic bundle $(E,\del_E,h,\theta^{\dagger})$
over $X^{\dagger}$,
we obtain the deformed holomorphic bundle
$(\nbige^{\dagger},d^{\twoprime\,\dagger})$
over $\nbigx^{\dagger}$.
Under the identification of $X$ and $X^{\dagger}$,
the underlying $C^{\infty}$-bundle
of $\nbige^{\dagger}$ is $p_{\mu}^{\dagger\,\ast}(E)$,
and the holomorphic structure $d^{\twoprime\dagger}$
is given by the operator
$\del_E+\mu\cdot\theta+\delbar_{\mu}$.

We put $\nbigx^{\dagger\,\mu}=\{\mu\}\times X^{\dagger}$
and $\nbigx^{\dagger\,\shikaku}=\cnum_{\mu}^{\ast}\times X$.
The restrictions of $(\nbige^{\dagger},d^{\twoprime\,\dagger})$
to $\nbigx^{\dagger\,\mu}$ and $\nbigx^{\dagger\,\shikaku}$
are denoted by $(\nbige^{\dagger\,\mu},d^{\twoprime\,\mu})$
and $(\nbige^{\dagger\,\shikaku},d^{\twoprime\,\shikaku})$.
The operator $d^{\twoprime\,\mu}$ is same as
$\del_E+\mu\cdot\theta$.


We have the associated $\mu$-connections,
which we denote by $\DD^{\dagger\,\mu}$.
Namely we have the following operator:
\[
 \DD^{\dagger\,\mu}:=
 \del_E+\theta^{\dagger}+\mu\cdot(\delbar_E+\theta):
 C^{\infty}(\nbigx^{\dagger\,\mu},\nbige^{\dagger\,\mu})
\lrarr
 C^{\infty}(\nbigx^{\dagger\,\mu},
    \nbige^{\dagger\,\mu}\otimes\Omega^{1}_{\nbigx^{\dagger\,\mu}}).
\]
We have the following operator, which we also call the $\mu$-connection:
\[
 \DD^{\dagger}:=\del_E+\theta^{\dagger}+\mu\cdot(\delbar_E+\theta):
 C^{\infty}\big(\nbigx^{\dagger},\nbige^{\dagger}\big)
\lrarr
 C^{\infty}\big(\nbigx^{\dagger},
    \nbige^{\dagger}\otimes
      p_{\mu}^{\dagger\,\ast}\Omega_{X^{\dagger}}^1\big).
\]
For any $\mu\neq 0$,
The associated flat connections are given as follows:
\[
 \DD^{\dagger\,\mu\,f}:=
 \del_E+\mu\cdot\theta+\delbar_E+\mu^{-1}\theta^{\dagger}:
 C^{\infty}(\nbigx^{\dagger\,\mu},\nbige^{\dagger\,\mu})
\lrarr
 C^{\infty}(\nbigx^{\dagger\,\mu},
    \nbige^{\dagger\,\mu}\otimes\Omega_{\nbigx^{\dagger\,\mu}}^1).
\]
We have the family of the flat connections:
\[
 \DD^{\dagger\,f}:=
 \del_E+\mu\cdot\theta+\delbar_E+\mu^{-1}\theta^{\dagger}:
 C^{\infty}(\nbigx^{\dagger\,\shikaku},\nbige^{\dagger\,\shikaku})
\lrarr
 C^{\infty}(\nbigx^{\dagger\,\shikaku},
   \nbige^{\dagger\,\shikaku}\otimes
       p_{\mu}^{\dagger\,\ast}\Omega_{X^{\dagger}}^1).
\]

The following lemma is clear from the definition.
\begin{lem}  \label{lem;12.1.6}
If $\lambda=\mu^{-1}$,
then we have $\DD^{\dagger\,\mu,f}=\DD^{\lambda,f}$
as the operator $C^{\infty}(X,E)\lrarr C^{\infty}(X,E\otimes\Omega^1_X)$.
Namely they give the same flat connection.
\hfill\qed
\end{lem}

We have the morphism
$\cnum_{\lambda}^{\ast}\lrarr \cnum_{\mu}^{\ast}$
by the correspondence $\mu=\lambda^{-1}$.
It induces the $C^{\infty}$-morphism
$\conjugate_X:\nbigx^{\shikaku}\lrarr\nbigx^{\dagger\,\shikaku}$.
Although $\conjugate_X$ is not holomorphic,
it is holomorphic in the direction of $\cnum_{\lambda}^{\ast}$.
The following lemma can be shown directly from the definitions.
\begin{lem}
Under the identification of
$\nbigx^{\shikaku}$ and $\nbigx^{\dagger\,\shikaku}$
by the morphism $\conjugate_X$ above,
we have the relation
$\DD^f=\DD^{\dagger,f}$.
\hfill\qed
\end{lem}

\subsubsection{Another relation between $\nbige$ and $\nbige^{\dagger}$}
\label{subsubsection;12.1.12}

Let $(E,\delbar_E,h)$ be a holomorphic vector bundle with a hermitian metric.
In general, a hermitian metric $h$ induces
an anti-linear morphism $\psi$ of $E$ to the dual $E^{\lor}$.
The morphism $\psi$ gives an anti-holomorphic isomorphism
of $E$ and $E^{\lor}$.

Let $\vecv=(v_1,\ldots,v_n)$ be a holomorphic frame of $E$.
We have the dual frame $\vecv^{\lor}$ of $E^{\lor}$.
Then we put $\vecv^{\dagger}=\psi^{-1}(\vecv^{\lor})$.
Namely we put as follows:
we have the matrix $(b_{i,j})=\overline{H(h,\vecv)}^{-1}$,
and frame $\vecv^{\dagger}=(v^{\dagger}_1,\ldots,v^{\dagger}_n)$ defined
by the relation
$v^{\dagger}_j:=\sum b_{i,j}\cdot v_i$,
that is,
$\vecv^{\dagger}=\vecv\cdot \overline{H(h,\vecv)}^{-1}$.
Then $\vecv^{\dagger}$ is an anti-holomorphic frame of $E$,
in other words,
$\vecv^{\dagger}$ is a holomorphic frame
of $(E,\del_E)$ over $X^{\dagger}$.

\vspace{.1in}

Let $\harmonicbundle$ be a harmonic bundle over $X$.
Then we have the deformed holomorphic bundle
$\nbige^{\lambda}$ over $\nbigx^{\lambda}$,
whose holomorphic structure is given by
$\delbar_E+\lambda\cdot\theta^{\dagger}$.
The anti-holomorphic structure is given by
$\del_E-\bar{\lambda}\cdot\theta$.
Thus $\vecv^{\dagger}$ gives a holomorphic frame of 
$\nbige^{\dagger\,(-\bar{\lambda})}$ over
$\nbigx^{\dagger\,(-\bar{\lambda})}$.
The following lemma can be checked by a direct calculation.
\begin{lem}
When we have the relation $\DD^{\lambda}\vecv=\vecv\cdot\nbiga$,
then we have the relation
$\DD^{\dagger\,(-\bar{\lambda})}\vecv^{\dagger}
  =\vecv^{\dagger}\cdot{}^t\overline{\nbiga}$.
\end{lem}
\pf
It follows from the following equality:
\[
 \bigl(\DD^{\lambda}v_i,v_j^{\dagger}\bigr)_h=
 \bigl((\lambda\cdot\del_E+\theta)v_i,v_j^{\dagger}\bigr)_h
 =\bigl(v_i,(-\bar{\lambda}\del_E+\theta^{\dagger})v_j^{\dagger}\bigr)_h
 =\bigl(v_i,\DD^{\dagger\,(-\lambdabar)}v_j\bigr)_h.
\]
\hfill\qed


We denote the metric connection of $(\nbige,d'',h)$
by $d''+d'$.
Then we have the holomorphic vector bundle
$(\nbige,d')$ over $\cnum_{\lambda}^{\dagger}\times X^{\dagger}$.
We have a holomorphic map $F:\cnum_{\mu}\lrarr \cnum_{\lambda}^{\dagger}$
defined by
$\lambdabar=-\mu$.
Then we have the naturally defined holomorphic map
$F:\nbigx^{\dagger}=\cnum_{\mu}\times X^{\dagger}
\lrarr \cnum_{\lambda}^{\dagger}\times X^{\dagger}$.
\begin{lem} \label{lem;12.1.10}
The holomorphic bundle
$F^{-1}(\nbige,d')$ is same as
$(\nbige^{\dagger},d^{\twoprime\,\dagger})$.
We also have $F^{-1}\DD=\DD^{\dagger}$.
\end{lem}
\pf
It can be checked by direct calculations.
\hfill\qed

Let $\vecv$ be a (not-necessarily holomorphic)
frame of $\nbige$.
Then $\vecv^{\dagger}$ can naturally be regarded
as a frame of $\nbige^{\dagger}$
in the sense of Lemma \ref{lem;12.1.10}.

The following lemmas can be checked by direct calculations.
\begin{lem} 
Let $\vecv_i$ $(i=1,2)$ be frames of $\nbige$
related by the matrices $B$,
that is,
$\vecv_1=\vecv_2\cdot B$.
Then we have the relation
\[
 \vecv_1^{\dagger}
=\vecv_2^{\dagger}\cdot F^{\ast}({}^t\bar{B}^{-1}).
\]
Here $F$ denotes the above morphism
$\cnum_{\mu}\times X^{\dagger}\lrarr 
 \cnum_{\lambda}^{\dagger}\times X^{\dagger}$
given by $\lambdabar=-\mu$.
\hfill\qed
\end{lem}

\begin{lem}
Let $\vecv$ be a holomorphic frame of $\nbige$
and $\nbiga$ be a $\lambda$-connection form of $\DD$
with respect to $\vecv$.
Then $F^{\ast}({}^t\overline{\nbiga})$
is the $\mu$-connection form of $\DD^{\dagger}$.
\hfill\qed
\end{lem}

\subsubsection{Functoriality}

Let $\harmonicbundle$ be a harmonic bundle over $X$.
We have a various kind of functorial construction.
We recall some of them to fix our notation.

\noindent
{\bf (Dual)}\\
We denote the dual bundle of $E$ by $E^{\lor}$.
The metric $h$ induces the naturally defined metric on $E^{\lor}$,
which we denote by $h^{\lor}$.
The Higgs field $\theta$ naturally induces
the holomorphic section $\theta^{\lor}=-\theta$
of $End(E^{\lor})\otimes\Omega^{1,0}_X\simeq
 End(E)\otimes\Omega^{1,0}_X$.
Then the tuple $(E^{\lor},\delbar_{E^{\lor}},\theta^{\lor},h^{\lor})$
gives a harmonic bundle.

Let $(\nbige,\DD)$ be a $\lambda$-connection.
We have the naturally induced $\lambda$-connection
$\DD$ on $\nbige^{\lor}$.
Let $f$ be a holomorphic section of $\nbige^{\lor}$
and $g$ be a holomorphic section of $\nbige$.
We denote the natural pairing by $(\cdot,\cdot)$.
We have the obvious $\lambda$-connection $\DD$ of $\nbigo_{\nbigx}$.
Then $\DD^{\lor}$ is determined by the following relation:
\[
 \DD(f,g)=(\DD^{\lor}f,g)+(f,\DD^{\lor}g).
\]
It is easy to see that
the deformed holomorphic bundle with $\lambda$-connection of
$(E^{\lor},\delbar_{E^{\lor}},\theta^{\lor},h^{\lor})$
is naturally isomorphic to
$(\nbige^{\lor},\DD^{\lor})$.

Similarly the conjugate deformed holomorphic bundle
with $\mu$-connection
of $\harmonicbundledual$ is naturally isomorphic
to the dual of $(\nbige^{\dagger},\DD^{\dagger})$.

\vspace{.1in}
\noindent
{\bf (Tensor product)}\\
Let $(E_i,\delbar_{E_i},\theta_i,h_i)$ $(i=1,2)$
are harmonic bundles over $X$.
We have the tensor product $E_1\otimes E_2$.
We have the naturally defined metric $\tilde{h}:=h_1\otimes h_2$
and the Higgs field
$\tilde{\theta}:=\theta_1\otimes id_{E_2}+id_{E_1}\otimes \theta_2$.
Here $id_{E}$ denotes the identity of $E$.
The tuple $(E_1\otimes E_2,\tilde{\theta},\tilde{h})$
gives a harmonic bundle over $X$.
It is denoted by $(E_1,\theta_1,h_1)\otimes (E_2,\theta_2,h_2)$.

It is clear that the deformed holomorphic bundle with $\lambda$-connection
of $(E_1,\theta_1,h_1)\otimes (E_2,\theta_2,h_2)$
is naturally isomorphic
to the tensor product
$(\nbige_1,\DD_1)\otimes (\nbige_2,\DD_2)$.
Here $(\nbige_i,\DD_i)$ denotes the deformed holomorphic bundle
of $(E_i,\theta_i,h_i)$.
We have a similar relation for the conjugate deformed holomorphic
bundles with $\mu$-connections.

\vspace{.1in}
\noindent
{\bf (A direct summand)}\\
Assume that $(E_1,\theta_1,h_1)\oplus (E_2,\theta_2,h_2)$ be
harmonic bundles over $X$.
Then it is easy to see that the direct summands
$(E_i,\theta_i,h_i)$ are also harmonic bundles.
We have the obvious relations for the deformed holomorphic bundles
and the conjugate deformed holomorphic bundles.

\vspace{.1in}
\noindent
{\bf (The invariant part with respect to the group action)}
 \\
Let $\harmonicbundle$ be a harmonic bundle over $X$.
Let $G$ be a finite group acting on the vector bundle $E$.
We assume that the action of $G$ on $X$ is trivial,
and preserves $\theta$ and the metric.
We have the decomposition of $E$ 
by the irreducible representations of $G$.
The direct summands are naturally harmonic bundles.
Such decomposition is compatible with 
the construction of the deformed holomorphic bundles
and the conjugate deformed holomorphic bundles.

\vspace{.1in}
\noindent{\bf (The symmetric products and the exterior products)}\\
We have the symmetric product and the exterior products
of the harmonic bundle $\harmonicbundle$.
They are characterized as the direct summands of 
the tensor products.

We denote the symmetric product by
$\Sym^l\harmonicbundle=(\Sym^l(E),\theta^{\sym\,l},h^{\Sym\,l})$.
The deformed holomorphic bundle is denoted by
$(\Sym^l\nbige,\DD^{\sym\,l},h^{\sym\,l})$.

We denote the exterior product by 
$\bigwedge^l\harmonicbundle=
 (\bigwedge^l(E),\theta^{\wedge\,l},h^{\wedge l})$.
The deformed holomorphic bundle
is denoted by
$(\bigwedge^l\nbige,\DD^{\wedge\,l},h^{\wedge\,l})$.

\vspace{.1in}
\noindent
{\bf (Determinant)}\\
As a special case of the exterior product,
we have the determinant line bundle $\det(E)$.
The metric $h$ naturally induces the metric of $\det(E)$,
which we denote by $\det(h)$.
We have the natural isomorphism
$End(\det(E),\det(E))\simeq\nbigo_X$,
and the trace morphism
$\tr:End(E,E)\lrarr \nbigo_X$.
The Higgs field $\theta$ naturally induces 
the Higgs field $\tr(\theta)\in \Gamma(X,\Omega^{1,0})$.
The tuple $(\det(E),\delbar_{\det(E)},\tr(\theta),\det(h))$
gives a harmonic bundle.
The deformed holomorphic bundle
is denoted by 
$(\det(\nbige),\tr(\DD),\det(h))$.

\vspace{.1in}
\noindent
{\bf (Hom)}\\
Let $(E_i,\theta_i,h_i)$ $(i=1,2)$ be harmonic bundles.
Then we have
$(Hom(E_1,E_2),\theta_1^{\lor}+\theta_2,h_1^{\lor}\otimes h_2)
\simeq
(E_1^{\lor},\theta_1^{\lor},h_1^{\lor})
\otimes
(E_2,\theta_1,h_2)$.
We have the obvious relations for the deformed holomorphic bundles
with $\lambda$-connections.

\vspace{.1in}
\noindent
{\bf (Pull back)}\\ 
Let $f:Y\lrarr X$ be a morphism of complex manifolds.
Let $\harmonicbundle$ be a harmonic bundle over $X$.
Then we have the pull back $(f^{\ast}E,f^{\ast}\theta,f^{\ast}h)$,
which gives a harmonic bundle over $Y$.
The deformed holomorphic bundle with $\lambda$-connection
of $f^{\ast}(E,\theta,h)$ is naturally isomorphic
to the pull back of $(\nbige,\DD)$.

\subsection{Model bundles} \label{subsection;11.27.27}

In this subsection, we recall some examples of harmonic bundles
over the punctured disc $\Delta^{\ast}$.

\subsubsection{line bundles $L(a)$}

We put $L=\nbigo_{\Delta^{\ast}}\cdot e$.
For any real number $a\in \real$,
we put $h_a(e,e)=|z|^{2a}$.
Then $h_a$ gives a hermitian metric of $L$.
The anti-holomorphic structure $\del_L$
is determined by
\[
 \del_L(e)=e\cdot a\frac{dz}{z}.
\]
Then $\delbar_L+\del_L$ gives a flat connection of $L$.

We have the trivial Higgs field $\theta=0$ of $L$,
and the conjugate $\theta^{\dagger}=0$.
Then the tuple $L(a)=(L,h_a,\theta)$ gives a harmonic bundle.

Let see the deformed holomorphic bundle $(\nbigl(a),\DD)$
over $\cnum_{\lambda}\times \Delta^{\ast}$.
The holomorphic structure $d''$ is $\delbar_L+\delbar_{\lambda}$.
Thus the natural $C^{\infty}$-frame $p_{\lambda}^{-1}(e)$
is holomorphic.
The $\lambda$-connection $\DD$ is $\delbar_L+\lambda\cdot\del_L$.
Thus we have the following relation:
\[
 \DD(e)=e\cdot
 \Bigl(
 \lambda a\frac{dz}{z}\Bigr).
\]
The associated family of the flat connections
is as follows:
\[
 \DD^f(e)=e\cdot 
 \Bigl(
 a\frac{dz}{z}\Bigr).
\]

\noindent
{\bf The parabolic structure
 and the eigenvalue of the residue}\\
We have the natural prolongment of $\nbigl(a)$ over
$\cnum_{\lambda}\times\Delta^{\ast}$
to the holomorphic line bundle $\nbigl(a)_a$
over $\cnum_{\lambda}\times \Delta$
by using the frame $v$ above.
Note that this prolongment is same as 
the prolongment by increasing order $a$
(See the subsection \ref{subsubsection;12.9.60}).
In fact, the norm $|v|$  is $|z|^a$.
The $\lambda$-connection $\DD$ is of log type
on $\nbigl(a)_a$.
Namely, for a holomorphic section $f$ of $\nbigl(a)_a$,
$\DD(f)$ is a holomorphic section
of $\nbigl(a)_a\otimes p_{\lambda}^{-1}\Omega^{1,0}(\log(O))$.
The residue of $\DD$ is
$\lambda\cdot a\in 
\Gamma\big(\cnum_{\lambda},
  End(\nbigl(a)_{a\,|\,\cnum_{\lambda}\times O})\big)
\simeq \Gamma(\cnum_{\lambda},\nbigo)$.

\vspace{.1in}
We have the conjugate $(\nbigl(a)^{\dagger},\DD^{\dagger})$.
The holomorphic structure is given by
$\del_L+\delbar_{\mu}$.
The section $e^{\dagger}=|z|^{-2a}\cdot e$ gives
a holomorphic frame.
The $\mu$-connection is $\DD^{\dagger}=\del_L+\mu\cdot\delbar_L$,
and we have the relation:
\[
 \DD^{\dagger}e^{\dagger}=
 e^{\dagger}\cdot \mu \cdot(-a) \frac{d\bar{z}}{\bar{z}}.
\]

The associated family of the flat connections
are as follows:
\[
 \DD^{\dagger,f}e^{\dagger}=
 e^{\dagger}\cdot (-a)\frac{d\bar{z}}{\bar{z}}.
\]

\subsubsection{line bundles $L(\alpha)$}

We put $L=\nbigo\cdot e$ and $h_0(e,e)=1$.
Then we have $\del_L(e)=0$.
For any complex number $\alpha\in\cnum$,
we have the Higgs field $\theta=\alpha\cdot dz/z$.
The conjugate $\theta^{\dagger}$
is $\bar{\alpha}\cdot d\bar{z}/\bar{z}$.
We have the following relations:
\[
 R(\del_L+\delbar_L)=0,
\quad\quad
 \theta\wedge\theta^{\dagger}+\theta^{\dagger}\wedge\theta=
 |\alpha|^2
 \bigl(
 dz\cdot d\bar{z}+d\bar{z}\cdot dz
 \bigr)=0.
\]
Thus $L(\alpha)=(L,h_0,\theta)$ is a harmonic bundle.

We have the deformed holomorphic bundle
$(\nbigl(\alpha),\DD)$.
The holomorphic structure is given by the operator
$d''=\delbar_L+\lambda\bar{\alpha}\cdot
d\bar{z}/\bar{z}+\delbar_{\lambda}$.
Thus we have the relation:
\[
 d''(e)=e\cdot (\lambda\bar{\alpha})\frac{d\bar{z}}{\bar{z}}.
\]
We put $v=\exp(-\lambda\baralpha\log|z|^2)\cdot e$.
Then $v$ gives a holomorphic frame of $\nbigl(\alpha)$.
We have the equality:
\[
 h(v,v)=\exp\Bigl(-(\lambda\baralpha+\barlambda\alpha)\log|z|^2\Bigr)
 =|z|^{-4Re(\lambda\alphabar)}.
\]
The $\lambda$-connection $\DD$ is as follows:
\[
 \DD(v)
=\bigl(
  \lambda\del_L+\theta
 \bigr)v
=\Bigl(\lambda(-\lambda\alphabar)+\alpha\Bigr)\cdot v\cdot
 \frac{dz}{z}
=v\cdot(-\lambda^2\alphabar+\alpha)\frac{dz}{z}.
\]
The associated family of the flat connection
is as follows:
\[
 \DD^{f}v=
 v\cdot(-\lambda\alphabar+\lambda^{-1}\alpha)\frac{dz}{z}.
\]

\noindent
{\bf The parabolic structure and the eigenvalues of
the residues}\\
Let fix a $\lambda\in\cnum_{\lambda}$.
We have the natural prolongment
of $\nbigl(\alpha)^{\lambda}$ over $\Delta^{\ast}$
to the line bundle $\nbigl(\alpha)^{\lambda}_{-2Re(\lambda\alphabar)}$
by using the holomorphic $v$ above.
Note that the prolongment is same as the prolongment
by the increasing order $-2Re(\lambda\alphabar)$.
In fact, the norm $|v|$ is $|z|^{-2Re(\lambda\alphabar)}$.
The $\lambda$-connection $\DD$ is of log type
on the $\nbigl(\alpha)_{-2Re(\lambda\alphabar)}$.
The residue is
$(-\lambda^2\alphabar+\alpha)\in
End\bigl(
 \nbigl(\alpha)_{-2Re(\lambda\alphabar)\,|\,(\lambda,O)}
 \bigr)\simeq \cnum$.

\subsubsection{line bundles $\nbigl(a,\alpha)$}

As is already seen,
we have the harmonic bundles $\nbigl(a)$ and $\nbigl(\alpha)$
for $a\in\real$ and $\alpha\in\cnum$.
The tensor product $\nbigl(a,\alpha):=\nbigl(a)\otimes\nbigl(\alpha)$
is also a harmonic bundle.
It is a tuple $(L,h_a,\alpha\!\cdot\!dz/z)$.
We denote the deformed holomorphic bundle
by $\nbigl(a,\alpha)$.
We have the natural frame
$v=\exp\big(-\lambda\alphabar \log|z|^2\big)\cdot e$
of $\nbigl(a,\alpha)$.
We have the following:
\[
 h_a(v,v)=|z|^{2(a-2Re(\lambda\alphabar))},
\quad
 \DD(v)=v\cdot \big(-\lambda^2\alphabar+\alpha+\lambda a\big)\frac{dz}{z},
\quad
 \DD^f(v)=v\cdot
 \big(-\lambda\alphabar+\lambda^{-1}\alpha+a\big)\frac{dz}{z}.
\]
Let fix a $\lambda$.
We have the prolongment of
$\nbigl(\alpha,a)^{\lambda}$
over $\Delta^{\ast}$
to the line bundle
$\bigl(
 \nbigl(\alpha,a)^{\lambda}
 \bigr)_{a-2Re(\lambda\alphabar)}$
over $\Delta$.
This is same as the prolongment by the increasing order
$a-2Re(\lambda\alphabar)$.
The $\lambda$-connection $\DD^{\lambda}$ is of log type.
The residue is $-\lambda^2 \alphabar+\alpha+\lambda a$.

\begin{lem}
Let take an arbitrary pair $(A,B)\in \real\times\cnum$.
Let fix a $\lambda$.
We can take a pair $(a,\alpha)\in \real\times\cnum$
such that
the parabolic structure of
$\nbigl(a,\alpha)^{\lambda}_{a-2Re(\lambda\alphabar)}$
is $A$,
and that
the residue of $\DD^{\lambda}$
is same as $B$.
In fact, we only have to put as follows:
\[
 a=
\frac{
 (1-|\lambda|^2)\cdot A+2Re(\lambdabar\cdot B)}
 {1+|\lambda|^2},
\quad
 \alpha=
  \frac{B-\lambda\cdot A}
 {1+|\lambda|^2}
\]
\end{lem}
\pf
We only have to check that
$A=a-2Re(\lambda\alphabar)$ and
$B=-\lambda^2\alphabar+\alpha+\lambda a$.
It can be checked by a direct calculation.
\hfill\qed

\subsubsection{$Mod(l,a,C)$}

\label{subsubsection;12.4.90}
We put $V_2=\cnum\cdot e_1\oplus\cnum\cdot e_{-1}$.
We have the nilpotent map $N_2$
defined by $N_2(e_1)=e_{-1}$ and $N_2(e_{-1})=0$.
We have the holomorphic vector bundle
$E_2:=V_2\otimes\nbigo_{\Delta^{\ast}}$,
which have the natural frame $\vece=(e_1,e_{-1})$.
We have the Higgs field $\theta_2=N_2\cdot dz/z$.
We have the following relation:
\[
 \theta_2(e_1,e_{-1})=(e_{-1},0)=
 (e_1,e_{-1})
 \left(
 \begin{array}{ll}
 0 & 0\\ 1 & 0
 \end{array}
 \right)\frac{dz}{z}.
\]

We put $y=-\log|z|^2$.
We take a metric $h_2$ as follows:
\[
 h_2(e_{-1},e_{-1})=y^{-1},
\quad
  h_2(e_{-1},e_1)=0,\quad
 h_2(e_1,e_1)=y.
\]
The conjugate $\theta^{\dagger}_2$ of $\theta_2$
with respect to the metric $h_2$ is as follows:
\[
 \theta^{\dagger}_2(e_1,e_{-1}) =
 (e_1,e_{-1})\cdot
 \left(
 \begin{array}{ll}
 0 & 1\\ 0 & 0
 \end{array}
 \right)
 \frac{d\zbar}{\zbar\cdot y^2}
\]

\begin{lem}
A tuple $Mod(2)=(E_2,\bardel_{E_2},h_2,\theta_2)$
gives a harmonic bundle.
\end{lem}
\pf
We denote
the hermitian matrix $H(h_2,\vece)$
by $H$.
The matrix $H$ and and the connection form $H^{-1}\del H$
is as follows:
\[
 H=
 \left(
 \begin{array}{cc}
 y & 0 \\ 0 & y^{-1}
 \end{array}
 \right),
\quad
 H^{-1}\del H=
 \left(
 \begin{array}{cc}
 -y & 0 \\ 0 & y^{-1}
 \end{array}
 \right)
 \frac{dz}{z\cdot y}.
\]
The curvature 
$R(\del_{E_2}+\delbar_{E_2})=\delbar\big(H^{-1}\del H\big)$
is as follows:
\[
 \delbar
 \bigl(H^{-1}\del H
 \bigr)
=\left(
 \begin{array}{cc}
 1 & 0 \\ 0 & -1
 \end{array}
 \right)\frac{dz\cdot d\zbar}{y^2\cdot |z|^2}.
\]
On the other hand, we have the following:
\[
 \theta_2\cdot \theta_2^{\dagger}+\theta_2^{\dagger}\cdot\theta_2
=
 \left(
 \left(
 \begin{array}{cc}
 0 & 0 \\ 0 & 1
 \end{array}
 \right)
+\left(
 \begin{array}{cc}
 -1 & 0 \\ 0 & 0
 \end{array}
 \right)
  \right)
\frac{dz\cdot d\barz}{|z|^2\cdot y^2}.
\]
Thus we obtain the equality
$R(\delbar_{E_2}+\del_{E_2})+[\theta_2,\theta_2^{\dagger}]=0$.
The other relations $\del_{E_2}(\theta_2)=0$ and 
$\delbar_{E_2}(\theta_2^{\dagger})$
are trivial in this case, because 
we have $\Omega_{\Delta^{\ast}}^{2,0}=\Omega_{\Delta^{\ast}}^{0,2}=0$.
\hfill\qed

Thus we obtain the harmonic bundle $Mod(2)$.
Note that we have the natural frame $\vece=(e_{1},e_{-1})$.

\vspace{.2in}

On the vector space $\bigotimes^{l-1}V_2$ $(l\geq 1)$,
we have the nilpotent map
\[
 N_l:=\sum_{a=0}^{l-2}
 \overbrace{1\otimes \cdots\otimes 1}^a\otimes N_2\otimes
 1\otimes \cdots\otimes 1.
\]
We have the action of $(l-1)$-th symmetric group $\gbigs_{l-1}$
on the $(l-1)$-th tensor product $V_2^{\otimes\,(l-1)}$
defined by the transposition of the components.
The $l$-dimensional vector space $V_l:=\Sym^{l-1}(V_2)$
is characterized as the invariant part of the action.
Since the nilpotent map $N_l$ commutes with the action of
$\gbigs_{l-1}$,
the map $N_l$ preserves the subspace $V_l$.
We denote the restriction of $N_l$ to $V_l$ also by $N_l$.

We have the harmonic bundle
$(E_2^{\otimes\,(l-1)},\theta_l,h_l)$,
where we put $\theta_l=N_ldz/z$
and $h_l=\bigotimes^{l-1}h_2$.
We have the natural inclusion
$E_l:=V_l\otimes\nbigo\subset E_2^{\otimes\,(l-1)}$,
which is the invariant part of the $\gbigs_{l-1}$-action.
The action of $\gbigs_{l-1}$
preserves $\delbar_{E_2^{\otimes(l-1)}}$, $\theta_l$ and $h_l$.
Thus the action preserves $\del_{E_2^{\otimes(l-1)}}$,
and $\theta_l^{\dagger}$.
Hence $\del_{E_2^{\otimes\,(l-1)}}$, $\theta_l$ and $\theta_l^{\dagger}$
preserves the subbundle $E_l$.
The adjoint of $\theta_l|_{E_l}$ with respect to the metric
$h_l|_{E_l}$
is same as the restriction $\theta_l^{\dagger}|_{E_l}$.
We denote the restrictions of $\theta_l$, $\theta_l^{\dagger}$
and $h_l$ to $E_l$ by the same notation.
From the consideration above,
it is easy to see the following.

\begin{lem}
The tuple $Mod(l)=(E_l,\theta_l,h_l) $ is a harmonic bundle.
\hfill\qed
\end{lem}

Thus we obtain the harmonic bundle $Mod(l)$.
Note that we have the characterization
of $Mod(l)$ as the $\gbigs_{l-1}$-fixed part
of $Mod(2)^{\otimes\,(l-1)}$.
Note also that we have the natural frame
$\vece=(e^{p,q}\,|\,p+q=l-1)$ defined as follows:
\[
 e^{p,q}:=e_{1}^p\cdot e_{-1}^q.
\]
Here the product is the symmetric product.
Clearly $\vece$ is an orthogonal holomorphic frame of $E_l$.

\vspace{.1in}
Let take a positive number $a>0$
and a complex number $C\in\cnum$ such that $|C|\leq 1$.
Then the following is easy.
\begin{lem}
The tuple $Mod(l,a,C)=(E_l,\theta_l,a\cdot h_l(C\cdot z))$
gives a harmonic bundle
for any $a>0$ and $0<|C|\leq 1$.
\end{lem}
\pf
We have no contribution of the positive constant $a$
to the condition for the harmonic bundles.
For any $C$ such that $0<|C|\leq 1$,
we have the morphism $f_c:\Delta^{\ast}\lrarr \Delta^{\ast}$
defined by $z\longmapsto C\cdot z$.
Then we have the isomorphism between
$Mod(l,a,C)$ and $f^{\ast}_cMod(l,a,1)$
by using the natural frame $\vece$.
Thus $Mod(l,a,C)$ is also harmonic.
\hfill\qed

We have the prolongment of $E_l$ over $\Delta^{\ast}$
to $(E_l)_0$ over $\Delta$,
by using the frame $\vece$.
Note that it is same as the prolongment
by the increasing order $0$.
The $\theta_l$ is of log type.
The residue $\bigl((E_l)_{0\,|\,O},\Res(\theta_l)\bigr)$
is naturally isomorphic to $(V_l,N_l)$.

Let $V$ be a finite dimensional vector bundle,
and $N$ be a nilpotent map on $V$.
We can take an isomorphism
$(V,N)\simeq \bigoplus_i (V_{l_i},N_{l_i})$.
Thus we have the following:
\begin{lem}
For any pair $(V,N)$ as above,
we have a harmonic bundle $\bigoplus_i Mod(l_i,a_i,C_i)$
such that the residue
is isomorphic to $(V,N)$.
\hfill\qed
\end{lem}

\noindent
{\bf Notation}\\
Although we do not have a canonical choice of such harmonic bundle,
we often denote such harmonic bundle
by $\modelbundle{V}{N}$ for simplicity.

\subsubsection{The deformed holomorphic bundle of $Mod(2,a,C)$}

First we see the deformed holomorphic bundle
$\modeldeform(2,a,C)$
of $Mod(2,a,C)$.
We put $y=-\log|z|^2$ and $c=-\log|C|^2\geq 0$.
We have the natural $C^{\infty}$-section 
$p_{\lambda}^{-1}e_{i}$ $(i=1,-1)$
of $\modeldeform(2,a,C)$.
For simplicity, we denote $p_{\lambda}^{-1}e_i$ by $e_i$.
We put as follows:
\[
 v^{1,0}=e_1,\quad
 v^{0,1}=e_{-1}-\frac{\lambda}{y+c}\cdot e_1.
\]
Then it can be directly checked that $v^{1,0}$ and $v^{0,1}$
are holomorphic.
Thus $\vecv=(v^{1,0},v^{0,1})$ gives a holomorphic frame
of $\modeldeform(2,a,C)$.
The $\lambda$-connection $\DD$ is as follows:
\[
 \DD(v^{1,0},v^{0,1})=
 \bigl(v^{0,1}\cdot dz/z,0\bigr)=
 \bigl(v^{1,0},v^{0,1}\bigr)\cdot
 \left(
 \begin{array}{cc}
 0 & 0 \\ 1 & 0
 \end{array}
 \right)
 \frac{dz}{z},
\quad
\mbox{i.e.,}
\quad
 \DD\vecv=\vecv\cdot
 \left(
 \begin{array}{cc}
 0 & 0 \\ 1 & 0
 \end{array}
 \right)
 \frac{dz}{z}.
\]
The associated family of the flat connections is as follows:
\[
 \DD^f\vecv=
 \vecv\cdot
  \left(
 \begin{array}{cc}
 0 & 0 \\ \lambda^{-1} & 0
 \end{array}
 \right)
 \frac{dz}{z}.
\]
To see the conjugate $\modeldeform(2,a,C)^{\dagger}$,
the frame $\vecv^{\dagger}=(v^{\dagger\,1,0},v^{\dagger,0,1})$
is given by the following formula:
\[
 \vecv^{\dagger}=\vecv\cdot \overline{H(h_2,\vecv)}^{-1}.
\]
Then $\vecv^{\dagger}$ gives a holomorphic frame
of $\modeldeform(2,a,C)^{\dagger}$,
and we have the following relation:
\[
 \DD^{\dagger}\vecv^{\dagger}=
 \vecv^{\dagger}\cdot
 \left(
 \begin{array}{cc}
 0 & 1 \\ 0 & 0
 \end{array}
 \right)
\frac{d\zbar}{\zbar},
\quad
 \DD^{\dagger\,f}\vecv^{\dagger}=
 \vecv^{\dagger}\cdot
  \left(
 \begin{array}{cc}
 0 & \mu^{-1} \\ 0 & 0
 \end{array}
 \right)
\frac{d\zbar}{\zbar}.
\]

We take a $C^{\infty}$-frame $\vecv'=(v^{\prime\,1,0},v^{\prime\,0,1})$
given as follows:
\[
 v^{\prime\,1,0}=y^{-1/2}\cdot v^{1,0},
\quad
 v^{\prime\,0,1}=y^{1/2}\cdot v^{0,1}.
\]
Let consider the transformation matrices
$\transmatrix_0$ and $\transmatrix_0'$
defined as follows:
\[
 \vecv=\vece\cdot \transmatrix_0=\vecv(0)\cdot\transmatrix_0,
\quad\quad
 \vecv'=\vecv'(0)\cdot\transmatrix_0'.
\]
Here $\vecv(0)$ denotes $p_{\lambda}^{-1}(\vecv|_{\lambda=0})$
precisely. Similar to $\vecv'(0)$.
Clearly $B_0$ and $B_0'$ are 
elements of the space
$C^{\infty}(\cnum_{\lambda}\times \Delta^{\ast},M(r))$.
In fact, they are given as follows:
\begin{equation} \label{eq;12.1.15}
 \transmatrix_0(\lambda,z)=
 \left(
 \begin{array}{cc}
 1 & -\lambda\cdot(y+c)^{-1}\\ 0 & 1
 \end{array}
 \right),
\quad
 \transmatrix_0'(\lambda,z)=
  \left(
 \begin{array}{cc}
 1 & -\lambda\cdot y\cdot(y+c)^{-1}\\ 0 & 1
 \end{array}
 \right).
\end{equation}

\begin{lem}
Let $\transmatrix_0$ and $\transmatrix_0'$ be as above.
\begin{enumerate}
\item
 Let $R_1$ and $R_2$ be real numbers 
satisfying $0<R_1$ and $0<R_2<1$.
The function $\transmatrix_0(\lambda,z)$ is bounded over
$\Delta_{\lambda}(R_1)\times \Delta_z^{\ast}(R_2)$.
The boundedness is independent of $C$ with $|C|\leq 1$.
\item
 Let $R_1$ be a real number
 satisfying $R_1>0$.
 The function $\transmatrix_0'(\lambda,z)$
is bounded over $\Delta_{\lambda}(R_1)\times \Delta^{\ast}$.
The boundedness is independent of $C$ with $|C|\leq 1$.

In particular,
the $C^{\infty}$-frame $\vecv'$ is adapted on the region.
\end{enumerate}
\end{lem}
\pf
The boundedness of the $\transmatrix_0$ and $\transmatrix_0'$
are clear.
It is easy to directly see that $\vecv'(0)$ is adapted.
Then the adaptedness of $\vecv'$ follows from
the boundedness of $\transmatrix_0'$.
\hfill\qed

For the conjugate,
the $C^{\infty}$-frame
$\vecv^{\prime\,\dagger}
  =(v^{\prime\,\dagger,1,0},v^{\prime\,\dagger,0,1})$
is given as follows:
\[
 v^{\prime\,\dagger,1,0}=y^{1/2}v^{\dagger\,1,0},
\quad
 v^{\prime\,\dagger,0,1}=y^{-1/2}v^{\dagger\,0,1}.
\]
Then we obtain the transformation matrices
$\transmatrix_0^{\dagger}$ and $\transmatrix_0^{\prime\,\dagger}$
satisfying
$\vecv^{\dagger}=\vecv^{\dagger}(0)\cdot \transmatrix_0^{\dagger}$
and
$\vecv^{\prime\dagger}=\vecv^{\prime\,\dagger}(0)
 \transmatrix_0^{\prime\dagger}$.
The formula of $\transmatrix_0^{\dagger}$
and $\transmatrix_0^{\prime\dagger}$ are given as follows:
\begin{equation} \label{eq;12.1.16}
 \transmatrix^{\dagger}_0(\mu,\zbar)=
 \left(
 \begin{array}{cc}
 1 & 0 \\ -\mu\cdot (y+c)^{-1} & 1
 \end{array}
 \right),
\quad
  \transmatrix^{\prime\dagger}_0(\mu,\zbar)=
 \left(
 \begin{array}{cc}
 1 & 0 \\ -\mu\cdot y\cdot(y+c)^{-1} & 1
 \end{array}
 \right).
\end{equation}

\subsubsection{The deformed holomorphic bundle of $Mod(l+1,a,C)$}

We describe the deformed holomorphic bundle
$\modeldeform(l+1,a,C)$ of $Mod(l+1,a,C)$.
We remark that the bundle
$\modeldeform(l+1,a,C)$ is the $\gbigs_{l}$-invariant
part of $\bigotimes^l\modeldeform(2,a,C)$.
Thus the following holomorphic sections
$\vecv=\bigl(v^{p,q}\,|\,p+q=l\bigr)$
of $\bigotimes^l\modeldeform(2,a,C)$ gives a holomorphic frame
of $\modeldeform(l+1,a,C)$:
\[
 v^{p,q}:=(v^{1,0})^p\cdot (v^{0,1})^q.
\]
Here the product is the symmetric product.
The $\lambda$-connection $\DD$ is described as follows:
\begin{equation}\label{eq;11.26.1}
 \DD (v^{p,q})=p\cdot v^{p-1,q+1}\frac{dz}{z}.
\end{equation}
By using the trivialization $\vecv$,
the $\modeldeform(l+1,a,C)$ is prolonged
to the holomorphic vector bundle
over $\cnum_{\lambda}\times \Delta$.
It is same as the prolongment by the increasing order $0$.
Thus the prolongment is denoted by
$\modeldeform(l+1,a,C)_0$.

The $\lambda$-connection $\DD$ is of log type
on $\modeldeform(l+1,a,C)_0$.
Thus we obtain the residue:
\[
 \Res(\DD)\in \Gamma\Bigl(\cnum_{\lambda}\times \{O\},
 End\bigl(\modeldeform(l+1,a,C)_{0\,|\,\cnum_{\lambda}\times \{O\}}\bigr)
 \Bigr).
\]
\begin{lem}
The conjugacy class of $\Res(\DD)$ of the model bundle
is independent of $\lambda$.
\end{lem}
\pf
Clear from the equation (\ref{eq;11.26.1}).
\hfill\qed

\begin{cor}\label{cor;12.15.110}
Let $V$ be a finite dimensional vector space
and $N$ be a nilpotent map on $V$.
Let $\modelbundle{V}{N}$ be a harmonic bundle
whose residue of the Higgs field is isomorphic to $(V,N)$.
Then the residue of the $\lambda$-connection
is also isomorphic to $(V,N)$.
\hfill\qed
\end{cor}

We have the $C^{\infty}$-frame
$\vecv'=\bigl(v^{\prime\,p,q}\,|\,p+q=l\bigr)$
of $\modeldeform(l+1,a,C)$, defined as follows:
\[
 v^{\prime\,p,q}:=
 v^{p,q}\cdot y^{\frac{q-p}{2}}.
\]
We have the transformation matrices
$\transmatrix_0=(\transmatrix_{0,i,p})$ and
$\transmatrix_0'=(\transmatrix'_{0,i,p})$
satisfying the equalities,
$\vecv=\vecv(0)\cdot \transmatrix_0$
and $\vecv'=\vecv'(0)\cdot \transmatrix_0'$,
or more precisely:
\[
 v^{p,l-p}(\lambda,z)=\sum_i \transmatrix_{0,i,p}(\lambda,z)\cdot
 v^{i,l-i}(0,z),
\quad
 v^{\prime\,p,l-p}(\lambda,z)=
 \sum_i \transmatrix'_{0,i,p}(\lambda,z)\cdot v^{\prime,i,l-i}(0,z).
\]
They are elements of
$C^{\infty}(\cnum_{\lambda}\times\Delta^{\ast},
 M(r))$.
\begin{lem}
We have the following  equalities:
\begin{equation} \label{eq;12.1.20}
 \transmatrix_{0,i,p}(\lambda,z)=
 \left\{
 \begin{array}{ll}
 0, & (i<p) \\
 c(l-p,i-p)\cdot \lambda^{i-p}\cdot (-y-c)^{-i+p}, &(i\geq p).
 \end{array}
 \right.
\end{equation}
Here $c(l,p)$ denotes the constant given as the
coefficient of $x^p$ in the polynomial $(1+x)^l$.

Similarly we have the following:
\begin{equation} \label{eq;11.26.2}
 \bigl(
  \transmatrix_0\bigr)^{-1}_{i,p}(\lambda,z)=
 \left\{
 \begin{array}{ll}
 0, & (i<p) \\
 c(l-p,i-p)\cdot \lambda^{i-p}\cdot (y+c)^{-i+p}, &(i\geq p).
 \end{array}
 \right.
\end{equation}

As a result, the matrices $\transmatrix_0$
and $\transmatrix^{-1}_0$ are bounded
on 
the region $\Delta_{\lambda}(R_1)\times \Delta^{\ast}_z(R_2)$
for any $0<R_1$ and $0<R_2<1$.
The boundedness is independent of the parameter $C\geq 0$.
\end{lem}
\pf
We put $\gamma=-\lambda\cdot(y+c)^{-1}$.
Then we have the following equalities:
\begin{multline}
 v^{p,l-p}(\lambda,z)
=e_1^p\cdot\bigl(e_{-1}(z)+\gamma\cdot e_1(z)\bigr)^{l-p}
=\sum_{j=0}^{l-p}
 e_1^p\cdot c(l-p,j)\cdot e_{-1}^j(z)\cdot
 e_1^{l-p-j}(z)\cdot\gamma^{l-p-j}
 \\
=\sum_{i=p}^{l-p}
 c(l-p,i-p)\gamma^{i-p}\cdot e_1^i(z)\cdot e_{-1}^{l-i}(z)
=\sum_{i=p}^{l}c(l-p,i-p)\cdot\gamma^{i-p}\cdot v^{i,l-i}(0,z)
\end{multline}
Thus we are done.
\hfill\qed

We have the following immediate corollary.
\begin{cor} \label{cor;11.26.5}
We have the following equalities:
\[
 \transmatrix'_{0,i,p}(\lambda,z)=
 \left\{
 \begin{array}{ll}
 0, & (i<p), \\
 c(l-p,i-p)\cdot\lambda^{i-p}\cdot(-y-c)^{-i+p}\cdot y^{i-p}, &(i\geq p).
 \end{array}
 \right.
\]
We also have the following equalities:
\[
 \bigl(
\transmatrix^{\prime\,-1}_0\bigr) _{i,p}(\lambda,z)=
 \left\{
 \begin{array}{ll}
 0, & (i<p), \\
 c(l-p,i-p)\cdot\lambda^{i-p}\cdot(y+c)^{-i+p}\cdot y^{i-p}, &(i\geq p).
 \end{array}
 \right.
\]
In particular,
$\transmatrix_0'(\lambda,z)$ and their inverses are bounded
on the region
$\Delta_{\lambda}(R_1)\times \Delta_z^{\ast}$
for any $0<R_1$.
The boundedness is independent of the parameter $C\geq 0$.
\end{cor}
\pf
The formula is obtained from (\ref{eq;11.26.2}).
The boundedness follows from the boundedness
of $y^{-1}(y+c)$.
\hfill\qed

\begin{cor}\mbox{{}} \label{cor;11.27.41}
\begin{enumerate}
\item
 We put $H_l:=H(h_l,\vecv)$.
Then $H_l(\lambda,z)$ and $H_l(0,z)$ are mutually bounded
over $\Delta_{\lambda}(R_1)\times \Delta^{\ast}_z(R_2)$
for any $0<R_1$ and $0<R_2<1$.
\item
The $C^{\infty}$-frame $\vecv'$ is adapted over
$\Delta_{\lambda}(R_1)\times \Delta_{z}^{\ast}$
for any $R_1>0$.
\end{enumerate}
\end{cor}
\pf
The first claim follows from the boundedness of
$\transmatrix_0$ and $\transmatrix_0^{-1}$.
The adaptedness of $\vecv'(0)$ is easy to see.
Thus the adaptedness of $\vecv'$ is obtained from
the boundedness of $\transmatrix_0'$.

Otherwise, it is also easy to see that
the adaptedness of $\vecv'$ follows from the adaptedness of $\vecv'$
in the case $l=2$.
\hfill\qed

\begin{cor} \label{cor;12.2.5}
When $|z|\to 0$,
$\transmatrix_0'(\lambda,z)$ converge to
$\transmatrix^{\heartsuit}_{0}$:
\[
  \transmatrix^{\heartsuit}_{0,j,p}(\lambda)=
 \left\{
 \begin{array}{ll}
 0, & (i<p), \\
 c(l-p,i-p)\cdot\lambda^{i-p}\cdot(-1)^{i-p}, &(i\geq p),
 \end{array}
 \right.
\quad
 \bigl(
 \transmatrix^{\heartsuit\,-1}_0\bigr)_{j,p}(\lambda)=
 \left\{
 \begin{array}{ll}
 0, & (i<p), \\
 c(l-p,i-p)\cdot\lambda^{i-p}, &(i\geq p),
 \end{array}
 \right.
\]
In particular,
we have
$\bigl(\transmatrix^{\heartsuit}_0\bigr)_{p,p}=
 \bigl(\transmatrix^{\heartsuit\,-1}_0\bigr)_{p,p}=1$
for any $0\leq p\leq l$.
We also have
$\bigl(
 \transmatrix^{\heartsuit}_0\bigr)_{l,0}=(-1)^l\cdot\lambda^l$
and
$\bigl(
 \transmatrix^{\heartsuit\,-1}_0
 \bigr)_{l,0}=\lambda^l$.
\end{cor}
\pf
This is a direct corollary of Corollary \ref{cor;11.26.5}
\hfill\qed

It will be convenient to replace the adaptedness for $\vecv'$
in Corollary \ref{cor;11.27.42}
with the adaptedness for some more flexible frame.
We have the residues $\Res(\DD)$,
which is a holomorphic section of $End(\modeldeform(l+1,a,C)_0)$
on $\cnum_{\lambda}\times \{O\}$.
Since the conjugacy classes are independent of $\lambda$,
the nilpotent maps $\Res(\DD)_{|(\lambda,O)}$
induces the weight filtration $\weightfilt$ of
$\modeldeform(l+1,a,C)|_{\cnum_{\lambda}\times O}$
by vector subbundles.
Let take a holomorphic frame $\vecw$ of $\modeldeform(l+1,a,C)_0$,
which is compatible with the filtration $\weightfilt$
on $\cnum_{\lambda}\times \{O\}$.
We put as follows:
\[
 k(w_i):=\frac{1}{2}\deg^{\weightfilt}(w_i).
\]
We have the $C^{\infty}$-frame $\vecw'=(w_1',\ldots,w_{l+1}')$
of $\modeldeform(l+1,a,C)$ over $\Delta^{\ast}$
defined as follows:
\[
 w_i':=y^{-k(w_i)}\cdot w_i.
\]
The following is an easy corollary of Corollary \ref{cor;11.27.41}.
\begin{cor} \label{cor;11.27.42}
Let $\vecw'$ be as above.
Then $\vecw'$ is adapted over $\Delta_{\lambda}(R)\times X$
for any $R>0$.
\end{cor}
\pf
We only have to see that the transformation matrices
of $\vecw'$ and $\vecv'$ and their inverses
are bounded.
The section $w_i$ is described by $\vecv$ as follows:
\[
 w_i=\sum_p C_{p,i}\cdot v^{p,l-p}.
\]
Here $C_{p,i}$ are holomorphic on
$\cnum_{\lambda}\times \Delta$.
Since both of $\vecw$ and $\vecv$ are compatible
with the $\weightfilt$,
the following holds:
\begin{quote}
If $\deg^{\weightfilt}(w_i)<\deg^{\weightfilt}(v^{p,l-p})$,
then $C_{p,i}(\lambda,O)=0$.
\end{quote}
We have the following relation:
\[
 w'_i=\sum_p C_{p,i}\cdot y^{-k(w_i)+k(v^{p,l-p})}\cdot v^{\prime\,p,l-p}.
\]
Then the function $C_{p,i}\cdot y^{-k(w_i)+k(v^{p,l-p})}$
is bounded over $\Delta_{\lambda}(R)\times X$
for any $R>0$.
Similarly the inverse of the transformation matrices
are bounded.
\hfill\qed

\vspace{.1in}

Similarly we can see the conjugate $\modeldeform(l+1,a,C)$.
We have the holomorphic frame
$\vecv^{\dagger}=\bigl(v^{\dagger,p,q}\,|\,p+q=l\bigr)$,
given as follows:
\[
 v^{\dagger,p,q}=
 \bigl(v^{\dagger,1,0}\bigr)^p\cdot
 \bigl(v^{\dagger,0,1}\bigr)^q.
\]
The $\mu$-connection $\DD^{\dagger}$
is described as follows:
\begin{equation} \label{eq;11.26.6}
 \DD^{\dagger}v^{\dagger,p,q}=
 q\cdot v^{\dagger,p+1,q-1}.
\end{equation}
The prolongment of $\modeldeform(l+1,a,C)^{\dagger}$
by the frame $\vecv^{\dagger}$
is same as the prolongment by increasing order $0$.
We denote the prolongment by $\modeldeform(l+1,a,C)^{\dagger}_0$.
The $\mu$-connection $\DD^{\dagger}$
is of log type on $\modeldeform(l+1,a,C)^{\dagger}_0$.
The residue can be obtained from (\ref{eq;11.26.1}).

We have the $C^{\infty}$-frame
$\vecv^{\prime\,\dagger}=
 \bigl(v^{\prime\,\dagger,p,q}\,|\,p+q=l\bigr)$,
given as follows:
\[
 v^{\prime\,\dagger,p,q}=y^{\frac{p-q}{2}}\cdot v^{\dagger,p,q}.
\]
Then we have the transformation matrices
$\transmatrix_0^{\dagger}=\bigl(\transmatrix^{\dagger}_{0\,i,q}\bigr)$
and $\transmatrix_0^{\prime\,\dagger}=
 \bigl(\transmatrix^{\prime\dagger}_{0\,i,q}\bigr)$
defined as follows:
\[
 v^{\dagger,l-q,q}=\sum_i \transmatrix^{\dagger}_{0,i,q}\cdot
 v^{\dagger,l-i,i},
\quad\quad
 v^{\prime\,\dagger,l-q,q}=
 \sum_i \transmatrix^{\prime\,\dagger}_{0,i,q}\cdot
 v^{\dagger,l-i,i}.
\]
\begin{rem}
Note that the rules of the subscripts
are slightly different in the cases
$\transmatrix^{\dagger}_0,\transmatrix^{\prime\dagger}_0$
and $\transmatrix_0$ and $\transmatrix_0'$.
Briefly speaking, the order is reversed.
\hfill\qed
\end{rem}

We have the following formulas by the same calculation:
\begin{equation} \label{eq;12.1.21}
 \transmatrix^{\dagger}_{0,i,q}:=
\left\{
 \begin{array}{ll}
 0 & (i<q)\\
 c(l-q,i-q)\cdot \mu^{i-q}\cdot (-y-c)^{-i+q} & (i\geq q).
 \end{array}
\right.
\end{equation}
And thus we have the following:
\[
 \transmatrix^{\prime\dagger}_{0,i,q}:=
\left\{
 \begin{array}{ll}
 0 & (i<q)\\
 c(l-q,i-q)\cdot \mu^{i-q}\cdot y^{i-q}\cdot(-y-c)^{-i+q} & (i\geq q),
 \end{array}
\right.
\]
and 
\[
  \transmatrix^{\prime\dagger\,-1}_{0,i,q}:=
\left\{
 \begin{array}{ll}
 0 & (i<q)\\
 c(l-q,i-q)\cdot \mu^{i-q}\cdot y^{i-q}(y+c)^{-i+q} & (i\geq q).
 \end{array}
\right.
\]
The limit $\transmatrix^{\heartsuit\dagger}_0$
of $\transmatrix^{\prime\,\dagger}_{0}$
is as follows:
\[
 \transmatrix^{\heartsuit\dagger}_{0,i,q}:=
\left\{
 \begin{array}{ll}
 0 & (i<q)\\
 c(l-q,i-q)\cdot \mu^{i-q} & (i\geq q),
 \end{array}
\right.
\quad\quad
 \transmatrix^{\heartsuit\dagger\,-1}_{0,i,q}:=
\left\{
 \begin{array}{ll}
 0 & (i<q)\\
 c(l-q,i-q)\cdot \mu^{i-q} & (i\geq q).
 \end{array}
\right.
\]

\subsubsection{Limit of $Mod(l+1,a,C)$}

Let $\psi_n$ denote the morphism $\Delta\lrarr\Delta$
defined by $\psi_n(z)=z^n$.
Later we will consider a `limit' of
harmonic bundles $\psi_n^{-1}(E,\delbar_E,\theta,h)$
for a harmonic bundle $(E,\delbar_E,\theta,h)$
on $\Delta^{\ast}$.
Here we see what happens in the case of model bundles
$Mod(l,a,C)=(E_l,\theta_l,a\cdot h_l(C\cdot z))$.
We can assume that $C\in\real$.
Recall that we have the natural frame
$\vece:=\vecv_{|\lambda=0}$ of $Mod(l,a,C)$,
that is we put $e^{p,q}(z):=v^{p,q}(0,z)$.

Let $F$ be a holomorphic vector bundle with a frame
$\vecu=\bigl(u^{p,q}\,|\,p+q=l\bigr)$
over $\Delta^{\ast}$.
We take a holomorphic isomorphism
of $\psi_n^{-1}(E_l)$ to $F$
by the following correspondence:
\[
 n^{-(l-2i)/2}\cdot \psi_n^{-1}(e^{i,l-i})\longmapsto u^{i,l-i}.
\]
Then we obtain the Higgs fields $\theta^{(n)}$
and the metrics $h_l^{(n)}$ on $F$
for any $n$.

We also have the natural frame $\vece:=\vecv|_{\lambda=0}$
of $Mod(l,a,C^{1/n})$.
By the frames $\vece$ and $\vecu$,
we have the isomorphism
of $F$ and $E_l$.
It induces the isomorphism of harmonic bundles
$(F,\theta^{(n)},h^{(n)})$
and $Mod(l,a,C^{1/n})$.

We take $\Theta^{(n)}\in
 \Gamma\bigl(\Delta^{\ast},M(r)\otimes\Omega^{1,0}_{\Delta^{\ast}}\bigr)$
and $H^{(n)}\in C^{\infty}\bigl(\Delta^{\ast},\nbigh(r)\bigr)$
defined as follows:
\[
 \theta^{(n)}\vecu=\vecu\cdot\Theta^{(n)},
\quad\quad
 H^{(n)}=H(h^{(n)},\vecu).
\]
Then it is easy to see that
the sequences  $\{\Theta^{(n)}\}$ and $\{H^{(n)}\}$
are convergent on any compact subset $K\subset \Delta^{\ast}$,
as the elements of
$\Gamma\bigl(\Delta^{\ast},M(r)\otimes\Omega^{1,0}_{\Delta^{\ast}}\bigr)$
or
$C^{\infty}\bigl(\Delta^{\ast},\nbigh(r)\bigr)$.
In fact, we have $\Theta^{(n)}=\Theta$
and $H^{(n)}(z)=a\cdot H(C^{1/n}\cdot z)$.
Here $\Theta$ and $H$ are given as follows:
\[
 \theta_l\vece=\vece\cdot\Theta,
\quad\quad
 H=H(h_l,\vece).
\]
Thus $\Theta^{(n)}$ converges to $\Theta$
and $H^{(n)}$ converges to $a\cdot H$.

\vspace{.1in}

In the discussion above,
we take a limit at $\lambda=0$.
Clearly, we can take a limit at any $\lambda$.
Namely, we put $\vecw(z):=\vecv(\lambda,z)$
for a fixed $\lambda$,
and consider the frames
$\vecw^{(n)}=(w^{(n)\,p,q})$ 
defined by $w^{(n),p,q}:=n^{-(p-q)/2}\psi_n^{-1}(w^{p,q})$.
We have the $\lambda$-connection form
$\psi_n^{-1}(\DD^{\lambda})\vecw^{(n)}=\vecw^{(n)}\cdot\nbiga^{(n)}$.
Then
$\bigl\{H(\psi_n^{-1}(h_l),\vecw^{(n)})\bigr\}$
and $\bigl\{\nbiga^{(n)}\bigr\}$
converges on any compact subset $K\subset \Delta^{\ast}$
as the sequences in $C^{\infty}(\Delta^{\ast},\nbigh(r))$
and $\Gamma(\Delta^{\ast},M(r)\otimes\Omega^{1,0}_{\Delta^{\ast}})$.

\subsection{A convergence of a sequence of harmonic metrics}

\label{subsection;12.4.10}

\subsubsection{A convergence at $\lambda=1$}

Let $X$ be $\Delta^{\ast\,l}\times \Delta^{d-l}$,
and $(E^{(n)},\theta^{(n)},h^{(n)})$ be a harmonic bundles on $X$
such that $\rank(E^{(n)})=r$.
Recall that
we have the deformed holomorphic bundles
$(\nbige^{1\,(n)},d^{\twoprime\,1\,(n)},\DD^{1\,(n)},h^{(n)})$
on $\nbigx^{1}=\{1\}\times X\subset \nbigx$.
In this subsection, the metric and the measure of $X$
are $\sum_{i=1}^ndz_i\cdot d\zbar_i $
and $\prod_{i=1}^n |dz_i\cdot d\zbar_i|$.

Assume that we have frames $\vecw^{(n)}=(w_1^{(n)},\ldots,w_r^{(n)})$
of $\nbige^{1\,(n)}$,
and the following holds:
\begin{condition}\mbox{{}} \label{condition;11.27.5}
\begin{enumerate}
\item
 We have the connection forms
 $A^{(n)}\in \Gamma(X,M(r)\otimes\Omega_X^{1,0})$
determined by $\DD^{1\,(n)}\vecw^{(n)}=\vecw^{(n)}A^{(n)}$.
Then the sequence $\bigl\{ A^{(n)} \bigr\}$ converges to
$A^{(\infty)}\in \Gamma(X,M(r)\otimes \Omega_X^{1,0})$
on any compact subset $K\subset X$.
\item \label{number;11.27.11}
We put $H^{(n)}:=H(h^{(n)},\vecw^{(n)})\in
 C^{\infty}(X,\nbigh(r))$.
On any compact subset $K\subset X$,
$H^{(n)}$ and $H^{(n)\,-1}$ are bounded
independently of $n$.
Namely we have a constant $C_K$ depending on $K$
such that $\big|H^{(n)}\big|<C_K$ and $\big|H^{(n)\,-1}\big|<C_K$.
\item
On any compact subset $K\subset X$,
the norms $|\theta^{(n)}|_{h^{(n)}}$ are bounded independently of $n$.
\end{enumerate}
\end{condition}

The element $\Theta^{(n)}\in C^{\infty}(X,M(r)\otimes\Omega^{1,0}_X)$
is determined by the following relation:
\[
 \theta^{(n)}\vecw=\vecw\cdot \Theta^{(n)}.
\]
The following lemma is easy to see.
\begin{lem}
Under the condition $2$,
the condition $3$ is equivalent to the following:
On any compact subset $K\subset X$,
$|\Theta^{(n)}|$ are bounded independently of $n$.
\hfill\qed
\end{lem}

\begin{prop} \label{prop;11.26.11}
Let $X$, $(E^{(n)},\theta^{(n)},h^{(n)})$
and $\vecw^{(n)}$ be as above.
Then we can take a subsequence $\{n_i\}$
of $\{n\}$ such that
the sequences $\big\{H^{(n_i)}\big\}$,
$\bigl\{H^{(n_i)\,-1}\bigr\}$
 and $\big\{\Theta^{(n_i)}\big\}$ 
are convergent on any compact subset $K$
in the $L^p_l$-sense.
Here $p$ is a sufficiently large real number,
and $l$ is an arbitrary large real number.
\end{prop}
\pf
We will use the following lemma without mention.
\begin{lem}
Let $K$ be a compact subset of $X$.
\begin{enumerate}
\item
Let $\{f_n\}$ be a bounded sequence in $L_l^p(X)$
such that all of $f_n$ vanish on $X-K$.

Assume that the sequence $\{\delbar(f_n)\}$
are bounded in $L^p_{l}(X)$.
Then $\{f_n\}$ is bounded in $L_{l+1}^p(X)$.

Similarly, when the sequence $\{\del(f_n)\}$ is bounded in $L_l^p(X)$,
then $\{f_n\}$ is bounded in $L_{l+1}^p(X)$.

\item
Let $\{f_n\}$ be a sequence of $L_l^p(X)$
satisfying that
any $f_n$ vanish on $X-K$, and that $f_n$ converges
to $g$ in $L_l^p(X)$.

Assume that $\{\delbar f_n\}$ converges to $\delbar g$ in $L_l^{p}(X)$.
Then $\{f_n\}$ converges to $g_n$ in $L_{l+1}^p(X)$.

Similarly when $\{\del f_n\}$ converges to
$\del g$ in $L_l^p(X)$,
then $\{f_n\}$ converges to $g_n$ in $L_{l+1}^p(X)$.
\end{enumerate}
\label{lem;11.26.12}
\end{lem}
\pf
Let see the first claim.
We only have to show the boundedness of $\{\del f_n\}$
in $L_l^p(X)$.
Let $D$ be a differential operator with constant coefficient
whose degree is less than $l$.
We have the following equality:
\[
 \int (\del D f_n, \del D f_n)=
 \int (\Delta D f_n, D f_n)=
 \int (\delbar D f_n,\delbar D f_n).
\]
Thus the $L_l^p$-boundedness of $\{\delbar f_n\}$
implies that the $L_l^p$-boundedness of $\{\del f_n\}$.

The second claim can be shown similarly.
Thus the proof of Lemma \ref{lem;11.26.12} is completed.
\hfill\qed

\vspace{.1in}
Let us return to the proof of Proposition \ref{prop;11.26.11}.
The elements $\Theta^{\dagger\,(n)}\in
 C^{\infty}(X,M(r)\otimes\Omega^{0,1}_X)$
are determined by the following condition:
\[
 \theta^{(n)\dagger}\vecw^{(n)}=
 \vecw^{(n)}\Theta^{(n)\dagger}.
\]
Then we have the following relations:
\begin{equation}\label{eq;11.26.18}
 \Theta^{(n)\,\dagger}=
 \overline{H}^{(n)\,-1}
\cdot {}^t\overline{\Theta}^{(n)}
\cdot \overline{H}^{(n)}.
\end{equation}
Hence we also have the boundedness of the sequence
$\{\Theta^{(n)\,\dagger}\}$
on any compact subset $K\subset X$.

Let $\nabla^{(n)}$ be the unitary connection
of $(\nbige^{1\,(n)},d^{\twoprime\,(n)},h^{(n)})$.
Then we have the relation
\[
 \nabla^{(n)}=
 \delbar_{E^{(n)}}+\theta^{\dagger\,(n)}
+\del_{E^{(n)}}-\theta^{(n)}
=\DD^{1\,(n)}-2\theta^{(n)}.
\]
We determine the element
$B^{(n)}\in C^{\infty}(X,M(r)\otimes\Omega^{1,0}_X)$
by the relation
$\nabla^{(n)}\vecw^{(n)}=\vecw^{(n)}\cdot B^{(n)}$.
Then we have the following equalities by definitions:
\begin{equation} \label{eq;11.27.14}
 B^{(n)}=A^{(n)}-2\Theta^{(n)},
\quad
\mbox{i.e.,}
\quad
 \Theta^{(n)}=\frac{1}{2}(A^{(n)}-B^{(n)}).
\end{equation}
On the other hand, we have the following relation:
\begin{equation} \label{eq;11.26.19}
 B^{(n)}=\bigl(H^{(n)}\bigr)^{-1} \del H^{(n)},
\quad
\mbox{i.e.,}
\quad
 \del H^{(n)}=H^{(n)}\cdot B^{(n)}.
\end{equation}

Let pick a compact subset $K\subset X$,
and let pick compact subsets $K_1$ and $K_2$ of $X$
such that 
$K$ is contained in the interior of $K_1$,
and that $K_1$ is contained in the interior of $K_2$.
We can pick an element $\varphi\in C^{\infty}(X,\real)$,
satisfying the following:
\[
 0\leq \varphi(x)\leq 1,\quad
 \varphi(x)=
 \left\{
 \begin{array}{ll}
 1 & (x\in K)\\
 0 & (x\not\in K_2)
 \end{array}
 \right.
\]
Let take a sufficiently large number $p$ 
satisfying the following:
\begin{quote}
 Let $l$ be any real number larger than $1$.
 For any elements $f,g\in L_l^p(X)$,
 the product $f\cdot g$ is contained in $L_l^p(X)$.
\end{quote}
Such $p$ can be taken, depending only on $\dim X$.

We need some equalities.
\begin{lem}
We have the following equality:
\begin{equation} \label{eq;11.26.20}
 \delbar(\varphi^m\Theta^{(n)})=
 \varphi^m\cdot
 \bigl[\Theta^{(n)},\Theta^{(n)\dagger}
 \bigr]
+m\cdot \delbar\varphi\cdot\varphi^{m-1}\cdot\Theta^{(n)}.
\end{equation}
\end{lem}
\pf
We have $\delbar_{E^{(n)}}(\theta^{(n)})=0$,
in other words,
$\bigl(d^{\twoprime\,1\,(n)}-\theta^{(n)\,\dagger}\bigr)
(\theta^{(n)})=0$.
(Recall that we have $d^{\twoprime\,1}=\delbar_E+\theta^{\dagger}$,
by definition for a harmonic bundle, in general.)
It is reworded as follows:
\[
 \delbar \Theta^{(n)}=
 \bigl[
 \Theta^{(n)\,\dagger},\Theta^{(n)}
 \bigr].
\]
The equality (\ref{eq;11.26.20}) follows immediately.
\hfill\qed

\begin{lem}
We have the following equality:
\begin{equation}\label{eq;11.26.21}
\del\bigl(
 \varphi^m\cdot H^{(n)}
 \bigr)=
 m\cdot(\del\varphi)\cdot\varphi^{m-1}\cdot H^{(n)}
+\varphi^m\cdot H^{(n)}
  \Bigl(A^{(n)}-2\Theta^{(n)}
  \Bigr).
\end{equation}
\end{lem}
\pf
The equality (\ref{eq;11.26.21})
follows from (\ref{eq;11.26.19})
and the equality $B^{(n)}=A^{(n)}-2\Theta^{(n)}$.
\hfill\qed

\begin{lem}
For non-negative integers $k$ and $l$,
we have the following equality:
\begin{equation} \label{eq;11.26.22}
 \del \bigl(
 \varphi^{k+l}\cdot H^{(n)\,-1}
 \bigr)=
 (k+2l)\cdot \del\varphi\cdot\varphi^{k+l-1}
 \cdot H^{(n)\,-1}
-\varphi^k\cdot H^{(n)\,-1}\cdot
 \del\bigl(\varphi^l\cdot H\bigr)
 \cdot H^{(n)-1}.
\end{equation}
\end{lem}
\pf
We have the following direct calculation.
We omit $(n)$ for simplicity, and we put $k+l=m$
\begin{multline}
\del(\varphi^{m}H^{-1})=
m\cdot\del\varphi\cdot\varphi^{m-1}\cdot H^{-1}
-\varphi^{m}\cdot H^{-1}\cdot\del(H)\cdot H^{-1}
\\
=m\cdot\del\varphi\cdot \varphi^{m-1}\cdot H^{-1}
-\varphi^k\cdot H^{-1}
 \Bigl(\del(\varphi^l\cdot H)-\del\varphi^l\cdot H\Bigr)H^{-1}\\
=(m+l)\cdot\del\varphi\cdot\varphi^{m-1}\cdot H^{-1}
-\varphi^k\cdot H^{-1}\del\bigl(\varphi^l\cdot H\bigr)H^{-1}.
\end{multline}
Thus we are done.
\hfill\qed

\begin{lem}
We have the following lemma:
\begin{equation} \label{eq;11.26.23}
 \varphi^{k+l+m}\Theta^{(n)\,\dagger}=
 \bigl(\varphi^k\overline{H}^{(n)\,-1}
 \bigr)
\cdot
 \bigl(
 \varphi^l\cdot{}^t\overline{\Theta}^{(n)}
 \bigr)
\cdot
 \bigl(
 \varphi^m\cdot \overline{H}^{(n)}
 \bigr).
\end{equation}
\end{lem}
\pf
The equality follows immediately from (\ref{eq;11.26.18}).
\hfill\qed

\vspace{.1in}
Let us return to the proof of Proposition \ref{prop;11.26.11}.
We use the standard boot strapping.
We divide the argument into some lemmas.

\begin{lem} \mbox{{}} \label{lem;11.26.30}
\begin{itemize}
\item
The sequences $\{\varphi\cdot\Theta^{(n)}\}$
and $\{\varphi^3\cdot\Theta^{\dagger\,(n)}\}$
are bounded in $L_1^p(X)$.
\item
The sequences $\{\varphi^2\cdot H^{(n)}\}$
and $\{\varphi^4\cdot H^{(n)\,-1}\}$
are bounded in $L_2^p(X)$.
\end{itemize}
\end{lem}
\pf
It is clear that $\varphi\cdot\Theta^{(n)}$
are bounded in $L^p(X)$ independently of $n$.
Then $\delbar(\varphi\cdot\Theta^{(n)})=
 \varphi\bigl[\Theta^{(n)\,\dagger},\Theta^{(n)}\bigr]
 +\delbar\varphi\cdot \Theta^{(n)}$
are bounded independently of $n$.
Thus $\varphi\cdot\Theta^{(n)}$
are bounded in $L^p_1(X)$.

It is easy to see that $\varphi\cdot H^{(n)}$ are bounded in $L^p(X)$.
Thus
$\del(\varphi\cdot H^{(n)})=\del\varphi\cdot H^{(n)}+
 \varphi\cdot H^{(n)}\bigl(A^{(n)}-2\Theta^{(n)}\bigr)$
are bounded in $L^p(X)$.
It implies that $\varphi\cdot H^{(n)}$ are bounded in $L^p_1(X)$.
Then
$\del(\varphi^2\cdot H^{(n)})=
  2\del\varphi\cdot H^{(n)}+
 \varphi\cdot H^{(n)}\cdot
  \bigl(\varphi \cdot A^{(n)}-2\varphi\cdot\Theta^{(n)}\bigr)$
are bounded in $L^p_1(X)$.
Thus $\varphi^2\cdot H^{(n)}$ are bounded in $L^p_2(X)$.

It is easy to see that $\varphi\cdot H^{(n)-1}$ are bounded in $L^p(X)$.
Then
$\del(\varphi\cdot H^{(n)\,-1})=
 \del\varphi\cdot H^{(n)\,-1}-
 \varphi\cdot H^{(n)\,-1}\del H^{(n)}\cdot H^{(n)\,-1}
 $ are bounded in $L^p(X)$.
Thus $\varphi\cdot H^{(n)\,-1}$ are bounded in $L_1^p(X)$.
We have
\[
 \del(\varphi^4\cdot H^{(n)\,-1})
=6\cdot\del\varphi\cdot \varphi^3\cdot H^{-1}
-(\varphi\cdot H^{(n)\,-1})
\cdot
 \del\bigl(\varphi^2\cdot H^{(n)}\bigr)
\cdot (\varphi\cdot H^{(n)-1}).
\]
Thus $\del(\varphi^4 H^{(n)\,-1})$ are bounded in $L_1^p(X)$,
and hence $\varphi^4 H^{(n)\,-1}$ are bounded in $L_2^p(X)$.

Clearly,
$\varphi^3\cdot\Theta^{\dagger\,(n)}=
 \bigl(\varphi \cdot\overline{H}^{(n)\,-1}\bigr)\cdot
 \bigl(\varphi\cdot {}^t\overline{\Theta}^{(n)}\bigr)
 \cdot
 \bigl(\varphi \cdot \overline{H}^{(n)}\bigr)$
are bounded in $L_1^p(X)$.
\hfill\qed

\begin{lem} \mbox{{}}
\begin{itemize}
\item
 The sequences $\{\varphi^4\cdot\Theta^{(n)}\}$ and
 $\{\varphi^{24}\Theta^{\dagger\,(n)}\}$
 are bounded in $L_2^p(X)$.
\item
 The sequences $\{\varphi^6\cdot H^{(n)}\}$ and
 $\{\varphi^{14}\cdot H^{(n)\,-1}\}$ are bounded in $L_3^p(X)$.
\end{itemize}
\end{lem}
\pf
By using the formulas (\ref{eq;11.26.20})
and the boundedness in Lemma \ref{lem;11.26.30},
we obtain the boundedness
of $\bardel(\varphi^4\Theta^{(n)})$ in $L_1^p(X)$.
Hence we obtain the boundedness
of $\varphi^4\Theta^{(n)}$ in $L_2^p(X)$.

Similarly, by using the formulas
(\ref{eq;11.26.21}), (\ref{eq;11.26.22})
and the boundedness in Lemma \ref{lem;11.26.30},
we obtain the boundedness
of $\del(\varphi^6 H^{(n)})$
and $\del(\varphi^{14}\cdot H^{(n)\,-1})$
in $L_3^p(X)$.

Lastly we obtain the boundedness
of $\varphi^{10}\cdot \Theta^{(n)\,\dagger}$
in $L_2^p(X)$.
And thus we obtain the boundedness
of $\varphi^{24}\cdot \Theta^{(n)\,\dagger}$
in $L_2^p(X)$.
\hfill\qed

\vspace{.1in}

Let take a subsequence $\{n_i\}$ of $\{n\}$
satisfying the following:
\begin{itemize}
\item
The sequences
$\{\varphi^4\cdot \Theta^{(n_i)}\}$
and 
$\{\varphi^{24}\cdot\Theta^{(n_i)\,\dagger}\}$
are convergent in $L_1^p(X)$.
The limit is denoted by $a$ and $d$ respectively.
\item
The sequences 
$\{\varphi^6\cdot H^{(n_i)}\}$
and $\{\varphi^{14}\cdot H^{(n_i)\,-1}\}$
are convergent in $L_2^p(X)$.
The limit is denoted by $b$ and $c$ respectively.
\end{itemize}
Such subsequence $(n_i)$ can be taken 
due to the Sobolev's embedding theorem.

\vspace{.1in}
We take the sequences of natural numbers
$\{\alpha_l\}$,
$\{\beta_l\}$,
$\{\gamma_l\}$ and
$\{\delta_l\}$ as follows:
\begin{itemize}
\item $\alpha_1=4$, $\beta_1=6$, $\gamma_1=14$ and $\delta_1=24$.
\item
 We have the relations:
 \[
 \begin{array}{ll}
  \alpha_l=\alpha_{l-1}+\delta_{l-1},
 & \beta_l=\beta_{l-1}+\alpha_l,\\
 \gamma_l=2\gamma_{l-1}+\beta_l,
 &\delta_l=\alpha_l+\beta_l+\gamma_l. 
\end{array}
 \]
Such relations determines the sequence of the numbers uniquely.
\end{itemize}

\begin{lem}\mbox{{}}
\begin{itemize}
\item
 The sequences $\bigl\{\varphi^{\alpha_l}\cdot\Theta^{(n_i)}\bigr\}$
and $\bigl\{\varphi^{\delta_l}\cdot \Theta^{(n_i)\,\dagger}\bigr\}$
are convergent in $L_l^p(X)$.
The limits are
$\varphi^{\alpha_l-2}\cdot a$
and $\varphi^{\delta_l-24}d$ respectively.
\item
 The sequences $\bigl\{\varphi^{\beta_l}\cdot H^{(n_i)} \bigr\}$
and $\bigl\{\varphi^{\gamma_l}\cdot H^{(n_i)\,-1}\bigr\}$
are convergent in $L_{l+1}^p(X)$.
The limits are $\varphi^{\beta_l-6}\cdot b$
and $\varphi^{\gamma_l-14}\cdot c$
respectively.
\end{itemize}
\end{lem}
\pf
We use the induction on $l$.
The assertions for $l=1$ hold because of our choice.
We assume that the assertions for $l-1$ hold,
and we will prove the assertions for $l$.

Due to the formula (\ref{eq;11.26.20}),
we have the following:
\[
 \bardel(\varphi^{\alpha_{l-1}+\delta_{l-1}}\cdot \Theta^{(n_i)})
=
 \Bigl[
 \varphi^{\alpha_{l-1}}\Theta^{(n_i)},
 \varphi^{\delta_{l-1}}\Theta^{(n_i)\,\dagger}
 \Bigr]+
(\alpha_{l-1}+\delta_{l-1})
\cdot\delbar\varphi\cdot\varphi^{\delta_{l-1}-1}
\cdot\bigl(
\varphi^{\alpha_{l-1}}\cdot \Theta^{(n_i)}\bigr).
\]
Thus we obtain the convergence
of $\delbar(\varphi^{\alpha_l}\cdot\Theta^{(n_i)})$
in $L_{l-1}^p(X)$.
Thus we obtain the convergence of $\{\varphi^{\alpha_l}\Theta^{(n_i)}\}$
in $L_l^p(X)$.

We have the following:
\[
 \del\bigl(
 \varphi^{\alpha_l+\beta_{l-1}}H^{(n_i)}
 \bigr)=
 (\alpha_l+\beta_{l-1})\cdot
 \del(\varphi)\cdot
 \varphi^{\alpha_l+\beta_{l-1}-1}\cdot H^{(n_i)}
+\bigl(\varphi^{\beta_{l-1}}H^{(n_i)}\bigr)
\cdot
 \bigl(
 \varphi^{\alpha_l}\cdot A^{(n_i)}-\varphi^{\alpha_l}\cdot\Theta^{(n_i)}
 \bigr)
\]
The convergence of the right hand side is in $L_l^p(X)$.
Thus we obtain the convergence of
$\{\varphi^{\beta_l}\cdot H^{(n_i)}\}$
in $L_{l+1}^p(X)$.

We have the following:
\begin{multline}
 \del\bigl(\varphi^{2\gamma_{l-1}+\beta_l}H^{(n_i)\,-1}\bigr)
=2(\gamma_{l-1}+\beta_{l})\cdot
 \del\varphi\cdot\varphi^{2\gamma_{l-1}+\beta_l-1}
 \cdot\bigl(
 \varphi^{\gamma_{l-1}}\cdot H^{(n_i)-1}\bigr)\\
-\varphi^{\gamma_{l-1}}\cdot H^{(n_i)\,-1}
\cdot\del\bigl(\varphi^{\beta_l}\cdot H^{(n_i)}\bigr)
\cdot \varphi^{\gamma_{l-1}}\cdot H^{(n_i)\,-1}.
\end{multline}
The convergence of the right hand side is in $L_l^p(X)$,
and thus the convergence
of $\{\varphi^{\gamma_l}\cdot H^{(n_i)\,-1}\}$ is
in $L_{l+1}^p(X)$.

Lastly we have the following:
\[
 \varphi^{\alpha_l+\beta_l+\gamma_l}\cdot \Theta^{(n)\,\dagger}
=\bigl(
 \varphi^{\gamma_l}\cdot \overline{H}^{(n_i)\,-1}
 \bigr)
\cdot
 \bigl(
 \varphi^{\alpha_l}\cdot {}^t\overline{\Theta}^{(n)}
 \bigr)
\cdot
 \bigl(
 \varphi^{\beta_l}\cdot \overline{H}^{(n_i)}
 \bigr).
\]
Thus the induction can proceed.
\hfill\qed

By our choice of $\varphi$,
we have $\varphi=1$ on $K$.
Thus we obtain the convergence of
$\{\Theta^{(n_i)}\}$,
$\{H^{(n_i)}\}$,
$\{H^{(n_i)\,-1}\}$
and $\{\Theta^{\dagger\,(n_i)}\}$
on $K$.
By using a standard diagonal argument,
we can take a sequence $\{n_i\}$
such that 
the sequences
$\{\Theta^{(n_i)}\}$,
$\{H^{(n_i)}\}$,
$\{H^{(n_i)\,-1}\}$
and $\{\Theta^{\dagger\,(n_i)}\}$
converge on any compact subset $K\subset X$.
Thus the proof of Proposition \ref{prop;11.26.11}
is completed.
\hfill\qed

\vspace{.1in}

Let take a holomorphic bundle
$F=\bigoplus \nbigo_X\cdot e_i$
with the frame $\vece=(e_i)$.
The frames $\vecw^{(n)}$ and $\vece$
induce the holomorphic isomorphism
$\Phi_n:(\nbige^{(n)\,1},\vecw^{(n)})
 \lrarr
 (F,\vece)$.
Then we obtain the sequence of the metrics
$\{h^{(n)}\}$,
(non-holomorphic) Higgs fields
$\{\theta^{(n)}\}$,
the conjugates
$\{\theta^{(n)\,\dagger}\}$ and
the holomorphic structures
$\{\delbar^{(n)}:=d''_F-\theta^{(n)\,\dagger}\}$
on $F$.
Let take the subsequence $\{n_i\}$
as in Proposition \ref{prop;11.26.11}.
Then the corresponding sequences converge,
and thus we obtain the limits
$h^{(\infty)}$,
$\theta^{(\infty)}$,
$\theta^{(\infty)\,\dagger}$,
$\delbar^{(\infty)}$
on $F$.

\begin{df}
The tuple of limits
$(F,\delbar^{(\infty)},\theta^{(\infty)},h^{(\infty)})$
will be called the limiting harmonic bundle
of the sequence $\{(E^{(n)},\theta^{(n)},h^{(n)})\}$,
although we do not have the uniqueness of the limit,
in general.
\hfill\qed
\end{df}

\begin{rem}
In this subsubsection,
we discussed the convergence at $\lambda=1$.
Clearly, the same argument works if $\lambda\neq 0$.
The key equality is {\rm (\ref{eq;11.27.14})}.

In the case $\lambda=0$, the conditions $1$ and $2$ in 
{\rm Condition \ref{condition;11.27.5}}
imply the condition $3$.
Hence the convergence of $\theta$ is obtained from the definition.
But we do not have the relation of the metric connection
and $\theta$ at $\lambda=0$.
Thus the argument to treat the metric breaks.
\hfill\qed
\end{rem}


\subsubsection{The dependence on the holomorphic frames $\vecw^{(n)}$}
\label{subsubsection;11.27.17}

Let $\dotvecw^{(n)}$ be holomorphic frames
of $\nbige^{(n)\,1}$ satisfying 
Condition \ref{condition;11.27.5}.
We take a holomorphic vector bundle
$\dotF=\bigoplus_{i=1}^r\nbigo_X\cdot \dote_i$.
The frames $\dotvecw^{(n)}$ and $\dotvece$
induces the isomorphism
$\dotPhi_n:(\nbige^{(n)\,1},\dotvecw^{(n)})\lrarr (\dotF,\dotvece)$.
Then we obtain the sequences
$\{\doth^{(n)}\}$, $\{\dottheta^{(n)}\}$,
$\{\dottheta^{(n)\,\dagger}\}$,
and $\{\dot{\delbar}^{(n)}\}$.
We take a subsequence $\{\dotn_i\}$ of $\{n\}$
such that the corresponding subsequences are convergent
to $\doth^{(\infty)}$, $\dottheta^{(\infty)}$,
$\dottheta^{(\infty)\,\dagger}$
and $\dot\delbar^{(\infty)}$ respectively.
We can assume that $\{\dotn_i\}$ is a subsequence of $\{n_i\}$.

\begin{lem} \label{lem;11.27.13}
We can take a subsequence $\{\dotn_i\}$ such that
the limiting harmonic bundles
$(\dotF,\dotdelbar^{(\infty)},\dottheta^{(\infty)},\doth^{(\infty)})$
is isomorphic to
$(F,\delbar^{(\infty)},\theta^{(\infty)},h^{(\infty)})$.
\end{lem}
\pf
We have the sequence of holomorphic frames
$\Phi_n(\dotvecw^{(n)})$ of $F$.
Due to the condition \ref{number;11.27.11}
in Condition \ref{condition;11.27.5}
for $\vecw^{(n)}$ and $\dotvecw^{(n)}$,
we can take a subsequence $\{n_{i_j}\}$ of $\{n_i\}$
such that $\Phi_{n_{i_j}}(\dotvecw^{(n_{i_j})})$ converges
to a holomorphic frame of $F$.
We can assume that $\{n_{i_j}\}$ is a subsequence
of $\{\dotn_i\}$.
We can replace $\{\dotn_i\}$ by $\{n_{i_j}\}$,
and thus we can assume that
$\{\Phi_{\dotn_i}(\dotvecw^{(\dotn_i)})\}$
are convergent.

The frames $\vecw^{(n)}$ and $\dotvecw^{(n)}$
induce the isomorphism $G^{(n)}$
of the harmonic bundles
$(F,\delbar^{(n)},\theta^{(n)},h^{(n)})$
and
$(\dotF,\dotdelbar^{(n)},\dottheta^{(n)},\doth^{(n)})$,
by our construction.
For the subsequence $\{\dotn_i\}$ above,
the sequence $G^{(\dotn_i)}$ converges to a holomorphic
isomorphism $G^{(\infty)}$.
Clearly, $G^{(\infty)}$ gives an isomorphism
of limiting harmonic bundles
$(F,\delbar^{(\infty)},\theta^{(\infty)},h^{(\infty)})$
and
$(\dotF,\dotdelbar^{(\infty)},\dottheta^{(\infty)},\doth^{(\infty)})$.
\hfill\qed

\vspace{.1in}
Lemma \ref{lem;11.27.13} is reworded as follows.
\begin{cor} 
Let $\vecw^{(n)}$ and $\dotvecw^{(n)}$ be two holomorphic
frames of $\nbige^{(n)\,1}$
satisfying Condition {\rm\ref{condition;11.27.5}}.
Then we can take a subsequence $\{n_i\}$ of $\{n\}$
such that we have the natural isomorphism
between the limiting harmonic bundles
for the frames $\vecw^{(n)}$ and $\dotvecw^{(n)}$.
\hfill\qed
\end{cor}


\subsubsection{Convergence of a sequence of the frames 
of the deformed holomorphic bundles}

Let continue to use the notation in the previous subsubsection.
From harmonic bundles $(F,\delbar^{(n)},h^{(n)},\theta^{(n)})$,
we obtain the deformed holomorphic bundles
$\nbigf^{(n)}=(p_{\lambda}^{-1}F,d^{\twoprime\,(n)})$
over $\nbigx$.
Note that $\nbigf^{(n)\,1}$ over $\nbigx^1$ is same as
$(F,d''_F)$ by our construction.
Let take a subsequence $\{n_i\}$ as above.
Since we have
$d^{\twoprime\,(n)}
 =\delbar^{(n)}+\lambda\cdot\theta^{(n)\,\dagger}+\delbar_{\lambda}
 =d''_F+(\lambda-1)\cdot\theta^{(n)\,\dagger}+\delbar_{\lambda}$,
and since the sequence of
the Higgs  fields $\{\theta^{(n_i)}\}$
and the metrics $\{h^{(n_i)}\}$ are convergent,
the sequence $\{d^{\twoprime\,(n_i)}\}$ converges to 
$d^{\twoprime\,(\infty)}
=d''_F+(\lambda-1)\theta^{(\infty)\,\dagger}
+\delbar_{\lambda}$.

Let $\{f^{(n)}\}$ be a sequence of holomorphic sections
of $\nbigf^{(n)}$ over $\nbigx$.
Namely $f^{(n)}$ are elements of
$C^{\infty}(\nbigx,p_{\lambda}^{-1}F)$
satisfying $d^{\twoprime\,(n)}(f^{(n)})=0$.
Assume the following:
\begin{itemize}
\item For any compact subset $K\subset \nbigx$,
we have a constant $C_K>0$ such that
$|h^{(n)}(f^{(n)},f^{(n)})|<C_K$ for any $n$.
\end{itemize}
This is equivalent to the following:
\begin{itemize}
\item
We have the $C^{\infty}$-functions $f_i^{(n)}$
determined by $f^{(n)}=\sum f^{(n)}_i\cdot p_{\lambda}^{-1}e_i$.
For any compact subset $K\subset \nbigx$,
we have a constant $C_K>0$ such that
$|f_i^{(n)}|<C_K$.
\end{itemize}

\begin{lem} \label{lem;11.27.1}
Let $f^{(n)}$ be as above.
We can take a subsequence $\{n_{i_j}\}$ of $\{n_i\}$
such that the sequence $\{f^{(n_{i_j})}\}$ is convergent
to a holomorphic section of $\nbigf^{(\infty)}$
on any compact subset $K$
in the $L^p_l$-sense.
Here $p$ is a sufficiently large real number,
and $l$ is an arbitrary large real number.
\end{lem}
\pf
We can use an argument similar to the proof
of Proposition \ref{prop;11.26.11}.
\hfill\qed

\vspace{.1in}

Let $\vecv^{(n)}=(v^{(n)}_1,\ldots,v^{(n)}_r)$ be frames of 
$\nbigf^{(n)}$ over $\nbigx$.
Assume that we have positive constants $C_{K,1}$ and $C_{K,2}$
for any compact subset satisfying the following:
\[
 h^{(n)}\bigl(v_i^{(n)},v_i^{(n)}\bigr)<C_{K,1},
\quad
 h^{(n)}\bigl(\volume(\vecv^{(n)}),\volume(\vecv^{(n)})\bigr)
 >C_{K,2}.
\]

\begin{lem} \label{lem;12.4.20}
We have a subsequence $\{n_{i_j}\}$
of $\{n_i\}$
such that 
the sequence $\{\vecv^{(n)}\}$
converges to a holomorphic frame $\{\vecv^{(\infty)}\}$
of $\nbigf^{(\infty)}$
on any compact subset $K$
in the $L^p_l$-sense.
\end{lem}
\pf
By the boundedness of $h^{(n)}(v^{(n)}_i,v^{(n)}_i)$
and Lemma \ref{lem;11.27.1},
we have a subsequence $\{n_{i_j}\}$
such that the sequence of sections
$\{v_i^{(n_{i_j})}\}$  converges to
a holomorphic section $v_i^{(\infty)}$.
We only have to show that $\vecv^{(\infty)}=(v^{(\infty)}_i)$
gives a holomorphic frame.

On any compact subset $K$,
we have the estimate
 $h^{(n)}\bigl(\volume(\vecv^{(n)}),\volume(\vecv^{(n)})\bigr)
 >C_{K,2}$.
Thus we obtain
$h^{(\infty)}
  \bigl(\volume(\vecv^{(\infty)}),\volume(\vecv^{(\infty)})\bigr)
 >C_{K,2}$.
Hence $\vecv^{(\infty)}$ gives a frame.
\hfill\qed

\begin{cor} \label{cor;11.27.2}
Let $\vecv^{(n)}$ be a holomorphic frame of $\nbigf^{(n)}$
satisfying the following:
\begin{itemize}
\item For any compact subset $K$,
the hermitian matrices
$H(h^{(n)},\vecv^{(n)})$ and the inverse $H(h^{(n)},\vecv^{(n)})^{-1}$
are bounded independently of $n$.
\end{itemize}
Then we can take a subsequence $\{n_{i_j}\}$
of $\{n_i\}$
such that the sequence $\{\vecv^{(n_{i_j})}\}$
converges to a holomorphic frame of $\nbigf^{(\infty)}$.
\hfill\qed
\end{cor}

Let $\vecv^{(n)}$ be a holomorphic frame of $\nbigf^{(n)}$.
Then we have the transformation matrices
$B^{(n)}\in C^{\infty}(\nbigx,M(r))$
determined by the relation
$\vecv^{(n)}=p_{\lambda}^{-1}\vece\cdot B^{(n)}$.
Corollary \ref{cor;11.27.2} can be reworded as follows:
\begin{cor}
Let $\vecv^{(n)}$ be a holomorphic frame of $\nbigf^{(n)}$.
Assume that the transformation matrices $B^{(n)}$
and their inverses $B^{(n)\,-1}$ are bounded
independently of $n$,
on any compact subset $K$.
Then we can take a subsequence $\{n_{i_j}\}$ of $\{n_i\}$ such that
$\{\vecv^{(n_{i_j})}\}$
converges to a holomorphic frame of $\nbigf^{(\infty)}$.
\hfill\qed
\end{cor}


\subsubsection{Convergence at any $\lambda$ and $\mu$}

We reword the result as follows:
Let $(E^{(n)},\delbar_{E^{(n)}},h^{(n)},\theta^{(n)})$
be a sequence of harmonic bundles over $X$.
We have the deformed holomorphic bundles
$(\nbige^{(n)},d^{\twoprime\,(n)},\DD^{(n)},h^{(n)})$
over $\nbigx$.
Let $U_1$ be an open subset of $\cnum_{\lambda}$.
We assume that $U_1$ contains $1\in\cnum_{\lambda}$
for simplicity.
We put $\nbigu:=U_1\times X\subset \nbigx$.

Let $\vecv^{(n)}$ be a holomorphic frame of $\nbige^{(n)}$
over $\nbigu$.
\begin{condition}\mbox{{}} \label{condition;11.27.16}
\begin{enumerate}
\item We have the $\lambda$-connection form
$\nbiga^{(n)}\in \Gamma(\nbigu,M(r)\otimes p_{\lambda}^{\ast}\Omega_{X}^{1,0})$,
such that
$\DD(\vecv^{(n)})=\vecv^{(n)}\cdot\nbiga^{(n)}$.
Then the sequence $\{\nbiga^{(n)}\}$
converges to
$\nbiga^{(\infty)}\in
  \Gamma(\nbigu,M(r)\otimes p_{\lambda}^{\ast}\Omega_X^{1,0})$
on any compact subset $K\subset \nbigu$.
\item
We put $H^{(n)}:=H(h^{(n)},\vecv^{(n)})\in C^{\infty}(\nbigu,\nbigh(r))$.
Then $H^{(n)}$ and $H^{(n)\,-1}$
are bounded independently of $n$,
on any compact subset $K\subset \nbigu$.
\item
On any compact subset $K$, we have a positive constant $C_K$
such that $|\theta^{(n)}|_{h^{(n)}}<C_K$.
\hfill\qed
\end{enumerate}

\end{condition}

The element
$\Theta^{(n)}\in
 C^{\infty}(\nbigu,
    M(r)\otimes p_{\lambda}^{\ast}\Omega^{1,0}_X)$
is determined by the relation
$\theta^{(n)}(\vecv^{(n)})=\vecv^{(n)}\cdot\Theta^{(n)}$.
\begin{prop} \label{prop;11.27.15}
We can take a subsequence $\{n_i\}$
such that the sequences $\{H^{(n_i)}\}$, $\{H^{(n_i)\,-1}\}$
and $\{\Theta^{(n_i)}\}$ are convergent
in $L^p_l$ on any compact subset $K$.
\end{prop}
\pf
When we take the restriction of $\vecv^{(n)}$
to $\lambda=1$,
then Condition \ref{condition;11.27.5} is satisfied.
Thus we can take $\{n_i\}$ as in Proposition \ref{prop;11.26.11}.
Then we take a subsequence $\{n_{i_j}\}$
as in Corollary \ref{cor;11.27.2}.
By construction,
it is clear that
the sequences
$\{H^{(n_{i_j})}\}$, $\{H^{(n_{i_j})\,-1}\}$
and $\{\Theta^{(n_{i_j})}\}$
are convergent.
We only have to replace $\{n_i\}$ with $\{n_{i_j}\}$.
\hfill\qed

Let pick an element $\lambda\in U$.
We denote the restriction of $\vecv^{(n)}$
to $\{\lambda\}\times X$
by $\vecv^{(n)}_{|\lambda}$.
We take a holomorphic vector bundle
$F=\bigoplus_{i=1}^r\nbigo_X\cdot e_i$.
By the frames $\vecv^{(n)}_{|\lambda}$
and $\vece=(e_i)$,
we have the holomorphic isomorphism
$\Phi_{n,\lambda}:
(\nbige^{(n)\,\lambda},\vecv^{(n)}_{|\lambda})
 \lrarr (F,\vece)$.
Then we obtain the sequence of the metrics
$\{h^{(n)}\}$, the (non-holomorphic) Higgs fields
$\{\theta^{(n)}\}$,
the conjugates
$\{\theta^{(n)\,\dagger}\}$
and the holomorphic structures
$\{\delbar^{(n)}:=d''_F-\lambda\theta^{(n)\,\dagger}\}$.
Let take a subsequence $\{n_i\}$
as in Proposition \ref{prop;11.27.15}.
Then we obtain the limits of the sequences
$h^{(\infty)}$, $\theta^{(\infty)}$, 
$\theta^{(\infty)\,\dagger}$
and $\delbar^{(\infty)}$.
The tuple of limits
$(F,\delbar^{(\infty)},\theta^{(\infty)},h^{(\infty)})$
is a harmonic bundle,
called a limiting harmonic bundle.

By our construction,
we can see that
the limit is independent of a choice of $\lambda$,
in the sense that we have the canonical isomorphisms
between the limits,
once we fix an appropriate subsequence.

Let $\{\vecv^{(n)}\}$ and $\{\dotvecv^{(n)}\}$
are two sequences of holomorphic frames
of $\nbige^{(n)}_{|\nbigu}$ 
satisfying Condition \ref{condition;11.27.16}.
By an argument similar to that
in the subsubsection \ref{subsubsection;11.27.17},
we can show the following.
\begin{prop} 
We have a subsequence $\{n_i\}$
of $\{n\}$
satisfying the following:
\begin{itemize}
\item We have the limiting harmonic bundles
for $\{\vecv^{(n_i)}\}$ and $\{\dotvecv^{(n_i)}\}$.
\item The limits are naturally isomorphic.
\hfill\qed
\end{itemize}
\end{prop}

\vspace{.1in}
Let consider the case $\nbigu=\nbigx$ for simplicity.
We have the frame $\vecv^{(n)\,\dagger}$ of
$\nbige^{(n)\,\dagger}$ over $\nbigx^{\dagger}$.
We denote the conjugate deformed holomorphic bundle
of $(F,\delbar^{(n)},\theta^{(n)},h^{(n)})$
by $\nbigf^{(n)\,\dagger}$.
The morphism $\Phi_n$ induces the holomorphic isomorphism
$\nbige^{(n)\,\dagger}\lrarr \nbigf^{(n)\,\dagger}$,
which we also denote by $\Phi_n$.
Thus we obtain the sequence of holomorphic frames
$\{\Phi_n(\vecv^{(n)\,\dagger})\}$.
The following lemma can be easily seen from our construction.
\begin{lem}\label{lem;11.27.20}
The sequence $\{\Phi_{n_i}(\vecv^{(n_i)\,\dagger})\}$
is convergent.
\hfill\qed
\end{lem}

Let pick $\mu\in\cnum_{\mu}$.
We denote the restriction of $\vecv^{(n)\dagger}$
to $\{\mu\}\times X^{\dagger}$ by $\vecv^{(n)\,\dagger}_{|\mu}$.
We take a holomorphic bundle
$F^{\dagger}=\bigoplus_{i=1}^r\nbigo_{X^{\dagger}}\cdot e_i^{\dagger}$.
Due to the frames
$\vecv^{(n)\,\dagger}_{|\mu}$
and $\vece^{\dagger}=(e^{\dagger}_{i})$,
we have the holomorphic isomorphism
$(\nbige^{(n)\,\dagger\,\mu},\vecv^{(n)\,\dagger}_{|\mu})
\lrarr
 (F^{\dagger},\vece^{\dagger})$.
Thus we obtain the sequence of the hermitian metrics
$\{h^{(n)}\}$,
the Higgs fields $\{\theta^{\dagger\,(n)}\}$
on $X^{\dagger}$,
the conjugates
$\{\theta^{(n)}\}$,
and the holomorphic structure
$\{\del^{(n)}=d''_{F^{\dagger}}-\mu\cdot\theta^{(n)}\}$.
Due to Lemma \ref{lem;11.27.20},
the corresponding subsequences for $\{n_i\}$
converge.
The limit is independent of a choice of $\mu$
in the sense that we have the canonical isomorphism
between the limits,
once we fix an appropriate subsequence.
Moreover the limit obtained from the frames
$\bigl\{\vecv^{(n)\,\dagger}_{|\mu}\bigr\}$
are canonically isomorphic to the limit
obtained from the frames $\bigl\{\vecv^{(n)}_{|\lambda}\bigr\}$.

\section{Norm estimate in one dimensional case and
prolongation of the deformed holomorphic bundle}
\label{section;12.15.71}
\subsection{Prolongation by increasing orders}

Let $X$ be an $n$-dimensional complex manifold,
and $D=\bigcup_{i\in \soeji} D_i$ be a simple normal crossing divisor.

\label{subsubsection;12.9.60}

\begin{df}\label{df;12.9.1}
Let $P$ be a point of $X$,
and $D_{i_j}$ $(j=1,\ldots,l)$ be components of $D$
containing $P$.
An admissible coordinate around $P$ is the tuple
$(\nbigu,\varphi)$:
\begin{itemize}
\item $\nbigu$ is an open subset of $X$ containing $P$.
\item $\varphi$ is a holomorphic isomorphism
$\nbigu\lrarr \Delta^n=\{(z_1,\ldots,z_n)\,|\,|z_i|<1\}$
such that $\varphi(P)=(0,\ldots,0)$,
and $\varphi(D_{i_j})=\{z_j=0\}$
for any $j=1,\ldots,l$.
\hfill\qed
\end{itemize}
\end{df}

Let $(E,\delbar_E,h)$ be a holomorphic bundle with a hermitian metric
defined over $X-D$.
Let $\vecalpha=(\alpha_i\,|\,i\in \soeji)\in\real^{\soeji}$
be a tuple of real numbers.

\begin{df}
Let $U$ be an open subset of $X$,
and $s$ be an element of $\Gamma(U-D,E)$.
We say that the increasing order of $s$
is less than $\vecalpha$,
if the following holds:
\begin{itemize}
\item
Let $P$ be a point of $U$.
Let take an admissible coordinate $(\nbigu,\varphi)$ around $P$.
Let $\epsilon$ be any positive real number.
Then we have a positive constant $C$
such that the following inequality holds on $\nbigu$:
\[
 |s|_h\leq C\cdot \prod_{j=1}^l |z_j|^{\alpha_{i_j}-\epsilon}.
\]
\end{itemize}
In the case,
it is described as $\ord(s)\leq \vecalpha$.
\hfill\qed
\end{df}

Let pick a tuple $\vecalpha$.
Then the $\nbigo_X$-sheaf $E_{\leq\vecalpha}$
is defined as follows:
For any open subset $U\subset X$,
\[
 \Gamma(U,E_{\leq\vecalpha}):=
 \{s\in\Gamma(U-U\cap D,E)\,|\,\ord(s)\leq \vecalpha\}.
\]
We often use the notation $E_{\vecalpha}$
instead of $E_{\leq\vecalpha}$.
The sheaf $E_{\vecalpha}$ is called the prolongment of
$E$ by an increasing order $\vecalpha$.

\vspace{.1in}
\noindent
{\bf Notation ($\prolong{E}$)}\\
In this paper,
we are mainly interested in the case $\alpha_i=0$ for any $i$.
In that case we will use the notation
$\prolong{E}$.

\vspace{.1in}

We also define the $\nbigo_X$-sheaf $E_{<\vecalpha}$
as follows:
\[
 \Gamma(U,E_{<\vecalpha}):=
 \{
 s\in\Gamma(U-U\cap D,E)\,|\,
 \exists \epsilon>0, \mbox{ such that }
 \ord(s)\leq \vecalpha-\epsilon\cdot\vecdelta
 \}
\]
Here we put $\vecdelta=(1,\ldots,1)$.


The Poincar\'e metrics $g_1$ on $\Delta$ and
$g_0$ on $\Delta^{\ast}$
are given by the following formulas
up to some minor modification:
\[
 g_1=\frac{2dz\cdot d\bar{z}}{(1-|z|^2)^2},
\quad
 g_0=\frac{2dz\cdot d\zbar}{|z|^2(-\log|z|^2)^2}
\]
The associated Kahler forms $\omega_1$ and $\omega_0$
are as follows:
\[
 \omega_1=
 \frac{\sqrt{-1}dz\cdot d\bar{z}}{(1-|z|^2)^2},
\quad
 \omega_0=\frac{\sqrt{-1}dz\cdot d\zbar}{|z|^2(-\log|z|^2)^2}
\]
As the metric on $\Delta^{\ast\,l}\times \Delta^{n-l}$,
we have the metric and the Kahler form:
\[
 g_{\Poin}:=
 \sum_{j=1}^l \pi_j^{\ast}g_0
+\sum_{j=l+1}^n\pi_{j}^{\ast}g_1,
\quad
 \omega_{\Poin}
=\sum_{j=1}^l \pi_j^{\ast}\omega_0
+\sum_{j=l+1}^n\pi_{j}^{\ast}\omega_1,
\]

Let $P$ be a point of $X$
and $(\nbigu,\varphi)$ be an admissible coordinate around $P$.
By the isomorphism
$\varphi:\nbigu-D\simeq \Delta^{\ast\,l}\times\Delta^{n-l}$,
we take the Poincar\'e metric $g_{\Poin}$ on $\nbigu-D$.
The metric $h$ of $E$ and the metric $g_{\Poin}$ on $T(\nbigu-D)$
induce the metric $(\cdot,\cdot)_{h,g_{\Poin}}$
of $End(E)\otimes \Omega^{p,q}$
over $\nbigu-D$.

\begin{df} \label{df;11.28.5}
We say that $(E,\delbar_E,h)$ is acceptable at $P$,
if the following holds:
\begin{itemize}
\item
 Let $(\nbigu,\varphi)$ be an admissible coordinate around $P$.
  The norms of the curvature $R(h)$ with respect to
 the metric $(\cdot,\cdot)_{h,g_{\poin}}$
 is bounded over $\nbigu-D$.
\end{itemize}
When $(E,\delbar_E,h)$ is acceptable at any point $P$,
then we say that it is acceptable.
\hfill\qed
\end{df}

\subsection{Tameness and nilpotentness}

\label{subsection;12.13.1}

Let $(E,\delbar_E,\theta,h)$ be a harmonic bundle of rank $r$
defined over $X-D$.

\begin{df}
Let $P$ be any point of $X$, and $(\nbigu,\varphi)$
be an admissible coordinate around $P$.
On $\nbigu$, we have the description:
\[
 \theta=
 \sum_{j=1}^l f_j\cdot\frac{dz_j}{z_j}
+\sum_{j=l+1}^{n}g_j\cdot dz_j.
\]
\begin{description}
\item[(Tameness)]
 Let $t$ be a formal variable.
 We have the polynomials
 $\det(t-f_j)$ and $\det(t-g_j)$ of $t$,
 whose coefficients are holomorphic functions
 defined over $\nbigu-\bigcup_{j=1}^l D_{i_j}$.
 When the functions are extended to the holomorphic functions
 over $\nbigu$,
 the harmonic bundle is called tame at $P$.
\item[(Nilpotentness)]
Assume that the harmonic bundle is tame at $P$.
When $\det(t-f_j)|_{\nbigu\cap D_{i_j}}=t^r$,
then the harmonic bundle is called nilpotent at $P$.
\end{description}
When $(E,\delbar_E,h,\theta)$ is a tame nilpotent at any point $P\in X$,
then it is called a tame nilpotent harmonic bundle.
\hfill\qed
\end{df}

For the tame nilpotent harmonic bundle,
we have the following estimate.
Note that the proposition is essentially contained in
Theorem 1 of \cite{s2}.
There, Simpson gave the proof when the residue is not nilpotent,
which is more difficult than the nilpotent case.
For completeness, here,
we give a proof of the nilpotent case,
i.e., the easier case,
by following the ideas of Simpson and Ahlfors (\cite{a}).
\begin{prop}[Simpson]\label{prop;1.25.1}.
Let $(E,\delbar_E,h,\theta)$ be a tame nilpotent harmonic bundle.
Let $P$, $f_j$, $g_j$ be as above.
We put $y_j=-\log|z_j|^2$ for $j=1,\ldots,l$.
Then there exists a positive constant $C>0$
satisfying the following:
\[
 \begin{array}{ll}
 |f_j|_h\leq C\cdot y_j^{-1}, & (j=1,\ldots,l)\\
\mbox{{}}\\
 |g_j|_h\leq C, & (j=l+1,\ldots,n).
 \end{array}
\]
\end{prop}
\pf
From the beginning,
we can assume that $\nbigu=\Delta^n$,
and we can assume that $D=\bigcup_{j=1}^l D_j$
for $l\leq n$,
where we put $D_j:=\{z_j=0\}$.
We only see that there exists a positive constant $C$
satisfying
$|f_j|\leq C\cdot y_j^{-1}$
independently of $(z_1,\ldots,z_n)$.
The estimate for $g_j$ can be shown similarly.
We can assume that $j=1$.
We denote $z_1$ by $z$, and
we put $f_0=z^{-1}\cdot f_1$.
For any positive number $r$, we put 
$\psi(r):=(-r\cdot\log|r|)^{-1}$.
Let $\pi_1$ denote the projection $\Delta^n\lrarr \Delta^{n-1}$,
omitting the first component.

We denote the dual of $f_0$ with respect the metric $h$
by $f_0^{\dagger}$.

First we note the following.
\begin{lem}
There exists positive constants $C'>0$ and $\epsilon>0$
satisfying the following:
\begin{itemize}
\item
 Let $a$ be one of the eigenvalues of
 the endomorphism $f_{0}(P)$ of the fiber $E_{|P}$.
Then the following inequality holds:
\begin{equation}
 |a|\leq C'\cdot |z(P)|^{-1+\epsilon}.
\end{equation}
\end{itemize}
\end{lem}
\pf
Since $\det(t-z\cdot f_0)$ is holomorphic over $\Delta^n$,
and since $\det(t-z\cdot f_0)|_{z=0}=t^n$,
the eigenvalues of $(z\cdot f_0)(P)$ is dominated by 
$C'\cdot |z(P)|^{\epsilon}$ for some $C'>0$ and $\epsilon>0$.
Then we obtain the inequality desired.
\hfill\qed

\vspace{.1in}
For an element $p\in \Delta^{\ast\,n-1}$,
the tuple $(E_{|\pi_1^{-1}(P)},f_0\cdot dz,h)$
gives a harmonic bundle.
Thus we obtain the following inequality due to Simpson
(See \cite{s1} and \cite{s2}):
\[
 \Delta \log|f_0|_h^2
 \leq
 -\frac{|[f_0,f_0^{\dagger}]|_h^2}{|f_0|_h^2}.
\]
Here $\Delta$ is the operator $-(\del_s^2+\del_t^2)$ for the
real coordinate $z=s+\sqrt{-1}t$.
We put $|f_0|_h^2$.

\begin{lem}
There are positive numbers
$C_1$ and $C_2$
such that either one of the following holds
for any $P\in \Delta^n-D$
(i)  $Q(P)\leq C_1\cdot \psi(|z(P)|)^2$
or 
(ii) $\Delta(\log Q)(P)\leq -C_2\cdot Q(P)$.
\end{lem}
\pf
Let $P$ be a point of $\Delta^{n}-D$.
We take an orthogonal base $\vece=(e_i)$
of the fiber $E_{|P}$ with respect to the metric $h$
such that the filtration
$\big\{F_i=\langle e_j\,|\,j\leq i\rangle\big\}$
is preserved by the map $f_{0\,|\,P}$.
We have the matrix representation $\Gamma$ of $f_{0\,|\,P}$
with respect to the base $e$,
that is,
$f_{0\,|\,P}\vece=\vece\cdot \Gamma$.
Then $\Gamma$ is a triangular matrix
by our choice of $\vece$.
We denote the diagonal part of $\Gamma$
by $\Gamma_0$
and we put $\Gamma_1:=\Gamma-\Gamma_0$.
Let $\Gamma^{\dagger}$ be the matrix
representation of $f_0^{\dagger}$ with respect to the base $\vece$.
It is the adjoint matrix of $\Gamma$,
that is,
$\Gamma^{\dagger}={}^t\overline{\Gamma}$.
We denote the diagonal part by $\Gamma_0^{\dagger}$
and put $\Gamma_1^{\dagger}:=\Gamma^{\dagger}-\Gamma_0^{\dagger}$.

We have the equality
$[\Gamma,\Gamma^{\dagger}]=
 [\Gamma_0,\Gamma_1^{\dagger}]+[\Gamma_1,\Gamma_0^{\dagger}]
+[\Gamma_1,\Gamma_1^{\dagger}]$.
We denote the diagonal part of $[\Gamma_1,\Gamma_1^{\dagger}]$
by $\Xi$.
It is easy to see that there is a positive number $c_1$
such that
$|\Xi|^2\geq c_1|\Gamma_1|^2$,
where $c_1$ depends only on the dimension of the vector space
$E_{|P}$.
Note that the diagonal parts
of $[\Gamma_0,\Gamma_1^{\dagger}]$
and $[\Gamma_1,\Gamma_0^{\dagger}]$
are $0$, so that we obtain the inequality
$|[\Gamma,\Gamma^{\dagger}]|\geq c_1|\Gamma_1|^2$.

On the other hand,
we have $|\Gamma_0|\leq c_2|z(P)|^{-1+\epsilon}$
where $c_2$ is a positive constant,
which follows from the estimates for the eigenvalues.
We have the equality $Q=|f_{0\,|\,P}|_h^2=|\Gamma_0|^2+|\Gamma_1|^2$.
Thus there are positive constants $C_1$ and $C_2$
such that
if $Q\geq C_1\cdot \psi(|z(P)|)^2$ then 
$|\Gamma_1|^4\geq 2^{-1}\cdot Q^2$, so
$\Delta\log Q\leq -C_2\cdot Q$.
\hfill\qed

\vspace{.1in}
We can assume that $C_2=-8$
by taking a multiplication of some constant.
Then we know that one of the following holds
for any point $P\in \Delta^{\ast\,l}\times \Delta^{n-l}$:
(i) $Q(P)\cdot|dz\cdot
d\bar{z}|\leq C_1\cdot\psi(|z(P)|)^2\cdot|dz\cdot d\bar{z}|$,
or (ii) $\Delta(\log Q^{1/2})(P)\leq -4 Q(P)$.
Note that we can make $C_1$ larger.
In particular, we can assume that $C_1>1$.
The (ii) means that the curvature given by the pseudo metric
$Q|dz\cdot d\bar{z}|$ is less than $-1$.

We use the notation $\Delta_x:=\{x\in\Delta\}$
and $\Delta^{\ast}_z:=\{z\in\Delta^{\ast}\}$.
We take a holomorphic covering map $\Delta_x\lrarr \Delta_z^{\ast}$.
We have the natural isomorphism $\Delta_z^{\ast}\simeq \pi_1^{-1}(p)$
for any $p\in \Delta^{\ast\,l-1}\times \Delta^{n-l}$.
Then we obtain the holomorphic covering map
$F:\Delta_x\lrarr \pi_1^{-1}(p)$.
The pull back $F^{-1}\bigl(\psi(|z|^2)|dzd\bar{z}|\bigr)$
is same as $d\sigma^2:=(1-|z|^2)^{-2}|dxd\bar{x}|$,
because they are the Poincar\'{e} metrics of $\Delta_x$.
We put $Q_x:=F^{-1}(Q)\cdot |dF/dx|$
and $u:=\log Q_x^{1/2}$.
Then one of the following holds for any $P\in \Delta_x$:
(i) $Q_x(P)|dxd\bar{x}|\leq C_1 d\sigma^2$,
or
(ii) $\Delta(u)(P)\leq -4e^{2u(P)}$.

For any $R<1$ which is sufficiently close to $1$,
we put $v_R(x):=\log R(R^2-|x|^2)^{-1}$
for $|x|<R$.
Note that the following equalities:
\[
 \Delta v_R=-4e^{2v_R},
\quad
d\sigma^2=e^{2v_1}|dxd\bar{x}|,
\quad
v_1\leq v_R.
\]
We only have to show that
$e^{2u}\leq C_3e^{2v_1}$
over $\Delta_x$
for some positive constant $C_3$,
independently of $p\in\Delta^{\ast\,l-1}\times \Delta^{n-l}$.

We have already known that
there exists a constant $C_4>0$ such that one of the following holds
for any $P\in \Delta_x$:
(i) $u(P)\leq v_1(P)+C_4$,
or
(ii) $\Delta(u)(P)\leq -4e^{2u(P)}$.
We will prove that (i) holds for any $P\in \Delta_x$, in fact,
in this stage.

We consider the region
$S(R):=\{P\in \Delta_x\,|\,u(P)>v_R(P)+C_4\}$.
On the region,
we have the inequality
$\Delta(u-v_R+C_4)\leq -4(e^{2u}-e^{2v_R})<0$.
On $\{P\in \Delta_x\,|\,|x(P)|=R\}$, we have $v_R=\infty$.
Thus the boundary of $S(R)$
does not intersect with $\{P\in \Delta_x\,|\,|x(P)|=R\}$.
Thus we have the inequality 
$u-v_R+C_4\leq 0$ on the boundary of $S(R)$,
which raises a contradiction.
Thus the region $S(R)$ is empty.

Taking a limit $R\to 1$,
we obtain the desired inequality $u\leq v_1+C_4$.
Thus we completed the proof of Proposition \ref{prop;1.25.1}.
\hfill\qed

The following corollary is just a reformulation.
\begin{cor} \label{cor;11.27.26}
Let $P$ be a point of $X$ and $(\nbigu,\varphi)$
be an admissible coordinate around $P$.
We have the metric $(\cdot,\cdot)_{h,g_{\poin}}$
of $End(E)\otimes \Omega^{p,q}_{\nbigu-D}$.
The norm of $\theta$ with respect to the metric
$(\cdot,\cdot)_{h,g_{\poin}}$
is bounded.
\hfill\qed
\end{cor}

Let $(E,\delbar_E,\theta,h)$ be a tame nilpotent harmonic bundle.
Then we have the deformed holomorphic bundle
$(\nbige,d'',h)$ over $\nbigx$,
and
$(\nbige^{\lambda},d^{\twoprime\,\lambda},h)$
over $\nbigx^{\lambda}$.
\begin{prop}\mbox{{}}
The hermitian holomorphic vector bundles 
$(\nbige^{\lambda},d^{\twoprime\,\lambda},h)$
and
$(\nbige,d'',h)$
is acceptable.
\end{prop}
\pf
The assertion for $(\nbige,d'',h)$ follows from the formula
(\ref{eq;11.27.25}) and Corollary \ref{cor;11.27.26}.
Similar to $(\nbige^{\lambda},d^{\twoprime\,\lambda},h)$.
\hfill\qed

\subsection{The tame nilpotent harmonic bundle over the punctured disc}

\label{subsection;12.15.10}

\subsubsection{Prolongation}

We recall some results of Simpson.
See the section 10 of \cite{s1}
and the sections 3, 4 and 5 of \cite{s2}.
Let $(E,\delbar_E,h)$ be a hermitian holomorphic bundle
over the punctured disc $\Delta^{\ast}$.
We denote the origin of $\Delta$ by $O$,
which gives a smooth divisor of $\Delta$.
Simpson showed the following:
\begin{lem}
Let $\alpha$ be a real number.
If $(E,\delbar_E,h)$ is acceptable,
then the prolongation $E_{\alpha}$ is coherent locally free.
\hfill\qed
\end{lem}

\begin{rem}\label{rem;11.27.32}
Our definition of acceptable is slightly different from Simpson's.
He showed the stronger result than that stated here.
If the curvature $R(h)$ is dominated by
$|z|^{-2}(-\log|z|^2)^{-2}+f$ for some $L^p$-function $f$
on $\Delta$,
then prolongments $E_{\alpha}$ are coherent locally free.
We will use his stronger result without mention.
\hfill\qed
\end{rem}

Let $\beta\geq\alpha$ be real numbers.
Then we have the naturally defined morphism
$E_{\beta}\lrarr E_{\alpha}$ of coherent sheaves.
We obtain the morphism
$E_{\beta|O}\lrarr E_{\alpha|O}$,
which gives a descending filtration
of the vector space $E_{\alpha|O}$.
We denote the image of $E_{\beta|O}$ by $F^{\beta}(E_{\alpha|O})$.
Similarly we have the morphism
$E_{<\beta}\lrarr E_{\alpha}$,
and thus $E_{<\beta|O}\lrarr E_{\alpha|O}$.
We put $Gr^{\alpha}=E_{\alpha|O}/E_{<\alpha|O}$.
Then the graduation of the filtration $F^{\cdot}(E_{\alpha|O})$
is $\bigoplus_{\alpha\leq\beta<\alpha+1}Gr^{\beta}$.

\begin{condition} \label{condition;11.28.1}
We mainly consider the case
that $\alpha=0$ and $\dim(Gr^{0})=\dim(E_{|O})$.
In the case, we say that 
the parabolic structure of $(E,\delbar_E,h)$
is trivial.
\hfill\qed
\end{condition}



We also refer the following.
\begin{lem}\mbox{{}} \label{lem;12.2.20}
Prolongation is compatible with the procedures
taking determinant, dual, and tensor products.
(See the papers {\rm\cite{s1}} and {\rm\cite{s2}}
of Simpson for a precise statement.)

Let assume the parabolic structures of all hermitian holomorphic
vector bundles are trivial, the the following holds:
\begin{itemize}
\item
For $(E,\delbar_E,h)$,
we have $\prolong{\det(E)}\simeq\det(\prolong{E})$
and $(\prolong{E})^{\lor}\simeq \prolong{(E^{\lor})}$.
\item
For $(E_i,\delbar_{E_i},h_i)$ $(i=1,2)$,
we have
$\prolong{(E_1\otimes E_2)}=\prolong{E_1}\otimes \prolong{E_2}$.
In particular,
we have $\prolong{Sym^{l}(E)}\simeq Sym^l(\prolong{E})$
and $\prolong{(\bigwedge^l E)}\simeq \bigwedge^l(\prolong{E})$.
\hfill\qed
\end{itemize}
\end{lem}


In the case of harmonic bundle,
we obtain the following corollary.

\begin{cor}
Let $(E,\delbar_E,h,\theta)$ be a tame nilpotent harmonic bundle
over $\Delta^{\ast}$.
\begin{itemize}
\item
$\nbige^{\lambda}_{\alpha}$ is coherent locally free
for any $\alpha$.
\item
Let $f$ be a holomorphic section of $\nbige^{\lambda}_{\alpha}$,
and then $\DD^{\lambda}(f)$ is a holomorphic section
of $\nbige^{\lambda}_{\alpha-1}\otimes\Omega^{1,0}_{\Delta}$.
\end{itemize}
\end{cor}
\pf
Since $(\nbige^{\lambda},h)$ is acceptable,
$\nbige^{\lambda}_{\alpha}$  is locally free coherent sheaf.
By the same argument as that in 737 page of \cite{s2},
we can obtain the estimate for
$d^{\prime\lambda}(f)$ for a holomorphic section $f$
of $\nbige^{\lambda}_{\alpha}$.
Here $d^{\prime\,\lambda}$ denotes the $(1,0)$-part of
of the metric connection
of the hermitian holomorphic bundle
$(\nbige^{\lambda},d^{\twoprime\,\lambda},h)$.
The second claim follows from such estimate.
\hfill\qed

For each $\lambda$, we obtain the residue $\Res(\DD^{\lambda})$:
\[
 \Res(\DD^{\lambda}):
 \nbige^{\lambda}_{\alpha|O}\lrarr \nbige^{\lambda}_{\alpha|O}.
\]
It preserves the parabolic filtration $F^{\alpha}$.
Thus we also obtain the elements of
$End(Gr^{\beta})$,
which we denote by $\Res(\DD^{\lambda})_{\beta}$.

\begin{lem}\mbox{{}}
\begin{itemize}
\item
Let $(E,\delbar_E,\theta,h)$ be a tame nilpotent harmonic bundle.
Assume that the parabolic structure of $(\nbige^{\lambda},h)$
is trivial.
Then we have
$\prolong{\det(\nbige^{\lambda})}=
 \det(\prolong{\nbige^{\lambda}})$,
and
$\prolong{(\nbige^{\lambda\,\lor})}=
 (\prolong{\nbige^{\lambda}})^{\lor}$.
Moreover, the induced residues are same.
\item
Let $(E_i,\delbar_{E_i},\theta_i,h_i)$ $(i=1,2)$
be tame nilpotent harmonic bundles.
Assume that the parabolic structure of $(\nbige^{\lambda}_i,h_i)$
are trivial.
Then we have
$\prolong{(\nbige_1^{\lambda}\otimes\nbige^{\lambda}_2)}=
 \prolong{\nbige_1^{\lambda}}\otimes \prolong{\nbige^{\lambda}_2}$.
The induced residues are isomorphic.
In particular, similar things hold
for symmetric products and exterior products.
\hfill\qed
\end{itemize}
\end{lem}

\subsubsection{Some inequalities}

We recall the $\lambda$-connection version
of the inequality due to Simpson (Lemma 4.1 of \cite{s2}).
In this subsubsection,
the metric of $\Delta^{\ast}$ is the standard 
metric given by $|dz\cdot d\barz|$.

Let $(E,d''_E)$ be a holomorphic bundle
over $\Delta^{\ast}$,
and $\DD^{\lambda}$ be a flat holomorphic $\lambda$-connection
on $(E,d''_E)$.
It is not necessarily obtained from a harmonic metric.
Let $h$ be a hermitian metric on $E$,
which is not necessarily harmonic.
We denote the $(1,0)$-part of the metric connection
of $d''_E$ with respect to $h$
by $d'_E$.
We denote the $(1,0)$-part of $\DD^{\lambda}$
by $\DD^{\lambda\,\prime}$,
that is $\DD^{\lambda}=d''_E+\DD^{\lambda\,\prime}$.

We put as follows:
\begin{equation} \label{eq;11.27.35}
 \theta:=\frac{1}{1+|\lambda|^2}
 (\DD^{\lambda\prime}-\lambda d_E')
\in C^{\infty}(X,End(E)\otimes\Omega^{1,0}_X).
\end{equation}
It is not necessarily holomorphic.
Here we start from $\DD^{\lambda}$ and the metric.
Thus we use the notation $\theta(\DD^{\lambda},h)$
if we emphasize the dependence of $\theta$
on $\DD^{\lambda}$ and $h$.

We denote the adjoint of $\theta$ with respect to $h$
by $\theta^{\dagger}$.
Then we put as follows:
\[
 \delbar_E:=d''_E-\lambda\cdot\theta^{\dagger},
\quad
 \del_E:=d'_E+\lambdabar\cdot\theta.
\]
Then we put as follows:
\[
 G(\DD^{\lambda},h):=\delbar_E(\theta).
\]
It is easy to see that
the tuple $(E,\delbar_E,\theta,h)$
is harmonic if and only if $G(\DD^{\lambda},h)=0$ .

The following lemma is just a $\lambda$-connection version
of Lemma 4.1 of \cite{s2}.
\begin{lem}
Assume that $\lambda\neq 0$.
Let $s$  be a section of $E$ such that $\DD^{\lambda}(s)=0$.
Then we obtain the following inequality:
\[
 \Delta\log|s|^2
 \leq 2\bigl(|\lambda|^{-1}+|\lambda|\bigr)
 \cdot |G(\DD^{\lambda},h)|_h.
\]
\end{lem}
\pf
We denote the curvature of $d'_E+d''_E$ by $R(h,d''_E)$.
By our assumption $\DD^{\lambda}(s)=0$,
the section $s$ is holomorphic.
Thus we have the following equality:
\[
 \del\delbar |s|_h^2=
(d'_Es,d'_Es)_h+(s,R(h,d_E'')s)_h.
\]
We also have the following:
\[
 R(h,d''_E)=
\bigl( \delbar_E+\lambda\cdot\theta^{\dagger}
\bigr)\cdot
\bigl( \del_E-\lambdabar\cdot\theta
\bigr)
+
\bigl( \del_E-\lambdabar\cdot\theta
\bigr)\cdot
\bigl( \delbar_E+\lambda\cdot\theta^{\dagger}
\bigr)
\]
By our assumption of the flatness of $\DD^{\lambda}$,
we have the following:
\[
0
=\DD^{\lambda}\circ \DD^{\lambda}
=(\delbar_E+\lambda\cdot\theta^{\dagger})\cdot
 (\lambda\del_E+\theta)
+(\lambda\del_E+\theta)\cdot
 (\delbar_E+\lambda\cdot\theta^{\dagger}).
\]
Then we obtain the following by a direct calculation:
\begin{equation} \label{eq;11.27.31}
 R(h,d''_E)=
-(1+|\lambda|^2)\cdot
\bigl(
 \theta\wedge\theta^{\dagger}
+\theta^{\dagger}\wedge\theta
\bigr)
-\lambda^{-1}\cdot(1+|\lambda|^2)\cdot G(\DD^{\lambda},h).
\end{equation}
By our assumption $\DD^{\lambda}(s)=0$,
we have  $\bigl(\lambda\del_E+\theta\bigr)s=0$.
Thus we obtain the following:
\begin{equation} \label{eq;11.27.30}
 d'_E(s)=\bigl(\del_E-\lambdabar\theta\bigr)s=
-(\lambda^{-1}+\lambdabar)\cdot\theta (s)=
-\lambda^{-1}(1+|\lambda|^2)\cdot \theta (s).
\end{equation}
Hence we obtain the following equality:
\begin{multline}
 \del\delbar|s|^2_h=
 \Bigl(
 |\lambda|^{-2}(1+|\lambda|^2)^2
-(1+|\lambda|^2)
 \Bigr)\cdot(\theta s,\theta s)_h
-\bigl(1+|\lambda|^2)\cdot(\theta^{\dagger}s,\theta^{\dagger}s\bigr)_h\\
-\lambdabar^{-1}(1+|\lambda|^2)
 \cdot
 \bigl(s,G(\DD^{\lambda},h)s\bigr)_h.
\end{multline}
Then we obtain the following:
\[
 \Delta''|s|^2_h=
-(1+|\lambda|^{-2})\cdot |\theta s|_h^2
-(1+|\lambda|^2)\cdot |\theta^{\dagger}s|_h^2
+\lambda^{-1}(1+|\lambda|^2)\cdot 
 \bigl(s,\sqrt{-1}\Lambda G(\DD^{\lambda},h)s\bigr)_h.
\]
Here $\Lambda$ denotes the operator
from $\Omega^{1,1}$ to $\cnum$,
such that $dz\cdot d\zbar\longmapsto -\sqrt{-1}$.

We have the following equality:
\[
 \Delta''\log|s|_h^2=
 \frac{\Delta''|s|_h^2}{|s|_h^2}
-\frac{|(s,d'_Es)_h|^2}{|s|_h^4}.
\]
Due to (\ref{eq;11.27.30}),
we have the following:
\[
 |(s,d'_Es)_h|^2=
 (|\lambda|^{-2}+1)\cdot |(s,\theta s)_h|^2
+(|1+\lambda|^2)\cdot |(s,\theta^{\dagger}s)_h|^2.
\]
Thus we obtain the inequality desired.
\hfill\qed

Assume that $|\theta|_h$ is dominated by
$|z|^{-2}(-\log|z|^2)^2$.
Then,
due to the equality (\ref{eq;11.27.31}),
the prolongments $E_{\alpha}$ are locally free
(See Remark \ref{rem;11.27.32}).
When we emphasize the dependence of the prolongments $E_{\alpha}$
to the metric $h$,
we use the notation $E_{\alpha}(h)$
instead of $E_{\alpha}$.

\begin{cor}\label{cor;11.27.37}
Assume that $\lambda\neq 0$.
Let $(E,\DD^{\lambda})$ be as above.
Let $h_1$ and $h_2$ satisfies the following conditions:
\begin{enumerate}
\item \label{number;11.27.35}
 For $i=1,2$,
the hermitian bundles $(E,h_i)$ are acceptable.
and  the functions
 $|G(\DD^{\lambda},h_i)|_{h_i}$ are $L^p$.
\item \label{number;11.27.36}
 We have  $E_{\alpha}(h_1)=E_{\alpha}(h_2)$ 
 and $E^{\lor}_{\alpha}(h_1)=E^{\lor}_{\alpha}(h_2)$
 for any $\alpha$.
 Note that it also implies the coincidence of the parabolic structures.
\end{enumerate}
Then the metrics $h_1$ and $h_2$ are mutually bounded.
\end{cor}
\pf
The argument is same as Corollary 4.2 and Corollary 4.3
of \cite{s2}.
\hfill\qed

\subsubsection{Norm estimate in one dimensional case}

In this subsubsection,
we give a norm estimate in one dimensional case.
In \cite{s2}, Simpson discussed the cases $\lambda=0$ and $\lambda=1$.
Clearly his argument works in general case.
We only have to indicate how to change.

We recall the argument of Theorem 4 of \cite{s2}.
Let $\harmonicbundle$ be a tame harmonic bundle over
a punctured disc.
We have the $\lambda$-connection
$(\nbige^{\lambda},\DD^{\lambda})$.
Let pick a real number $\alpha\in\real$.
Then we have the residues of $\DD^{\lambda}$:
\[
 \bigl(Gr^{\beta},\Res(\DD^{\lambda})_{\beta}\bigr)
\quad
 (\alpha\leq\beta<\alpha+1).
\]
We decompose $Gr^{\beta}$ into the generalized
eigenspaces of $\Res(\DD^{\lambda})_{\beta}$.
\[
 (Gr^{\beta},\resddlambda_{\beta})=
 \bigoplus_{\omega}
 (Gr^{\beta}_{\omega},\resddlambda_{\beta,\omega}).
\]
Here $\omega$ runs through the set of eigenvalues
of $\resddlambda_{\beta}$.
The pair $(Gr^{\beta}_{\omega},\resddlambda_{\beta,\omega})$
of the vector space and the endomorphism is called
the $(\beta,\omega)$-part of
$(\nbige^{\lambda}_{\alpha},\DD^{\lambda},h)$.

For $(\beta,\omega)$, we denote the nilpotent part
of $\resddlambda_{\beta,\omega}$
by $N(\beta,\omega)$,
and we put $V(\beta,\omega):=Gr^{\beta}_{\omega}$.

Let consider the following harmonic bundle
(See the subsection \ref{subsection;11.27.27}):
\[
 (E_0,\delbar_{E_0},\theta_0,h_0):=
 \bigoplus_{(\beta,\omega)}
 \modelbundle{V(\beta,\omega)}{N(\beta,\omega)}
 \otimes L(C_1,C_2).
\]
Here $C_1$ and $C_2$ are real numbers given as follows:
\[
  C_1=
\frac{
 (1-|\lambda|^2)\cdot \beta+2Re(\lambdabar\cdot \omega)}
 {1+|\lambda|^2},
\quad
 C_2=
  \frac{\beta-\lambda\cdot \omega}
 {1+|\lambda|^2}
\]
We have the corresponding deformed holomorphic bundle
$\nbige_0^{\lambda}$ on $\nbigx^{\lambda}$,
and the $\lambda$-connection $\DD^{\lambda}_0$.
By our construction,
the $(\beta,\omega)$-part of
$(\nbige^{\lambda}_{0,\alpha},\DD^{\lambda}_0,h_0)$
is isomorphic to the $(\beta,\omega)$-part of
$(\nbige^{\lambda}_{\alpha},\DD^{\lambda},h)$.

\begin{lem} \label{lem;11.27.41}
There exists a holomorphic isomorphism
$f:\nbige^{\lambda}_0\lrarr \nbige^{\lambda}$
satisfying the following:
\begin{enumerate}
\item \label{number;11.27.28}
We put
$ g_1:=f_{|O}\circ \Res(\DD^{\lambda}_0)-\resddlambda\circ f_{|O}$
and 
$g_2:=f^{-1}_{|O}\circ \resddlambda-\Res(\DD^{\lambda}_0)\circ f^{-1}_{|O}$.
Then $g_1(F^{\beta})\subset F^{<\beta}$
and $g_2(F^{\beta})\subset F^{<\beta}$.
\end{enumerate}
\end{lem}
\pf
First of all, we take an isomorphism
$f_{|O}:\nbige^{\lambda}_{0,\alpha|O}\lrarr \nbige^{\lambda}_{\alpha|O}$
such that the condition \ref{number;11.27.28}.
It is possible because the graded parts
of the endomorphisms
$\resddlambda$ and $\Res(\DD^{\lambda}_0)$ are isomorphic.
And then we only have to extend $f_{|O}$
to a holomorphic map $f$ over $\Delta$.
\hfill\qed

\vspace{.1in}
By the isomorphism $f$,
we identify $\nbige^{\lambda}$ and $\nbige^{\lambda}_0$.
Thus the metric $h_0$ and the $\lambda$-connection $\DD_0^{\lambda}$
induces the metric and the $\lambda$-connection
on $\nbige^{\lambda}$.
Let compare $h$ and $h_0$.
It is clear by our construction
that the condition 2 in Corollary \ref{cor;11.27.37}
is satisfied for $h$ and $h_0$.
Recall that we obtain the non-holomorphic Higgs field
from the $\lambda$-connection and the metric
(See (\ref{eq;11.27.35})).
Here we have the following:
\[
 \theta(h,\DD^{\lambda})=\theta,
\quad
  \theta(h_0,\DD^{\lambda}_0)=\theta_0
\]
They are not same in general.
We also have $\theta_1=\theta(h_0,\DD^{\lambda})$,
which is not same as both of them above.

\begin{lem}
The condition $1$ in Corollary {\rm \ref{cor;11.27.37}}
is satisfied for the metrics $h$ and $h_0$ and the $\lambda$-connection
$\DD^{\lambda}$.
\end{lem}
\pf
Precisely we have to show the following:
\begin{itemize}
\item
 $|G(h,\DD^{\lambda})|_{h}$ is $L^p$.
\item
 $|G(h_0,\DD^{\lambda})|_{h_0}$ is $L^p$.
\end{itemize}
Since $(E,\delbar_E,\theta,h)$ is harmonic,
we know $G(\DD^{\lambda},h)=\delbar_E(\theta)=0$.
Thus we only have to care $\theta_1=\theta(h_0,\DD^{\lambda})$
and $G(h_0,\DD^{\lambda})$.
Let $\theta_1^{\dagger}$ denote the conjugate of $\theta_1$
with respect to $h_0$.

Let $\delta'$ denote the $(1,0)$-part of the metric connection
of $\nbige^{\lambda}$ with respect to the metric $h_0$.
By definition, we have the following:
\[
 \theta_1-\theta_0=
\frac{1}{|\lambda|^2+1}
 \Bigl(
 \DD^{\lambda}-\lambda\cdot\delta'
-\DD^{\lambda}_0+\lambda\cdot\delta'
 \Bigr)=
\frac{1}{|\lambda|^2+1}
(\DD^{\lambda}-\DD^{\lambda}_0).
\]
We only have to recall the argument in page 747 of \cite{s2}:
We put $A=\DD^{\lambda}-\DD^{\lambda}_0\in 
\Gamma(\Delta^{\ast},End(\nbige^{\lambda}))$.
By our choice of $f$,
the order of $A$ is less than $-1+\epsilon$
for some positive $\epsilon>0$.
By (\ref{eq;11.27.31}),
we have the following:
\[
\lambda^{-1}\cdot(1+|\lambda|^2)\cdot G(\DD^{\lambda},h_0)
=
-(1+|\lambda|^2)\cdot
\bigl(
 \theta_1\wedge\theta_1^{\dagger}
+\theta_1^{\dagger}\wedge\theta_1
\bigr)
- R(h_0,d''_{\nbige^{\lambda}_0}).
\]
Since $h_0$  and $\DD^{\lambda}_0$ is obtained from harmonic bundles,
we have the following:
\[
0=
\lambda^{-1}\cdot(1+|\lambda|^2)\cdot G(\DD^{\lambda}_0,h_0)
=
-(1+|\lambda|^2)\cdot
\bigl(
 \theta_0\wedge\theta_0^{\dagger}
+\theta_0^{\dagger}\wedge\theta_0
\bigr)
- R(h_0,d''_{\nbige^{\lambda}_0}).
\]
Hence we obtain the following:
\[
 \lambda^{-1}G(\DD^{\lambda},h_0)
=
\bigl(
 \theta_0\wedge\theta_0^{\dagger}
+\theta_0^{\dagger}\wedge\theta_0
- \theta_1\wedge\theta_1^{\dagger}
+\theta_1^{\dagger}\wedge\theta_1
\bigr)
=\frac{-1}{|\lambda|^2+1}
 \Bigl(
[A,A^{\dagger}]+[A,\theta^{\dagger}]+[\theta,A^{\dagger}]
 \Bigr)
\]
Hence we obtain the estimate
$G(\DD^{\lambda},h_0)\leq C\cdot |z|^{-2+\epsilon}$
for some $\epsilon>0$.
Hence $|G(\DD^{\lambda},h_0)|_{h_0}$ is $L_p$
with respect to the measure $|dz\cdot d\zbar|$.
\hfill\qed

We have a direct corollary.
\begin{cor} \label{cor;11.27.40}
The metrics $h$ and $h_0$ above are mutually bounded.
\hfill\qed
\end{cor}

We can reword Corollary \ref{cor;11.27.40}.
The nilpotent part of $\resddlambda_{(\beta,\omega)}$
induces the weight filtration $W_{(\beta,\omega)}$
on $Gr^{\beta}_{\omega}$.
The filtrations $\{W_{(\beta,\omega)}\,|\,\omega\,\,\mbox{eigenvalue}\}$
give the filtration of $Gr^{\beta}=\bigoplus_{\omega}Gr^{\beta}_{\omega}$.
Let $W$ denote the filtration of $\bigoplus_{\beta} Gr^{\beta}$.
We put $F_{\beta}=F^{-\beta}\subset \nbige^{\lambda}_{\alpha|O}$.
Then $F_{\beta}$ gives an ascending filtration.
Then we obtain the sequence of the filtrations
$(F_{\cdot},W_{\cdot})$.
Let take a holomorphic frame $\vecv=(v_1,\ldots,v_r)$
of $\nbige^{\lambda}_{\alpha}$
over $\Delta$ satisfying the following:
\begin{itemize}
\item
 $\vecv_{|O}$ is compatible with the ascending filtration
     $\{F_{\beta}\}$.
\item
 The induced base $\vecv^{(1)}$ of $\bigoplus Gr_{\beta}$
 is compatible with the filtration $W$.
\end{itemize}
Such frame $\vecv$ is called compatible with $(F_{\cdot},W_{\cdot})$.
The element $v_{i\,|\,O}$ induces an element of $\bigoplus Gr_{\beta}$,
which we denote by $v_i^{(1)}$.
We put as follows:
\[
 \alpha_i=\deg^{F_{\cdot}}(v_{i\,|\,O}),
\quad
 k_i=\frac{1}{2}\deg^{W}(v_i^{(1)}).
\]
We have the $C^{\infty}$-frame $\vecv'=(v_1',\ldots,v_r')$
of $\nbige^{\lambda}$ over $\Delta^{\ast}$,
defined as follows:
\[
 v_i'=|z|^{\alpha_i}\cdot (-\log|z|)^{-k_i}\cdot v_i.
\]

\begin{cor}
The frame $\vecv'$ is adapted.
\end{cor}
\pf
The claim follows 
from Corollary \ref{cor;11.27.42}
and Corollary \ref{cor;11.27.40}.
\hfill\qed

Let $f$ be a holomorphic section of $\nbige^{\lambda}_{\alpha}$
over $\Delta$.
We have the number
$\alpha(f):=\deg^{F}(f(O))$
and $k(f):=2^{-1}\deg^{\weightfilt}(f^{(1)}(O))$.
Here $F$ is the ascending filtration above.
\begin{cor}
There exists positive constants $C_1$ and $C_2$
such that the following holds over $\Delta^{\ast}$:
\[
 0<C_1<
 (-\log|z|)^{-k(f)}\cdot |z|^{\alpha(f)}\cdot |f|_h
<C_2.
\]
\hfill\qed
\end{cor}

\subsubsection{Finiteness of some norms and some consequences}

In this section,
the metric of $\Delta^{\ast}$ is the standard one
given by $|dz\cdot d\zbar|$.
Let $\harmonicbundle$ be a tame nilpotent harmonic bundle
over $\Delta^{\ast}$.
We have a prolongment $\prolong{E}$ by increasing order $0$.
Then we have the vector spaces
We take a model bundle $(E_0,\theta_0,h_0)$
and $f:\prolong{E_0}\lrarr \prolong{E}$
as in Lemma \ref{lem;11.27.41}.
We will identify $E_0$ and $E$ over $\Delta^{\ast}$
by the morphism $f$.
We have already known that the metrics $h$ and $h_0$
are mutually bounded.

We have the deformed holomorphic bundle
$(\nbige^{\lambda},\DD^{\lambda},h)$
and $(\nbige^{\lambda}_0,\DD^{\lambda}_0,h_0)$
over $\{\lambda\}\times \Delta^{\ast}$,
which are induced by
$(E,\theta,h)$ and $(E_0,h_0,\theta_0)$ respectively.

We have 
the parabolic filtration and the weight filtration
for $\prolong{\nbige^{\lambda}}$
and $\prolong{\nbige^{\lambda}_0}$.
We denote them by $(F,\weightfilt)$ and
$(F_0,\weightfilt_0)$ respectively.

Let take a holomorphic frame $\vecv$
of $\prolong{\nbige^{\lambda}}$
such that it is compatible with
the sequence of the filtrations $(F,\weightfilt)$.
We put as follows:
\[
 \alpha(v_i):=\deg^{F}(v_i(O)),
\quad
 k(v_i):=\frac{1}{2}\deg^{\weightfilt}(v_i^{(1)}).
\]
Here $v^{(1)}_i$ denotes the induced element of $Gr^{F}$.
We also take a holomorphic frame $\vecv_0$
of $\prolong{\nbige^{\lambda}_0}$
such that it is compatible with the sequence of
the filtrations $(F_0,\weightfilt_0)$.
Similarly we obtain the numbers
$\alpha(v_{0,i})$ and $k(v_{0,i})$.

Since the underlying $C^{\infty}$-vector bundles
of $\nbige^{\lambda}$ and $\nbige^{\lambda}_0$
are naturally identified with $E$ over $\Delta^{\ast}$,
$\vecv$ and $\vecv_0$ give the $C^{\infty}$-frames of $E$
over $\Delta^{\ast}$.
Let $\nbigi$ denote the $C^{\infty}$-isomorphism.
Then we have the functions
$I_{i\,j}\in C^{\infty}(\Delta^{\ast},\cnum)$ determined as follows:
\[
 \nbigi(v_{0\,j})=\sum_i I_{i\,j}\cdot v_i.
\]
\begin{rem}
Note that our rule of subscription is different from
the rule in {\rm\cite{s2}}.
There, our $I_{i\,j}$ is denoted by $I_{j\,i}$ 
in the section {\rm 7} of {\rm\cite{s2}}.
We also note that the choice of the signature of $\alpha(v_i)$
is opposite to that of Simpson.
\hfill\qed
\end{rem}

Simpson considered the following norms
$||\cdot||_Z$ and $||\cdot||_W$
for a function $f$ on $\Delta^{\ast}(C)$:
\[
  ||f||_{Z,C}=
 \int _{\Delta^{\ast}(C)} |f|\cdot \frac{dr\cdot d\alpha}{r\cdot (-\log r)},
\quad\quad
 ||f||_{W,C}=
 \int_{\Delta^{\ast}(C)} |f|\cdot
  \frac{dr\cdot d\alpha}{r\cdot(-\log r)\cdot\log(-\log r)}.
\]
Here the real coordinate $z=r\cdot\exp(2\pi\sqrt{-1}\alpha)$
is used.
Simpson showed the following lemma
to show that the conjugacy classes of the resides
are invariant
(See section 7 in \cite{s2}).
\begin{lem}\label{lem;1.14.1}
For any $i$ and $j$, we have the following finiteness:
\begin{equation} \label{eq;1.14.1}
 \Bigl|\Bigl|
 \zbar\cdot\delbar_{z}I_{i\,j\,|\Delta^{\ast}(C)}
 \cdot
 |\log r|^{k(v_i)-k(v_{0,j})+1}
 \cdot
 r^{-\alpha(v_i)+\alpha(v_{0,j})}
  \Bigr|\Bigr|_{Z,C}<\infty.
\quad\quad
\end{equation}
For any $i$ and $j$ such that $\alpha(v_i)-\alpha(v_{0,j})\neq 0$,
there exists a positive constant $C>0$ satisfying the following:
\begin{equation} \label{eq;11.27.45}
 \Bigl|
 I_{i\,j}
 \Bigr|
 \cdot
 (-\log r)^{k(v_i)-k(v_{0,j})}
 \cdot
 r^{-\alpha(v_i)+\alpha(v_{0,j})}
\leq C(-\log|z|)^{-1}.
\end{equation}
For any $i$ and $j$ such that 
$\alpha(v_i)-\alpha(v_{0,j})=0$ and that
$k(v_i)-k(v^0_j)\neq 0$,
we have the following finiteness:
\begin{equation} \label{eq;1.14.2}
  \Bigl|\Bigl| I_{i\,j}|_{\{\lambda\}\times B^{\ast}(C)}
   |\log r|^{(k(v_i)-k(v_{0,j}))}
 \Bigr|\Bigr|_{W,C}<\infty.
\end{equation}
\end{lem}
\pf
We have
$\delbar_z(\nbigi)=\lambda\cdot(\theta^{\dagger}-\theta^{\dagger}_0)$
by definition of the holomorphic structures
of $\nbige^{\lambda}_0$ and $\nbige^{\lambda}$.
Simpson showed the following inequality (Lemma 7.7 in \cite{s2}):
\[
 \int |\theta^{\dagger}-\theta_0^{\dagger}|_{h}\cdot (-\log |z|)\cdot
 |dz\cdot d\bar{z}|<\infty.
\]
Here the metric $h$ is used. Since $h_0$ and $h$ are mutually bounded,
we can also use $h_0$. In fact, we can take any metric
mutually bounded to $h_0$ and $h$.
Thus we do not have to care a choice of the metric in the following.

Recall that the frames $\vecv'=(v_i)$ and $\vecv'_0=(v_{0,i}')$
are adapted, if we put as follows:
\[
 v'_i:=|z|^{\alpha(v_i)}\cdot(-\log|z|)^{-k(v_i)}\cdot v_i,
\quad
 v'_{0,j}:=
  |z|^{\alpha(v_{0,j})}\cdot(-\log|z|)^{-k(v_{0,j})}\cdot v_{0,j}.
\]
Since we know that $|\theta^{\dagger}|_h$ and $|\theta_0^{\dagger}|_h$
are bounded by $(-|z|\cdot\log |z|)^{-1}$,
we have the inequality:
\[
 \Big|\bar{\del}I_{i\,j}\cdot |z|^{-\alpha(v_i)+\alpha(v_{0,j})}
 \cdot (-\log|z|)^{k(v_i)-k(v_{0,j})}\Big|
 \leq C\cdot|\lambda|\times (-|z|\cdot\log |z|)^{-1}.
\]
Thus any components
$\bar{\del}_zI_{i\,j}$
satisfy the inequality (\ref{eq;1.14.1}).

We obtain the inequalities 
(\ref{eq;11.27.45}) and (\ref{eq;1.14.2})
by using Lemma 7.1, Lemma 7.8 and the argument
of corollary 7.10 in \cite{s2}.
\hfill\qed

Note that many of the arguments in section 7 in \cite{s2}
are not needed in our case, for we assumed that the residue
of $\harmonicbundle$ are nilpotent.

\begin{cor} \label{cor;6.25.10}
The conjugacy classes of the residues
$(Gr^{\beta}_{\omega},\Res(\DD^{\lambda})_{(\beta,\omega)})$
are independent of $\lambda$.
\end{cor}
\pf
Simpson showed that
the conjugacy classes of
$\Res(\DD^{\lambda})_{|(1,O)}$ and $\Res(\DD^{\lambda})_{|(0,O)}$
are same in \cite{s2}
by using the inequality (\ref{eq;1.14.2}).
By the same method and Lemma \ref{lem;1.14.1},
we can show that the conjugacy classes of
$\Res(\DD^{\lambda})_{(\beta,\omega)\,|\,(\lambda,O)}$
and $\Res(\DD^0)_{(\beta,\omega)|\,(0,O)}$ are same
for any $\lambda\neq 0$.
\hfill\qed

In particular, we obtain the following.
\begin{cor} \label{cor;11.27.50}
Let $(E,\theta,h)$ be a tame nilpotent harmonic bundle
over $\Delta^{\ast}$.
We also assume that the parabolic structure of
the prolongment $\prolong{E}$ is trivial.
Then the following holds:
\begin{itemize}
\item
The parabolic structure of the prolongment
$\prolong{\nbige^{\lambda}}$ is trivial
for each $\lambda$.
\item
All of the eigenvalues of $\Res(\DD^{\lambda})$ 
on $\prolong{\nbige^{\lambda}_{|O}}$ are $0$.
\item
The conjugacy classes of
$(\prolong{\nbige^{\lambda}}_{|O},\resddlambda)$ 
are independent of $\lambda$.\hfill\qed
\end{itemize}
\end{cor}

\vspace{.1in}
Let $\harmonicbundle$ be a tame nilpotent harmonic bundle
over $\Delta^{\ast}$.
Assume that the parabolic structure of
$\prolong{E}$ is trivial.
When $\lambda\neq 0$,
we have the flat holomorphic bundle
$(\nbige^{\lambda},\DD^{\lambda,f})$.
We have a holomorphic frame $\vecw$
of $\nbige^{\lambda}$ over $\Delta^{\ast}$
satisfying the following:
\begin{itemize}
\item
We have the flat connection form
$A\in \Gamma(\Delta,M(r)\otimes\Omega_{\Delta}(\log O))$
determined by the relation
$\DD^{\lambda,f}\vecw=\vecw\cdot A$.
Then $A$ is of the form $A_0dz/z$ for some constant matrix $A_0$.
\item
All of the eigenvalues of $A_0$ are $0$.
\end{itemize}
Such frame is called a normalizing frame in this paper.
We have the prolongment of $\nbige^{\lambda}$
by a normalizing frame.
\begin{lem} \label{lem;12.15.200}
The prolongment by a normalizing frame is same as
the prolongment $\prolong{\nbige^{\lambda}}$
by an increasing order $0$.

In particular,
a normalizing frame naturally gives a frame of
the prolongation $\prolong{\nbige^{\lambda}}$.
\end{lem}
\pf
The claim follows from the uniqueness of
the prolongation, for which the holomorphic connection
is of log type.
(See \cite{d}. In particular, II section 5.)
\hfill\qed

\subsection{Definition of trivial parabolic structure}
\label{subsection;12.15.11}

Let $X$ be an $n$-dimensional complex manifold,
and $D$ be a normal crossing divisor of $X$.
Let $(E,\delbar_E,h)$ be a hermitian holomorphic bundle over $X-D$.
Let $C$ be a curve contained in $X$, transversal with $D$.
Then we obtain the hermitian holomorphic bundle
$(E,\delbar_E,h)_{|C-C\cap D}$.
\begin{df} \label{df;12.9.2}
We say that the parabolic structure of
the hermitian holomorphic bundle
$(E,\delbar_E,h)$ over $X-D$ is trivial,
if the following holds:
\begin{quote}
For any curve $C\subset X$ transversal with $D$,
the hermitian holomorphic bundle
$(E,\delbar_E,h)_{|C-C\cap D}$ over $C-C\cap D$ is trivial
in the sense of Condition {\rm \ref{condition;11.28.1}}.
\hfill\qed
\end{quote}
\end{df}

As an example, let consider the case
that $X=\Delta^n$ and $D=\bigcup_{j=1}^l D_j$,
where $D_j=\{z_j=0\}$.
The projection
$\Delta^{\ast\,l}\times \Delta^{n-l}
\lrarr \Delta^{\ast\,l-1}\times \Delta^{n-l}$,
omitting the $j$-th component,
is denoted by $\pi_j$.
For any element $a\in \Delta^{\ast\,l-1}\times \Delta^{n-l}$,
we obtain the curve
$\pi_j^{-1}(a)\subset \Delta^{\ast\,l}\times\Delta^{n-l}$.
Let $(E,\delbar_E,h)$ be a hermitian holomorphic bundle over
$X-D=\Delta^{\ast\,l}\times\Delta^{n-l}$.
Then the parabolic structure is trivial if and only if
the parabolic structure of $(E,\delbar_E,h)_{|\pi_j^{-1}(a)}$
is trivial for any
$a\in \Delta^{\ast\,l-1}\times\Delta^{n-l}$ and $j=1,\ldots,l$.

\begin{cor}
Let $\harmonicbundle$ be a tame nilpotent harmonic bundle
over $X-D$.
Assume that the parabolic structure
of the hermitian holomorphic bundle $(E,\delbar_E,h)$ is trivial.
\begin{itemize}
\item
The parabolic structure of
$(\nbige^{\lambda},d^{\twoprime\,\lambda},h)$ is trivial
for any $\lambda$.
\item
All of the eigenvalues of $\resddlambda$ are $0$
for any $\lambda$.
\item
If $\lambda\neq 0$, then we have the flat holomorphic bundle
$(\nbige^{\lambda},\DD^{\lambda})$.
All of the eigenvalues of the monodromies around the divisor $D$
are $1$.
\end{itemize}
\end{cor}
\pf
The claims follow immediately from Corollary \ref{cor;11.27.50}.
\hfill\qed

\subsubsection{Rank 1}

Let see the easy case, that is, the rank 1 case.
\begin{lem}
Let $\harmonicbundle$ be a tame nilpotent harmonic bundle
of rank 1 over $\Delta^{\ast\,l}\times \Delta^{n-l}$.
Assume that the parabolic structure is trivial.
Then it is naturally extended to the harmonic bundle
over $\Delta^n$.
The deformed holomorphic bundle
is also extended to that over $\cnum_{\lambda}\times \Delta^{n}$.
\end{lem}
\pf
Let consider
the holomorphic bundle
$(E,\delbar_E+\theta^{\dagger})$
with the flat connection $\DD^1=\delbar_E+\del_E+\theta^{\dagger}+\theta$.
Since the eigenvalue of the monodromy is $1$,
we can take a holomorphic frame
$e$ of $E$ over $X-D$ satisfying the following:
\[
 \DD^1(e)
=\bigl(\delbar_E+\theta^{\dagger}\bigr)e
=\bigl(\del_E+\theta\bigr)e=0.
\]
The $(1,0)$-part of the metric connection
of $\delbar_E+\theta^{\dagger}$ with respect to $h$
is given by $\del_E-\theta$.
We put $h_0=h(e,e)\in C^{\infty}(X-D,\real)$.
Then we have the following equation:
\[
 (\del h_0)\cdot h_0^{-1}=\del(\log h_0)=-2\theta.
\]
Since the rank of $E$ is $1$,
the sheaf $End(E)$ is naturally isomorphic to $\nbigo$.
The tameness and the nilpotentness
of $\theta$ implies that $\theta$ is, in fact,
a holomorphic section of $\Omega_{\Delta^n}^{1,0}$.
We also have the equality $\del(\theta)=0$,
because $(E,\delbar_E,\theta,h)$ is harmonic.
Thus we have a holomorphic function $f$
such that $\del(f)=\theta$.

We have the following equality:
\[
 \del\bigl(\log (h_0)+2f\bigr)=0.
\]
Note that $\log(h_0)$ is $\real$-valued.
Thus we can conclude that
$\log(h_0)=-4Re(f)+C$ for some constant $C\in\real$,
that is, we obtain the following:
\[
 h_0=\exp(-4Re(f)+C).
\]
It means that the increasing order of $e$ is $0$,
i.e.,
$e$ is a holomorphic frame of $\prolong{\nbige^{1}}$,
and $h$ induces the $C^{\infty}$-metric of $\prolong{\nbige^{1}}$.

We put
$v:=\exp\bigl((1-\lambda)\bar{f}\bigr)\cdot e$.
We have the equality $\theta^{\dagger}=\delbar(\fbar)$.
Then we obtain the following equality:
\[
 \bigl(
 \delbar_E+\lambda\theta^{\dagger}
 \bigr)\cdot v
=\bigl(
 \delbar_E+\delbar(\fbar)
 \bigr)\cdot
 \Bigl(
 \exp\bigl((1-\lambda)\fbar\bigr)\cdot e\Bigr)
+(\lambda-1)\cdot\delbar(\fbar)\cdot v
=0.
\]
We have $h(v,v)=\exp\Bigl(-2Re\bigl((1+\lambda)\fbar\bigr)\Bigr)$.
Thus $v$ gives a holomorphic frame of
$\prolong{\nbige}$ over $\cnum_{\lambda}\times \Delta$.
\hfill\qed

\begin{cor}
Let $X$ be a complex manifold, and $D$ be a normal crossing divisor.
Let $\harmonicbundle$ be a tame nilpotent harmonic bundle
with trivial parabolic structure
over $X-D$.
Then it is naturally extended to the harmonic bundle
over $X$.
The deformed holomorphic bundle $(\nbige,h)$
over $\nbigx-\nbigd$ is extended to that over $\nbigx$.
\hfill\qed
\end{cor}


\subsection{Some preliminary}
\label{subsection;12.15.15}

Recall something from \cite{cg}.
\subsubsection{Some results of Andreotti-Vesentini}

We recall some results of Andreotti-Vesentini in \cite{av}.
Let $(Y,g)$ be a complete Kahler manifold, not necessarily
compact.
We denote the natural volume form by $\vol$.
Let $(E,\delbar_E,h)$ be a hermitian holomorphic bundle
over $Y$.
The hermitian metric $h$ and the Kahler metric $g$
induces the fiberwise hermitian metric
of $E\otimes\Omega^{p,q}_Y$,
which we denote by $(\cdot,\cdot)_{h,g}$.
The space of $(p,q)$-forms with compact support
is denoted by $A_c^{p,q}(E)$.
For any $\eta_1,\eta_2\in A_c^{p,q}(E)$,
we put as follows:
\[
 \langle
 \eta_1,\eta_2
 \rangle_h=
 \int (\eta_1,\eta_2)_{h,g}\cdot\vol,
\quad
 ||\eta||^2_h=\langle \eta,\eta\rangle_h.
\]
The completion of $A_c^{p,q}$
with respect to the norm $||\cdot||_h$
is denoted by $A_h^{p,q}$.

We have the operator
$\delbar_E:A_c^{p,q}(E)\lrarr A_c^{p,q+1}(E)$,
and the formal adjoint
$\delbar_E^{\ast}:A_c^{p,q}(E)\lrarr A_c^{p,q-1}(E)$.
We use the notation
$\laplacian=\delbar_E^{\ast}\delbar_E+\delbar_E\delbar_E^{\ast}$.
We have the maximal closed extensions
$\delbar_E:A_h^{p,q}(E)\lrarr A_h^{p,q+1}(E)$
and 
$\delbar_E^{\ast}:A_h^{p,q}(E)\lrarr A_h^{p,q-1}(E)$.
We denote the domains of $\delbar_E$
and $\delbar_E^{\ast}$ by $Dom(\delbar_E)$
and $Dom(\delbar_E^{\ast})$ respectively.

\begin{prop}[Proposition 5 of \cite{av}] \label{prop;11.28.6}
In $W^{p,q}:=Dom(\delbar_E)\cap Dom(\delbar_E^{\ast})$,
the space $A_c^{p,q}(E)$ is dense
with respect to the the graph norm:
$||\eta||_h^2+||\delbar_E\eta||^2_h+||\delbar_E^{\ast}\eta||^2_h$.
(See also {\rm\cite{cg}}).
\hfill\qed
\end{prop}


\begin{prop}[Theorem 21 of \cite{av}] \label{prop;12.11.1}
Assume that 
there exists a positive number $c>0$ satisfying the following:
\begin{quote}
Then, for any $\eta\in W^{p,q}$,
we have $||\delbar_E\eta||_h^2+||\delbar_E^{\ast}\eta||_h^2
 \geq c\cdot ||\eta||_{h}^2$.
\end{quote}
For any $C^{\infty}$-element $\eta\in A_h^{p,q}(E)$
such that $\delbar_E(\eta)=0$,
we have a $C^{\infty}$-solution $\rho\in A_h^{p,q-1}(E)$
satisfying the equation $\delbar_E(\rho)=\eta$.
\hfill\qed
\end{prop}

\subsubsection{Kodaira identity}

For the Kahler manifold $Y$,
we have the operator $\Lambda:\Omega^{p,q}\lrarr \Omega^{p-1,q-1}$
(see  62 page of \cite{koba}).
For a section $f$ of $End(E)\otimes \Omega^{p_0,q_0}_Y$,
we have the natural morphism $A_c^{p,q}(E)\lrarr A_c^{p+p_0,q+q_0}(E)$,
defined by $\eta\longmapsto f\wedge \eta$.
We denote the morphism by $e(f)$.

We have the metric connection of $E$
induced by the holomorphic structure $\delbar_E$
and the hermitian metric $h$.
We denote the curvature by $R(h)$.
We have the Levi-Civita connection of the tangent bundle
of $Y$.
It induces the connection of $E\otimes\Omega^{0,1}$:
\[
\nabla:
 A_c^{0,0}(E\otimes\Omega^{0,1})\lrarr
 A_c^{0,1}(E\otimes\Omega^{0,1})\oplus
 A_c^{1,0}(E\otimes \Omega^{0,1}).
\]
We denote the $(0,1)$-part of $\nabla$ by $\nabla''$
to distinguish with $\delbar_E:A_c^{0,1}(E)\lrarr A_c^{0,2}$.
The $(1,0)$-part of $\nabla$ is same as
$\del$ of $E\otimes\Omega^{0,1}$.
We denote the curvature of $\nabla$
by $R(\nabla)$.

We denote the Ricci curvature of the Kahler metric $g$
by $\Ric(g)$.
We can naturally
regard $\Ric(g)$ as a section of $End(E)\otimes \Omega^{1,1}$,
by the natural diagonal inclusion $\cnum\lrarr End(E)$.

Let $f$ be a section of $End(E)\otimes\Omega^{1,1}_Y$,
and $\eta$ be an element of $A_c^{0,1}(E)$.
Then we put as follows:
\[
 \doublelangle
 f,\eta
 \doublerangle_h
:=
 -\sqrt{-1}
\bigl(
\xi,\eta
\bigr)_h,
\quad
 \xi:=
 \Bigl(
 \Lambda\circ e(f)
 -e\bigl(\Lambda(f)\bigr)
 \Bigr)
 (\eta)
=
 \Lambda
 \bigl(
 f\cdot\eta
 \bigr)
-\Lambda(f)\cdot\eta
\]

Let $\varphi_i$ be an $C^{\infty}$-orthogonal local coframe 
$(\varphi_i)$ of the tangent bundle of $Y$,
that is $g=\sum \varphi_i\cdot\varphibar_i$.
We also take a $C^{\infty}$-orthogonal local frame
$(e_i)$ of $E$.
We denote the dual frame by $(e_i^{\lor})$.
We have the local description:
\[
\eta=\sum \eta_{\mu,i}\cdot e_{\mu}\otimes\varphibar_i,
\quad
 f=\sum f_{\mu,\nu,i,\bar{j}}\cdot
 e_{\mu}^{\lor}\otimes e_{\nu}
 \otimes
 \bigl(
 \varphi_i\wedge\varphibar_j
 \bigr).
\]
Then we have the following local description
((9.1) in \cite{cg}):
\begin{equation} \label{eq;11.28.8}
 \doublelangle
 f,\eta
 \doublerangle_h:=
 \sum f_{\mu,\nu,i,\bar{j}}
\cdot
 \eta_{\mu,i}
\cdot
 \bar{\eta}_{\nu,j}.
\end{equation}

We recall the identity which is called Kodaira identity in \cite{cg}.
We only need the following special case.

\begin{prop}[Kodaira \cite{ko1}, Cornalba-Griffiths \cite{cg}]
Let $\eta$ be an element of $A_c^{0,1}(E)$.
We have the following equality:
\[
 ||\delbar_E(\eta)||_h^2
+||\delbar_E^{\ast}(\eta)||_h^2
=
 ||\nabla''\eta||^2+
\int\!\!
\big\langle\big\langle
R(h)+Ric(g),
\eta
\big\rangle\big\rangle_h\!
\vol.
\]
\end{prop}
\pf
(See \cite{ss} for some formulas of Nakano type.)
We have the equality:
\[
 ||\delbar_E(\eta)||_h^2+||\delbar_E^{\ast}(\eta)||_h^2=
 \langle\laplacian\eta,\eta\rangle_h.
\]
On the $(0,1)$-forms, we also have the equality:
\[
 \laplacian=
 \del\del^{\ast}+\del^{\ast}\del
-\sqrt{-1}\Lambda\circ e(R(h)).
\]
On $A^{0,1}_c(E)$, we have $\del^{\ast}=0$ for an obvious reason.
Thus we obtain the following:
\[
 \langle\laplacian\eta,\eta\rangle_h
=\langle\del \eta,\del \eta\rangle_h
-\sqrt{-1}
\Bigl\langle
 \Lambda \bigl(R(h)\cdot\eta\bigr),
\eta
\Bigr\rangle_h.
\]
We also have the operator on $A_c^{0,0}(E\otimes\Omega^{0,1})$:
\[
 \laplacian_1=
 \nabla''\nabla^{\twoprime\,\ast}+
 \nabla^{\twoprime\,\ast}\nabla''
=
 \del\del^{\ast}+\del^{\ast}\del
-\sqrt{-1}\Lambda\circ e(R(\nabla)).
\]
In this case, we have
$\Lambda\circ e(R(\nabla))(\eta)=
 \Lambda(R\nabla)\cdot \eta$.
Thus we obtain the following:
\[
 \langle\laplacian\eta,\eta\rangle_h=
 ||\nabla''\eta||_h^2
-\sqrt{-1}
\Bigl\langle
 \Lambda\bigl(R(h)\cdot\eta\bigr),\eta
\Bigr\rangle_h
+\sqrt{-1}
\Bigl\langle
 \Lambda\bigl(R(\nabla)\bigr)\cdot\eta,\eta
\Bigr\rangle_h.
\]
We have $R(\nabla)=R(h)+R(\Omega^{0,1})$.
It can be checked 
that $\sqrt{-1}(\Lambda(R(\Omega^{0,1})\cdot\eta),\eta)_h$
is same as
$\doublelangle \Ric(g),\eta\doublerangle$,
by a direct calculation and the coincidence
of the Ricci curvature and the mean curvature of the Kahler metric
(See (7.23) in page 28 of \cite{koba}.)
Thus we obtain the following:
\[
 \langle\laplacian\eta,\eta\rangle_h=
 ||\nabla''\eta||_h^2
+\int
 \left[
 -\sqrt{-1}
\Bigl(
 \Lambda\bigl(R(h)\cdot\eta\bigr)
 -\Lambda(R(h))\cdot \eta ,
 \eta
\Bigr)_h
+\doublelangle
 \Ric(g),\eta
 \doublerangle_h\right]\vol.
\]
Thus we are done.
\hfill\qed

\begin{cor}\label{cor;12.11.2}
Let $\eta$ be an element of
$Dom(\delbar_E)\cap Dom(\delbar_E^{\ast})$
in $A^{0,1}_h(E)$.
Then we have the following inequality:
\[
 ||\delbar_E(\eta)||^2_h
+||\delbar_E^{\ast}(\eta)||_h^2
\geq
\int\!\!
 \doublelangle
 R(h)+\Ric(g),\eta
 \doublerangle_h
\vol.
\]
\end{cor}
\pf
We only have to note the density of $A_c^{0,1}(E)$
in $Dom(\delbar_E)\cap Dom(\delbar^{\ast}_E)$
(Proposition \ref{prop;11.28.6}).
\hfill\qed

\subsubsection{The modification of acceptable metrics}

Let consider the Poincare metric $g_{\poin}$
on $Y=\Delta^{\ast\,l}\times \Delta^{n-l}$,
and the acceptable hermitian holomorphic vector bundle
$(E,\delbar_E,h)$
(Definition \ref{df;11.28.5}).

Let $\veca=(a_1,\ldots,a_l)$ be a tuple of real numbers.
Let $N$ be a real number.
Then we put as follows:
\[
 \tau(\veca,N):=
 \sum_{i=1}^l a_i\log|z_i|^{2}
+
 N\cdot
\left(
 \sum_{i=1}^l
 \log\bigl(-\log|z_i|^2\bigr)
+\sum_{i=l+1}^n
 \log(1-|z_i|^2)
\right).
\]
Recall the following formulas:
\[
 \delbar\del \log|z|^2=0,
\quad\quad
 \delbar\del\log(-\log|z|^2)
=\frac{dz\cdot d\zbar}
 {(-\log|z|^2)^2\cdot |z|^2},
\quad\quad
 \delbar\del\log(1-|z|^2)
=\frac{dz\cdot d\zbar}{(1-|z_i|^2)^2}.
\]
For a metric $h$, we put as follows:
\[
 h_{\veca,N}:=
 h\cdot\exp(-\tau(\veca,N))
=h\times
 \prod_{i=1}^l |z_i|^{-2a_i}(-\log|z_i|)^{-N}
 \times
\prod_{i=l+1}^n
 (1-|z_i|^2)^{-N}.
\]
We use the notation
$|\cdot|_{\veca,N}$,
$||\cdot||_{\veca,N}$,
$(\cdot,\cdot)_{\veca,N}$ and
$\doublelangle\cdot,\cdot\doublerangle_{\veca,N}$
instead of
$|\cdot|_{h_{\veca,N}}$,
$||\cdot||_{h_{\veca,N}}$,
$(\cdot,\cdot)_{h_{\veca,N}}$ and
$\doublelangle\cdot,\cdot\doublerangle_{h_{\veca,N}}$
for simplicity.
We also use the notation $A_{\veca,N}^{p,q}(E)$
instead of $A_{h_{\veca,N}}^{p,q}(E)$.

The $(1,0)$-part of the metric connection for the $h_{\veca,N}$
is denoted by $\del_{\veca,N}$.
The curvature is denoted by $R(h_{\veca,N})$.
We have the following formula for
$R(h_{\veca,N})$:
\begin{equation} \label{eq;11.28.7}
 R(h_{\veca,N})=
 R(h)-\delbar\del\tau(\veca,N)
=R(h)-
 N\cdot
 \left(
 \sum_{i=1}^l
 \frac{dz_i\cdot d\zbar_i}
 {(-\log|z_i|^2)^2|z_i|^2}
+\sum_{i=l+1}^n
\frac{dz_i\cdot d\zbar_i}{(1-|z_i|^2)^2}
 \right).
\end{equation}
In particular, we have the equality
$\sqrt{-1}\Lambda (R(h_{\veca,N}))
=\sqrt{-1}\Lambda(R(h))-n\cdot N$.

We will use the following later.
\begin{lem}\label{lem;12.11.3}
When $N$ is sufficiently smaller than $0$,
then the following inequality holds for any 
$(0,1)$-form $\eta$:
\[
 \doublelangle
 R(h_{\veca,N})+\Ric(g),\eta
 \doublerangle_{\veca,N}
\geq |\eta|_{\veca,N}^2.
\]
\end{lem}
\pf
It immediately follows from the equality (\ref{eq;11.28.7})
and the local description (\ref{eq;11.28.8}).
\hfill\qed

\begin{rem}
The real number $N$ has no effect to the increasing order
of holomorphic sections at the divisors.
Namely,
if we have a holomorphic section $s$,
the increasing order with respect to $h_{\veca,N}$
is independent of $N$.
On the contrary,
$\veca$ has an effect.
\hfill\qed
\end{rem}

\subsubsection{When $N$ is sufficiently larger than $0$.}
Let see the case $N$ is sufficiently larger than $0$.
The function $\tau(\veca,N)$ diverges at the boundary of $Y$.
Thus we consider
$Y(C):=\Delta^{\ast}(C)^l\times \Delta(C)^{n-l}$
for a number $0<C<1$.
The projection of
$\Delta^{\ast}(C)^l\times \Delta(C)^{n-l}
\lrarr
 \Delta^{\ast}(C)^{l-1}\times \Delta(C)^{n-l}$,
forgetting $j$-th component $(1\leq j\leq l)$,
is denoted by $\pi_j$.
For any element $p\in\Delta^{\ast}(C)^{l-1}\times \Delta(C)^{n-l}$,
we obtain the curve $\pi_j^{-1}(p)\subset Y(C)$,
which is isomorphic to the punctured disc.

Let $s$ be a holomorphic section of $E$ of $Y(C)$.
Let consider the restriction of $F$ to
the curve $\pi_j^{-1}(p)$.
We denote the restriction by $F_{\pi_j^{-1}(p)}$.
On $\pi_j^{-1}(p)$,
we have the following formula:
\begin{multline}
 \Delta''_{\poin}|F_{|\pi_j^{-1}(p)}|_{\veca,N}
=-\big|\del_{\veca,N}F_{|\pi_j^{-1}(p)}\big|^2_{\veca,N}
+\big(F_{|\pi_j^{-1}(p)},
 \sqrt{-1}\Lambda R(h_{\veca,N})\cdot
 F_{|\pi_j^{-1}(p)}
 \big)_{\veca,N}\\
=
-\big|\del_{\veca,N}
 F_{|\pi_j^{-1}(p)}
 \big|^2_{\veca,N}
+\big(F_{|\pi_j^{-1}(p)},
 \sqrt{-1}\Lambda R(h)\cdot 
 F_{|\pi_j^{-1}(p)}
 \big)_{\veca,N}
-n\cdot N\cdot |F_{\pi_j^{-1}(p)}|^2_{\veca,N}.
\end{multline}
Here $\Delta''_{\poin}$ denotes the Laplacian
for the punctured disc $\pi_j^{-1}(p)$ with the Poincare metric.
Then we obtain the following inequality:
\[
 \Delta''_{\poin} \log|F_{\pi_j^{-1}(p)}|_{\veca,N}^2
\leq
 \frac{\Delta''_{\poin}|F_{|\pi_j^{-1}(p)}|^2_{\veca,N}}
  {|F_{|\pi_j^{-1}(p)}|^2_{\veca,N}}
+\frac{|(\del_{\veca,N}
 F_{|\pi_j^{-1}(p)},
 F_{|\pi_j^{-1}(p)})_{\veca,N}|^2}
 {|F_{|\pi_j^{-1}(p)}|^4_{\veca,N}}
\leq
 -n\cdot N+ |R|_{h,g_{\poin}}
\]
If $(E,\delbar_E,h)$ is acceptable,
then $|R|_{h,g_{\poin}}$ is bounded by definition,
and thus we obtain the following inequality
for any holomorphic section $F$ and for any sufficiently large $N$:
\[
 \Delta''_{\poin}\log|F_{|\pi_j^{-1}(p)}|^2_{\veca,N}<0.
\]
We can take such $N$ independently of $p$, $j$ and $F$.
We can replace $\Delta''_{\poin}$
with the standard Laplacian $-(\del_s^2+\del_t^2)$
for the real coordinate $z_j=s+\sqrt{-1}t$.

As an example, we obtain the following corollary.
It says, we obtain the increasing order of
a holomorphic section over $\Delta^{\ast\,l}\times \Delta^{n-l}$
from the increasing order of the restriction
to the curves.
\begin{cor} \label{cor;11.28.15}
Let $F$ be a holomorphic section of $\nbige^{\lambda}$.
Let $a_j$ and $k_j$ be real numbers $(j=1,\ldots,l)$.

For any point $p\in\Delta^{\ast}(C)^{l-1}\times \Delta(C)^{n-l}$
and any $1\leq j\leq l$,
we assume that we are given numbers
$C_1(p,j)$, $C_2(p,j)$, $a(p,j)$ and $k(p,j)$
satisfying the following:
\begin{enumerate}
\item $C_1(p,j)$ and $C_2(p,j)$ are positive numbers.
\item
 $a(p,j)$ and $k(p,j)$ are real numbers satisfying
$a(p,j)\leq a_j$ and $k(p,j)\leq k_j$.
\item
The following inequality holds on $\pi_j^{-1}(p)$:
\[
 0<C_1(p,j)\leq |F_{|\pi_j^{-1}(p)}|_h\cdot |z_j|^{-a(p,j)}\cdot
 (-\log|z_j|)^{-k(p,j)}
 \leq C_2(p,j).
\]
\item $C_1(p,j)$, $C_2(p,j)$, $a(p,j)$ and $k(p,j)$ may depend on $p$ and $j$.
\end{enumerate}
Then there exists a positive constant $C_3$ and a large number $M$,
satisfying the following:
\begin{itemize}
\item
 The inequality
 $|s|_h\leq C_3\cdot \prod_{j=1}^l|z_j|^{a_j}(-\log |z_j|)^M$
 holds over $Y(C)$.
\item
$C_3$ depends only on the values of $|s|_h$
at $\{(z_1,\ldots,z_n)\,|\,|z_j|=C,\,\,j=1,\ldots,l\}$.
\end{itemize}
\end{cor}
\pf
We put $\veca=(0,\ldots,0)$.
Let pick $p\in \Delta^{\ast}(C)^{l-1}\times \Delta(C)^{n-l}$
and $1\leq j\leq l$.
We note that the following holds:
By our assumption, we have the following for any number $N$:
\[
 \lim_{|z_j|\to 0}
 \frac{\log|F_{|\pi_j^{-1}(p)}|_{\veca,N}-a\log|z_j| }
 {\log|z_j|}=0.
\]
If $N$ is sufficiently large,
then we obtain the following inequality
on $\pi_j^{-1}(a)$:
\[
 \Delta(\log|F_{|\pi_j^{-1}(p)}|_{\veca,N}-a\log|z_j|)\leq 0
\]
Here $\Delta$ denotes the standard Laplacian
$-(\del_s^2+\del_t^2)$ on $\pi^{-1}_j(p)$.
Then the values of
$\log|F_{|\pi_j^{-1}(p)}|-a\log|z_j|$ on $\pi_j^{-1}(p)$
are dominated by the values at $\{|z_j|=C\}\cap \pi_j^{-1}(p)$.
(See \cite{s2}, in particular,
Lemma 2.2 and the proof of Corollary 4.2).
Thus we obtain the following inequality:
\begin{equation} \label{eq;11.28.10}
 |F_{|\pi_j^{-1}(p)}|_{\veca,N}
\leq C(p,j)\cdot |z_j|^a.
\end{equation}
Here $C(p,j)$ denotes a constant depending only 
on the values $F$ 
at $\{|z_j|=C\}\cap \pi_j^{-1}(p)$:
\[
 C(p,j)=\max\Bigl\{|F_{|\pi_j^{-1}(p)}(z_j)|_{\veca,N}\,
 \Big|\,z_j\in\cnum,
 \,|z_j|=C\Bigr\}.
\]
From the inequality (\ref{eq;11.28.10}),
we obtain the following on $\pi_j^{-1}(p)$:
\[
 |F_{|\pi_j^{-1}(p)}|_h
\leq C(p,j)\cdot |z_j|^{a_j}\cdot (-\log|z_j|)^{N}.
\]
Then we obtain the result by using an induction on $l$.
\hfill\qed

\vspace{.1in}
We also have the following:
\begin{lem} \label{lem;12.8.15}
Let $(E,h)$ be a hermitian holomorphic bundle
over $\Delta^{\ast\,l}\times \Delta^{n-l}$.
Let $F$ be a holomorphic section of $E$
over $\Delta^{\ast\,l}\times \Delta^{n-l}$.
Assume that $||\,F\,||_{\veca,N}<\infty$.
If $N$ is sufficiently large,
then $||\,F_{|\pi_j^{-1}(p)}||_{\veca,N}<\infty$
for any $j$ and $p\in \Delta^{\ast\,l-1}\times \Delta^{n-l}$.
\end{lem}
\pf
We only have to use the subharmonicity of 
the function $|\,F_{|\pi_j^{-1}(a)}|_{\veca,N}$.
\hfill\qed

\subsection{Some reductions toward 
the prolongation of the deformed holomorphic bundles}
\label{subsection;12.15.16}

\subsubsection{The prolongation of $\nbige^{\lambda}$ when $\lambda\neq 0$}
\label{subsubsection;12.8.6}

Let $\harmonicbundle$ be a tame nilpotent harmonic bundle
with trivial parabolic structure over
$\Delta^{\ast\,l}\times\Delta^{n-l}$.
Let consider the prolongation of $\nbige^{\lambda}$
over $\Delta^{\ast\,l}\times \Delta^{n-l}$
by increasing order $(0,\ldots,0)$.

Since the eigenvalues of the monodromies are $1$,
we can take a normalizing frame $\vecw$ of $\nbige^{\lambda}$
over $\Delta^{\ast\,l}\times \Delta^{n-l}$,
namely,
the connection form $A$ of $\DD^{\lambda,f}$
with respect to $\vecw=(w_1,\ldots,w_r)$ 
satisfies the following:
\begin{itemize}
\item $A$ is of the form $A=\sum_{j=1}^l A_j\cdot dz_j/z_j$
 for some constants matrices $A_j\in M(r)$.
\item
 All of the eigenvalues of $A_j$ are $0$.
\end{itemize}
Note that the restriction $\vecw_{|\pi_j^{-1}(p)}$
of $\vecw$ to the curves $\pi_j^{-1}(p)$
gives a normalizing frame
of $(\nbige^{\lambda}_{|\pi_j^{-1}(p)},\DD^{\lambda}_{|\pi_j^{-1}(p)})$.
Thus $\vecw_{|\pi_j^{-1}(p)}$ is a frame of
$\prolong{\nbige^{\lambda}_{|\pi_j^{-1}(p)}}$,
due to Lemma \ref{lem;12.15.200}.
In particular,
the restrictions $w_{i\,|\,\pi_j^{-1}(p)}$
satisfies the conditions of Corollary \ref{cor;11.28.15}.
Thus we can conclude that 
the increasing order of $w_{i}$
is less than $(0,\ldots,0)$,
i.e.,
$w_{i}$ is a holomorphic section of
$\prolong{\nbige^{\lambda}}$.

On the other hand,
let consider the following section of
$\det(\nbige^{\lambda})$ over $\Delta^{\ast\,l}\times \Delta^{n-l}$:
\[
 \volume(\vecw)=w_1\wedge\cdots\wedge w_r.
\]
By our choice, we know that
$\DD^{\lambda}(\volume(\vecw))=0$.
Thus it naturally gives a holomorphic frame
of $\prolong{\det(\nbige^{\lambda})}$
over $\Delta^n$.

\begin{lem} \label{lem;11.28.20}
Let $\vecs=(s_1,\ldots,s_r)$ be a tuple of holomorphic sections
of the sheaf $\prolong{\nbige^{\lambda}}$.
\begin{enumerate}
\item
The section $\Omega(\vecs)=s_1\wedge\cdots\wedge s_r$ over
$\Delta^{\ast\,l}\times \Delta^{n-l}$ naturally induces
the section of $\prolong{\det(\nbige^{\lambda})}$.
We denote the section by the same notation.
\item
Assume that $\Omega(\vecs)(P)\neq 0$ in
$\prolong{\det(\nbige^{\lambda})}_{|P}$.
Then $\prolong{\nbige^{\lambda}}$ is locally free around $P$,
and $\vecs$ gives a holomorphic frame around $P$.
\end{enumerate}
\end{lem}
\pf
The first claim is clear, for we only have to see
the increasing order of $\Omega(\vecs)$.

Let assume that  $\Omega(\vecs)(P)\neq 0$.
Then we have an open subset $U$ of $\Delta^n$
such that $\Omega(\vecs)$ gives a frame of
$\prolong{\det(\nbige^{\lambda})}$
over $U$.
On the open set $U'=U\cap (\Delta^{\ast\, l}\times \Delta^{n-l})$,
the tuple $\vecs$ gives a holomorphic frame of $\nbige^{\lambda}$.
Let $f$ be an element of $\Gamma(U,\prolong{\nbige^{\lambda}})$.
Since we have already known that $\vecs$ gives a frame on $U'$,
we have a holomorphic functions $f_i$
defined over $U'$ satisfying the following over $U'$:
\[
 f=\sum f_i\cdot s_i.
\]
Let see $f_1$.
We know that
$f\wedge s_2\wedge\cdots s_r=
 f_1\cdot\Omega(\vecs)$ gives a section
of $\prolong{\det(\nbige^{\lambda})}$ on $U$.
We also know that $\Omega(\vecs)$ gives a frame of
$\prolong{\det(\nbige^{\lambda})}$ on $U$.
Thus we can conclude that $f_1$ is in fact a holomorphic function
over $U$.
Similarly we can show that the other $f_i$ are also
holomorphic over $U$.
\hfill\qed

As a corollary, we obtain the following:
\begin{cor}
When $\lambda\neq 0$,
then $\prolong{\nbige^{\lambda}}$ is locally free coherent sheaf,
and the normalizing frame gives a frame.
\hfill\qed
\end{cor}

\begin{cor}
The $\lambda$-connection $\DD^{\lambda}$
and the flat connection $\DD^{\lambda,f}$
is of log type on $\prolong{\nbige^{\lambda}}$.
\hfill\qed
\end{cor}

\begin{rem}
The argument of Lemma {\rm\ref{lem;11.28.20}},
which we learn at {\rm\cite{cg}},
will be used without mention.
\hfill\qed
\end{rem}

Let $\vecv$ be a holomorphic frame of $\prolong{\nbige^1}$
over $\Delta^n$.
Let consider the $\nbigh(r)$-valued function
$H(h,\vecv)$.
\begin{lem} \label{lem;12.10.30}
There exists a large number $M$ and positive numbers $C_i$ $(i=1,2)$
satisfying the following:
\begin{equation} \label{eq;12.8.1}
 C_1\cdot
 \prod_{i=1}^l(-\log|z_i|)^{-M}
\leq
 H(h,\vecv)
\leq
 C_2\cdot
 \prod_{i=1}^l(-\log|z_i|)^M.
\end{equation}
\end{lem}
\pf
We only have to consider the case that $\vecv$
is a normalizing frame.
We have already known the following estimate
for some positive number $C_3$ and $M_1$:
\begin{equation} \label{eq;12.8.2}
 |v_i|\leq C_3\cdot \prod_{i=1}^l(-\log|z_i|)^{M_1}
\end{equation}
Thus we obtain the right inequality in (\ref{eq;12.8.1}).
To see the left inequality,
we use the dual.
Namely let $\vecv^{\lor}$ denote the dual frame
of $\vecv$.
Then $\vecv^{\lor}$ is also a normalizing frame
of $\nbige^{1\,\lor}$.
Hence it satisfies the inequality 
similar to (\ref{eq;12.8.2}).
Hence we obtain the right inequality for $H(h^{\lor},\vecv^{\lor})$.
Since we have $H(h^{\lor},\vecv^{\lor})\cdot H(h,\vecv)=1$,
we obtain the result.
\hfill\qed

\begin{rem}
The argument of the proof of Lemma {\rm\ref{lem;12.10.30}}
works for any $\lambda\neq 0$.
If $\lambda=0$, we do not have a normalizing frame, in general.
However, the argument works once we know that the dual frame
$\vecv^{\lor}$ of $\nbige^{0\,\lor}$ gives a frame
of $\prolong{\nbige^{0\,\lor}}$.
\hfill\qed
\end{rem}

\subsubsection{The prolongation of $\nbige^{\shikaku}$}
\label{subsubsection;12.15.17}

Let $\nbigo\bigl(-\{\lambda_0\}\bigr)$ denote
the sheaf of the holomorphic functions
which vanish on the divisor $\{\lambda=\lambda_0\}$.
The sheaf $\nbigo\bigl(-\{\lambda_0\}\bigr)$ is a line bundle.
We have the natural inclusion:
\[
 \nbigo\bigl(-\{\lambda_0\}\bigr)\subset \nbigo,
\mbox{ and thus }
\nbige\otimes\nbigo\bigl(-\{\lambda_0\}\bigr)\subset  \nbige.
\]

\begin{prop}\mbox{{}} \label{prop;11.28.21}
Let $\harmonicbundle$ be a tame nilpotent harmonic bundle
with trivial parabolic structure
over $\Delta^{\ast\,l}\times \Delta^{n-l}$.
\begin{enumerate}
\item
The prolongment of the deformed holomorphic bundle
$\prolong{\nbige^{\shikaku}}$ over $\cnum_{\lambda}^{\ast}\times\Delta^{n}$
is locally free.
\item
If $\prolong{\nbige^0}$ over $\Delta^n$
is locally free,
then $\prolong{\nbige}$ over $\cnum_{\lambda}\times\Delta^n$
is locally free.
\end{enumerate}
\end{prop}
\pf
Let pick $\lambda_0\in \cnum_{\lambda}$.
We assume that $\lambda_0\neq 0$ or that we know that $\nbige^{0}$ is
locally free.
We put
$\Delta_0:=\{\lambda\in\cnum_{\lambda}\,|\,
 |\lambda-\lambda_0|<1 \}$.
We can naturally identify $\Delta_0$ with $\Delta$.
We denote the projection of
$\Delta_0\times \Delta^{\ast\,l}\times \Delta^{n-l}\lrarr 
 \Delta^{\ast\,l}\times \Delta^{n-l}$ by $p_{\lambda}$.
We denote the restriction of $\nbige$ to
$\Delta_0\times\Delta^{\ast\,l}\times \Delta^{n-l}$
by the same notation $\nbige$.

Let take a holomorphic frame $\vecv=(v_1,\ldots,v_n)$
of $\prolong{\nbige^{\lambda_0}}$ over $\Delta^n$.
The restriction of $\vecv$
to $\Delta^{\ast\,n}$ is denoted by $\vect=(t_1,\ldots,t_n)$.
We have the $C^{\infty}$-frame $p_{\lambda}^{-1}(\vect)$
of $\nbige$ over $\Delta_0\times\Delta^{\ast\,l}\times \Delta^{n-l}$.

Let $d''$ denotes the holomorphic structure of $\nbige$.
We have the following:
\[
 d''(p_{\lambda}^{-1}t_i)=
 (\lambda-\lambda_0)\cdot
 p_{\lambda}^{-1}
\bigl(
 \theta^{\dagger}(t_i)
\bigr)=:\eta_i.
\]

Let take a sufficiently negative number $N$.
Let $\epsilon$ be any sufficiently small positive number.
Let $\vecdelta$ denotes the tuple
$\overbrace{(1,\ldots,1)}^l$.
Due to the estimate of $\theta^{\dagger}$
and the assumption on the increasing order of $t_i$,
the $(0,1)$-form $\eta_i$ is 
an element of the following space:
\[
 A^{0,1}_{-\epsilon\cdot\vecdelta,N}
 \Bigl(
 \nbige\otimes \nbigo\bigl(\{-\lambda_0\}\bigr)
 \Bigr).
\]
Since $N$ is sufficiently negative,
We can find an element $\rho_i$ satisfying the following
(Proposition \ref{prop;11.28.6},
 Proposition \ref{prop;12.11.1},
 Corollary \ref{cor;12.11.2},
 Lemma \ref{lem;12.11.3}):
\[
 \rho_i\in
 A^{0,0}_{-\epsilon\cdot\vecdelta,N}
 \Bigl(
 \nbige\otimes\nbigo\bigl(-\{\lambda_0\}\bigr)
 \Bigr),
\quad
 \delbar_E\rho_i=\eta_i.
\]
Then we put $\tilde{t}_i:=p_{\lambda}^{-1}(t_i)-\rho_i$.
By our construction,
it is easy to check the following:
\begin{itemize}
\item
$\tilde{t}_i$ is holomorphic.
\item
The restriction to
$\{\lambda_0\}\times \Delta^{\ast\,l}\times\Delta^{n-l}$
is same as $t_i$.
\item
The increasing order of $\tilde{t}_i$
is less than $-\epsilon\cdot\vecdelta$.
\end{itemize}

\begin{lem}
The increasing order of $\tilde{t}_i$ is less than
$(0,\ldots,0)$.
\end{lem}
\pf
Let consider the restriction of
$\{\lambda\}\times\tilde{t}_i$
to $\pi_j^{-1}(p)$ for any $p\in \Delta^{\ast\,l-1}\times \Delta^{n-l}$
and for any $1\leq j\leq l$.
We know that the parabolic structure of
$\prolong{\nbige^{\lambda}_{|\pi_j^{-1}(p)}}$
is trivial.
Since $\epsilon$ is sufficiently small,
we obtain that the increasing order of
the restriction of $\tilde{t}_i$ to $\{\lambda\}\times\pi_j^{-1}(a)$
is less than $0$.

Then we can apply Corollary \ref{cor;11.28.15}.
(We put $a_j=0$ in applying.)
Then we obtain that the increasing order of $\tilde{t}_i$
is less than $(0,\ldots,0)$.
\hfill\qed

Thus $\tilde{\vect}=(\tilde{t}_i)$ are holomorphic sections
of $\prolong{\nbige}$ around $(\lambda_0,O)$.
Here $O$ denotes the origin $(0,\ldots,0)$ of $\Delta^n$.
By our choice,
$\Omega(\tilde{\vect})(\lambda_0,O)=
 \Omega(\vect)(\lambda_0,O)$
is not $0$ in $\prolong{\det(\nbige)}_{|(\lambda_0,O)}$.
By the same argument to the proof of Lemma \ref{lem;11.28.20},
we obtain that $\nbige$ is locally free
on a neighborhood $U$ of $(\lambda_0,O)$,
and that the $\tilde{\vect}$ gives a holomorphic frame
on $U$.
Thus the proof of Proposition \ref{prop;11.28.21}
is completed.
\hfill\qed

\subsubsection{One dimensional case, revisited}
\label{subsubsection;12.15.18}

\begin{cor}
Let $\harmonicbundle$ be a tame nilpotent harmonic bundle
with trivial parabolic structure over
the punctured disc $\Delta^{\ast}$.
Then the prolongment $\prolong{\nbige}$
of the deformed holomorphic bundle is locally free.
The $\lambda$-connection $\DD^{\lambda}$ is of log type.
\end{cor}
\pf
We have already known that
the prolongment of $\nbige^0$ is coherent locally free.
Thus $\prolong{\nbige}$ is locally free due to Proposition
\ref{prop;11.28.21}.

Let $f$ be a holomorphic section of $\prolong{\nbige}$.
Let consider $\DD(f)$.
It is holomorphic over $\nbigx-\nbigd$.
Since a normalizing frame gives a frame of $\prolong{\nbige^{\shikaku}}$
over $\nbigx^{\shikaku}$,
$\DD(f)$ is holomorphic as a section of
$\prolong{\nbige}^{\shikaku}\otimes
 p_{\lambda}^{\ast}\Omega_{X}^{1,0}$
over $\nbigx^{\shikaku}$.
Then $\DD(f)$ is a holomorphic section
of locally free sheaf
$\prolong{\nbige}\otimes p_{\lambda}^{\ast}\Omega_X(\log D)$,
outside the codimension two subset $\nbigd^0$.
Thus $\DD(f)$ is a holomorphic section of
$\prolong{\nbige}\otimes p_{\lambda}^{\ast}\Omega_X(\log D)$,
over $\nbigx$,
namely $\DD$ is of log type.
\hfill\qed

We have already known that
the conjugacy classes of $\resddlambda$ is independent of $\lambda$.
Thus the weight filtrations of $\resddlambda$ form
the filtration of $\prolong{\nbige}_{|\cnum_{\lambda}\times O}$.
We denote it by $\weightfilt$.
The associated graded vector bundle over $\cnum_{\lambda}\times\{O\}$
is denoted by $\graded=\bigoplus_h \graded_h$.
We denote the primitive part of $\graded_h$
by $P_{h+2a}\graded_h$.
\begin{lem}\label{lem;12.1.25}
There exists a holomorphic frame $\vecv$ of
$\prolong{\nbige}$ over $\cnum_{\lambda}\times \Delta$,
satisfying the following:
\begin{itemize}
\item
 We have the $M(r)$-valued holomorphic function $J(\lambda)$
 over $\cnum_{\lambda}\times \{O\}$,
 such that
 $\resddlambda\cdot\vecv_{|\cnum_{\lambda}\times O}=
 \vecv_{|\cnum_{\lambda}\times O}\cdot J(\lambda)$.
 Then $J(\lambda)$ is, in fact, a constant matrix $J$.
\end{itemize}
\end{lem}
\pf
For any $h\geq 0$,
we have the surjective morphism of vector bundles:
\[
 \pi_h:
 \Ker\bigl(\resddlambda^{h+1}\bigr)\lrarr
 P_{h}\graded_h.
\]
We take a splitting $\phi_h$ of $\pi_h$.
We denote the image of $\phi_h$
by $P_h\graded'_h$,
which is a subbundle of $\prolong{\nbige}_{\cnum_{\lambda}\times O}$.
For any $0\leq a\leq h$,
we put $P_h\graded'_{h-2a}:=\resddlambda^{a}(P_h\graded'_h)$.
Note that $\resddlambda^{h+1}(P_h\graded'_h)=0$
due to our choice of $P_h\graded'_h$.
Then we obtain a decomposition:
\[
 \prolong{\nbige}_{|\cnum_{\lambda}\times O}=
 \bigoplus_h\bigoplus_{a\geq 0}
 P_{h+2a}\graded'_{h}.
\]
We take a holomorphic frame $\vecu_h$ of
$P_h\graded'_h$.
Then we obtain the holomorphic frame of
$\prolong{\nbige}_{|\cnum_{\lambda}\times O}$
given as follows:
\[
 \bigcup_h 
\left(
 \bigcup_{a=0}^h
 \resddlambda^{a}(\vecu_h)\right).
\]
We extend the frame to a holomorphic frame of
$\prolong{\nbige}$ over $\cnum_{\lambda}\times \Delta$,
which is a desired frame.
\hfill\qed

Let $\vecv$ be a holomorphic frame of $\prolong{\nbige}$
over $\cnum_{\lambda}\times \Delta$,
compatible with the filtration $\weightfilt$.
We put $2\cdot k(v_i):=\deg^{W}(v_i)$.
Then we put $v_i':=v_i\cdot (-\log|z|)^{-k(v_i)}$,
and we obtain the $C^{\infty}$-frame
$\vecv'=(v_i')$ over $\cnum_{\lambda}\times \Delta^{\ast}$.

\begin{prop} \label{prop;12.3.70}
For any compact subset $K$ of $\cnum_{\lambda}$,
the frame $\vecv'$ is adapted over $K\times \Delta^{\ast}$.
\end{prop}
\pf
We can assume that $\vecv$ satisfies the condition
in Lemma \ref{lem;12.1.25}.
We take a model bundle $E(V,J)=(E_0,\theta_0,h_0)$.
We denote the deformed holomorphic bundle
by $(\nbige_0,\DD_0,h_0)$.
We have the canonical frame $\vecv_0$ of $\prolong{\nbige_0}$,
such that $\DD_0\cdot\vecv_0=\vecv_0\cdot J\cdot dz/z$.
Then the frames $\vecv$ and $\vecv_0$
induce the holomorphic isomorphism
$\Phi:\prolong{\nbige_0}\lrarr \prolong{\nbige}$.
Then we have the following:
\begin{equation} \label{eq;12.1.26}
 \Phi\circ\Res(\DD_0)-\Res(\DD)\circ\Phi=0,
\quad
 \mbox{ on }
 \cnum_{\lambda}\times \{O\}.
\end{equation}
We only have to show the boundedness of $\Phi$
and $\Phi^{-1}$
over $K\times \Delta^{\ast}$
for any compact subset $K\subset \cnum_{\lambda}$.

Note that $\Phi$ is a holomorphic section
of $Hom(\nbige_0,\nbige)$,
which is a deformed holomorphic bundle of
$(Hom(E_0,E),-\theta_0\otimes 1+1\otimes\theta,h_0\otimes h)$.
Due to the Lemma \ref{lem;12.8.5} below,
we have the inequality
$\Delta''_{z}|\Phi|_{h_0\otimes h}\leq
 |\DD(\Phi)|_{h_0\otimes h}$.
Due to (\ref{eq;12.1.26}),
$\DD(\Phi)$ is holomorphic over $\cnum_{\lambda}\times\Delta$,
Thus we have the following inequality over $K\times \Delta^{\ast}$
for any compact subset $K\subset \cnum_{\lambda}$:
\[
 \Delta''_{z}|\Phi|_{h_0\otimes h}\leq
 |\DD(\Phi)|_{h_0\otimes h}
 \leq C\cdot r^{-1+\epsilon}.
\]
Here $\Delta''_z=-(\del_t^2+\del_s^2)$ for a real coordinate
$z=s+\sqrt{-1}t$.
Hence we obtain the estimate,
which is locally uniform for $\lambda$.
\hfill\qed

\begin{lem}\label{lem;12.8.5}
Let $\harmonicbundle$ be a tame nilpotent harmonic bundle
over $\Delta^{\ast}$.
Let $f$ be a holomorphic section of $\nbige^{\lambda}$.
Then we obtain the following inequality:
\[
 \Delta''_z|f|_h\leq |\DD^{\lambda}f|_h.
\]
Here $\Delta''_z=-(\del_t^2+\del_s^2)$ for a real coordinate
$z=s+\sqrt{-1}t$.
\end{lem}
\pf
We denote the metric connection of $\nbige^{\lambda}$
by $d^{\twoprime\,\lambda}+d^{\prime\,\lambda}$.
We have the following equality:
\[
 \del\delbar|f|_h^2=
 (d^{\prime\,\lambda}f,d^{\prime\,\lambda}f)_h
+(f,R(d^{\twoprime\,\lambda})f)_h.
\]
We also have the following equality:
\[
 R(d^{\twoprime\,\lambda})
=\delbar_E\del_E+\del_E\delbar_E
-|\lambda|^2
 \cdot(\theta^{\dagger}\cdot\theta+\theta\cdot\theta^{\dagger})
=-(1+|\lambda|^2)\cdot
 \bigl(\theta^{\dagger}\cdot\theta+\theta\cdot\theta^{\dagger}\bigr).
\]
Thus we obtain the following:
\[
 \del\delbar|f|_h^2=
 (d^{\prime\,\lambda}f,d^{\prime\,\lambda}f)_h
-\bigl(1+|\lambda|^2\bigr)
 \cdot
 \bigl(\theta f,\theta f\bigr)_h
-\bigl(1+|\lambda|^2\bigr)
 \cdot
 \bigl(\theta^{\dagger}f,\theta^{\dagger}f
 \bigr)_h.
\]
We also have the following:
\[
 \bigl(
 d^{\prime\,\lambda}f,d^{\prime\,\lambda}f
 \bigr)_h
+(\DD^{\lambda}f,\DD^{\lambda}f)_h
-(1+|\lambda|^2)\cdot(\theta f,\theta f)_h
=-(1+|\lambda|^2)\cdot(\del_Ef,\del_Ef)_h.
\]
Thus we obtain the following:
\[
 \del\delbar|f|_h^2=
 (1+|\lambda|^2)\cdot (\del_Ef,\del_Ef)_h
-(\DD^{\lambda}f,\DD^{\lambda}f)_h
-(1+|\lambda|^2)\cdot(\theta^{\dagger}f,\theta^{\dagger}f)_h.
\]
Then we obtain the following:
\[
 \Delta''|f|^2_h
=-(1+|\lambda|^2)\cdot |\del_Ef|_h^2
+|\DD^{\lambda}f|^2_h
-(1+|\lambda|^2)\cdot |\theta^{\dagger}f|_h^2
\leq |\DD^{\lambda}f|_h^2.
\]
Thus we are done.
\hfill\qed

We have a $\lambda$-family version of Lemma \ref{lem;1.14.1}.
Let $\vecv$ be a holomorphic frame of $\prolong{\nbige}$
over $\cnum_{\lambda}\times \Delta$,
compatible with the filtration $\weightfilt$.
We have the vector space $V=\prolong{\nbige}_{|(0,O)}$
and the endomorphism $N=\Res(\DD)_{|(0,O)}$.
Then we obtain the model bundle $E(V,N)=(E_0,\theta_0,h_0)$.
We denote the deformed holomorphic bundle by
$(\nbige_0,\DD_0,h_0)$. We have the canonical frame $\vecv_0$.

By the frames $\vecv_{|\{0\}\times \Delta^{\ast}}$
and $\vecv_{0|\{0\}\times \Delta^{\ast}}$,
we obtain the holomorphic isomorphism
$\nbigi:\nbige^0_0\lrarr \nbige^0$
over $\{0\}\times \Delta$.
Then we obtain the $C^{\infty}$-isomorphism
$\nbigi:\nbige_0\lrarr \nbige$ defined over
$\cnum_{\lambda}\times \Delta^{\ast}$.
Then we obtain the elements
$I_{i\,j}(\lambda,z)\in
 C^{\infty}(\cnum_{\lambda}\times \Delta_z^{\ast})$,
determined as follows:
\[
 \nbigi(v_{0\,j})=\sum I_{i\,j}(\lambda,z)\cdot v_i,
\quad
 \nbigi(\vecv_0)=\vecv\cdot I.
\]
Note that $I_{i\,j}$ are holomorphic along the direction of $\lambda$,
although it is not holomorphic along the direction of $z$.
The following lemma is clear from the proof
of Lemma \ref{lem;1.14.1}.
\begin{lem} \label{lem;12.2.10}
The finiteness in Lemma {\rm\ref{lem;1.14.1}} is locally uniform
on $\cnum_{\lambda}$.
Namely 
the values of integrals {\rm(\ref{eq;1.14.1})}
and {\rm(\ref{eq;1.14.2})} are bounded on
any compact subset $K\subset \cnum_{\lambda}$.
\hfill\qed
\end{lem}

\subsection{Extension of holomorphic sections over hyperplanes}

\label{subsection;12.15.19}

\subsubsection{Preliminary}
Let $\harmonicbundle$ be a tame nilpotent harmonic bundle
with the trivial parabolic structure
over $\Delta_{\zeta}\times \Delta_z^{\ast n-1}=
 \{(\zeta,z_1,\ldots,z_{n-1})\}$.
We would like to show that the prolongment of
$(E,h)=(\nbige^{0},h)$ is locally free.
We note that we have already known that
the prolongment $\prolong{\nbige}^1$ is locally free
(the subsubsection \ref{subsubsection;12.8.6}).

We put $X-D:=\Delta_{\zeta}\times \Delta_{z}^{\ast\,n-1}$,
$X:=\Delta_{\zeta}\times \Delta_z^{n-1}$,
$X_0-D_0:=\{0\}\times \Delta_{z}^{n-1\,\ast}$,
and $X_0:=\{0\}\times \Delta_z^{n-1}$.

We assume the following as a hypothesis of the induction.
\begin{assumption} \label{assumption;11.30.2}
\begin{enumerate}
\item
 The prolongment of
 the restriction $\nbige^0_{|X_0-D_0}$
 is locally free on $X_0$.
\hfill\qed
\end{enumerate}
\end{assumption}
Recall that we have already known the following
(Lemma \ref{lem;12.10.30}):
\begin{condition}
\label{number;11.30.1}
 Let $\vecv$ be a frame of $\prolong{\nbige^1}$
 over $\Delta_{\zeta}\times\Delta_z^{\ast\,n-1}$.
 The components of the hermitian matrices
 $H(h,\vecv)$ and the inverse $H(h,\vecv)^{-1}$ 
 are dominated by polynomials of $-\log |z_i|$ for $i=1,\ldots,n-1$.
\hfill\qed
\end{condition}

By the relations
$ \theta\cdot\vecv=
 \vecv\cdot\Theta$ or
$\theta^{\dagger}\cdot\vecv=
 \vecv\cdot\Theta^{\dagger}$,
we have the following elements:
\[
\begin{array}{l}
 \Theta=
 \Theta^{\zeta}\cdot d\zeta
+\sum\Theta^k\cdot dz_k/z_k
 \in C^{\infty}(X,M(r)\otimes \Omega_{X-D}^{1,0}),\\
\mbox{{}}\\
\Theta^{\dagger}
 =\Theta^{\zeta\,\dagger}\cdot d\zetabar
+\sum\Theta^{k\,\dagger}\cdot d\zbar_k/\zbar_k
\in C^{\infty}(X,M(r)\otimes\Omega_{X-D}^{0,1}) .
\end{array}
\]
By the estimate for $\theta$ (Proposition \ref{prop;1.25.1}),
the absolute values of the components
of $\Theta^{\zeta}$, $\Theta^k$,
$\Theta^{\zeta\,\dagger}$ and $\Theta^{k\,\dagger}$
are bounded by the polynomials of
$-\log|z_i|$ for $i=1,\ldots,n-1$.

We have the flat connection form $A$ of $\DD^1$
with respect to the frame $\vecv$,
that is,
$\DD^1\vecv=(\del_E+\theta)\vecv=\vecv\cdot A$.
The form $A$ is a holomorphic section of
$M(r)\otimes \Omega^{1,0}_{X}(\log D)$.

We have the relation $\del_E(\theta)=0$.
It is translated as follows:
\begin{equation} \label{eq;12.8.10}
 \del\Theta^{\dagger}+
\bigl[
 A-\Theta,\Theta^{\dagger}
\bigr]=0.
\end{equation}
\begin{lem}\mbox{{}} \label{lem;12.14.100}
\begin{itemize}
\item
The components of
$\del_{\zeta}\Theta^{\dagger}$ is of the form
$\sum_k C_k\cdot d\zbar_k/\zbar_k\cdot d\zeta+
 C_{\zeta}\cdot d\zetabar\cdot d\zeta$.
Then the absolute values of $C_k$ and $C_{\zeta}$
are dominated by a polynomial of
$(-\log|z_i|)$ for $i=1,\ldots,n-1$.
\item
The components of $\del_k\Theta^{\zeta\,\dagger}$ is of the form
$C\cdot d\barz_k/\zbar_k\cdot d\zetabar$.
Then the absolute values of $C$ is dominated
by a polynomials of $(-\log|z_i|)$ for $i=1,\ldots,n-1$.
\item
Hence, for any $z=(z_1,\ldots,z_{n-1})\in X_0$,
$\Theta(\zeta,z)$ is $L_k^p$
as a function of $\zeta$.
And the $L_k^p$-norm is dominated by 
the polynomials of $(-\log|z_i|)$ for $i=1,\ldots,n-1$.
\end{itemize}
\end{lem}
\pf
The claims immediately follow from (\ref{eq;12.8.10}).
\hfill\qed

\subsubsection{A construction of a $(0,1)$-cocycle and the extension
   property I}
\label{subsubsection;12.15.201}

The restriction of the frame $\vecv$ to $\nbige^1$
over $X-D$ can be regarded as a $C^{\infty}$-frame
of $\nbige^0$.
Let $f$ be a holomorphic section of $\prolong{(\nbige^0_{|X_0-D_0})}$.
We have $C^{\infty}$-functions $f_i$ on $X_0-D_0$
determined by the relation:
\[
 f(z)=\sum_i f_i(z)\cdot v_i(0,z).
\]
Since an increasing order of $f$ is less than $0$,
we have some inequalities
$|f|_h\leq C_{\epsilon}\cdot
 \bigl(
 \prod_{j=1}^l |z_j|\bigr)^{-\epsilon}$
for any $\epsilon>0$.
Due to the condition \ref{number;11.30.1},
we have some inequalities
$|f_i|_h\leq C_{\epsilon}\cdot
 \bigl(
 \prod_{j=1}^l |z_j|\bigr)^{-\epsilon}$
for any $\epsilon>0$.
Since $f$ is holomorphic,
we have the following equality:
\[
 0=\delbar f(z)=
\sum
 \Bigl(
  \delbar f_i(z)
 +\sum_{k}\Theta^{k\dagger}_{i\,j}(0,z)\cdot \etabar_k f_j(z)
 \Bigr)\cdot v_i(0,z).
\]
Here we use the notation $\etabar_i=d\zbar_i/\zbar_i$.

We put as follows:
\[
 \begin{array}{l}
 f^1:=
 \sum_i f_i(z)\cdot v_i(\zeta,z) ,\quad\quad
 f^2:= \zetabar\cdot g^2,\\
\mbox{{}}\\
g^2:=
\sum_{i,j}
 \Theta^{\zeta\,\dagger}_{i\,j}(0,z)
 \cdot f_{j}(z)\cdot v_{i}(\zeta,z).
\end{array}
\]
We have the following equality:
\begin{multline}
\delbar f^1=
\sum
 \Bigl(
  \delbar f_i(z)+
  \Bigl(
    \sum _k \Theta^{k\,\dagger}_{i\,j}(\zeta,z)\etabar_k
    +\Theta^{\zeta\,\dagger}_{i\,j}(\zeta,z)d\zetabar
  \Bigr)\cdot f_j
 \Bigr)v_i(\zeta,z) \\
=
\sum
 \Bigl(
   \sum_k
      \Bigl(
    \Theta^{k\,\dagger}_{i\,j}(\zeta,z)-
    \Theta^{k\,\dagger}_{i\,j}(0,z)
      \Bigr)\cdot\etabar_k\
 \Bigr)\cdot f_j(z)\cdot v_{i}(\zeta,z)
+
 \sum \Theta^{\zeta\,\dagger}_{i\,j}(\zeta,z)\cdot d\zetabar
 \cdot f_j(z)\cdot v_i(\zeta,z).
\end{multline}
We also have the following equality:
\[
 \delbar f^2=
 \sum \Theta^{\zeta\,\dagger}_{i\,j}(0,z)
 \cdot d\zetabar
 \cdot f_j(z)
 \cdot v_i(\zeta,z)
+\zetabar\cdot \delbar g^2.
\]
Thus we obtain the equality
$\delbar(f_1-f_2)=\zeta\cdot F$,
when we put as follows:
\[
 F:=
 \zeta^{-1}\cdot
 \left[
 \sum
 \Bigl(
 \sum_k
  \Bigl(
    \Theta^{k\,\dagger}_{i\,j}(\zeta,z)-\Theta^{k\,\dagger}_{i\,j}(0,z)
  \Bigr)d\zbar_k
+\Bigl(
   \Theta^{\zeta\,\dagger}_{i\,j}(\zeta,z)
  -\Theta^{\zeta\,\dagger}_{i\,j}(0,z)
 \Bigr)d\zetabar
 \Bigr)\cdot f_j\cdot v_i(\zeta,z)
+\zetabar\cdot \delbar g^2
\right].
\]
By our construction,
we have the equality $\delbar F=0$,
and the norm $|F|_{h,g_{\poin}}$  is dominated
by a polynomial of $(-\log|z_i|)$ for $i=1,\ldots,l$.
(Use Lemma \ref{lem;12.14.100}.)
Thus we have the finiteness
for any $\epsilon>0$ and any number $M$:
\[
 \int |F|_h^2\cdot
 \prod_{j=1}^{n-1} 
 \Bigl(|z_j|^{\epsilon}\cdot (-\log |z_j|)^M
 \Bigr)\vol<\infty
\]
Here $\vol$ is obtained from the Poincar\'{e} metric.

Let take a small positive number $\epsilon$
and a sufficiently large real number $M$.
Then we can find a solution
$G$ satisfying the equation
$\delbar(G)=F$ and the finiteness
$\int |G|_h^2\cdot
 \prod_{j=1}^{n-1} 
 \bigl(|z_j|^{\epsilon}\cdot (-\log |z_j|)^M
 \bigr)\vol<\infty$,
due to Proposition \ref{prop;11.28.6},
 Proposition \ref{prop;12.11.1},
 Corollary \ref{cor;12.11.2},
 Lemma \ref{lem;12.11.3}).

We put $f_3=f_1-f_2-\zeta\cdot G$.
Then by our construction, it satisfies the following:
\[
 \delbar(f_3)=0,
\quad
 f_{3\,|\,X_0-D_0}=f,
\quad
 \int |f|_h^2\cdot
 \prod_{j=1}^{n-1} 
 \bigl(|z_j|^{\epsilon}\cdot (-\log |z_j|)^M
 \bigr)\vol<\infty.
\]
By lemma \ref{lem;12.8.15},
we obtain the finiteness of the integral
over $\pi_j^{-1}(a)$
for any $a\in \Delta^{\ast\,l-1}\times \Delta^{n-l}$
and any $1\leq j\leq l$.
Then the increasing order of $f$ is dominated by $-\epsilon$.
By the triviality of the parabolic structures,
the increasing order of $f_3$ is, in fact, $0$.
Thus $f_3$ is a holomorphic section of $\prolong{\nbige^0}$
satisfying $f_{3\,|\,X_0}=f$,
due to Corollary \ref{cor;11.28.15}.

By using the extension result above,
we obtain the following immediate corollary.
\begin{cor}
The sheaf $\prolong{\nbige^0}$ is locally free
in codimension $1$,
under the assumption {\rm\ref{assumption;11.30.2}}.
\hfill\qed
\end{cor}

\begin{rem}{\rm
It is interesting that we can derive the fact
$\prolong{\nbige^{0}}$ is locally free
if the dimension of the base complex manifold is $2$.
Let consider the case
$X=\Delta^2=\{(z_1,z_2)\}$ and $D=\{z_1\cdot z_2=0\}$.
We put $O=(0,0)$.
We have already known that 
$\prolong{\nbige^{0}}$ is locally free
over $X-\{O\}$.
We can show that we have a sections $f_1,\ldots, f_N$
of $\prolong{\nbige^0}$ such that
they generate $\prolong{\nbige^0}$ over $X-\{O\}$.
Use, for example, an argument in page 35 in \cite{cg}.
Then we obtain the morphism
$\varphi:\nbigo^{\oplus\,N}\lrarr \prolong{\nbige^0}$.
The image of $\varphi$ is denoted by $\Image(\varphi)$,
which is coherent.
On $X-\{O\}$, we have $\Image(\varphi)=\prolong{\nbige^0}$.
Let consider the double dual of $\Image(\varphi)$,
namely,
we put $\nbigf:=Hom(\Image(\varphi,\nbigo),\nbigo)$.
Then $\nbigf$ is coherent, and locally free.
Since $\nbigf$ is locally free,
we have $\prolong{\nbige^0}\subset\nbigf$
due to Hartogs theorem.
On the other hand, we can show that 
the increasing orders any holomorphic sections
of $\nbigf$ are less than $0$,
due to Corollary \ref{cor;11.28.15}.
Thus we have $\nbigf=\prolong{\nbige^0}$.

The argument works in higher dimensional case,
i.e.,
we can show that $\prolong{\nbige^0}$ is coherent and reflexive
by the same argument.
But the double dual of coherent sheaf
is just reflexive, and not locally free in general.
\hfill\qed
}
\end{rem}

\subsubsection{The metrics and the curvatures of $\nbigo(i)$ on $\proj^1$}

To resolve the difficulty in the intersection of the divisors,
we need some preparation.
We denote the one dimensional projective space over $\cnum$
by $\proj^{1}=\{[t_0:t_1]\}$.
The points $[0:1]$ and $[1:0]$ are denoted by $0$ and $\infty$
respectively.
We use the coordinates $t=t_0/t_1$ and $s=t_1/t_0$.
We have a line bundle $\nbigo(i)$ over $\proj^1$.
The coordinates of $\nbigo(i)$ is given as follows:
$(t,\zeta_1)$ over $\proj^1-\{\infty\}$,
and $(s,\zeta_2)$ over $\proj^1-\{0\}$.
The relations are given by
$s=t^{-1}$ and $t^{-i}\cdot \zeta_1=\zeta_2$.

Recall that we have the smooth metric $h_i$ of $\nbigo(i)$.
Let $\xi=(t,\zeta_1)=(s,\zeta_2)$ be an element of $\nbigo(i)$.
\[
 h_i(\xi,\xi):=
 |\zeta_1|^2\cdot \bigl(1+|t|^2\bigr)^{-i}
=|\zeta_2|^2\cdot \bigl( 1+|s|^2\bigr)^{-i}.
\]
For any real numbers $a$ and $b$,
we have the possibly singular metrics $h_{i,(a,b)}$
of $\nbigo(i)$:
Let $\xi=(t,\zeta_1)=(s,\zeta_2)$ be an element of $\nbigo(i)$.
\[
 h_{i,(a,b)}(\xi,\xi):=
 h_i(\xi,\xi)\cdot \bigl(1+|t|^{-2}\bigr)^a\cdot \bigl(1+|t|^2\bigr)^b
=h_i(\xi,\xi)\cdot \bigl(1+|s|^{2}\bigr)^a\cdot \bigl(1+|s|^{-2}\bigr)^b.
\]
Around $|t|=0$,
the order of $h_{0,(a,b)}$ is equivalent to $|t|^{-2a}$.
Around $|s|=0$,
the order of $h_{0,(a,b)}$ is equivalent to $|s|^{-2b}$.
The curvature $R(h_{i,a,b})$ is as follows:
\[
 R(h_{i,a,b})=
 (-a-b+i)\cdot
 \frac{dt\cdot d\tbar}{\bigl(1+|t|\bigr)^2}.
\]
When $i=-1$ and $a=b=-1/2$,
then we have the following:
 \[
 R(h_{-1,a,b})=0.
\]
Let take a point $P\in \proj^1$.
Then we obtain a morphism $\nbigo(i)\lrarr \nbigo(i+1)$
of coherent sheaves.
The morphism is bounded with respect to the metrics
$h_{i,(a,b)}$ and $h_{i+1,(a,b)}$.

\vspace{.1in}
We are mainly interested in the case $i=-1$.
We regard $\nbigo(-1)$ as a complex manifold.
The open submanifold
$Y$ is defined to be
$\bigl\{\xi\in\nbigo(-1)\,\big|\,h_{-1,(0,0)}(\xi,\xi)<1\bigr\}$.
We denote the naturally defined projection of $Y$
onto $\proj^1$ by $\pi$.
We denote the image of the $0$-section $\proj^1\lrarr Y$
by $\proj^1$.
Then we have the normal crossing divisor
$D'=\proj^1\cup \pi^{-1}(0)\cup \pi^{-1}(\infty)$
of $Y$.
The manifold $Y-D'$ is same as 
$\bigl\{(t,x)\in \cnum^{\ast\,2}\,\big|\,
 |x|^2\bigl(1+|t|^2\bigr)<1\bigr\}$.
We have the complete Kahler metric $g:=g_1+g_2+g_3$ of $Y-D'$
given as follows:
As a contribution of the $0$-section $\proj^1$,
we put $\tau_1=-\log\Bigl[\bigl(1+|t|^2\bigr)\cdot |x|^2\Bigr]$,
and as follows:
\[
 g_1:=\frac{1}{\tau_1^2}
 \left(
 \frac{\tbar\cdot dt}{1+|t|^2}
+\frac{dx}{x}
 \right)
\cdot
 \left(
 \frac{t\cdot d\tbar}{1+|t|^2}
+\frac{d\xbar}{\xbar}
 \right)
+\frac{1}{\tau_1}
 \frac{dt\cdot d\tbar}{\bigl(1+|t|^2\bigr)^2}.
\]
Note that $g_1$ gives the complete Kahler metric of 
$Y-\proj^1$.

As a contribution of $\pi^{-1}(\infty)$,
we put $\tau_2=\log\bigl(1+|t|^2\bigr)$, and as follows:
\[
 g_2=\frac{1}{\tau_2}
 \Bigl(
 -1+\frac{|t|^2}{\tau_2}
 \Bigr)
 \cdot
 \frac{dt\cdot d\tbar}{\bigl(1+|t|^2\bigr)}.
\]
Note that we have $-1+|t|^2\cdot \tau_2>1$.
Around $|t|=\infty$, or equivalently, around $|s|=0$,
the $g_2$ is similar to $(-|s|\log|s|)^{-2}ds\cdot d\sbar$.
Around $|t|=0$,
we have $g_2=\bigl(2^{-1}+o(|t|^2)\bigr)\cdot dt\cdot d\tbar$.

As the contribution of the divisor $\pi^{-1}(0)$,
we put $\tau_3:=\log(1+|t|^2)-\log|t|^2=\log(1+|s|^2)$,
where we use $s=t^{-1}$.
And we put as follows:
\[
 g_3=\frac{1}{\tau_3}
 \cdot
 \Bigl(
 -1+\frac{|s|^2}{\tau_3}
 \Bigr)
 \frac{ds\cdot d\sbar}{\bigl(1+|s|^2\bigr)^2}.
\]
By the symmetry, the behaviour of $g_3$ is similar to $g_2$.

The following lemma is easy to see.
\begin{lem}
$g$ gives  the complete Kahler metric of $Y-D'$.
Around the $0$-section $\proj^1$, 
$\pi^{-1}(0)$ and $\pi^{-1}(\infty)$,
the behaviours of the metric
are equivalent to the Poincar\'{e} metric.
\hfill\qed
\end{lem}

We note the following formulas:
\[
 \begin{array}{l}
 {\displaystyle
 \delbar\del\log \tau_1=
 \frac{1}{\tau_1^2}
 \Bigl(
  \frac{\tbar\cdot dt}{1+|t|^2}+\frac{dx}{x}
 \Bigr)
 \wedge
 \Bigl(
  \frac{t\cdot d\tbar}{1+|t|^2}+\frac{d\xbar}{\xbar}
 \Bigr)
 +\frac{1}{\tau_1}
 \frac{dt\wedge d\tbar}{\bigl(1+|t|^2\bigr)^2}=:\omega_1,}\\
\mbox{{}}\\
{\displaystyle
 \delbar\del\log \tau_2=
 \frac{1}{\tau_2}
 \Bigl(
 -1+ \frac{|t|^2}{\tau_2}
 \Bigr)
\cdot
 \frac{dt\wedge d\tbar}{\bigl(1+|t|^2\bigr)^2}=:\omega_2,}\quad
{\displaystyle
 \delbar\del\log \tau_3=
 \frac{1}{\tau_3}
 \Bigl(
 -1+ \frac{|s|^2}{\tau_3}
 \Bigr)
\cdot
 \frac{ds\wedge d\sbar}{\bigl(1+|s|^2\bigr)^2}=:\omega_3.}
 \end{array}
\]
We put $\omega=\omega_1+\omega_2+\omega_3$.
We put as follows:
\[
 H_0=
 \frac{1}{\tau_1}
+\frac{1}{\tau_2}
 \Bigl(
   -1+\frac{|t|^2}{\tau_2}
 \Bigr)
+\frac{1}{\tau_3}
 \Bigl(
  -1+\frac{|s|^2}{\tau_3}
 \Bigr)>0.
\]
Then we have the following:
\[
 \omega^2=\det(g)\cdot dt\wedge d\tbar \wedge dx\wedge d\xbar
= \Bigl(
 \frac{1}{\tau_1^2\cdot |x|^2\cdot \bigl(1+|t|^2\bigr)^2}
 \times H_0
 \Bigr)
 \cdot dt\wedge d\tbar \wedge dx\wedge d\xbar.
\]
We put as follows:
\[
 H_1:=\frac{H_0}{\bigl(1+|t|^2\bigr)\cdot \bigl(1+|s|^2\bigr)}
\]
Recall that we have $\Ric(g)=\delbar\del(\det(g))$.
\begin{lem}\mbox{{}}
\begin{itemize}
\item
Let $C$ be a number such that $0<C<1$.
On the domain
$\bigl\{(t,x)\in\cnum^{\ast\,2}\,\big|\,
 |x|^2\cdot\bigl(1+|t|^2\bigr)\leq C\bigr\}$,
we have the following similarity of the behaviour:
\[
 H_1\sim (\log|t|)^{-2},\quad(|t|\to\infty, \mbox{\rm or, } |t|\to 0).
\]
\item
We have the equality:
$Ric(g)-\delbar\del\log(H_1)=-\delbar\del\log\tau_1^2$.
\hfill\qed
\end{itemize}
\end{lem}

\subsubsection{Extension property II}

Let consider the blow up of $\Delta^2=\{(z_1,z_2)\}$ at $O=(0,0)$
which we denote by $\widetilde{\Delta}^2$.
We denote the proper birational projective
morphism $\widetilde{\Delta}^2\lrarr \Delta^2$
by $\phi'$.
We have the divisor $\{z_i=0\}$ of $\Delta^2$.
The proper transformation of $\{z_i=0\}$  in $\widetilde{\Delta}^2$
is denoted by $D_i$.
We can take a holomorphic embedding $\iota$
of $Y$, given in the previous subsubsection,
to $\widetilde{\Delta}^2$
satisfying the following:
\begin{itemize}
\item The image of the $0$-section $\proj^1$ is the exceptional curve.
\item
We have $\iota^{-1}(D_1)=\pi^{-1}(\infty)$
and $\iota^{-1}(D_2)=\pi^{-1}(0)$.
\end{itemize}
We denote the composite $\phi'\circ\iota:Y\lrarr \Delta^2$
by $\phi$.
The restriction $\phi'_{|Y-D'}$ is an open embedding
into $\Delta^2-D_1\cup D_2$.

Let consider the one dimensional subvariety $C'(1,1)={(z,z)\in\Delta^2}$
of $\Delta^2$.
We denote the closure of 
the inverse image of $C'(1,1)-\{O\}$ via $\phi$ in $Y$
by $C(1,1)$.
Clearly the natural morphism $\phi_{|C(1,1)}$ gives an open embedding
of $C(1,1)$ into $C'(1,1)$.

Let consider the product with $\Delta^{n-2}$.
We obtain the naturally defined morphism
$\phi:Y\times \Delta^{n-2}\lrarr \Delta^2\times\Delta^{n-2}$.
The restriction $\phi_{|(Y-D')\times \Delta^{\ast\,n-2}}$
gives an open embedding of $(Y-D')\times \Delta^{\ast\,n-2}$
into $\Delta^{2\,\ast}\times \Delta^{\ast\,n-2}$.
We have the divisors
$C'(1,1)\times \Delta^{n-2}\subset \Delta^2\times\Delta^{n-2}$,
and
$C(1,1)\times \Delta^{n-2}\subset Y\times \Delta^{n-2}$.

Let consider the tame nilpotent harmonic bundle
$\harmonicbundle$ over $\Delta^{\ast\,n}$.
Then we obtain the harmonic bundle
$\phi^{\ast}(E,\delbar_E,h)$ 
over $(Y-D')\times\Delta^{\ast\,n-2}$.
We also have the deformed holomorphic bundles
$\phi^{\ast}\nbige^{1}$ over $\{1\}\times
(Y-D')\times \Delta^{\ast\,n-2}$.

We will use the coordinate $(t,x)$ of $Y-D'$ and 
$(z_2,\ldots,z_{n})$ of $\Delta^{\ast,n-2}$.

Note the following:
\begin{itemize}
\item The prolongment
 $\prolong{\phi^{\ast}\nbige^{1}}$ is locally free.
\item
 We have the holomorphic frame $\vecv$
 of $\prolong{\phi^{\ast}\nbige^{1}}$.
 Around $P=[1:1]$ in the $0$-section $\proj^1$,
 the components of the hermitian matrices
 $H(\phi^{\ast}h,\vecv)$ and $H(\phi^{\ast}h,\vecv)^{-1}$
 are dominated by the polynomials of
 $(-\log |x|)$ and  $(-\log|z_i|)$ for $i=2,\ldots,n$.
\item
We denote the curvature of $(E,\delbar_E,h)$
by $R(h)$.
It is dominated by $\phi^{\ast}\delbar\del\log\tau(\veca,N)$.
\end{itemize}

Let $\epsilon$ be sufficiently small positive number.
Let $a$ and $b$ be positive numbers
such that $a=b=-1/2$.
The metric $\tilde{h}$ of $\phi^{\ast}(E)\otimes\nbigo_{\proj}(-1)$
is defined as follows:
\[
 \tilde{h}:=
 \phi^{\ast}h_{\veca,N}\cdot
 h_{-1,a,b}\cdot
 H_1^{-1}\cdot
  \tau_1^{2+\epsilon}
 \bigl(
 \tau_2\cdot\tau_3
 \bigr)^{\epsilon}.
\]

\begin{lem}
When $N$ is sufficiently smaller than $0$,
then the following inequality holds for any
$\eta\in A^{0,1}_c(\phi^{\ast}E)$:
\[
 \doublelangle
 R(\tilde{h})+\Ric(g),\eta
 \doublerangle_{\tilde{h}}
\geq \epsilon||\eta||_{\tilde{h}}^2.
\]
\end{lem}
\pf
We have the following equality:
\begin{multline}
 R(\tilde{h})+\Ric(g)
=R(\phi^{\ast}h_{\veca,N})+R(h_{-1,a,b})
-\delbar\del\log H_1
+(2+\epsilon)\cdot\delbar\del\log \tau_1
+\epsilon \cdot\delbar\del(\log \tau_2+\log\tau_3)
+\Ric(g)\\
= R(\phi^{\ast}h_{\veca,N})
+\epsilon(\omega_1+\omega_2+\omega_3).
\end{multline}
By taking sufficiently negative $N$,
we can assume the following inequality
for any $\eta\in A_c^{0,1}(E)$:
\[
 \doublelangle
 R(h_{\veca,N}),\eta
 \doublerangle_{\veca,N}\geq\,0.
\]
Then, by a fiberwise linear algebraic argument,
it is easy to see that the following inequality holds
for any $\eta\in A_c^{0,1}(\phi^{\ast}(E))$:
\[
 \doublelangle
 \phi^{\ast}R(h_{\veca,N}),\eta
 \doublerangle_{\tilde{h}}\geq\,0.
\]
On the other hand, we obtain
$\big\langle\big\langle
 \omega_1+\omega_2+\omega_3,\eta
 \big\rangle\big\rangle_{\tilde{h}}
\geq \epsilon\cdot||\eta||_{\tilde{h}}$,
which can be checked directly from definition.
Thus we are done.
\hfill\qed

\begin{lem}
Let $f$ be a holomorphic section of
$\prolong{\nbige^0_{|C'(1,1)\times \Delta^{n-2}}}$.
Then we have a holomorphic section $\tilde{f}$
of $\prolong{\nbige^0}$ 
such that $\tilde{f}_{|C'(1,1)\times\Delta^{n-2}}=f$.
\end{lem}
\pf
We identify $C(1,1)\times \Delta^{n-2}$ and
$C'(1,1)\times \Delta^{n-2}$.
By the same argument as that
in the subsubsection \ref{subsubsection;12.15.201},
we can show that we have an extension of $f$, say
$f_3$, to the section of $\prolong{\phi^{\ast}\nbige^0}$,
by using the metric $\tilde{h}$.
Let consider the norm of $f_3$
with respect to our original metric $h$.
We consider the restriction of  $f_3$ to
$\pi_j^{-1}(a)$, for any $a\in\Delta^{\ast\,l-1}\times\Delta^{n-l}$
and for $1\leq j\leq l$.
First we consider the case $j=1$.
On $\pi_1^{-1}(a)$,
the metric $\tilde{h}$ is equivalent
to $h\cdot |z_1|^{1/2+\epsilon}\times Q$,
where $Q$ denotes a polynomial of $\log|z_1|$.
Thus the increasing order of $|f_{3\,|\,\pi_1^{-1}(a)}|_{h}$
is less than $-1/2-\epsilon$.
Since the parabolic structure is trivial,
we can conclude that the increasing order
is, in fact, less than $0$.
Similarly we can show that
the increasing order of $f_{3\,|\,\pi_2^{-1}(a)}$
with respect to the metric $h$ is $0$.
If $j>2$, then the metrics $h$ and $\tilde{h}$
are equivalent on $\pi_j^{-1}(a)$.
Thus the increasing order of $f_{3\,\pi_j^{-1}(a)}$
with respect to the metric $h$, is less than $0$,
also in this case.
Thus we obtain that $f_3$ is, in fact,
a section of $\prolong{\nbige^0}$,
due to  Corollary \ref{cor;11.28.15}.
\hfill\qed

\begin{cor}
The prolongments $\prolong{\nbige^0}$ and $\prolong{\nbige}$
are locally free.
\end{cor}
\pf
Assume that the prolongment
of the restriction to the hyperplane is locally free.
Let $\vecv=(v_i)$ be a holomorphic section of
$\prolong{(\nbige^0_{|C(1,1)\times\Delta^{n-2}})}$.
We extend $v_i$ to a section $\tilde{v}_i$ of $\prolong{\nbige^0}$
as in the previous lemma.
Then $\tilde{\vecv}=(\tilde{v}_i)$ gives a holomorphic frame
of $\nbige^0$ around $(0,\ldots,0)$,
for $\volume(\tilde{\vecv})_{(0,\ldots,0)}$ is not $0$
in $\prolong{\det(\nbige)}_{(0,\ldots,0)}$.
Thus the induction can proceed.
\hfill\qed

\begin{prop}\mbox{{}}
\begin{itemize}
\item
For any $\lambda\in\cnum$, $\DD^{\lambda}$ is of log type.
Namely if $f$ is a holomorphic section of $\prolong{\nbige^{\lambda}}$,
then $\DD^{\lambda}(f)$ is a holomorphic section of
 $\prolong{\nbige^{\lambda}}
 \otimes
 \Omega^{1,0}_{\nbigx^{\lambda}}
 (\log \nbigd^{\lambda})$.
\item $\DD$ is of log type.
\end{itemize}
\end{prop}
\pf
If $\lambda\neq 0$, a normalizing frame gives 
a trivialization of $\prolong{\nbige^{\lambda}}$.
It implies that $\DD^{\lambda}$ is of log type.
If $\lambda=0$, we know the estimate of the norm of $\theta$.
It implies that $\theta$ is of log type.

Let $f$ be a holomorphic section of $\prolong{\nbige}$.
The section $\DD(f)$ is holomorphic over
$\nbigx^{\shikaku}\cup (\nbigx-\nbigd)=\nbigx-\nbigd^0$.
Thus it naturally gives a holomorphic section
of $\prolong{\nbige}\otimes p_{\lambda}^{\ast}\Omega_{X}(\log D)$.
\hfill\qed

\begin{rem}
Note that we also obtain the extendability of the holomorphic sections
over the curves.
We will use the fact without mention.
\hfill\qed
\end{rem}

\subsection{Functoriality}
\label{subsection;12.15.20}

\begin{prop} We have the following functoriality of
the prolongment.
\[
\begin{array}{lll}
 \det(\prolong{\nbige})\simeq \prolong(\det(\nbige)),
 &
\prolong{\nbige_1}\otimes\prolong{\nbige_2}\simeq
 \prolong{(\nbige_1\otimes\nbige_2)},
 &
\prolong{(\nbige^{\lor})}
 \simeq (\prolong{\nbige})^{\lor}.
\\
Hom(\prolong{\nbige_1},\prolong{\nbige_2})\simeq
 \prolong{Hom(\nbige_1,\nbige_2)}
&
\Sym^l(\prolong{\nbige})\simeq \prolong{\Sym^l\nbige},
 &
\bigwedge^l(\prolong{\nbige})\simeq \prolong{\bigwedge^l\nbige}.
\end{array}
\]
In each case, the $\lambda$-connection $\DD$ is also isomorphic.

We also a similar isomorphisms of the conjugate deformed holomorphic 
bundles with $\mu$-connections.
\end{prop}
\pf
Let see the case of determinant bundles.
We have the natural morphism
$\det(\prolong{\nbige})\lrarr \prolong{\det(\nbige)}$.
By using the result of Simpson
(see \cite{s2}, or Lemma \ref{lem;12.2.20} in this paper)
and the extension property of the holomorphic sections,
we can show that the morphism is surjective, and thus isomorphic.
The other cases are similar.
\hfill\qed

Let consider the morphism $f:\Delta_z^m\lrarr \Delta_{\zeta}^n$
given as follows:
\[
 f^{\ast}(\zeta_i)=
 \prod_{j}z_j^{a_{i\,j}}.
\]
Let $\harmonicbundle$ be a tame nilpotent harmonic bundle
with trivial parabolic structure over $\Delta_{\zeta}^{\ast\,n}$.
Then the deformed holomorphic bundle
of $f^{\ast}\harmonicbundle$ is naturally isomorphic to
the pull back of $(\nbige,\DD,h)$.

\label{subsection;12.10.1}

\begin{prop}\label{prop;12.10.2}
We have the natural isomorphism of 
$f^{\ast}(\prolong{\nbige})$ and $\prolong{(f^{\ast}\nbige)}$.
\end{prop}
\pf
Since the prolongment is characterized by increasing orders,
we have the natural morphism
$f^{\ast}(\prolong{\nbige})\lrarr\prolong{(f^{\ast}\nbige)}$.
We have to show that it is isomorphic.
When rank of $E$ is $1$,
then the claim is obvious.
Let consider the general case.
Let $\vecv$ be a holomorphic frame of $\prolong{\nbige}$.
Then $f^{\ast}\Omega(\vecv)$ gives a holomorphic frame
of $f^{\ast}\det(\prolong{\nbige})\simeq
\det(\prolong{f^{\ast}\nbige})$.
Thus $f^{\ast}(\vecv)$ gives naturally a frame
of $\prolong{f^{\ast}\nbige}$.
\hfill\qed

We see the relation of prolongments of $\nbige$ and $\nbige^{\dagger}$.
We put $X=\Delta^n$ and $D=\{\prod_{i=1}^n z_i=0\}$.
Let $\vecv$ be a holomorphic frame of $\prolong{\nbige}$ over $\nbigx$.
We have the restriction $\vecv_{|\nbigx-\nbigd}$
of $\vecv$ to $\nbigx-\nbigd$,
which is a frame of $\nbige$ over $\nbigx-\nbigd$.
Then we obtain the conjugate frame
$(\vecv_{|\nbigx-\nbigd})^{\dagger}$
of $\nbige^{\dagger}$ over $\nbigx^{\dagger}-\nbigd^{\dagger}$.

\begin{prop}
$(\vecv_{|\nbigx-\nbigd})$ naturally induces the frame
of $\prolong{\nbige^{\dagger}}$ over $\nbigx^{\dagger}$.
\end{prop}
\pf
When the rank of $E$ is $1$,
the claim is obvious.
Let consider the general case.
By using the result in the case of curves,
we can show that the increasing order of
$v_{i\,|\,\nbigx-\nbigd}^{\dagger}$ is less than $0$,
and thus
it gives a holomorphic section of $\prolong{\nbige^{\dagger}}$.
Let consider the section
$\Omega(\vecv^{\dagger}_{\nbigx-\nbigd})$.
By functoriality of our construction,
it is same as $\Omega(\vecv)^{\dagger}_{\nbigx-\nbigd}$,
which gives a holomorphic frame of $\det(\prolong{\nbige^{\dagger}})$.
Thus $\vecv^{\dagger}_{\nbigx-\nbigd}$ gives a frame
of $\prolong{\nbige^{\dagger}}$.
\hfill\qed

\vspace{.1in}
We will denote the induced frame by $\vecv_{|\nbigx-\nbigd}^{\dagger}$
simply by $\vecv^{\dagger}$.


\section{Limiting mixed twistor structure}
\label{section;12.15.75}

\subsection{Preliminary}
\label{subsection;12.15.21}
\subsubsection{Normalizing frame}

Let $(E,\DD^1)$ be a holomorphic vector bundle with a holomorphic
flat connection
over $\Delta^{\ast\,n}$.
Assume that the monodromies are unipotent.
Let fix the coordinate $(z_1,\ldots,z_n)$ 
of $\Delta^{\ast\,n}$.
Recall that we have a holomorphic frame
such that the corresponding connection form is
of the form $\sum_iA_idz/z_i$ for $A_i\in M(r)$.
Such a frame is called the normalizing frame in this paper.
Let us recall the construction of such a frame.

Let $\hyperh$ denote the upper half plain
$\{\zeta\in\cnum\,|\,Im(\zeta)>0\}$.
Let take the universal covering
$\pi:\hyperh^n\lrarr\Delta^{\ast\,n}$,
given by $z_i=\exp(2\pi\sqrt{-1}\zeta_i)$.
Let $P$ be a point of $\Delta^{\ast\,n}$
and $\tilde{P}$ be a point of $\hyperh^n$
such that $\pi(\tilde{P})=P$.
Let $\gamma_i$ denote the following path:
\[
 \{t\in\real\,|\,0\leq t\leq 1\}\ni t\longmapsto
 \bigl(z_1,\ldots,z_{i-1},\exp(2\pi\sqrt{-1}t)\cdot z_i,
 z_{i+1},\ldots,z_n\bigr)
 \in \Delta^{\ast\,n}.
\]
We denote the monodromy along $\gamma_i$
by $M(\gamma_i)$.

Let pick a base $\vece(P)$ of $E_{|P}$.
Let $M_i$ be the matrix representing
the endomorphism $M(\gamma_i)\in End(E_{|P})$
with respect to the frame $\vecw(P)$,
that is,
$M(\gamma_i)\vecw=\vece\cdot M_i$.
By our assumption, $M_i$ is unipotent.
We denote the logarithm of $M_i$
by $N_i$,
i.e.,
\[
 N_i=\log(M_i):=\sum_{m=l}^{\infty}\frac{1}{m}\bigl(M_i-1\bigr)^m,
\quad
 \exp\bigl(N_i\bigr)=M_i.
\]

The frame $\vece(P)$ of $E_{|P}$
naturally induces the base $\vece(\tilde{P})$
of $\pi^{\ast}E_{|\tilde{P}}$.
We can take the flat frame $\vecu$ of $\pi^{\ast}E$
such that $\vecu_{|\tilde{P}}=\vece(\tilde{P})$.
Then we put as follows:
\[
 \vecw=\vecu\cdot
 \exp\Bigl(
 -\sum_{i=1}^n
 \bigl(\zeta_i-\zeta_i(\tilde{P})
 \bigr)
  \cdot N_i\Bigr).
\]
Then $\vecw$ is, in fact, the pull back of
some holomorphic frame of $E$ on $\Delta^{\ast\,n}$.
We denote the frame over $\Delta^{\ast\,n}$
by the same notation $\vecw$.
We have the following equality:
\begin{equation} \label{eq;12.1.1}
 \DD^{1}\vecw=\vecv\cdot
 \Bigl( -\sum_{i=1}^n N_i\cdot d\zeta_i
 \Bigr)
=\vecw\cdot
 \Bigl(-\sum_{i=1}^nN_i\cdot \frac{dz_i}{2\pi\sqrt{-1}z_i}
 \Bigr).
\end{equation}
Note that such normalizing frame is determined
once we pick a point $P$ and the frame $\vece(P)$,
when we fix a coordinate $(z_1,\ldots,z_n)$.
Namely, we have the following lemma.
\begin{lem}
Let $\vecw_1$  and $\vecw_2$ be normalizing frames
such that $\vecw_1(P)=\vecw_2(P)$.
Then we have $\vecw_1=\vecw_2$.
\end{lem}
\pf
On $\hyperh^n$,
the frames $\vecw_i$ are solutions of the equation (\ref{eq;12.1.1}),
such that $\vecw_1(\tilde{P})=\vecw_2(\tilde{P})$.
Such solution is uniquely determined.
\hfill\qed

Let see the dependence of a normalizing frame 
$\vecw$ on a choice of $P$, $\vece(P)$ and the coordinate.
Let take another frame $\vece_1(P)$ of $E_{|P}$.
We denote the corresponding normalizing frame by $\vecw_1$.
We have an element $g\in GL(r)$ such that
$\vece_1(P)=\vece(P)\cdot g$.
Then we have $\vecw_1(P)=\vecw(P)\cdot g$.
We also have
$\DD^1\vecw\cdot g=\vecw\cdot A\cdot g=
 \vecw\cdot g\cdot (g^{-1}\cdot A\cdot g)$.
Thus we obtain that $\vecw_1=\vecw\cdot g$.

Let $Q$ be another point of $\Delta^{\ast\,n}$,
and $\vece(Q)$ be a base of $E_{|Q}$.
We denote the corresponding normalizing frame
by $\vecw_2$.
We have an element $g\in GL(r)$
such that $\vecw_2(P)=\vecw(P)\cdot g$.
Then we know that $\vecw_2=\vecw\cdot g$.

Let $(z_1',\ldots,z_n')$ denotes another holomorphic coordinate,
satisfying $\{z_i'=0\}=\{z_i=0\}=:D_i$ for any $i$.
Note that $(z_i')$ be a coordinate of an open subset $U$
of $\cnum^n$, which is not necessarily same as $\Delta^n$ above.
Let take a covering $\pi':\hyperh^n\lrarr U-\bigcup_i D_i$.
We use the coordinate $(\zeta_i')$ for $\hyperh^n$ in this case.
We denote the corresponding normalizing frame $\vecw'$.

We have a holomorphic function $g_i$ defined over $\Delta^n$
satisfying
$(\zeta_i-\zeta_i')=g_i$,
or $z_i/z_i'=\exp(2\pi\sqrt{-1}g_i)$.
Then we have the following relation:
\begin{multline} \label{eq;12.1.2}
 \vecw'=
 \vecu\cdot
 \exp\Bigl(
 -\sum 
    \bigl(
      \zeta'_i-\zeta'_i(P)
    \bigr)
  \cdot N_i
 \Bigr)\\
=\vecu\cdot
 \exp\Bigl(
 -\sum 
    \bigl(
      \zeta_i-\zeta_i(P)
    \bigr)
  \cdot N_i
+\sum_i
 \bigl(
  g_i+  (\zeta_i'(P)-\zeta_i(P))
 \bigr)\cdot N_i
 \Bigr)\\
=\vecw\cdot
 \exp\Bigl[
  \sum_i 
 \Bigl(
 g_i+
  \bigl(
 \zeta_i'(P)-\zeta_i(P)
 \bigr)
\Bigr)
\cdot N_i
 \Bigr].
\end{multline}


\subsubsection{Isomorphisms} \label{subsubsection;12.1.5}

Let $\vecw=(w_1,\ldots,w_r)$ be a normalizing frame of
a holomorphic flat connection $(E,\DD^1)$.
We regard it as a frame of the prolongment $\prolong{E}$
over $\Delta^n$.
Let $P$ and $Q$ be points of $\Delta^n$.
Let $v$ be an element of $\prolong{E}_{|P}$.
Then $v$ can be uniquely described as a linear combination
of $w_{i\,|\,P}$,
that is,
$v=\sum_i v_i\cdot w_{i\,|\,P}$ for $v_i\in\cnum$.
Then we obtain the element
$\Phi_{P,Q}(v)=\sum _i v_i\cdot w_{i\,|\,Q}$
of $\prolong{E}_{|Q}$.
Then we obtain the linear isomorphism
$\Phi_{(P,Q)}:\prolong{E}_{|P}\lrarr \prolong{E}_{|Q}$.

\begin{lem}
When we fix a coordinate of $\Delta^n$,
then the morphism $\Phi_{(P,Q)}$ is independent of a choice of
normalizing frame.
\end{lem}
\pf
Let $\vecw_1$ and $\vecw_2$ be normalizing frames.
Then we have the element $g\in GL(r,\cnum)$
such that $\vecw_1=\vecw_2\cdot g$.
Thus we have the relation
$\vecw_{1\,|\,P}=\vecw_{2\,|\,P}\cdot g$
and $\vecw_{1\,|\,Q}=\vecw_{2\,|\,Q}\cdot g$.
It implies that $\Phi_{(P,Q)}$ is independent of a choice of
a normalizing frame.
\hfill\qed

Let $P$ be a point contained in $\Delta^{\ast\,n}$.
Then we have the action of monodromies $M(\gamma_i)$
and the logarithms $N(\gamma_i)$ on $E_{|P}$.
Let $P$ be a point contained in $D_i=\{z_i=0\}$.
Then we have the endomorphism $\Res_{D_i}(\DD^1)_{|P}$
of $\prolong{E}_{|P}$.
We denote it by $N(\gamma_i)$ by abbreviation.

\begin{lem}
Let $P$ and $Q$ be points of $\Delta^{n}$.
Then the isomorphism $\Phi_{(P,Q)}$
is compatible with $N(\gamma_i)$.
\end{lem}
\pf
Both of
the endomorphisms $N(\gamma_i)_{|P}$ and $N(\gamma_i)_{|Q}$
are represented by $N_i$ with respect to the frames
$\vecw_{|P}$ and $\vecw_{|Q}$.
Thus the isomorphism $\Phi_{(P,Q)}$ are compatible with
$N(\gamma_i)$.
\hfill\qed

We see the dependence of $\Phi_{(P,Q)}$
on a choice of the coordinate.
Let $(z'_i)$ be another coordinate.
We denote the corresponding morphism
by $\Phi'_{(P,Q)}$.
\begin{lem}
We have the relation
\[
 \Phi_{(P,Q)}=
 \exp\Bigl(
 \sum_i
 \bigl(g_i(P)-g_i(Q)\bigr)
 \cdot N(\gamma_i)
 \Bigr)\circ
 \Phi_{(P,Q)}'.
\]
\end{lem}
\pf
We have the relation (\ref{eq;12.1.2})
between $\vecw$ and $\vecw'$.
It implies the relation of $\Phi_{(P,Q)}$
and $\Phi'_{(P,Q)}$.
\hfill\qed

\subsection{The vector bundle of Simpson}

\label{subsection;12.15.22}

\subsubsection{Normalizing frames}

We put $X=\Delta^n$, $D_i=\{z_i=0\}$ and $D=\bigcup_{i=1}^l D_i$
for $l\leq n$.
Let $\harmonicbundle$ be a tame nilpotent harmonic bundle
with trivial parabolic structure over $X-D$.
We have the prolongment of the deformed holomorphic bundle
$\prolong{\nbige}$ over $\nbigx$,
and the $\lambda$-connection $\DD^{\lambda}$ of log type.
The following lemma is easy from the construction of
normalizing frames.
\begin{lem}
We can take a holomorphic frame
$\vecw$ of $\nbige^{\shikaku}$ over
$\nbigx^{\shikaku}$
such that the restrictions $\vecw_{|\nbigx^{\lambda}}$
are normalizing frames for $\DD^{\lambda,f}$.
\hfill\qed
\end{lem}

In that case, the following holds:
Let
$\nbiga^{f}\in
 \Gamma\Bigl(\nbigx^{\shikaku},
 M(r)\otimes p_{\lambda}^{\ast}\Omega_X^{1,0}\Bigr)$
be a holomorphic flat connection form
of $\DD^{f}$ with respect to $\vecw$,
that is,
$\DD^{f}\cdot \vecw
=\vecw\cdot \nbiga^f$ over $\nbigx^{\shikaku}$.
Then $\nbiga^f$ is of the following form:
\begin{equation} \label{eq;12.1.3}
 \sum_iC_i(\lambda)\cdot \frac{dz_i}{z_i}.
\end{equation}
Here $C_i(\lambda)$ are $M(r)$-valued holomorphic functions
of $\lambda$, but they are independent of $z_i$.
Let consider the $\lambda$-connection form
$\nbiga$ of $\DD$ with respect to the frame $\vecw$.
Then we have the relation
$\nbiga^{f}=\lambda^{-1}\cdot\nbiga$.
Thus $\nbiga$ is also of the form as in (\ref{eq;12.1.3}).
A frame such as $\vecw$ above is also called 
a normalizing frame for $\DD$.

Note that there does not exist a normalizing frame
for Higgs field.
Namely the normalizing frame is defined
only over $\nbigx^{\shikaku}$,
in general.

\subsubsection{Construction}
\label{subsubsection;12.9.10}

Let $P$ and $Q$ be points of $X$.
As in the subsubsection \ref{subsubsection;12.1.5},
we obtain the following isomorphisms
by using the normalizing frame for any $\lambda\neq 0$:
\[
 \Phi_{\lambda,P,Q}:
 \prolong{\nbige}_{|(\lambda,P)}
\lrarr
 \prolong{\nbige}_{|(\lambda,Q)}.
\]
Since $\vecw$ is holomorphic with respect to $\lambda$,
the dependence of the morphisms $\Phi_{\lambda,P,Q}$
on $\lambda$ is holomorphic.
Thus we obtain the holomorphic isomorphism
of the following locally free sheaves 
over $\cnum^{\ast}_{\lambda}$
for any $P,Q\in X$:
\[
 \Phi_{\lambda,P,Q}:
 \prolong{\nbige^{\shikaku}}_{|\cnum_{\lambda}^{\ast}\times \{P\}}
\lrarr
\prolong{\nbige^{\shikaku}}_{|\cnum_{\lambda}^{\ast}\times \{P\}}.
\]

Similarly we obtain the holomorphic isomorphism
for any $P,Q\in X^{\dagger}$:
\[
 \Phi^{\dagger}_{\mu,P,Q}:
 \prolong{\nbige^{\dagger\,\shikaku}}_{|\cnum_{\mu}^{\ast}\times \{P\}}
\lrarr
\prolong{\nbige^{\dagger\,\shikaku}}_{|\cnum_{\mu}^{\ast}\times \{P\}}.
\]

We have the isomorphism
$\cnum^{\ast}_{\lambda}
 \simeq
 \cnum^{\ast}_{\mu}$
given by $\mu=\lambda^{-1}$.
Recall that we have the following identification
of $C^{\infty}$-bundles with the families of the flat connections
(Lemma \ref{lem;12.1.6}.):
\[
 (\nbige^{\shikaku},\DD^{f})
=
 (\nbige^{\dagger\,\shikaku},\DD^{\dagger\,f}),
\mbox{ over }
 \nbigx^{\shikaku}-\nbigd^{\shikaku}
=\nbigx^{\dagger\,\shikaku}-\nbigd^{\dagger\,\shikaku}.
\]
Note that they are identification of the $C^{\infty}$-vector bundles,
and
the identification is holomorphic with respect to $\lambda$
and $\mu$.
Let $P$ be a point of $X-D$.
Then we have the following identification of {\em holomorphic} bundles
over $\cnum_{\lambda}^{\ast}=\cnum_{\mu}^{\ast}$:
\[
 \nbige^{\shikaku}_{|\cnum^{\ast}_{\lambda}\times \{P\}}
=\nbige^{\dagger\,\shikaku}_{|\cnum^{\ast}_{\mu}\times \{P\}}.
\]

Let $Q_1,Q_2$ be points of $X$,
and $P$ be a point of $X-D$.
Then we obtain the following sequence of isomorphisms:
\[
\begin{CD}
 \prolong{\nbige^{\shikaku}}_{|\cnum_{\lambda}^{\ast}\times\{Q_1\}}
@>{\Phi_{\lambda,Q_1,P}}>>
 \nbige^{\shikaku}_{|\cnum_{\lambda}^{\ast}\times \{P\}}
@=
 \nbige^{\dagger\,\shikaku}_{|\cnum_{\mu}^{\ast}\times \{P\}}
@>{\Phi^{\dagger}_{\mu,P,Q_2}}>>
 \prolong{\nbige^{\dagger\,\shikaku}}_{\cnum_{\mu}^{\ast}\times \{Q_2\}}.
\end{CD}
\]
We denote the composite by $\Psi(Q_1,Q_2,P)$.
We have the vector bundles
$\nbige_{|\cnum_{\lambda}\times\{Q_1\}}$ over $\cnum_{\lambda}$
and 
$\nbige^{\dagger}_{|\cnum_{\mu}\times \{Q_2\}}$ over $\cnum_{\mu}$.
By the gluing $\Psi(Q_1,Q_2,P)$,
we obtain the holomorphic vector bundle
over $\proj^1=\cnum_{\lambda}\cup\cnum_{\mu}$.
The vector bundle is denoted by 
$S(Q_1,Q_2,P,(z_i))$.
When we fix the coordinate, then we often omit to denote $(z_i)$.
When $Q_1=Q_2=Q$, it is denoted by $S(Q,P)$.
When we distinguish the harmonic bundle $\harmonicbundle$,
it is denoted by $S(Q,P,\harmonicbundle)$.

\subsubsection{Nilpotent maps}

Assume that $Q_1$ and $Q_2$ are contained in $D_i$.
Then we have the following morphisms of coherent sheaves:
\[
 \residuei\in
End\bigl(\prolong{\nbige}_{|\cnum_{\lambda}\times \{Q_1\}} \bigr),
\quad
 \residueidagger
 \in
End\bigl(\prolong{\nbige^{\dagger}}_{|\cnum_{\mu}\times \{Q_2\}} \bigr).
\]

\begin{lem}
The nilpotent morphisms
$\residuei$ and $-\residueidagger$
induce the morphism of the coherent sheaves:
\[
 \nbign_i^{\sankaku}:
 S(Q_1,Q_2,P)\lrarr S(Q_1,Q_2,P)\otimes\nbigo_{\proj^1}(2).
\]
\end{lem}
\pf
Then we have the following morphisms of the coherent sheaves:
\[
 \residueiflat\in
 End\bigl(\prolong{\nbige}_{|\cnum_{\lambda}\times \{Q_1\}} \bigr),
\quad
 \residueiflatdagger
 \in
 End\bigl(\prolong{\nbige^{\dagger}}_{|\cnum_{\mu}\times \{Q_2\}} \bigr).
\]
The relation of $\residueiflat$ and $\residuei$,
(resp. $\residueiflatdagger$ and $\residueidagger$)
is given as follows:
\[
 \residueiflat=\lambda^{-1}\cdot \residuei,
\quad
\bigl(\mbox{resp.}\,\,
 \residueiflatdagger=\mu^{-1}\cdot \residueidagger.
\bigr)
\]
On the other hand, the logarithms of
the monodromies of the flat connections $\DD$
along the path $\gamma_i$,
induces the following morphism of the coherent sheaves:
\[
 N(\gamma_i)\in 
 End\bigl(\prolong{\nbige}_{|\cnum_{\lambda}\times \{P\}} \bigr).
\]
By the morphisms $\Phi_{\lambda,Q_1,P}$,
the $\residueiflat$ is mapped to
$N(\gamma_i)$.
On the other hand,
the morphism $\Phi^{\dagger}_{\mu,Q_2,P}$
map the $\residueiflatdagger$ to
$-N(\gamma_i)$.
Thus we are done.
\hfill\qed

\subsubsection{Gluing matrices}

Let $\vecv$ be a holomorphic frame of $\prolong{\nbige}$
over $\nbigx$.
Then $\vecv^{\dagger}$ is naturally a holomorphic
frame of $\prolong{\nbige^{\dagger}}$
over $\nbigx^{\dagger}$.
We obtain the frames
$\vecv_{|\cnum_{\lambda}\times Q_1}$
and $\vecv^{\dagger}_{|\cnum_{\mu}\times Q_2}$
of $\prolong{\nbige}_{|\cnum_{\lambda}\times Q_1}$
and $\prolong{\nbige^{\dagger}}_{\cnum_{\mu}\times Q_2}$
respectively.
We would like to describe the gluing morphism $\Psi(Q_1,Q_2,P)$
with respect to the frames.

Let $\vecw$ be a normalizing frame of $\nbige^{\shikaku}$
for the flat connections $\DD^{f}$.
Then $\vecw^{\dagger}$ is naturally a normalizing frame
of $\nbige^{\dagger\,\shikaku}$
for the flat connection $\DD^{\dagger\,f}$.
Let $\Gamma\in C^{\infty}(\nbigx^{\shikaku},GL(r))$ be
the transformation matrices of $\vecv$ and $\vecw$
that is,
$\vecw=\vecv\cdot \Gamma$.
Then the transformation matrices of $\vecv^{\dagger}$
and $\vecw^{\dagger}$
is given by $\Gamma^{\dagger}=F^{\ast}({}^t\overline{\Gamma}^{-1})$,
that is,
$\vecw^{\dagger}=\vecv^{\dagger}\cdot \Gamma^{\dagger}$.
Here $F$ denotes the morphism 
$\cnum_{\mu}\lrarr \cnum_{\lambda}^{\dagger}$
given by $\mu=-\lambdabar$
(See subsubsection \ref{subsubsection;12.1.12}).

Let $\vece$ be a frame of $\nbige^0$
over $\nbigx^0-\nbigd^0$.
Then $\vece^{\dagger}$ gives a frame of $\nbige^{\dagger\,0}$.
The relation of $\vece$ and $\vece^{\dagger}$
is given by $\vece^{\dagger}=\vece\cdot \overline{H(h,\vece)}^{-1}$.

Let $B\in C^{\infty}(\nbigx-\nbigd,GL(r))$
be a transformation matrices of
$p_{\lambda}^{-1}(\vece)$ and $\vecv$,
that is,
$\vecv=p_{\lambda}^{-1}(\vece)\cdot B$.
Then the transformation matrices of
$\vecv^{\dagger}$ and $p_{\mu}^{\dagger\,-1}(\vece^{\dagger})$
is given by $B^{\dagger}=F^{\ast}({}^t\overline{B}^{-1})$,
that is,
$\vecv^{\dagger}=p_{\mu}^{\dagger\,-1}(\vece^{\dagger})\cdot B^{\dagger}$.

Then the gluing morphism $\Psi(Q_1,Q_2,P)$ for $S(Q_1,Q_2,P)$
is represented by the following $GL(r)$-valued function
$A(\lambda,P)$,
with respect to the frames
$\vecv_{|\cnum_{\lambda}^{\ast}\times Q_1}$
and $\vecv^{\dagger}_{|\cnum_{\mu}^{\ast}\times Q_2}$:
\begin{equation} \label{eq;12.2.1}
 \begin{array}{l}
 \vecv^{\dagger}_{|\cnum_{\mu}^{\ast}\times Q_2}=
 \vecv_{|\cnum_{\lambda}^{\ast}\times Q_1}
 \cdot
 A(\lambda,P),\\
\mbox{{}}\\
 A(\lambda,P):=
 \Gamma(\lambda,Q_1)\cdot
 \Gamma(\lambda,P)^{-1}\cdot
 B(\lambda,P)^{-1}\cdot
 \overline{H(h,\vece)}^{-1}(P)\cdot
 B^{\dagger}(\mu,P)\cdot
 \Gamma^{\dagger}(\mu,P)\cdot
 \Gamma^{\dagger}(\mu,Q_2)^{-1}.
\end{array}
\end{equation}
Here $\mu=\lambda^{-1}$.

\subsubsection{Functoriality}

Let see the functoriality of
the construction of the vector bundle $S(O,P)$.
To distinguish the harmonic bundle $\harmonicbundle$,
we use the notation
$S(O,P,\harmonicbundle)$.

\begin{prop}
We fix the coordinate of $\Delta^{n}$.
Then we have the following natural isomorphisms:
\[
\begin{array}{l}
 \det\Bigl(S\bigl(O,P,\harmonicbundle\bigr)\Bigr)\simeq
 S\bigl(O,P,\det\harmonicbundle\bigr),\\
\mbox{{}}\\
 S\bigl(O,P,
 (E_1,\theta_1,h_1)
 \bigr)
 \otimes
 S\bigl(O,P,(E_2,\theta_2,h_2)
  \bigr)
\simeq
 S\Bigl(
  O,P, (E_1,\theta_1,h_1)\otimes (E_2,\theta_2,h_2)
  \Bigr)\\
\mbox{{}}\\
 S\bigl(O,P,\harmonicbundle\bigr)^{\lor}\simeq
 S\bigl(O,P,\harmonicbundledual\bigr),\\
\mbox{{}}\\
 Hom\Bigl(
  S\bigl(O,P,(E_1,\theta_1,h_1)\bigr),
 S\bigl(O,P,(E_2,\theta_2,h_2)\bigr)
 \Bigr)
\simeq
 S\Bigl(O,P,
 Hom\bigl(
  (E_1,\theta_1,h_1),(E_2,\theta_2,h_2)
 \bigr)\Bigr),\\
\mbox{{}}\\
 S\bigl(O,P,\Sym^l\harmonicbundle\bigr)\simeq
 \Sym^lS\bigl(O,P,\harmonicbundle\bigr),\\
\mbox{{}}\\
 S\bigl(O,P,\bigwedge^l\harmonicbundle\bigr)\simeq
 \bigwedge^lS\bigl(O,P,\harmonicbundle\bigr),\\
\end{array}
\]
In each isomorphism,
the nilpotent morphisms $N_i^{\sankaku}$ are also isomorphic,
if $O$ is contained in $D_i=\{z_i=0\}$.
\end{prop}
\pf
Let see the case of the tensor product.
Let $\vecw_a=(w_{a\,i})$ be normalizing frames of $\prolong{\nbige_a}$
over $\cnum_{\lambda}\times \Delta^n$ for $a=1,2$.
Then  $\vecw_1\otimes\vecw_2=(w_{1\,i}\otimes w_{2\,j})$
gives a normalizing frame of $\prolong{\nbige_1\otimes\nbige_2}$.
Then it is easy to check that our construction is functorial
with respect to the tensor product.
The other cases are similar.
\hfill\qed

Let consider the morphism $f:\Delta_z^m\lrarr \Delta_{\zeta}^n$
considered in the subsection \ref{subsection;12.10.1},
given as follows:
\[
 f^{\ast}(\zeta_i)=
 \prod_{j}z_j^{a_{i\,j}}.
\]
Let $\harmonicbundle$ be a tame nilpotent harmonic bundle
with trivial parabolic structure over
$\Delta_{\zeta}^{\ast\,l}\times \Delta_{\zeta}^{n-l}$.

\begin{prop}
We have the natural isomorphism of 
\[
 S\Bigl(Q_1,Q_2,P,f^{\ast}\harmonicbundle\Bigr)\simeq
 S\Bigl(f(Q_1),f(Q_2),f(P),\harmonicbundle\Bigr).
\]
\end{prop}
\pf
Let $\vecw$ be a normalizing frame of
$\prolong{\nbige^{\shikaku}}$.
Then $f^{\ast}(\vecw)$ is a normalizing frame
of $f^{\ast}\prolong{\nbige^{\shikaku}}$.
Then the proposition follows from the functoriality of
the prolongation (Proposition \ref{prop;12.10.2}).
\hfill\qed
\subsection{A limiting mixed twistor theorem}
\label{subsection;12.15.23}

\subsubsection{The filtration}
Let consider the case that the base manifold is one dimensional,
i.e.,
$X=\Delta$ and $D=\{O\}$.
Due to the result of Simpson,
the conjugacy classes of $\nbign^{\sankaku}(\lambda,O)=\resddlambda$
is independent of a choice of $\lambda$.
Thus the weight filtrations $\weightfilt(\resddlambda)$
induced by $\resddlambda$
form the filtrations of $S(O,P)$ by subvector bundles.
We denote the filtration by $\weightfilt^{\sankaku}$.
Thus we obtain the filtered vector bundle
$(S(O,P),\weightfilt^{\sankaku})$.
\begin{quote}
{\em
Ideally, we hope that the filtered vector bundle is mixed twistor,
namely,
we hope that
the $l$-th graded part $\graded_l$ is a direct sum of
$\nbigo_{\proj^1}(l)$.
}
\end{quote}
We call it a conjecture of Simpson.

\subsubsection{Example (the case of $Mod(2,a,C)$)}

As an example,
let see what kind of filtered vector bundle
is obtained in the case $Mod(2,a,C)$.
This is an extremely easy example.
We have the canonical frame $\vecv=(v^{1,0},v^{0,1})$
of $\nbige$ over $\cnum_{\lambda}\times \Delta$.
The restriction
$\vecv_{|\cnum_{\lambda}\times O}$
gives a frame of $\modeldeform(2,a,C)_{|\cnum_{\lambda}\times O}$.
For simplicity, we omit to denote
the notation $|\cnum_{\lambda}\times O$.
The weight filtration $\weightfilt^{\sankaku}$ on $\cnum_{\lambda}$
is obviously given as follows:
\[
 \weightfilt^{\sankaku}_{-1}=
 \langle
  v^{0,1}
 \rangle,
\quad
 \weightfilt^{\sankaku}_{1}=
 \langle
 v^{1,0},
 v^{0,1}
 \rangle.
\]
On the other hand,
$\vecv^{\dagger}=(v^{\dagger\,1,0},v^{\dagger\,0,1})$
gives a holomorphic frame of $\modeldeform(2,a,C)^{\dagger}$.
Then the filtration $\weightfilt^{\sankaku}$ on $\cnum_{\mu}$
is as follows:
\[
 \weightfilt^{\sankaku}_{-1}=
 \langle v^{\dagger\,1,0} \rangle,
\quad
 \weightfilt^{\sankaku}_1=
 \langle
 v^{\dagger\,1,0},
 v^{\dagger\,0,1}
 \rangle
\]

Let consider the gluing matrices.
Luckily $\vecv$ and $\vecv^{\dagger}$ are normalizing frames.
We have the relations (\ref{eq;12.1.15}) 
between $\vecv$ and
$p_{\lambda}^{-1}\vecv_{|\nbigx^0-\nbigd_0}=
 p_{\lambda}^{-1}\vece$,
and the relation (\ref{eq;12.1.16})
between $\vecv^{\dagger}$
and
$p_{\mu}^{\dagger\,-1}\vecv_{|\nbigx^{\dagger\,0}-\nbigd^{\dagger\,0}}
=p_{\mu}^{\dagger\,-1}\vece^{\dagger}$.
Then the gluing matrices $A(\lambda,P)$ is given as follows:
\begin{equation} \label{eq;12.2.2}
 \left(
 \begin{array}{cc}
 1 & \lambda\cdot (y+c)^{-1}\\
 0 & 1
 \end{array}
 \right)
\cdot
 \left(
 \begin{array}{cc}
 (y+c)^{-1} & 0 \\
 0 & y+c
 \end{array}
 \right)
\cdot
 \left(
 \begin{array}{cc}
 1 & 0 \\
 -\mu\cdot(y+c)^{-1} & 1
 \end{array}
 \right)
=
 \left(
 \begin{array}{cc}
 0 & \lambda \\
 -\mu & y+c
 \end{array}
 \right).
\end{equation}
Here $y$ denotes $-\log|z(P)|$, and we used the relation
$\lambda\cdot\mu=1$.
Namely, we have the following relation:
\[
 (v^{\dagger\,1,0},v^{\dagger\,0,1})
=(v^{1,0},v^{0,1})\cdot
  \left(
 \begin{array}{cc}
 0 & \lambda \\
 -\mu & y+c
 \end{array}
 \right)=
\bigl(
 -\mu\cdot v^{0,1},
 \lambda\cdot v^{1,0}+
 (y+c)\cdot v^{0,1}
\bigr).
\]
We have the graded vector space $\graded=\graded_{-1}\oplus\graded_1$
associated with the filtration $\weightfilt^{\sankaku}$.
The line bundle $\graded_{-1}$ is obtained by the relation
$v^{\dagger\,1,0}=-\mu\cdot v^{0,1}$.
Thus the first Chern class is $-1$.
On the other hand,
the line bundle $\graded_1$ is obtained by the relation
$v^{\dagger\,0,1}=\lambda\cdot v^{1,0}$.
Thus the first Chern class is $1$.
Namely the filtered vector bundle
$(S(O,P),\weightfilt^{\sankaku})$ is a mixed twistor,
in this case.

\subsubsection{Example (the case of $Mod(l+1,a,C)$)}

Let see the case of $Mod(l+1,a,C)$ for general $l$.
Since we know that $Mod(l+1,a,C)=\Sym^lMod(2,a,C)$,
we obtain a mixed twistor structure also in this case.
We see the gluing matrices.
We have the canonical frame $\vecv=(v^{p,q}\,|\,p+q=l)$
of $\modeldeform(l+1,a,C)$ over $\cnum_{\lambda}\times \Delta$.
It is also a normalizing frame.
The filtration $\weightfilt^{\sankaku}$ is given as follows:
\[
 \weightfilt^{\sankaku}_h=
 \langle
 v^{i,l-i}\,|\,2i-l\leq h
 \rangle.
\]
We have the conjugate frame
$\vecv^{\dagger}=(v^{\dagger\,p,q}\,|\,p+q=l)$.
It is also a normalizing frame.
The filtration $\weightfilt^{\sankaku}_h$
is given as follows:
\[
 \weightfilt^{\sankaku}_h
=\langle
 v^{\dagger\,i,l-i}\,|\,2i-l\geq -h
 \rangle.
\]
We would like to show that the Chern class of
the $h$-th graded part $\graded_h$ is $h$.
For that purpose,
we would like to calculate the transformation
matrices $A$ of $\vecv$ and $\vecv^{\dagger}$
determined by the following relation:
\[
v^{\dagger\,p,l-p}=
\sum_i A_{i,p}\cdot v^{i,l-i}.
\]
However it is not so easy
to directly calculate whole of the matrices $A$.
We have already known that
the filtration $\weightfilt^{\sankaku}$ have to be
preserved by the gluing.
Then we can derive the following information.
\begin{itemize}
\item
It implies that $A_{i,p}=0$ if $i>l-p$.
\item
We also know that $A$ have to be an element of $GL(l+1)$,
so that $A_{i,l-i}\neq 0$.
\end{itemize}
On the other hand,
we can show the following by a direct calculation:
\begin{itemize}
\item
$A_{i,l-i}$ is of the form
$c_i\cdot\lambda^{2i-l}$ for some complex number $c_i$.
\end{itemize}
(We only have to use
the relations(\ref{eq;11.26.2}), (\ref{eq;12.1.21})
and the orthogonality of the frames
$\vecv_{|\{0\}\times \Delta}$
and $\vecv^{\dagger}_{|\{0\}\times \Delta^{\dagger}}$.)

Thus the the Chern class of $\graded_h$ is $h$,
namely
the filtered vector bundle $(S(O,P))$ is a mixed twistor in this case.

\subsubsection{Some modification of the frames of $\deformoverorigin$
 and $\deformoverorigindagger$}

The gluing matrices (\ref{eq;12.2.1}) are divergent
when $P\to O$.
(For example, see (\ref{eq;12.2.2}).)
We would like to modify our choice of frames
of $\deformoverorigin$ and $\deformoverorigindagger$,
so that the gluing matrices are convergent
for some sequence of points $P_i\to O$.

Let $\vecv$ be a holomorphic frame of $\prolong{\nbige}$,
compatible with the weight filtration
$\weightfilt^{\sankaku}$,
over $\cnum_{\lambda}\times\Delta$.
Then the weight filtration $\weightfilt^{\sankaku}$
over $\cnum_{\lambda}\times O$ is obviously given as follows:
\[
 \weightfilt^{\sankaku}_h=
 \big\langle
 v_{i\,|\,\cnum_{\lambda}\times O}
 \,\big|\,\deg^{\weightfilt^{\sankaku}}(v_i)\leq h
 \big\rangle
\]
We have the conjugate frame $\vecv^{\dagger}$
of $\prolong{\nbige^{\dagger}}$ over $\cnum_{\mu}\times \Delta$.
It is also compatible with the filtration $\weightfilt^{\sankaku}$.
Note that
$\deg^{\weightfilt^{\sankaku}}(v_i)=
-\deg^{\weightfilt^{\sankaku}}(v_i^{\dagger})$.
The weight filtration $\weightfilt^{\sankaku}$
over $\cnum_{\mu}\times O$ is given as follows:
\[
 \weightfilt^{\sankaku}_h=
 \big\langle
 v^{\dagger}_{i\,|\,\cnum_{\lambda}\times O}
 \,\big|\,\deg^{\weightfilt^{\sankaku}}(v_i^{\dagger})\leq h
 \big\rangle
=\big\langle
 v^{\dagger}_{i\,|\,\cnum_{\lambda}\times O}
 \,\big|\,\deg^{\weightfilt^{\sankaku}}(v_i)\geq  -h
 \big\rangle
\]

Let take holomorphic vector bundles
$F:=\bigoplus_{i=1}^r \nbigo_{\cnum_{\lambda}}\vecu_i$
over $\cnum_{\lambda}$.
We give the filtration $\nbigw_F$ to
the bundle $F$, defined as follows:
\[
 \nbigw_{F,h}:=
 \big\langle
 u_i\,\big|\,\deg^{\weightfilt^{\sankaku}}(v_i)\leq h
 \big\rangle.
\]
Similarly we put
$F^{\dagger}:=\bigoplus_{i=1}^r\nbigo_{\cnum_{\mu}}\vecu_i^{\dagger}$.
We give the filtration $\nbigw_{F^{\dagger}}$ defined as follows:
\[
 \nbigw_{F^{\dagger},h}
:=\big\langle
 u_i^{\dagger}\,\big|\,\deg^{\weightfilt^{\sankaku}}(v_i^{\dagger})\leq h
 \big\rangle
=\big\langle
 u_i^{\dagger}\,\big|\,\deg^{\weightfilt^{\sankaku}}(v_i)\geq -h
 \big\rangle
\]
We also put as follows:
\[
 U_h:=
 \langle u_i\,|\, \deg^{\nbigw^{\sankaku}}(v_i)=h
 \rangle,
 \quad
 U_h^{\dagger}:=
 \langle u_i^{\dagger}\,|\,\deg^{\nbigw^{\sankaku}}(v_i^{\dagger})=h
 \rangle.
\]
Then we have $\nbigw_{F,h}=\bigoplus_{l\leq h} U_l$
and $\nbigw_{F^{\dagger},h}=\bigoplus_{l\leq h}U_{l}^{\dagger}$.
Note the following obvious lemma.
\begin{lem}
Let $\vecc=(c_i\,|\,i=1,\ldots, r)$ be a tuple of non-zero complex numbers.
We have the morphism
$\eta_{\vecc}: \deformoverorigin\lrarr F$ defined by
$v_i\longmapsto c_i\cdot u_i$ for $1\leq i\leq r$.
Then the morphism $\eta_{\vecc}$ preserves the filtrations
$\weightfilt^{\sankaku}$ and $\nbigw_F$.

Similarly we have the morphism 
$\eta^{\dagger}_{\vecc}:\deformoverorigindagger\lrarr F^{\dagger}$
defined by
$v_i^{\dagger}\longmapsto c_i\cdot u_i$ for $1\leq i\leq r$.
Then the morphism $\eta^{\dagger}$ preserves the filtrations
$\weightfilt^{\sankaku}$ and $\nbigw_{F^{\dagger}}$.
\end{lem}
\pf
Clear from our choices of $\nbigw_F$ and $\nbigw_{F^{\dagger}}$.
\hfill\qed

For any $P\in\Delta^{\ast}$,
we take the tuple of complex numbers
$\vecc(P)=\bigl(c_i(P)\bigr)$
and $\vecc^{\dagger}(P)=\bigl(c_i^{\dagger}(P)\bigr)$,
given as follows:
\[
 c_i(P):=(-\log|z(P)|)^{k(v_i)},
\quad
 c_i^{\dagger}(P):=(-\log|z(P)|)^{k(v_i^{\dagger})}
=(-\log|z(P)|)^{-k(v_i)}.
\]
Here we put $2\cdot k(v_i):=\deg^{\weightfilt^{\sankaku}}(v_i)$
and
$2\cdot k(v_i^{\dagger}):=\deg^{\weightfilt^{\sankaku}}(v_i^{\dagger})$.
Then we obtain the following isomorphisms
preserving the filtrations:
\[
 \eta_{\vecc(P)}:
 \bigl(
 \prolong{\nbige}_{|\cnum_{\lambda}\times O},\weightfilt
 \bigr)
 \simeq
 \bigl(
 F,\nbigw_{F}
 \bigr),
\quad
 \eta_{\vecc^{\dagger}(P)}:
 \bigl(
 \prolong{\nbige^{\dagger}}_{|\cnum_{\mu}\times O},\weightfilt
 \bigr)
 \simeq
 \bigl(
 F^{\dagger},\nbigw_{F^{\dagger}}
 \bigr).
\]
Then the gluing morphism $\Psi(O,P)$ for the vector bundle $S(O,P)$
induces 
the isomorphism
$g(P):
 (F,\nbigw_{F})_{|\cnum_{\lambda}^{\ast}}
\simeq 
 (F^{\dagger},\nbigw_{F^{\dagger}})_{|\cnum_{\mu}^{\ast}}$.
The filtered vector bundle $(S(O,P),\weightfilt^{\sankaku})$
is naturally isomorphic to
the filtered vector bundle obtained from
$(F,\nbigw_{F})$, $(F^{\dagger},\nbigw_{F^{\dagger}})$
and the gluing $g(P)$.

We have the diagonal matrix $L(P)$ whose $(i,i)$-component
is $c_i(P)=(-\log|z(P)|)^{k(v_i)}$.
We also have the diagonal matrix
$L^{\dagger}(P)$ whose $(i,i)$-component
is $c_i^{\dagger}(P)=(-\log|z(P)|)^{k(v_i^{\dagger})}=
 (-\log|z(P))|^{-k(v_i)}$.
The matrix $L^{\dagger}(P)$ is the inverse of $L(P)$.
Let $A(\lambda,P)$ is an element of
$\Gamma(\cnum_{\lambda}^{\ast},GL(r))$
given in (\ref{eq;12.2.1}).

\begin{lem}
The gluing $g(P)$ is represented by
$L(P)\cdot A(\lambda,P)\cdot L^{\dagger}(P)^{-1}$
with respect to the frames $\vecu$ and $\vecu^{\dagger}$.
\end{lem}
\pf
Since $A(\lambda,P)$ represents $\Psi(O,P)$
with respect to the frames
$\vecv_{|\cnum_{\lambda}\times O}$
and
$\vecv_{|\cnum_{\mu}\times O}$,
the claim is clear from our choices of 
$\eta_{\vecc}$ and $\eta^{\dagger}_{\vecc^{\dagger}}$.
\hfill\qed

In the following,
we assume that we take $\vece=\vecv_{|\{0\}\times \Delta^{\ast}}$
for the construction of $A(\lambda)$.

\subsubsection{The statement and an outline of a proof}
\label{subsection;6.25.1}

We will prove the following theorem.
\begin{thm}[A limiting mixed twistor theorem]\label{thm;1.8.1}
Let $\harmonicbundle$ be a tame nilpotent harmonic bundle
with parabolic structure over $\Delta^{\ast}$.
For any open set $U$ containing $O$,
there is a point $P\in U\cap\Delta^{\ast}$
such that the filtered vector bundle $(S(O,P),\weightfilt^{\sankaku})$
gives a mixed twistor structure.
\end{thm}
An outline of the proof of the theorem is as follows:
\begin{enumerate}
\item
We can take some sequence of the points $\{P_i\}$
converging to $O$, such that
the corresponding sequence of the gluing functions $\{g(P_i)\}$
converge to the gluing function $g_{\infty}$
(the subsubsection \ref{subsubsection;12.2.3}).
\item
We have the filtered vector bundle $(S_{\infty},\nbigw)$
obtained from $(F,\nbigw_F)$, $(F^{\dagger},\nbigw_{F^{\dagger}})$
and $g_{\infty}$.
We will see that
$(S_{\infty},\nbigw)$ is a mixed twistor.
(the subsubsection \ref{subsubsection;12.10.5})
\item
Let consider
the condition that a filtered vector bundle is a mixed twistor.
We can see that the condition is open.
Since the sequence of 
the filtered vector bundles $(S(O,P_i),\weightfilt^{\sankaku})$
converge to $(S_{\infty},\nbigw)$ in a sense,
we can conclude that
they are mixed twistors for sufficiently large $i$
(the subsubsection \ref{subsubsection;12.10.6}).
\end{enumerate}

Before entering the proof,
we state an obvious corollary.
\begin{cor}
Let $P$ be any point of $\Delta^{\ast}$.
The Chern class $c_1(\graded^{\sankaku}_h)$ is
$h\cdot\rank(\graded^{\sankaku}_h)$.
\end{cor}
\pf
If $(S(O,P),\weightfilt^{\sankaku})$ is mixed twistor,
the claim is clear.
Since the Chern class is an topological invariant,
we obtain the result.
\hfill\qed

\subsubsection{Convergency of the gluing functions}
\label{subsubsection;12.2.3}

We only have to investigate the convergency of 
the sequence $\{L(P_l)\cdot A(P_l)\cdot L^{\dagger}(P_l)^{-1}\}$
for some sequence $\{P_l\}$.
We decompose $L(P)\cdot A(P)\cdot L^{\dagger}(P)^{-1}$
as follows:
\begin{multline}
 \Bigl(
 L(P)\cdot
  \Gamma(\lambda,O)\cdot
 \Gamma(\lambda,P)^{-1}
\cdot
 L(P)^{-1}
\Bigr)
 \cdot
 \Bigl(
 L(P)\cdot
 B(\lambda,P)^{-1}
 L(P)^{-1}
\Bigr)\\
\times
 \Bigl(
 L(P)\cdot
 \overline{H(h,\vece)}^{-1}(P) 
\cdot
 L^{\dagger}(P)^{-1}
 \Bigr)\\
\times
 \Bigl(
 L^{\dagger}(P)\cdot
 B^{\dagger}(\mu,P)\cdot
 L^{\dagger}(P)^{-1}
 \Bigr)
\cdot
 \Bigl(
 L^{\dagger}(P)\cdot
 \Gamma^{\dagger}(\mu,P)\cdot
 \Gamma^{\dagger}(\mu,O)^{-1}\cdot
 L^{\dagger}(P)^{-1}
 \Bigr).
\end{multline}

\begin{lem}
Let $R$ be a real number such that $0<R<1$.
On the region
$T(R)=\bigl\{\lambda\in\cnum_{\lambda}\,\big|\,
  R<|\lambda|<R^{-1}\bigr\}$,
the sequence
$\{L(P_l)\cdot \Gamma(\lambda,O)\cdot\Gamma(\lambda,P_l)^{-1}
 \cdot L(P_l)^{-1}\}$
converges to the identity matrix, with respect to the sup norms,
for any sequence $\{P_l\}$ converging to $O$.
\end{lem}
\pf
The sequence
$\bigl\{ \Gamma(\lambda,O)\cdot\Gamma(\lambda,P_l)^{-1}\bigr\}$
converges to the identity matrix for any sequence $\{P_l\}$.
In fact, we have an equality
$\bigl|\Gamma(\lambda,O)\cdot \Gamma(\lambda,P)^{-1}-1\bigr|
\leq C\cdot |z(P)|$ over the region $T(R)$,
for some positive constant $C$.
\hfill\qed

\begin{lem}
Let $\bigl\{P_l\bigr\}$ be a sequence of points of $\Delta^{\ast}$,
converging to $O$.
Then we can take a subsequence $\bigl\{P_{l_i}\bigr\}$ 
such that
the corresponding sequence
$\bigl\{  L(P_{l_i})\cdot
 \overline{H(h,\vece)}^{-1}(P_{l_i}) 
\cdot
 L^{\dagger}(P_{l_i})^{-1}\bigr\}$ converges
to a positive definite hermitian matrix $M$.
\end{lem}
\pf
We put $e_i':=(-\log|z|)^{-k(e_i)}\cdot e_i$ for $i=1,\ldots,r$.
Then the frame $\vece'=(e_1',\ldots,e_r')$ of $\nbige^0$
over $\Delta^{\ast}$ is adapted.
We have the following equality:
\[
 H(h,\vece')=
 L(P)^{-1}\cdot H(h,\vece)\cdot L(P)^{-1}
=L^{\dagger}(P)\cdot H(h,\vece)\cdot L(P)^{-1}.
\]
Thus we are done.
\hfill\qed

\vspace{.1in}
Let us investigate the convergency of 
$\bigl\{L(P)^{-1}\cdot B(\lambda,P)\cdot L(P)\bigr\}$.
For that purpose, we take a model bundle.
We put $V=\prolong{\nbige^0}_{|O}$ and 
$N=\Res(\theta)$.
Then we have a model bundle
$(E_0,\theta_0,h_0)=E(V,N)$ over $\Delta^{\ast}$.
We denote the deformed holomorphic bundle of $(E_0,\theta_0,h_0)$ 
by $(\nbige_0,\DD_0,h_0)$.
We have the canonical frame $\vecv_0$ of $\prolong{\nbige_0}$
over $\cnum_{\lambda}\times \Delta$.
We can assume $\deg^{\weightfilt}(v_{0\,i})=\deg^{\weightfilt}(v_i)$.

We put $\vece_0:=\vecv_{0\,|\{0\}\times \Delta}$.
Then the frames $\vece_0$ and $\vece$ induce
the isomorphism $F$ of the holomorphic bundles
$\prolong{\nbige^0_0}$ and $\prolong{\nbige_0}$,
such that 
$F_{|O}\circ \Res(\theta_0)=\Res(\theta)\circ F_{|O}$.
The morphism $F$ induces the holomorphic isomorphism
$F:\nbige_0^0\lrarr \nbige_0$ over $\Delta^{\ast}$,
such that $F$ and the inverse $F^{-1}$
are bounded with respect to the metrics $h_0$ and $h$.

Since we have
$\nbige_0=p_{\lambda}^{-1}\nbige_0^0$ and
$\nbige=p_{\lambda}^{-1}\nbige_0$
by our definition,
we obtain a $C^{\infty}$-isomorphism
$\nbigi:\nbige_0\lrarr\nbige$ over $\cnum_{\lambda}\times \Delta^{\ast}$.
The morphism $\nbigi$ and the inverse $\nbigi^{-1}$
are bounded with respect to the metrics
$h$ and $h_0$ by our construction.

We have the $GL(r)$-valued function $I$
determined by
$\nbigi(\vecv_0)=\vecv\cdot I$.
We have the $GL(r)$-valued function $B_0$
determined by
$\vecv_0=p_{\lambda}^{-1}(\vece_0)\cdot B_0$.
\begin{lem}
We have the relation
$B=B_0\cdot I^{-1}$.
\end{lem}
\pf
We have $\vecv=p_{\lambda}^{-1}(\vece)\cdot B$.
On the other hand, we have the following:
\[
 \vecv=\nbigi(\vecv_0)\cdot I^{-1}=
 \nbigi(p_{\lambda}^{-1}\vece_0)\cdot B_0\cdot I^{-1}.
\]
By our choice of $\nbigi$,
we have
$\nbigi(p_{\lambda}^{-1}\vece_0)=p_{\lambda}^{-1}\vece$.
Thus we are done.
\hfill\qed

Hence we have the following:
\[
 L(P)\cdot B(\lambda,P)\cdot L(P)^{-1}
=\Bigl(
 L(P)\cdot B_0(\lambda,P)\cdot L(P)^{-1}
 \Bigr)
\times
 \Bigl(
 L(P)\cdot I^{-1}(\lambda,P)\cdot L(P)^{-1}
 \Bigr).
\]

For the model bundle,
we have already known the behaviour of the transformation matrices
$B_0(P,\lambda)$ when $P\to 0$.
\begin{lem} \label{lem;1.27.1}
When $P\to O$,
the sequence $\bigl\{L(P)\cdot B_0(\lambda,P)\cdot L(P)^{-1}\bigr\}$
converges to
a $GL(r)$-valued holomorphic function
$\bar{B}(\lambda)=\bigl(\bar{B}_{i\,j}(\lambda)\bigr)$
over $\cnum_{\lambda}$ satisfying the following:
\begin{itemize}
\item We put
$k(v_i)
:=2^{-1}\cdot\deg^{\weightfilt}(v_i)
=2^{-1}\cdot\deg^{\weightfilt}(v_{0\,i})$.
Then we have the following:
\[
\bar{B}_{i\,j}(\lambda)=
\left\{
\begin{array}{ll}
0 & \mbox{\rm (if $k(v_i)<k(v_j)$ or if $k(v_j)-k(v_i)$ is not an integer),}\\
b_{i\,j}\cdot\lambda^{k(v_i)-k(v_j) },\,(b_{i\,j}\in \cnum)
 &
\mbox{\rm (if $k(v_j)-k(v_i)$ is an integer).}
\end{array}
\right.
\]
\item
The matrices $\bar{B}(h,h):=(b_{i\,j}\,|\,
  \deg^{\weightfilt}(v_i)=\deg^{\weightfilt}(v_j)=h)$
are invertible for any $h$.
\item
The matrices
$\bar{B}(h,-h):=(b_{i\,j}\,|\,
  \deg^{\weightfilt}(v_i)=h,\deg^{\weightfilt}(v_j)=-h )$ 
are invertible for any $h\leq 0$.
\end{itemize}

A similar convergence holds for
the sequence $L^{\dagger}(P)B_0^{\dagger}(P)L^{\dagger}(P)^{-1}$,
when $P\to O$.
\end{lem}
\pf
The claims can be checked directly.
(See Corollary \ref{cor;12.2.5}.)
\hfill\qed

\vspace{.1in}

We put $C(P):=L(P)\cdot I(\lambda,P)\cdot L(P)^{-1}$,
which is an element of
$C^{\infty}(\cnum_{\lambda}\times \Delta^{\ast},GL(r))$.
\begin{lem}
Let $K$ be a compact subset of $\cnum_{\lambda}$.
Then $C(P)$ and the inverse $C(P)^{-1}$ are 
bounded over $K\times \Delta^{\ast}$.
\end{lem}
\pf
We put $v_i':=(-\log|z|)^{-k(v_i)}\cdot v_i$.
Then we obtain $C^{\infty}$-frame $\vecv'=(v_i')$
of $\nbige$ over $\cnum_{\lambda}\times \Delta^{\ast}$.
We also put $v_{0\,i}'=(-\log|z|)^{-k(v_{0\,i})}\cdot v_{0\,i}$,
and then we obtain the $C^{\infty}$-frame $\vecv'_0$
of $\nbige_0$.
For any compact subset $K\subset \cnum_{\lambda}$,
the frames $\vecv'$ and $\vecv'_0$ are adapted
over $K\times \Delta^{\ast}$,
for the metrics $h$ and $h_0$ respectively.
Since $\nbigi$ and the inverse $\nbigi^{-1}$
are bounded over $K\times \Delta^{\ast}$,
the $C^{\infty}$-frame
$\nbigi(\vecv_0')$ is adapted for the metric $h$.
Then the transformation matrices
between $\vecv'$ and $\nbigi(\vecv_0')$ are bounded.
Now we have the relation
$\nbigi(\vecv'_0)=\vecv'\cdot C$
by our construction.
Thus $C$ and the inverse $C^{-1}$
are bounded over $K\times \Delta^{\ast}$.
\hfill\qed

We have the following immediate corollary.
\begin{cor}
Let $\{P_l\}$ be a sequence of points in $\Delta^{\ast}$
converging to $O$.
Then we can take a subsequence $\{P_{l_i}\}$
such that the corresponding sequence
$\{C(P_l)\}$ converges to 
a $GL(r)$-valued holomorphic function $\bar{C}$
with respect to the sup norm,
over the region $\Delta_{\lambda}(R)$ for any  $0<R<1$.
\hfill\qed
\end{cor}

Then we have already known the existence of the sequence
$\{P_l\}$ such that 
$\bigl\{L(P_l)\cdot A(\lambda,P_l)\cdot L^{\dagger}(P_l)^{-1}\bigr\}$
converges to a $GL(r)$-valued holomorphic functions
over the region $T(R)=\bigl\{\lambda\,|\,R<|\lambda|<R^{-1}\bigr\}$.
However, we need a better sequence.

\begin{lem}
There is a sequence $\{P_l\}$
satisfying the following:
\begin{itemize}
\item
The sequence $\bigl\{C(P_l)\bigr\}$
converges to $GL(r)$-valued holomorphic function $\bar{C}$
over $\cnum_{\lambda}$.
\item
The $(i,j)$-components $\bar{C}_{i\,j}$ are $0$
if $\deg(v_i)\neq \deg(v_{0\,j})$.
In other words,
$\bar{C}$ preserves the decomposition $F=\bigoplus_h U_h$.
\end{itemize}
\end{lem}
\pf
We only have to see that 
we can take a sequence $\{P_l\}$
satisfying  $C_{i\,j}(P_l)\to 0$
for any pair $(i,j)$ such that $\deg(v_{0\,j})\neq \deg(v_i)$.

Assume that $\deg(v_{0\,j})\neq \deg(v_i)$.
Due to Lemma \ref{lem;1.14.1} and 
Lemma \ref{lem;12.2.10},
we obtain the following finiteness:
\[
 \Bigl|\Bigl|
 \int_{|\lambda|<R}|C_{i\,j}(\lambda,z)|
 \cdot |d\lambda \cdot d\bar{\lambda}|
 \Bigr|\Bigr|_W
=\int _{|\lambda|<R} ||C_{i\,j}(\lambda,z)||_W
 \cdot |d\lambda \cdot d\bar{\lambda}|
 <\infty.
\]
Thus we can take a sequence $\{P_l\}$
such that the sequence
$\bigl\{
 \int_{|\lambda|<R}|C_{i\,j}(\lambda,P_l)|\cdot |d\lambda \cdot d\bar{\lambda}|
 \bigr\}$
converges to $0$ for any $R>0$
Since the components $C_{i\,j}(\lambda,P_l)$ depend on $\lambda$
holomorphically,
the sequence $\{C_{i\,j}(\lambda,P_l)\}$ converges to $0$
on $\Delta_{\lambda}(R)$ for any $R>0$.
\hfill\qed

In all,
we obtain the following:
\begin{prop}
There exists a sequence $\{P_l\}$ of points in $\Delta^{\ast}$
converging to $O$ satisfying the following:
\begin{itemize}
\item The corresponding sequence
 $\bigl\{L(P)\cdot A(\lambda,P)\cdot L^{\dagger}(P)^{-1}\bigr\}$
converges to a holomorphic $GL(r)$-valued function $\bar{A}$
on $\cnum_{\lambda}^{\ast}$.
The convergence is with respect to the sup norms
on $T(R):=\{\lambda\in\cnum_{\lambda}\,|\,R<|\lambda|<R^{-1}\}$
for any $0<R<1$.
\item
We have the decomposition of $\bar{A}$
into the product:
\[
 \bar{A}= \bar{C}\cdot\bar{B}_0^{-1}\cdot M
 \cdot \bar{B}_0^{\dagger}\cdot \bar{C}^{\dagger}.
\]
\item
$M$ is a positive definite hermitian matrix.
It is independent of $\lambda$ and $\mu$.
\item
$\bar{B}$ is a holomorphic $GL(r)$-valued function defined over
$\cnum_{\lambda}$.
It is given in Lemma {\rm \ref{lem;1.27.1}}.
\item
$\bar{C}$ is a holomorphic $GL(r)$-valued function defined over
$\cnum_{\lambda}$.
It preserves the decomposition $F=\bigoplus_h U_h$.
\item
$\bar{B}^{\dagger}$ and $\bar{C}^{\dagger}$
are $GL(r)$-valued holomorphic functions
defined over $\cnum_{\mu}$.
They are given as follows:
\[
 \bar{B}^{\dagger}(\mu)={}^t\overline{\bar{B}(-\lambdabar)}^{-1},
\quad
 \bar{C}^{\dagger}(\mu)={}^t\overline{\bar{C}(-\lambdabar)}^{-1}.
\]
\item
Similar things hold for the exterior products
$\bigwedge^l(E,\theta,h)$.
\end{itemize}
\end{prop}
\pf
We have already seen the first five properties.
The last two properties are clear
from our construction.
Note that the transformation matrices
for $S\bigl(O,P,\bigwedge^l(E,\theta,h)\bigr)$
is obtained from the transformation matrices of
$S\bigl(O,P,(E,\theta,h)\bigr)$,
by some standard linear algebraic procedures.
\hfill\qed

\subsubsection{The property of $g_{\infty}$}

\label{subsubsection;12.10.5}
Let $g_{\infty}$ denote the gluing given by the matrices
$\bar{A}$.
From the filtered vector bundles
$(F,\nbigw_{F})$,
$(F^{\dagger},\nbigw_{F^{\dagger}})$
and the gluing $g_{\infty}$,
we obtain the filtered vector bundle.
We denote it by $(S^{(\infty)},\nbigw^{(\infty)})$.
We denote the associated graded vector bundle
by $\bigoplus_{h}\graded_h$.
\begin{prop} \label{prop;12.2.15}
The filtered vector bundle $(S^{(\infty)},\nbigw^{(\infty)})$
is a mixed twistor.
Namely $\graded_h$ is isomorphic to a direct sum
of $\nbigo_{\proj^1}(h)$.
\end{prop}
\pf
First we see the following:
\begin{lem}
We only have to care
the transformation given by $\bar{B}^{-1}\cdot M\cdot
 \bar{B}^{\dagger}$.
\end{lem}
\pf
The matrix valued function $\bar{C}$
and the inverse are holomorphic for the variable $\lambda$
over $\cnum_{\lambda}$.
Moreover they preserve the decomposition
$\bigoplus U_h$.
In particular, they preserve the filtration
$\nbigw_F$ over $\cnum_{\lambda}$.
Thus we can ignore $\bar{C}$ to see
$\graded_h$.
Similarly we can ignore $\bar{C}^{\dagger}$.
\hfill\qed

We denote the associated graded vector bundle
of $\nbigw_F$ and $\nbigw_{F^{\dagger}}$
by $\graded_F$ and $\graded_{F^{\dagger}}$.
Since $\bar{B}^{-1}\cdot M\cdot \bar{B}^{\dagger}$
preserves the filtrations
$\nbigw_F$ and $\nbigw_{F^{\dagger}}$,
we obtain the morphisms
of the graded parts.
We denote them as follows:
\[
 Gr_l(\bar{B}^{-1}\cdot M\cdot\bar{B}^{\dagger}):
 \graded^{\nbigw_{F}}_l\lrarr
 \graded^{\nbigw_F^{\dagger}}_l.
\]

Let $b(W)$
be the bottom number of the filtration $\nbigw^{(\infty)}$.
\begin{lem} \label{lem;12.16.1}
The morphism
$Gr_{b(W)}(\bar{B}^{-1}M\bar{B}^{\dagger}):
 \graded_{b(W)}^{\nbigw_F}\lrarr \graded_{b(W)}^{\nbigw_{F^{\dagger}}}$
is of the form $\lambda^{-b(W)}\times \Phi_{b(W)}(\mu)$,
where $\Phi_{b(W)}$ is holomorphic and invertible on $\cnum_{\mu}$.
In particular,
$\graded_{b(W)}$ is a pure twistor of weight $b(W)$.
\end{lem}
\pf
We can regard the matrices $\bar{B}^{-1}$ and $\bar{B}^{\dagger}$
as the endomorphisms of
$F$ and $F^{\dagger}$.
And the hermitian matrix $M$ can be regarded as the
transformation of $\vecu=(u_i)$ to $\vecu^{\dagger}=(u_i^{\dagger})$.
We denote the submatrices
$\bigl((\bar{B}^{-1})_{i\,j}\,|\,\deg(u_i)=x,\deg(u_j)=y\bigr)$
by $\bar{B}^{-1}(x,y)$.
Similarly we put as follows:
\[
\begin{array}{l}
 \bar{B}^{\dagger}(x,y):=
 \bigl(\bar{B}_{i\,j}^{\dagger}\,|\,
 \deg(u_i^{\dagger})=x,\deg(u_j^{\dagger})=y\bigr),\\
\mbox{{}}\\
 M(x,y):=\bigl(M_{i\,j}\,|\,\deg(u_i)=x,\deg(u_j^{\dagger})=y\bigr),\\
\mbox{{}}\\
 (\bar{B}^{-1}\cdot M\cdot\bar{B}^{\dagger})(x,y):=
 \bigl((\bar{B}^{-1}\cdot M\cdot\bar{B}^{\dagger})_{i\,j}\,|\,
 \deg(u_i)=x,\deg(u_j^{\dagger})=y\bigr).
\end{array}
\]
We will calculate $(\bar{B}^{-1}\cdot M\cdot \bar{B}^{\dagger})(b(W),b(W))$.
Note the following obvious lemma.
\begin{lem}
We have $\bar{B}^{\dagger}(x,y)=0$ for any $x>y$.
In particular $\bar{B}^{\dagger}(x,b(W))=0$
if $x\neq b(W)$.
\hfill\qed
\end{lem}
Then we obtain the following equalities:
\begin{multline}
 (\bar{B}^{-1}\cdot M\cdot\bar{B}^{\dagger})(b(W),b(W))
=\sum_y \bar{B}^{-1}(b(W),y)\cdot M(y,b(W))\cdot\bar{B}^{\dagger}(b(W),b(W))\\
=\bar{B}^{-1}(b(W),-b(W))\cdot M(-b(W),b(W))\cdot
 \bar{B}^{\dagger}(b(W),b(W))\\
+
 \sum_{y<-b(W)}\bar{B}^{-1}(b(W),y)\cdot M(y,b(W))\cdot
 \bar{B}^{\dagger}(b(W),b(W)).
\end{multline}
If follows from Lemma \ref{lem;1.27.1} that
the function $\bar{B}^{-1}(b(W),-b(W))$ is of the form $\lambda^{-b(W)}Q_1$
where $Q_1$ is an element of $GL_{r'}(\cnum)$.
Here $r'$ denotes $\dim\graded_{b(W)}$.
The matrix $M(-b(W),b(W))$ is positive hermitian.
Note that we have the equality $\deg(u_i)=-\deg(u_i^{\dagger})$.
Moreover $\bar{B}^{\dagger}(b(W),b(W))$
is an element of $GL_{r'}(\cnum)$.
Thus the first term in the right hand side
$\bar{B}^{-1}(b(W),-b(W))\cdot M(-b(W),b(W))\cdot 
 \bar{B}^{\dagger}(b(W),b(W))$
is of the form $\lambda^{-b(W)}Q_2$ where $Q_2$
is an element of $GL_{r'}(\cnum)$.

If follows from Lemma \ref{lem;1.27.1} that
the function $\bar{B}^{-1}(b(W),y)$\, $(y<-b(W))$
is of the form $\lambda^{ (-b(W)+y)/2}S_y$
where $S_y$ denotes a matrix with $\cnum$-coefficient
and $S_y$ is $0$ unless $b(W)-y$ is even.
Also $M(y,b(W))$ is a matrix with $\cnum$-coefficient.
Thus we have the following equality:
\[
\begin{array}{l}
 (\bar{B}^{-1}M\bar{B}^{\dagger})(b(W),b(W))
=\lambda^{-b(W)} \times \Phi_{b(W)}(\mu)\\
\mbox{}\\
\Phi_{b(W)}(\mu)=
 Q_2+
 \sum_{i>0}\mu^iS'_i.
\end{array}
\]
Here $S_i'$ $(i>0)$ denote the matrices with $\cnum$-coefficient.
The function $\Phi_{b(W)}(\mu)$ is holomorphic and invertible
over $\cnum_{\mu}$.

In particular, $\Phi_{b(W)}(\mu)$ is just a change of the trivialization
of $\graded^{\nbigw_{F^{\dagger}}}_{b(W)}$ over $\cnum_{\mu}$.
As a result, we can conclude that
$\graded_{b(W)}$ is a pure twistor of weight $b(W)$.
Thus the proof of Lemma \ref{lem;12.16.1} is completed.
\hfill\qed

\vspace{.1in}
Then we will show the following claim by using an induction on $h$.
\begin{description}
\item[($\boldsymbol {P_h}$)]
For any $l<h$, the graded part
$\graded_l$ is a pure twistor of weight $l$.
\end{description}

We assume that $P_{h-1}$ holds
and show that $P_h$ holds.
By assumption,
$\graded_l$ are pure twistor of weight $l$.
In particular the first Chern class of $\det(\graded_l)$
is $l\cdot\rank(\graded_l)$.

We put $R=\rank \nbigw^{(\infty)}_{h-1}+1$.
Note that we have already seen that the bottom part
is a pure twistor of an appropriate weight
for any tame nilpotent harmonic bundle with nilpotent residues.
In particular, 
the bottom part of the filtered vector bundle
$\bigwedge^R(S^{(\infty)},\nbigw^{(\infty)})$ is
pure twistor of the following weight:
\[
 b_0:=
 \sum_{l=1}^{h-1} l\cdot \rank \graded_l
 +h=
\sum_{l=1}^{h-1}c_1(\det \graded_l)+h
=c_1(\det\nbigw_{h-1}\superinfty)+h.
\]
There is the natural isomorphism:
\[
 \graded_{b_0}
 \Bigl( \bigwedge^R\bigl(S^{(\infty)},\nbigw^{(\infty)}\bigr)
 \Bigr)
\simeq
 \det(\nbigw^{(\infty)}_{h-1})
\otimes
 \graded_h.
\]
Thus we can conclude that $\graded_h$ is
a pure twistor of weight $h$.

Hence the proof of Proposition \ref{prop;12.2.15}
is completed.
\hfill\qed

\subsubsection{The end of the proof of a limiting mixed twistor theorem}
\label{subsubsection;12.10.6}

For a positive number $R$ such that $1<R<\infty$,
we put $T(R):=\{z\in\cnum\,|\,R^{-1}<|z|<R\}$.
For any holomorphic function $g:T(R)\lrarr GL(n)$,
we can naturally associate the holomorphic vector bundle
of rank $n$,
which we denote by $V(g)$.
We denote the set of holomorphic function
$g:T(R)\lrarr GL(n)$ such that $c_1(V(g))=0$
by $\nbigc_0$.
We have the subset
$\nbigc_{triv}=\{g\in\nbigc_0\,|\,V(g)\mbox{ trivial }\}$.
\begin{lem} \label{lem;1.27.2}
$\nbigc_{triv}$ is open with respect to the topology
given by the sup norm,
i.e.,
if a sequence $\{g_i\in \nbigc_0\}$ 
converges to an element $g\in \nbigc_{tirv}$,
then $g_i\in \nbigc_{triv}$ for any sufficiently large $i$.
\end{lem}
\pf
We can translate the sequence of the gluing $g_i$
to the sequence of holomorphic structures $\bar{\del}_i$
on a $C^{\infty}$-vector bundle $E$ such that $c_1(E)=0$.
The sequence $\bar{\del}_i$ converges to 
a holomorphic structure $\bar{\del}$
in $L_l^p$ for any $l$,
such that $(E,\bar{\del})$ is isomorphic to the trivial bundle.
Then we have to show that $\bar{\del}_i$ gives
the holomorphic structure which is holomorphically trivial
for any sufficiently large $i$.
It is a consequence of the vanishing $H^1(\proj^1,\nbigo)=0$.
\hfill\qed

By Lemma \ref{lem;1.27.2},
we obtain the open-ness  of the condition
that a filtered vector bundle is a mixed twistor.
We have already known that
the sequence of the gluings $\{g(P_l)\}$
converges to a gluing $g(P)$, which gives a mixed twistor.
Then we can conclude that $g(P_l)$ gives
a mixed twistor if $l$ is sufficiently large.
Thus the proof of Theorem \ref{thm;1.8.1} is completed.
\hfill\qed

\subsection{Higher dimensional case of a limiting mixed twistor theorem}

\subsubsection{The morphisms  induced by the residues}

\label{subsubsection;12.15.24}

We put $X=\Delta^n$, $D_i=\{z_i=0\}$,
and $D=\bigcup_{i=1}^lD_i$ for some $l\leq n$.

Let consider the tame nilpotent harmonic bundle $\harmonicbundle$
with trivial parabolic structure
over $X-D$.
Let take a point $Q\in D$. We put
$I=\{i\in\lbar\,|\,Q\in D_i\}$.
Let take a point $P$ and consider $S(O,P)$
over $\proj^1$ for $\harmonicbundle$.
We have the nilpotent maps induced by the residues
of $\DD$ and $\DD^{\dagger}$ at $\nbigd_i$ and $\nbigd_i^{\dagger}$
respectively:
\[
 N_i^{\sankaku}:S(Q,P)\lrarr S(Q,P)\otimes\nbigo_{\proj}(2),
\quad
 (i\in I).
\]
For a tuple $\veca=(a_i\,|\,i\in I)\in\cnum^I$,
we put $N^{\sankaku}(\veca):=\sum_{i\in I} a_i\cdot N_i^{\sankaku}$.

\begin{lem}
Assume that all of $a_i$ are positive integers.
\begin{itemize}
\item
The conjugacy classes of $N^{\sankaku}(\veca)_{|\lambda}$
are independent of the choice of $\lambda\in\proj^1$.
\item
Let $W^{\sankaku}(\veca)$ denote
the weight filtration of $N^{\sankaku}(\veca)$.
For an appropriate point $P$,
the filtered vector bundle $(S(Q,P),W^{\sankaku}(\veca))$
is a mixed twistor.
\end{itemize}
\end{lem}
\pf
Let consider the embedding $\varphi:\Delta\lrarr X$
given as follows:
\begin{equation} \label{eq;12.10.10}
 z_i\bigl(\varphi(t)\bigr)=
 \left\{
 \begin{array}{ll}
 t^{a_i} & (i\in I)\\
 z_i(Q)\neq 0  & (i\not\in I)
 \end{array}
 \right.
\end{equation}
We denote the origin of $\Delta$ by $O$.
Let take a point $\blowup{P}\in\Delta$
such that $\varphi(\blowup{P})=P$.

We obtain the harmonic bundle
$\varphi^{\ast}\harmonicbundle$
over $\Delta^{\ast}$.
By our construction and the functoriality of the prolongation,
the residue $\Res( \varphi^{\ast}\DD)$ is
isomorphic to $N^{\sankaku}(\veca)_{|\cnum_{\lambda}}$.
Thus the conjugacy classes of $N^{\sankaku}(\veca)_{|\cnum_{\lambda}}$
are independent of a choice of $\lambda\in\cnum_{\lambda}$.
Similarly the conjugacy class of $N^{\sankaku}(\veca)_{|\cnum_{\mu}}$
are independent of a choice of $\mu\in\cnum_{\mu}$.
Moreover, the filtration $W^{\sankaku}(\veca)$
gives a mixed twistor if we take an appropriate point $\blowup{P}$
of $\Delta^{\ast}$,
which is a consequence of Theorem \ref{thm;1.8.1}.
\hfill\qed

\vspace{.1in}

We can take a general $\veca$,
in the sense of \ref{df;12.3.40},
from $\rnum_{>0}^{I}$.
Let take an appropriate point $P$
such that $(S(Q,P),\weightfilt^{\sankaku}(\veca))$
is a mixed twistor.
We denote the associated graded vector bundle
by $\graded^{\sankaku}(\veca)$.

Let consider the morphism
$N_i^{\sankaku}:S(Q,P)\lrarr S(Q,P)\otimes\nbigo_{\proj^1}(2)$
for $i\in I$.
Due to lemma \ref{lem;12.3.41},
we obtain the following morphisms for $h$ and $i\in I$:
\begin{equation}\label{eq;12.3.42}
 \widetilde{N}_{i,h}^{\sankaku}:
 \graded^{\sankaku}(\veca)_{h}\lrarr 
 \graded^{\sankaku}(\veca)_{h-1}\otimes\nbigo(2).
\end{equation}
\begin{lem} \label{lem;12.3.44}
Let consider the case  $\lambda=0\in\proj^1$ and $Q$ as above.
We have the following implication:
\[
 N^{\sankaku}_{i\,|\,0}
 \Bigl(
 W^{\sankaku}(\veca)_{h\,|\,0}
 \Bigr)
\subset
 W^{\sankaku}(\veca)_{h-2\,|\,0}.
\]
\end{lem}
\pf
Let describe $\theta$ as
\[
 \theta
=\sum_{j=1}^l f_j\cdot \frac{dz_j}{z_j}
+\sum_{j=l+1}^n g_j\cdot dz_j. 
\]
We know that $|f_i|_{h}\leq C\cdot(-\log|z_i|)^{-1}$.
Thus we have
\[
 |\varphi^{\ast}(f_i)|_{\varphi^{\ast}h}
 \leq C'\cdot(-\log|t|)^{-1}.
\]
Here $\varphi$ denotes the morphism given in
(\ref{eq;12.10.10}) for $\veca$.
Let consider the section $s$ of
$\prolong{\varphi^{\ast}E}$ over $\Delta$.
Let $\deg(s)$ be the degree of $s(O)$ with respect to
the weight filtration of $\Res(\varphi^{\ast}(\theta))$.
Then we have the norm estimate:
\[
 |s|_{\varphi^{\ast}(h)}\sim (-\log|t|)^{\deg(s)/2}.
\]
It implies the following:
\[
 |\varphi^{\ast}(f_i)\cdot s|_{\varphi^{\ast}(h)}
 \leq  C''\cdot (-\log|t|)^{(\deg(s)-2)/2}.
\]
By the norm estimate for the sections on the punctured disc,
we can conclude that the degree of $\varphi^{\ast}f_i\cdot s(O)$
is less than $\deg(s)-2$.
Thus we are done.
\hfill\qed

We have the following immediate corollary.
\begin{cor}
For $i\in I$,
let $\widetilde{N}_{i,h}^{\sankaku}$
be the morphism in {\rm (\ref{eq;12.3.42})}.
Then we have $\widetilde{N}^{\sankaku}_{i,h\,|\,0}=0$.
We also have $\widetilde{N}^{\sankaku}_{i,h\,|\,\infty}=0$.
\end{cor}
\pf
The first claim is obvious from Lemma \ref{lem;12.3.44}.
The second claim is obtained by applying
Lemma \ref{lem;12.3.44} to
the tame nilpotent harmonic bundle $(E^{\dagger},\theta^{\dagger},h)$.
\hfill\qed

\begin{cor} \label{cor;12.3.45}
Let $\widetilde{N}_{i,h}^{\sankaku}$ be
the morphism in {\rm (\ref{eq;12.3.42})} for $i\in I$.
Then $\widetilde{N}_{i,h}^{\sankaku}$ is, in fact, $0$.
\end{cor}
\pf
Since $(S(Q,P),\weightfilt^{\sankaku})$ is a mixed twistor,
the vector bundles $\graded_h(\veca)$ and
$\graded_{h-1}(\veca)\otimes\nbigo_{\proj^1}(2)$
are isomorphic to direct sums
of $\nbigo(h)$ and $\nbigo(h+1)$ respectively.
Thus $\widetilde{N}_{i,h}^{\sankaku}$ is a section of
the vector bundle isomorphic to
a direct sum of $\nbigo(1)$.

On the other hand, we know that $\widetilde{N}_{i,h}^{\sankaku}$
vanishes at two points $\{0,\infty\}$ due to
Corollary \ref{cor;12.3.45}.
Thus we can conclude that $\widetilde{N}_{i,h}^{\sankaku}$
vanishes over $\proj^1$.
\hfill\qed

As a direct corollary,
we obtain the following important theorem.
\begin{thm} \label{thm;12.3.50}
The morphism
$N^{\sankaku}_i:(S(Q,P),W^{\sankaku}(\veca))
 \lrarr (S(Q,P),W^{\sankaku}(\veca))$
$(i\in I)$
give morphisms of mixed twistor.
\end{thm}
\pf
The claim is equivalent to
$N^{\sankaku}_i\cdot W^{\sankaku}(\veca)_h
 \subset W^{\sankaku}(\veca)_{h-2}$.
It is proved in Corollary \ref{cor;12.3.45}.
\hfill\qed

\subsubsection{Some consequences}

\label{subsubsection;12.15.25}

The theorem \ref{thm;12.3.50} implies the following,
for example.
\begin{cor} \label{cor;12.3.51}
The conjugacy classes of $N^{\sankaku}_{i\,|\,\lambda}$
are independent of a choice of $\lambda\in\proj^1$ for each $i\in I$.
Moreover the conjugacy classes of $N^{\sankaku}(\veca)_{|\lambda}$ are
independent of a choice of $\lambda\in \proj$
for each $\veca\in\cnum^I$.
\hfill\qed
\end{cor}

The claim of the corollary \ref{cor;12.3.51} for $\lambda\neq 0,\infty$
is rather obvious.
However the fact that
the conjugacy classes of $N^{\sankaku}_{j}$ does not degenerate
at $\lambda=0,\infty$ 
is not trivial.

In each point $Q\in D_i$,
the nilpotent map $N_{i\,|\,(\lambda,Q)}$ induces
the weight filtration $W_{i\,|\,(\lambda,Q)}$.
\begin{cor} \label{cor;12.4.110}
The conjugacy classes of $N_{i\,|\,(\lambda,Q)}$
are independent of $(\lambda,Q)\in \nbigd_i$.
As a result,
the filtration $\{W_{i|\,(\lambda,Q)}\,|\,(\lambda,Q)\in\nbigd_i\}$
forms the filtration of $\prolong{\nbige}_{\nbigd_i}$
by vector subbundles.
\end{cor}
\pf
Let fix $\lambda\neq 0$.
Then it is easy to see that
the conjugacy classes of $N_{i\,|\,(\lambda,Q)}$
are independent of a choice of $Q\in D_i$.
To see it,
we only have to use a normalizing frame, for example.

Let fix $Q\in D_i$.
Then we know that the conjugacy classes of $N_{i\,|\,(\lambda,Q)}$
are independent of $\lambda$.
Thus we obtain our result.
\hfill\qed

For any subset $I\subset\lbar$, we put $\nbigd_I=\bigcap_{i\in I}\nbigd_i$.
On $\nbigd_I$, we have nilpotent maps,
$N(\veca)=\sum_{i\in I}a_i\cdot N_{i\,|\,\nbigd_I}$
for any $\veca\in \cnum^I$.
\begin{cor} \label{cor;12.11.31}
The weight filtrations $W(\veca)_{|(\lambda,Q)}$
of $N(\veca)_{|(\lambda,Q)}$ form the filtration
of $\prolong{\nbige}_{|\nbigd_I}$ by vector subbundles.
\end{cor}
\pf
Similar to Corollary \ref{cor;12.4.110}.
\hfill\qed

On $D_{\jbar}$, we have the residues $N_{i\,|\,\nbigd_{\jbar}}$
for $i\leq j$.
We put $N(\jbar)=\sum_{i\leq j}N_{i\,|\,\nbigd_{\jbar}}$.
We have the weight filtration $W(\jbar)$
of $N(\jbar)$,
which is a filtration of $\prolong{\nbige}_{\nbigd_{\jbar}}$.
In particular, we have the filtrations
$W(\jbar)$ on $\prolong{\nbige}_{\nbigd_{\mbar}}$ for any $j\leq m$.

\begin{lem}\label{lem;12.11.30}
Let $\vech=(h_1,\ldots,h_m)$ be an $m$-tuple of integers.
The intersections $\bigcap_{j=1}^m W(\jbar)_{h_j}$
form a vector subbundle of $\prolong{\nbige}_{\nbigd_{\mbar}}$.
\end{lem}
\pf
Let fix $\lambda\neq 0$.
Then it is easy to see that
the rank of $\bigcap_{j=1}^m W(\jbar)_{h_j\,|\,(\lambda,Q)}$
is independent of $Q\in D_{\mbar}$.
We only have to use a normalizing frame, for example.

Let fix $Q\in D_{\mbar}$.
Let $I$ denote the set $\{i\in\lbar\,|\,Q\in D_i \}$.
Let $\veca\in\rnum_{>0}^I$ be an element
such that $N(\veca)_{|(\lambda,Q)}$ is general for any $\lambda$.
We pick an appropriate point $P$
such that the vector bundle $S(Q,P)$ with
the filtration $W^{\sankaku}(\veca)$ 
is a mixed twistor.
Then we know that
$\bigcap_{j=1}^m W^{\sankaku}(\jbar)_{h_j}$
is a sub mixed twistor of $S(Q,P)$.
In particular, we obtain that
the rank of $\bigcap_{j=1}^m W(\jbar)_{h_j\,|\,(\lambda,Q)}$
is independent of $\lambda\in\cnum_{\lambda}$.
Thus we obtain our result.
\hfill\qed

Let $\graded^{\sankaku\,(1)}$
denote the associated graded vector bundle
to $W^{\sankaku}(\itibar)$.
We remark the following.
\begin{lem}
We have
$c_1(\graded^{\sankaku\,(1)}_h)=
 h\cdot \rank(\graded^{\sankaku\,(1)}_h)$.
Here $c_1(\nbigf)$ denotes the first Chern class
of a coherent sheaf $\nbigf$ on $\proj^1$.
\end{lem}
\pf
If $Q\in D_1$ is contained in $D_1-\bigcup_{i=2}^lD_1\cap D_i$,
then we only have to consider the restriction
of the harmonic bundle to a curve which 
transversally intersects with $D_1$ at $Q$.
In the general case,
we use the topological invariance of the Chern class.
\hfill\qed

\subsubsection{The graded part}

\label{subsubsection;12.10.40}

Let $Q$ be a point of $D_{\mbar}$.
We put $I=\{i\,|\,Q\in D_i\}$.
Let $\veca$ be an element of $\rnum_{>0}^I$
such that $N^{\sankaku}(\veca)$ is general.
For any $j\leq m$, we put
$N^{\sankaku}(\jbar):=\sum_{i\leq j}N^{\sankaku}_i$.
Let $W^{\sankaku}(\jbar)$ denote
the weight filtration of $N^{\sankaku}(\jbar)$.
In particular,
we obtain the associated graded vector bundle
$\graded^{\sankaku\,(1)}$ of $W^{\sankaku}(\itibar)$.
When $(S(W,P),W^{\sankaku}(\veca))$ is a mixed twistor,
we have the naturally induced mixed twistor structure
on $(\graded^{\sankaku\,(1)},W^{\sankaku\,(1)}(\veca))$.

We have the induced morphisms
$N^{\sankaku\,(1)}(\jbar):
 \graded^{\sankaku\,(1)}\lrarr
 \graded^{\sankaku\,(1)}\otimes\nbigo_{\proj^1}(2)$.
They are again the morphisms of mixed twistors.
In particular, the conjugacy classes
of $N^{\sankaku(1)}(\jbar)_{|\lambda}$
are not independent of $\lambda\in\proj^1$.

We will use the following special case later.
\begin{lem}
The morphism $N^{\sankaku(1)}(\nibar)$ induces the filtration
$W(N^{\sankaku(1)}(\nibar))$ 
on $\graded^{\sankaku(1)}=\bigoplus \graded^{\sankaku(1)}_{h}$
by vector subbundles.
Thus we obtain the graded vector bundle
$\graded^{W(N^{\sankaku(1)}(\nibar))}_{(h_1,h_2)}:=
Gr^{W(N^{\sankaku(1)}(\nibar))}_{h_2}
 \Bigl(
 \graded^{\sankaku(1)}_{h_1}\Bigr)$.
\hfill\qed
\end{lem}

\vspace{.1in}
On $\nbigd_{\itibar}$,
we obtain the graded vector bundle $\graded^{(1)}$ of $W(\itibar)$.
For $(\lambda,Q)\in\nbigd_{\mbar}$,
we have the induced filtrations
$W^{(1)}(\mbar)_{|(\lambda,Q)}$
on $\graded^{(1)}_{|(\lambda,Q)}$.
We also obtain the morphism
$N^{(1)}(\mbar)_{|(\lambda,Q)}
\in End(\graded^{(1)}_{|(\lambda,Q)})$.
\begin{lem} \mbox{{}}
\begin{itemize}
\item
The filtrations
$\{W^{(1)}(\mbar)_{|(\lambda,Q)}\,|\,(\lambda,Q)\in\nbigd_{\mbar}\}$
form a filtration of $\prolong{\nbige}_{|\nbigd_{\mbar}}$
by vector subbundles.
\item
The conjugacy classes of $N^{(1)}(\mbar)_{|(\lambda,Q)}$ are independent of
a choice of $(\lambda,Q)\in\nbigd_{\mbar}$.
\end{itemize}
\end{lem}
\pf
Let fix a $\lambda\neq 0$.
Then the claims can be checked by a normalizing frame.

Let fix $Q\in D_{\mbar}$.
Let $I$ denote the set $\{i\in\lbar\,|\,Q\in D_i \}$.
Let $\veca\in\rnum_{>0}^I$ be an element
such that $N(\veca)_{|(\lambda,Q)}$ is general for any $\lambda$.
We pick an appropriate point $P$
such that the vector bundle $S(Q,P)$ with
the filtration $W^{\sankaku}(\veca)$ 
is a mixed twistor.

Then we know that
$W^{\sankaku(1)}(\mbar)$ is a sub mixed twistor of
$(\graded^{\sankaku(1)},W^{\sankaku}(\veca))$.
Thus the rank of $W^{(1)}(\mbar)_h$ is independent of $\lambda$.
Since $N^{\sankaku(1)}(\mbar)$ is a morphism of mixed twistor,
the conjugacy classes are independent of a choice of $\lambda$.
Thus we obtain the result.
\hfill\qed

\begin{lem}
The rank and the first Chern class of 
$\graded^{W(N^{\sankaku(1)}(\nibar))}_{(h_1,h_2)}$
is independent of a choice of $Q$ and $P$.
\end{lem}
\pf
We have already seen the independence of the rank.
Since the dependence of the vector bundle on $Q$ and $P$
is continuous,
the Chern class is invariant.
\hfill\qed


\subsubsection{The weak constantness of the filtrations}

\label{subsubsection;12.15.26}

We continue to use the notation in the previous subsubsections.
By using the theorem \ref{thm;12.3.50},
we can show the following
weak constantness of the filtrations on the positive cones.
\begin{prop} \label{prop;12.13.10}
Let $\veca_1,\veca_2\in\cnum^I$ be general.
Then $W^{\sankaku}(\veca_1)=W^{\sankaku}(\veca_2)$.
\end{prop}
\pf
Let $\veca_1$ be general.
We have already known that
$N^{\sankaku}(\veca)$ drops the degree with respect to 
the filtration $W^{\sankaku}(\veca_1)$ by $2$,
for any $\veca$.
For non-negative $k\geq 0$
let consider the following morphism:
\[
 \bigl(N^{\sankaku}(\veca)\bigr)^k:
 \graded^{\sankaku}_k
 \lrarr
 \graded^{\sankaku}_{-k}\otimes\nbigo_{\proj^1}(2k).
\]
It is isomorphic when $\veca=\veca_1$.
Thus there is a Zariski open subset $U_k$ of $\cnum^n$
such that $\bigl(N^{\sankaku}(\veca)\bigr)^k$
is isomorphic for any $\veca\in O_k$.
Then we know that there is a Zariski open subset $U$,
such that $W^{\sankaku}(\veca)=W^{\sankaku}(\veca_1)$
for any $\veca\in U$.
Then we obtain  $W^{\sankaku}(\veca)=W^{\sankaku}(\veca_1)$
if $\veca$ is general.
\hfill\qed


\section{Limiting harmonic bundle in one direction}

\label{section;12.15.76}

\subsection{The method of Comparison}
\label{subsection;12.4.1}

We put $X=\Delta^n=\{(\zeta_1,\ldots,\zeta_n)\in\Delta^n\}$,
$D_i:=\{\zeta_i=0\}$, and $D=\bigcup_{i=1}^lD_i$
for $l\leq n$.
Let $\harmonicbundle$ be a tame nilpotent harmonic bundle
with trivial parabolic structure over $X-D$.
We have the deformed holomorphic bundle
$(\nbige,\DD,h)$ on $\nbigx-\nbigd$,
and the prolongment $\prolong{\nbige}$.
We have the residue $N_i:=\Res_{\nbigd_i}(\DD)$.

We have already known that the conjugacy classes
of $N_{1\,|\,(\lambda,Q)}$ are independent of a choice of
$(\lambda,Q)\in\nbigd_1$
(Corollary \ref{cor;12.4.110}).
We have the weight filtration $W(\itibar)$ of $N_1$.
For any $k\geq 0$ and $h\in\seisuu$,
we have the number $d(k,h):=\dim P_k Gr_{h}^{(1)}$
determined by the conjugacy class of $N_1$.

Let take a holomorphic frame $\vecv$ of $\prolong{\nbige}$
over $\nbigx$ satisfying the following:
\begin{itemize}
\item
 $\vecv=
 \Bigl(v_{k,h,\eta}\,\Big|\,
   k\geq 0,\,\,h\in\seisuu,\,\,
   \eta=1,\ldots,d(k,h)
 \Bigr)$.
\item
  $N_1(v_{k,h,\eta})=v_{k,h-2,\eta}$ if $h>-k$,
and $N_1(v_{k,-k,\eta})=0$.
\end{itemize}
Note that $N_1$ is represented by a constant matrix
with respect to the frame $\vecv_{|\nbigd_1}$.

We put $\blowup{X}:=\Delta^n=\{(z_1,\ldots,z_n)\in\Delta^n\}$,
$\blowup{D}_i=\{z_i=0\}$,
and $\blowup{D}=\bigcup_{i=1}^l\blowup{D}_i$.
We have the morphism $\blowup{X}\lrarr X$
defined as follows:
\begin{equation} \label{eq;12.10.20}
 \phi^{\ast}(\zeta_i)=
 \left\{
 \begin{array}{ll}
  {\displaystyle
 \prod_{j=i}^lz_j,} & (i\leq l)\\
 \mbox{{}}\\
 z_i, &(i>l).
 \end{array}
 \right.
\end{equation}
We obtain the tame nilpotent harmonic bundle
$\phi^{\ast}\harmonicbundle$,
and the deformed holomorphic bundle
$\phi^{\ast}(\nbige,\DD,h)$.
We put $\blowup{\nbige}:=\phi^{\ast}\nbige$.
We also have the prolongment
$\prolong{\blowup{\nbige}}=\phi^{\ast}\prolong{\nbige}$.

\label{subsection;12.14.101}

We have the projection
$\projection_1:\Omega^{1,0}_X\lrarr q_1^{\ast}\Omega^{1,0}_{\Delta}$,
and
$\projection_1:\Omega^{1,0}_X(\log D)\lrarr
 q_1^{\ast}\Omega^{1,0}_{\Delta}(\log O)$:
\[
 \projection_1\Bigl(\sum_{i=1}^n f_i\cdot dz_i\Bigr):=f_1\cdot dz_1.
\]
From the $\lambda$-connection
$\blowup{\DD}=\phi^{\ast}\DD:
\Gamma(\blowup{\nbigx},\prolong{\phi^{\ast}\nbige})
 \lrarr
 \Gamma(\blowup{\nbigx},
  \prolong{\phi^{\ast}\nbige}\otimes
       p_{\lambda}^{\ast}\Omega_{X}^{1,0}(\log D))$,
we obtain the family of $\lambda$-connections along the $z_1$-direction:
\[
 \projection_1(\blowup{\DD}):
 \Gamma(\blowup{\nbigx},\prolong{\blowup{\nbige}})
 \lrarr \Gamma(\blowup{\nbigx},\prolong{\blowup{\nbige}}\otimes
 p_{\lambda}^{\ast}q_1^{\ast}\Omega_{\Delta}(\log O) ).
\]
The residue $\Res_{\blowup{\nbigd}_1}(\projection_1(\blowup{\DD}))$
is same as $\phi^{\ast}N_1$.

We have the holomorphic frame $\blowup{\vecv}:=\phi^{\ast}\vecv$.
We have the $\lambda$-connection form
$A\in
 \Gamma\Bigl(\blowup{\nbigx},
   M(r)\otimes\nbigo_{\nbigx}\Bigr)$
of
$\blowup{\DD}$ with respect to $\blowup{\vecv}$,
that is,
$\blowup{\DD}\blowup{\vecv}=\blowup{\vecv}\cdot A\cdot dz_1/z_1$.

\begin{lem}
The restrictions $A_{|\blowup{\nbigd}_i}$ $(i=1,\ldots,l)$
are constant,
say $A_0$.
\end{lem}
\pf
Clear from our construction.
Note that $\phi(\blowup{\nbigd}_i)\subset \nbigd_1$.
\hfill\qed

Let $V$ be $\prolong{\nbige}_{|\,(0,O)}$ and 
$N$ be the residue $\Res(\DD)_{|\,(0,O)}$.
From the pair $(V,N)$,
we have a model bundle $E(V,N)=(E_0,\theta_0,h_0)$
on $\Delta^{\ast}$.
We denote the deformed holomorphic bundle
by $(\nbige_0,\DD_0,h_0)$.
We have the canonical frame $\vecv_0$
such that $\DD_0\vecv_0=\vecv_0\cdot A_0\cdot dz/z$.

Let $q_1$ denote the projection $\Delta^n\lrarr \Delta$
onto the first component.
We put $\blowup{\nbige}_0:=q_1^{\ast}\nbige_0$.
We have the $\lambda$-connection $\blowup{\DD}_0:=q_1^{\ast}\DD_0$
along the $z_1$-direction.
We also put $\blowup{\vecv}_0:=q_1^{\ast}\vecv_0$.

Due to the frames $\blowup{\vecv}$
and $\blowup{\vecv}_0$,
we obtain the holomorphic isomorphism
$\Phi:\prolong{\blowup{\nbige}}_0\lrarr \prolong{\blowup{\nbige}}$.

\begin{lem} \label{lem;12.3.60}
Note the following.
\begin{itemize}
\item
We have $\Phi\circ \blowup{\DD}_0-\blowup{\DD}\circ\Phi=0$
on $\blowup{\nbigd}_i$
for $i=2,\ldots,l$.
Similarly we have
$\Phi^{-1}\circ\blowup{\DD}-\blowup{\DD}_0\circ\Phi^{-1}=0$
on $\blowup{\nbigd}_i$
for $i=2,\ldots,l$.
\item
We have $\Res(\Phi\circ\blowup{\DD}_0-\blowup{\DD}\circ\Phi)=0$
on $\blowup{\nbigd}_1$.
Similarly we have
$\Res(\Phi^{-1}\circ\blowup{\DD}-\blowup{\DD}_0\circ\Phi^{-1})=0$
on $\blowup{\nbigd}_1$.
\end{itemize}
\end{lem}
\pf
Clear from our construction.
\hfill\qed

We have holomorphic bundles
$Hom(\prolong{\blowup{\nbige}_0},\prolong{\blowup{\nbige}})$
and
$Hom(\prolong{\blowup{\nbige}},\prolong{\blowup{\nbige}_0})$.
We have the naturally defined family of the $\lambda$-connections
along the $z_1$-direction.
induced by $\projection_1(\blowup{\DD})$ and $\blowup{\DD}_0$.
We denote them by $\DD_1$ and $\DD_2$.

The morphism $\Phi$ and $\Phi^{-1}$
can be regarded as the sections of
$Hom(\prolong{\blowup{\nbige}_0},\prolong{\blowup{\nbige}})$
and
$Hom(\prolong{\blowup{\nbige}},\prolong{\blowup{\nbige}_0})$
respectively.
Then Lemma \ref{lem;12.3.60} can be reworded as follows:
\begin{lem}\mbox{{}}
\begin{itemize}
\item
$\DD_1\Phi$ is holomorphic section of
$Hom(\prolong{\blowup{\nbige}_0},\prolong{\blowup{\nbige}})
 \otimes p_{\lambda}^{\ast}q_1^{\ast}\Omega_{\Delta}(\log O)$.
It vanishes on $\bigcup_{i=2}^l\blowup{\nbigd}_i$.
\item
$\DD_2\Phi$ is holomorphic section of
$Hom(\prolong{\blowup{\nbige}},\prolong{\blowup{\nbige}_0})
 \otimes p_{\lambda}^{\ast}q_1^{\ast}\Omega_{\Delta}(\log O)$.
It vanishes on $\bigcup_{i=2}^l\blowup{\nbigd}_i$.
\hfill\qed
\end{itemize}
\end{lem}

\vspace{.1in}
We explain our method to obtain some estimate of norms.
We have the metric $q_1^{\ast}(h_0)$ of $\blowup{\nbige}_0$.
Let $a$ and $b$ be functions as follows:
\begin{itemize}
\item $a$ is a positive function defined over
      $\Delta^{\ast\,n-1}$.
 For simplicity, we assume that
 $a(z_2,\ldots,z_n)$ is a polynomial
 of $-\log|z_2|,\ldots,-\log|z_n|$.
\item
 $b$ is a holomorphic function defined over $\Delta^{n-1}$,
 such that $|b(z_2,\ldots,z_n)|\leq 1$ for any
 $(z_2,\ldots,z_n)\in \Delta^{n-1}$.
\end{itemize}
If we are given such functions $a$ and $b$,
we put as follows:
\[
 \blowup{h}_0(\lambda,z_1,z_2,\ldots,z_n):=
 a(z_2,\ldots,z_n)\times
 h_0\bigl(\lambda,b(z_2,\ldots,z_n)\cdot z_1\bigr)
\]
From the metrics $\blowup{h}$ and $\blowup{h}_0$,
we obtain the metrics of
$Hom(\blowup{\nbige}_0,\blowup{\nbige})$
and $Hom(\blowup{\nbige},\blowup{\nbige}_0)$.
We denote them by $|\cdot|_{\blowup{h},\blowup{h}_0}$.

\begin{lem}\mbox{{}}
Let $C$ be a real number such that $0<C<1$.
Let $R$ be a positive number.
\begin{itemize} \label{lem;12.3.71}
\item
Assume that $\Phi$ is bounded with respect to 
the metrics $\tilde{h}$ and $\blowup{h}_0$ on the boundary:
\[
 \Delta_{\lambda}(R)
 \times 
 \bigl\{(z_1,\ldots,z_n)\in \blowup{X}-\blowup{D}
 \,|\,|z_1|=C\bigr\}.
\]
Then $\Phi$ is bounded over the following region:
\[
 \Delta_{\lambda}(R)\times
 \bigl\{(z_1,\ldots,z_n)\in\blowup{X}-\blowup{D}\,\big|\,
 |z_1|\leq C\bigr\}.
\]
\item
Assume that $\Phi^{-1}$ is bounded with respect to 
the metrics $\tilde{h}_0$ and $\blowup{h}_0$ on the boundary
$\Delta_{\lambda}(R)\times
 \bigl\{(z_1,\ldots,z_n)\in\blowup{X}-\blowup{D},\big|\, |z_1|=C\bigr\}$.
Then $\Phi^{-1}$ is bounded over
$\Delta_{\lambda}(R)\times
 \bigl\{(z_1,\ldots,z_n)\in\blowup{X}-\blowup{D}\,\big|\,
 |z_1|\leq C\bigr\}$.
\end{itemize}
\end{lem}
\pf
Since $\DD_1\Phi$ is holomorphic,
and since $\DD_1\Phi$ vanishes on $\bigcup_{i=1}^l\blowup{\nbigd}_i$,
we obtain the inequality
$|\DD_1\Phi|_{\blowup{h},\blowup{h}_0}<C_{\epsilon}\cdot |z_1|^{-\epsilon}$
over $\Delta_{\lambda}(R)\times (\blowup{X}-\blowup{D})$
for any $0<\epsilon<1$.
Then we can dominate the values
$|\Phi|_{\blowup{h},\blowup{h}_0}$
by the boundary values,
due to the same argument as
that in Proposition \ref{prop;12.3.70}.
Thus we obtain the result.
\hfill\qed

\vspace{.1in}

Let use the method.
Let $M$ be an integer.
We put as follows:
\[
 \blowup{h}_M(\lambda,z_1,\ldots,z_n)
:=\prod_{i=2}^l\bigl(-\log |z_i|\bigr)^M
 \cdot h_0(\lambda,z_1).
\]
\begin{lem}
If $M$ is sufficiently larger than $0$,
then $\Phi$ is bounded.
If $M$ is sufficiently smaller than $0$,
then $\Phi^{-1}$ is bounded.
\end{lem}
\pf
The claims are consequences of Lemma \ref{lem;12.3.71}
and Lemma \ref{lem;12.10.30}.
Note that Lemma \ref{lem;12.10.30} is stated in the case $\lambda=1$.
However the argument works for any $\lambda$.
\hfill\qed

\vspace{.1in}

We reword the lemma as follows:
Let $\vecv_1$ be a holomorphic frame of $\prolong{\nbige}$
over $\nbigx$, compatible with the filtration $\weightfilt(\itibar)$
on $\nbigd_1$.
Then we obtain the frame $\vecv_2=\phi^{\ast}\vecv_1$
of $\prolong{\blowup{\nbige}}$ over $\blowup{\nbigx}$.
It is compatible with the filtration on $\blowup{\nbigd}_1$.
We put $2\cdot k(v_{2,i}):=\deg^{\weightfilt(\itibar)}(v_{2,i})$.
We have the $C^{\infty}$-frame $\vecv'_2$ of
$\blowup{\nbige}$ over $\blowup{\nbigx}-\blowup{\nbigd}_1$,
given by $v'_{2\,i}:=(-\log|z_1|)^{-k(v_{2\,i})}\cdot v_{2\,i}$.

\begin{cor} \label{cor;12.4.2}
Let $\epsilon$ and $R$ be any positive numbers.
Let consider the following region:
\[
 \Delta_{\lambda}(R)\times
 \big\{(z_1,\ldots,z_n)\in\blowup{X}-\blowup{D}\,\big|\,\,
 \epsilon<|z_j|,\,\,(j=2,\ldots,l)
 \big\}.
\]
On the region, the frame $\vecv'_2$ is adapted.

We also have the following estimate over
$\Delta_{\lambda}(R)\times (\blowup{X}-\blowup{D})$:
\[
 C_1\prod_{i=2}^l(-\log|z_i|)^{-M}
\leq
 |v_{2\,i}|_{\blowup{h}}\cdot
 (-\log|z_1|)^{-k(v_{2\,i})}
\leq
 C_2\prod_{i=2}^l(-\log|z_i|)^{M}.
\]
Here $C_1$ and $C_2$ denote some positive constant,
and $M$ denotes a sufficiently large number.
\hfill\qed
\end{cor}

\subsection{Taking limit} \label{subsection;12.4.70}

\subsubsection{Replacement of Notation}

We essentially use the setting in the subsection
\ref{subsection;12.4.1}.
For simplicity of the notation,
we replace $\blowup{X}$ with $X$, and make the same replacement for others.
More precisely, we consider as follows:

We put $X=\Delta^n$, $D_i=\{z_i=0\}$ and $D=\bigcup_{i=1}^lD_i$
for $l\leq n$.
Let $\harmonicbundle$ be a tame nilpotent harmonic bundle
with trivial parabolic structure over $X-D$.
As usual,
$(\nbige,\DD)$ denote the deformed holomorphic bundle
with the $\lambda$-connection.

Let $\phi:\Delta^n\lrarr\Delta^n$ be the morphism
considered in the subsection \ref{subsection;12.4.1}.
It gives a morphism $\phi:X-D\lrarr X-D$.
\begin{assumption} \label{assumption;12.4.120}
Assume the following:
\begin{itemize}
\item We have a tame nilpotent harmonic bundle $(E_1,\theta_1,h_1)$
with trivial parabolic structure
over $X-D$,
and $(E,\theta,h)$ is $\phi^{\ast}(E_1,\theta_1,h_1)$.
\hfill\qed
\end{itemize}
\end{assumption}

Let $(\nbige_1,\DD_1)$ be the deformed holomorphic bundle 
with $\lambda$-connection of $(E_1,\theta_1,h_1)$.
On $\nbigd_1$, we have the weight filtration $\weightfilt(\itibar)$
of $\prolong{\nbige}_1$,
induced by the residue $\Res_{\nbigd_1}(\DD_1)$.
Let $\vecv_1$ be a holomorphic frame
of the prolongment $\prolong{\nbige}_1$
compatible with the filtration $\weightfilt(\itibar)$.

Then we have the holomorphic frame $\vecv=\phi^{\ast}\vecv_1$
of $\prolong{\nbige}$ over $\nbigx$.
It is compatible with the weight filtration
$\weightfilt(\itibar)$ induced by $\Res_{\nbigd_1}(\DD)$,
which is same as the pull back of the filtration above.

We put $2\cdot k(v_{i}):=\deg^{\weightfilt(\itibar)}(v_{i})$.
We have the $C^{\infty}$-frame $\vecv'$ of
$\nbige$ over $\nbigx-\nbigd_1$,
given by $v'_{i}:=(-\log|z_1|)^{-k(v_{i})}\cdot v_{i}$.

The following lemma is completely same as Corollary \ref{cor;12.4.2}.
\begin{lem} \label{lem;12.4.3}
Let $\epsilon$ and $R$ be any positive numbers.
Let consider the following region:
\[
 \big\{(\lambda,z_1,\ldots,z_n)\in \nbigx-\nbigd\,\big|\,\,
 |\lambda|<R,\,\,
 \epsilon<|z_i|,\,\,(i=2,\ldots,l)
 \big\}.
\]
On the region, the frame $\vecv'$ is adapted.

We also have the following estimate
over the region $\Delta_{\lambda}(R)\times (X-D)$:
\[
 C_1\prod_{i=2}^l(-\log|z_i|)^{-M}
\leq
 |v_{i}|_{\blowup{h}}\cdot
 (-\log|z_1|)^{-k(v_i)}
\leq
 C_2\prod_{i=2}^l(-\log|z_i|)^{M}.
\]
Here $C_1$ and $C_2$ denote some positive constant,
and $M$ denotes a sufficiently large number.
\hfill\qed
\end{lem}

By the frame $\vecv$, we decompose $\prolong{\nbige}$
as follows:
\[
 \prolong{\nbige}=
 \bigoplus_h U_h,
\quad\quad
 U_h:=
 \big\langle
  v_i\,\big|\,\deg^{\weightfilt(\itibar)}(v_i)=h
 \big\rangle.
\]

\subsubsection{Pull backs}

Let $m$ be a non-negative integers.
We have the morphism $\psi_{m,\itibar}:X\lrarr X$ 
or $X-D\lrarr X-D$, defined as follows:
\[
 \psi_{m,\itibar}(z_1,z_2,\ldots,z_n):=(z_1^m,z_2,\ldots,z_n).
\]
Then we obtain the harmonic bundles
$\psi_{m,\itibar}^{\ast}\harmonicbundle$ on $X-D$.
We also obtain the deformed holomorphic bundles
and the $\lambda$-connections
$(\psi_{m,\itibar}^{\ast}\prolong{\nbige},\psi_{m,\itibar}^{\ast}\DD)$.
We obtain the holomorphic frame $\vecv^{(m)}$
of $\psi_{m,\itibar}^{\ast}\prolong{\nbige}$,
defined as follows:
\[
 v^{(m)}_i:=
 \psi_{m,\itibar}^{-1}(v_i)\cdot m^{-k(v_i)}.
\]
We put $H^{(m)}:=H(\psi_{m,\itibar}^{\ast}(h),\vecv^{(m)})$,
which is an $\nbigh(r)$-valued function.

\begin{lem}\label{lem;12.5.10}
On any compact subset $K\subset X-D$,
the $\nbigh(r)$-valued functions
$\{H^{(m)}\}$ and $\{H^{(m)\,-1}\}$
are bounded independently of $m$.
\end{lem}
\pf
We put $H^{\prime}:=H(h,\vecv')$.
Let $K'$ be a compact subset of $\Delta^{\ast\,l-1}\times \Delta^{n-l}$.
Then $\Delta_1^{\ast}\times K'$ naturally gives a subset of $X-D$.
We have $\psi_{m,\itibar}(\Delta_1^{\ast}\times K)
\subset \Delta_1^{\ast}\times K$.
Due to Lemma \ref{lem;12.4.3},
$H^{\prime}$ and $H^{\prime\,-1}$ are bounded
over the region $\Delta_1^{\ast}\times K\subset X$.

We put as follows:
\[
 \vecv^{\prime\,(m)}=\psi_{m,\itibar}^{\ast}(\vecv'),
\quad
 H^{\prime\,(m)}:=H(\psi_{m,\itibar}^{\ast}(h),\vecv^{\prime\,m})
=\psi_{m,\itibar}^{\ast}H^{\prime}.
\]
Then $H^{\prime\,(m)}$ and $H^{\prime\,(m)\,-1}$
are bounded over the region $\Delta_1^{\ast}\times K\subset X$,
independently of $m$.

Let $L$ denote the diagonal matrix whose $(i,i)$-component
is $(-\log|z_1|)^{k(v_i)}$.
It is easy to check the following relation:
\[
 L\cdot H^{\prime\,(m)}\cdot L=H^{(m)}.
\]
Thus we obtain our result.
\hfill\qed

We have the $\lambda$-connection form
$\nbiga\in
    \Gamma(\nbigx,
    M(r)\otimes p_{\lambda}^{\ast}\Omega_{X}^{1,0}(\log X))$
of $\DD$ with respect to the frame $\vecv$.
We decompose $\nbiga$ as follows:
\[
 \nbiga=\sum_j \nbiga_j,
\quad\quad
 \nbiga_j
 \in \Gamma(\nbigx,
   M(r)\otimes
    p_{\lambda}^{\ast}q_i^{\ast}
       \Omega_{\Delta}(\log O)).
\]
We obtain $A_j\in M(r)\otimes\nbigo_{\nbigx}$ satisfying the following:
\[
 \nbiga_j=
 \left\{
 \begin{array}{ll}
{\displaystyle A_1\cdot \frac{dz_1}{z_1},} & (j=1),\\
\mbox{{}}\\
 A_j\cdot dz_j,     & (j\neq 1).
 \end{array}
 \right.
\]
Let $f_{A_j}$ denote the section of
$End(\prolong{\nbige})$ over $\nbigx$,
determined by the frame $\vecv$ and $A_j$,
that is
$f_{A_j}(\vecv)=\vecv\cdot A_j$.
The decomposition $\prolong{\nbige}=\bigoplus U_h$
induce the decomposition of $f_j$ as follows:
\[
 f_{A_j}=\sum_{h,k} f_{A_j\,(h,k)},
\quad\quad
 f_{A_j\,(h,k)}(U_k)\subset U_h.
\]
The section $f_{A_j\,(h,k)}$
induces the section of $M(r)\otimes\nbigo_{\nbigx}$
by the relation
$f_{A_j\,(h,k)}(\vecv)=\vecv\cdot A_{j\,(h,k)}$.
Thus we obtain the decomposition
$\nbiga_j=\sum_{h,k}A_{j\,(h,k)}\cdot dz_j$.

\label{subsubsection;12.5.10}

\begin{lem}
We have the following vanishing.
\begin{itemize}
\item
If $h>k-2$,
we have $A_{1\,(h,k)\,|\,\nbigd_1}=0$.
\item
If $h>k$,
we have $A_{j\,(h,k)\,|\,\nbigd_1}=0$
for $j=2,\ldots,n$.
\end{itemize}
\end{lem}
\pf
Since $\nbiga$ is the flat $\lambda$-connection form,
we obtain the following equality:
\[
 \lambda\cdot d\nbiga+\nbiga\wedge\nbiga=0.
\]
Note that $A_{1\,|\,\nbigd_1}$ is constant, by our construction,
i.e.,
$d(A_{1\,|\,\nbigd_1})=0$.
It implies that $[A_{1\,|\nbigd_1},A_{j\,|\,\nbigd_1}]=0$
for any $j$,
in other words,
$f_{A_1|\nbigd_1}$ and $f_{A_j|\nbigd_1}$ are commutative.
Hence the sections $f_{A_j|\nbigd_1}$ preserves the filtration
$\weightfilt(\itibar)$ on $\nbigd_1$.
Thus we obtain the second claim.
Moreover $f_{A_1|\nbigd_1}$ drops the degree by $2$.
Thus we obtain the first claim.
\hfill\qed

\label{subsection;12.4.130}
We have the $\lambda$-connection form $\nbiga^{(m)}$
of $\psi_{m,\itibar}^{\ast}\DD$ with respect to the frame $\vecv^{(m)}$.
We decompose $\nbiga^{(m)}$ into $\sum_j\nbiga^{(m)}_j$
as in the case of $\nbiga$.
\begin{lem} \label{lem;12.4.51}
We have the following equalities:
\[
 \nbiga^{(m)}_j=
 \left\{
 \begin{array}{ll}
 \sum_{h,k}\psi_{m,\itibar}^{\ast}A_{1\,(h,k)}\cdot m^{(h-k+2)/2}\cdot
 {\displaystyle\frac{dz_1}{z_1}
 }
 & (j=1)\\
\mbox{{}}\\
 \sum_{h,k}\psi_{m,\itibar}^{\ast}A_{j\,(h,k)}\cdot m^{(h-k)/2}\cdot dz_j
 &(j\neq 1).
 \end{array}
 \right.
\]
In particular, the sequences $\{\nbiga^{(m)}\}$ converges
to $\nbiga^{(\infty)}=\sum_j\nbiga_j^{(\infty)}$ given as follows:
\[
 \nbiga_j^{(\infty)}=
 \left\{
 \begin{array}{ll}
 \sum_h \pi_1^{\ast}A_{1\,(h,h+2)}\cdot 
 {\displaystyle  \frac{dz_1}{z_1}  } & (j=1)\\
 \mbox{{}}\\
 \sum_h \pi_1^{\ast}A_{j\,(h,h)}\cdot dz_j &(j\neq 1).
 \end{array}
 \right.
\]
Here $\pi_1$ denote the projection $\nbigx\lrarr \nbigd_1$,
omitting $z_1$.
\end{lem}
\pf
The $\lambda$-connection form of $\psi_{m,\itibar}^{\ast}\DD^{\ast}$
with respect to the frame $\psi_{m,\itibar}^{\ast}\vecv$
is represented by $\psi_{m,\itibar}^{\ast}\nbiga$.
Thus we obtain our result by a direct calculation.
\hfill\qed

\vspace{.1in}

We decompose $\theta$ into
$\sum_{i=1}^l f_i\cdot dz_i/z_i
 +\sum_{i=l+1}^n g_i\cdot dz_i$.
Then we have
$\psi_{m,\itibar}^{\ast}\theta=
 m\cdot\psi_{m,\itibar}^{\ast}(f_1)\cdot dz_1/z_1
+\sum_{i=2}^l \psi_{m,\itibar}^{\ast}(f_i)\cdot dz_i/z_i
+\sum_{i=l+1}^n \psi_{m,\itibar}^{\ast}(g_i)\cdot dz_i$.
\begin{lem}\label{lem;12.13.31}
We have the following inequality independent of $m$:
\[
\begin{array}{l}
 |m\cdot\psi_{m,\itibar}^{\ast}(f_1)|_{\psi_{m,\itibar}^{\ast}(h)}
 \leq C\cdot(-\log|z_1|)^{-1},\quad
 |\psi_{m,\itibar}^{\ast}(f_i)|_{\psi_{m,\itibar}^{\ast}(h)}\leq
 C\cdot (-\log|z_i|)^{-1},\,\,(2\leq i\leq l),\\
\mbox{{}}\\
 |\psi_{m,\itibar}^{\ast}(g_i)|_{\psi_{m,\itibar}^{\ast}(h)}\leq C,\,\,
(l+1\leq i\leq n).
\end{array}
\]
\end{lem}
\pf
We have the equality
$\psi_{m,\itibar}^{\ast}\bigl((-\log|z_1|)^{-1}\bigr)=
 m^{-1}(-\log|z_1|)^{-1}$.
Thus we obtain the result
due to Proposition \ref{prop;1.25.1}
\hfill\qed

In all we obtain the following.
\begin{lem}
The sequence $\{\vecv^{(m)}\}$
satisfies Condition {\rm \ref{condition;11.27.16}}.
\hfill\qed
\end{lem}

\subsubsection{Limit}

\label{subsubsection;12.14.2}

We can apply the result in subsection \ref{subsection;12.4.10}.
Let $F=\bigoplus_{i=1}^r\nbigo_{X-D}\cdot u_i$ be a holomorphic bundle
with the frame $\vecu=(u_i)$,
over $X-D$.
We put $\vece^{(m)}:=\vecv^{(m)}_{|\nbigx^0}$.
It is the frame of $\prolong{\psi_{m,\itibar}^{\ast}\nbige^{0}}$
over $\nbigx^0=\{0\}\times X$.
The frames $\vece^{(m)}$ and $\vecu$ give the holomorphic isomorphism
$\Phi_m:\nbige^0\lrarr F$
over $X-D$.

The morphism $\Phi_m$ induces the structure of harmonic bundle
on $F$.
Namely,
we obtain the metric $h^{(m)}$, the holomorphic Higgs fields
$\{\theta^{(m)}\}$ defined over $X-D$,
which are the image of $\psi_{m,\itibar}^{\ast}(h)$ and $\psi_{m,\itibar}^{\ast}\theta$
via the morphism $\Phi_m$.
The tuple $(F,\delbar_F,\theta^{(m)},h^{(m)})$
gives a harmonic bundle for each $m$.

We obtain the deformed holomorphic bundles with $\lambda$-connection
$(\nbigf^{(m)},\DD^{(m)})$.
The morphism $\Phi_m$ induces the holomorphic isomorphism
$\psi_{m,\itibar}^{\ast}\nbige\lrarr \nbigf^{(m)}$.
Thus we obtain the frame $\psi_{m,\itibar}(\vecv^{(m)})$
of $\nbigf^{(m)}$ over $\nbigx-\nbigd$.

\begin{lem}
We can pick a subsequence $\{m_i\}$ of $\{m\}$
satisfying the following:
\begin{itemize}
\item
We have the holomorphic Higgs field $\theta^{(\infty)}$
and the metric $h^{(\infty)}$.
The sequence $\{\theta^{(m_i)}\}$
and $\{h^{(m_i)}\}$ converges to $\theta^{(\infty)}$
and $h^{(\infty)}$ respectively,
on any compact subset $K\subset X-D$.
We denote the deformed holomorphic bundle of
$(F,\theta^{(\infty)},h^{(\infty)})$
by $(\nbigf^{(\infty)},\DD^{(\infty)})$
\item
We have the holomorphic frame of $\vecv^{(\infty)}$
of $\nbigf^{(\infty)}$ over $\nbigx-\nbigd$.
The sequence $\{\Phi_m(\vecv^{(m_i)})\}$
converges on any compact subset $K\subset\nbigx-\nbigd$.
\end{itemize}
\end{lem}
\pf
We only have to apply Proposition \ref{prop;11.27.15}
and Lemma \ref{lem;12.4.20}.
\hfill\qed

\begin{lem}
Let $N$ be a sufficiently large number.
We have the following estimate over $\nbigx-\nbigd$,
for the frame $\vecv^{(\infty)}=(v^{(\infty)}_i)$.
\[
 |v^{(\infty)}_i|_{h^{(\infty)}}
 <C\cdot(-\log|z_1|)^{k(v_i)}\prod_{j=2}^{l}(-\log|z_j|)^{N}
\]
We also have the following estimate over $\nbigx-\nbigd$,
for
$\Omega(\vecv^{(\infty)})=
 v_1^{(\infty)}\wedge \cdots\wedge v_r^{(\infty)}$:
\[
 0<C_1<|\Omega(\vecv^{(\infty)})|_{h^{(\infty)}}
\]
\end{lem}
\pf
By our construction of $\vecv^{(m)}$,
we obtain the estimate:
\[
 |v^{(m)}_i|_{\psi_{m,\itibar}^{\ast}(h)}
 <C\cdot (-\log|z_1|)^{k(v_i)}\prod_{j=2}^{l}(-\log|z_j|)^{N}
\]
Hence we obtain the inequality in the limit.

We have the inequality $0<C_1<|\Omega(\vecv)|_h$.
It induces the following:
\[
 0<C_1<|\Omega(\vecv^{(m)})|_{\psi_{m,\itibar}^{\ast}(h)}.
\]
Hence we obtain the inequality in the limit.
\hfill\qed

\begin{cor} \label{cor;12.4.50}
The frame $\vecv^{(\infty)}$ of $\nbigf^{(\infty)}$
over $\nbigx-\nbigd$ naturally induces the frame of
the prolongment
$\prolong{\nbigf^{(\infty)}}$
over $\nbigx$.

The $\lambda$-connection form of $\DD^{(\infty)}$
with respect to the frame $\vecv^{(\infty)}$
is given by $\nbiga^{(\infty)}$.
\hfill\qed
\end{cor}

\subsection{The decomposition in Limit}

\label{subsection;12.15.30}

\subsubsection{Construction}

Let $\directsummand_h$ denote
the vector subbundle of $\prolong{\nbigf^{(\infty)}}$
generated by
$\bigl\{v^{(\infty)}_i\,\big|\,\deg^{\nbigw(\itibar)}(v_i)=h\bigr\}$.

\begin{lem} \label{lem;12.7.10}
The vector bundle $\directsummand_h$ does not depend on a choice of
the original frame $\vecv$ of $\prolong{\nbige}$
compatible with the filtration $\weightfilt(\itibar)$.
\end{lem}
\pf
Let $\dotvecv$ be another frame of $\prolong{\nbige}$
over $\nbigx$ compatible with the filtration $\weightfilt(\itibar)$.
We have the holomorphic functions $b_{j\,i}$ satisfying
$v_i=\sum_j b_{j\,i}\dotv_j$.
Since $\vecv$ and $\dotvecv$ are compatible with
the filtration $\weightfilt(\itibar)$,
we have the following vanishing:
\[
 b_{j\,i\,|\,\nbigd_1}=0,
 \quad
 \mbox{if }
 \deg^{\weightfilt(\itibar)}(\dotv_j)
 >
 \deg^{\weightfilt(\itibar)}(v_i).
\]
We obtain the relation of $\vecv^{(m)}$ and $\dotvecv^{(m)}$
for each $m$:
\[
 v_i^{(m)}=
 \sum_j \psi_m^{\ast}(b_{j\,i})\cdot m^{k(\dotv_j)-k(v_i)}
 \cdot \dotv_j^{(m)}.
\]
Here we have $2k(\dotv_j)=\deg^{\weightfilt(\itibar)}(\dotv_j)$.
Then we obtain the relation of $\vecv^{(\infty)}$
and $\dotvecv^{(\infty)}$ as follows:
\[
 v_i^{(\infty)}=
 \sum_{k(\dotv_j)=k(v_i)}
 \pi_1^{\ast}(b_{j\,i})\cdot \dotv_j^{(\infty)}.
\]
Here $\pi_1$ denotes the projection $\nbigx\lrarr\nbigd_1$,
omitting the first component.
Thus $\directsummand_h$ does not depend on a choice of the original frame.
\hfill\qed

Thus we obtain the decomposition
$\nbigf^{(\infty)}=\bigoplus_h \directsummand_h$.

We have the $\lambda$-connection form
$\nbiga^{(\infty)}=\sum_{j=1}^lA^{(\infty)}_1\cdot dz_j/z_j+
 \sum_{j=l+1}^n A^{(\infty)}_j\cdot dz_j$
of $\DD^{(\infty)}$
with respect to the frame $\vecv^{(\infty)}$.
Let $f_{A^{(\infty)}_j}$ denote
the sections of $End(\nbigf^{(\infty)})$
over $\nbigx$ determined by $A^{(\infty)}_j$
and the frame $\vecv^{(\infty)}$.

\begin{lem} \mbox{{}} \label{lem;12.7.11}
\begin{itemize}
\item
We have $f_{A^{(\infty)}_1}\Bigl(
 \directsummand_h\Bigr)\subset \directsummand_{h-2}$.
\item
When $j\neq 1$, the morphisms
$f_{A^{(\infty)}_j}$ preserve the decomposition
$\nbigf^{(\infty)}=\bigoplus_h \directsummand_h$.
\end{itemize}
\end{lem}
\pf
The claims immediately follow from
Lemma \ref{lem;12.4.51} and Corollary \ref{cor;12.4.50}.
\hfill\qed

\subsubsection{Orthogonality}

We put
$\directsummand^{\lambda}_h:=\directsummand_{h\,|\nbigx^{\lambda}}$.
\begin{thm} \label{thm;12.4.100}
If $h\neq h'$,
then $\directsummand^{0}_{h}$ and $\directsummand^0_{h'}$ are
orthogonal.
\end{thm}
\pf
We only have to consider the case $X=\Delta$ and $D=\{O\}$.
We use the notation $\psi_m$ instead of $\psi_{m,\itibar}$
for simplicity.
We put $V=\prolong{\nbige^0}_{|(0,O)}$ and
$N=\Res(\DD)_{|(0,O)}$.
We have a model bundle $E(V,N)=(E_0,\theta_0,h_0)$.
Let $(\nbige_0,\DD_0)$ denote the deformed holomorphic bundle
with the $\lambda$-connection.
We have the canonical frame $\vecv$
of the prolongment $\prolong{\nbige_0}$,
such that $\DD_0\cdot\vecv=\vecv\cdot N\cdot dz/z$.

We put $\vece_0:=\vecv_{0\,|\,\nbigx^0}$.
Due to the frames $\vece_0$ and $\vece$,
we obtain the holomorphic isomorphism
$\nbigi:\prolong{\nbige^0_0}\lrarr \prolong{\nbige^0}$
such that
$\Res(\nbigi\circ\theta_0-\theta\circ\nbigi)=0$.
We obtain
$\psi_m^{\ast}\nbigi:
 \psi_m^{\ast}\nbige^0_0\lrarr \psi_m^{\ast}\nbige^0$.

We take a limit of $(\nbige_0,\DD_0,h_0)$
by using the frame $\vecv_0$,
as in the subsection \ref{subsection;12.4.70}.
Namely we have the frames $\vecv_0^{(m)}$
of $\psi_m^{\ast}\nbige_0$ defined as follows:
\[
 v_{0\,i}^{(m)}:=
 \psi_m^{-1}(v_{0\,i})\cdot m^{-k(v_{0\,i})}.
\]
We have the holomorphic isomorphism
$\Phi_{0\,m}:\nbige^{0}_0\lrarr F$
given by the frames $\vece_0^{(m)}$ and $\vecu$.
Note that we have $\Phi_{0\,m}=\Phi_m\circ\psi_m^{\ast}\nbigi$.

Then $\Phi_{0\,m}$ induces the holomorphic Higgs field
$\theta_0^{(m)}$ and the metric $h_0^{(m)}$
on $F$.
The tuple $(F,\delbar_F,\theta_0^{(m)},h_0^{(m)})$
gives a harmonic bundle.
We denote the deformed holomorphic bundles
by $(\nbigf_0^{(m)},\DD_0^{(m)})$.
We also obtain the holomorphic frames
$\Phi_{0\,m}(\vecv_0^{(m)})$ of $\nbigf_0^{(m)}$.

For the subsequence $\{m_i\}$,
we have the limits
$\theta_0^{(\infty)}$,
$h_0^{(\infty)}$
and $\vecv_0^{(\infty)}$ of
the sequences $\{\theta_0^{(m_i)}\}$,
$\{h_0^{(m_i)}\}$ and $\{\Phi_{0\,m_i}(\vecv_0^{(m_i)})\}$
respectively.
By our construction,
we have the equality $\theta_0\superinfty=\theta\superinfty$.

From the two harmonic bundles
$(F,\theta^{(\infty)},h^{(\infty)})$
and $(F,\theta^{(\infty)}_0,h^{(\infty)}_0)$,
we obtain the two deformed holomorphic bundles
$\nbigf\superinfty$ and $\nbigf_0\superinfty$.
Since the underlying $C^{\infty}$-vector bundles of them
is same as $p_{\lambda}^{-1}(F)$,
we have the natural $C^{\infty}$-isomorphism
$\nbigi\superinfty:\nbigf_0\superinfty\lrarr\nbigf\superinfty$.

\begin{prop} \label{prop;12.4.80}\mbox{{}}
\begin{itemize}
\item
The morphism $\nbigi\superinfty$ is holomorphic.
\item
$\nbigi\superinfty$ naturally induces the isomorphism
of the prolongments
$\prolong{\nbige_0}\lrarr\prolong{\nbige}$.
\item
$\nbigi\superinfty$ preserves the weight filtrations
of the residues $\Res(\DD^0)$ and $\Res(\DD)$.
\end{itemize}
\end{prop}
\pf
The holomorphic map $\nbigi:\nbige_0^0\lrarr \nbige^0$
induces the $C^{\infty}$-isomorphism
$\nbigi:\nbige_0\lrarr\nbige$ defined over $\nbigx-\nbigd$.
We have the elements $I_{i\,j}(\lambda,z)\in C^{\infty}(\nbigx-\nbigd)$
determined as follows:
\[
 \nbigi(v_{0\,j})=
 \sum_i I_{i\,j}\cdot v_i,
\quad
\mbox{or equivalently,}
\quad
 \Phi_{0\,0}(v^{(0)}_{0\,j})=
 \sum_i I_{i\,j}\cdot \Phi_0(v^{(0)}_{i}).
\]
Similarly we obtain the functions
$I^{(m)}_{i\,j}\in C^{\infty}(\nbigx-\nbigd)$
determined as follows:
\[
  \Phi_{0\,m}(v^{(m)}_{0\,j})=
 \sum_i I^{(m)}_{i\,j}\cdot \Phi_m(v^{(m)}_{i}).
\]
Since $\{\Phi_{0\,m_i}(\vecv_0^{(m_i)})\}$
and $\{\Phi_{m_i}(\vecv^{(m_i)})\}$ converge
to $\vecv_0\superinfty$ and $\vecv\superinfty$
in $L_k^p$ for any large $k$ and for a sufficiently large $p$
over any compact subset $K\subset \nbigx-\nbigd$,
the functions $\{I_{i\,j}^{(m_i)}\}$ converges similarly.
We denote the limit by $I_{i\,j}\superinfty$.
We obtain the following relation:
\[
 v\superinfty_{0\,j}=\sum_i I\superinfty_{i\,j}\cdot v_i\superinfty.
\]
We will show that
$I\superinfty_{i\,j}$ are holomorphic,
which implies that $\nbigi\superinfty$ is holomorphic.
In fact,
$I_{i\,j}^{(m)}$ are holomorphic along the direction
of $\lambda$.
Thus we will check that $I\superinfty_{i\,j}$
are holomorphic along the direction of $z$.

We use the following lemma.
\begin{lem}
Note that we have the equality
$||\,\psi_m^{-1}(f)\,||_{Z,C}=||\,f\,||_{Z,C^n}$
for any $C^{\infty}$-function $f$ on $\Delta^{\ast}$
such that $||f||_{Z,C}<\infty$.
\end{lem}
\pf
By using the real coordinate $z=r\cdot\exp(\sqrt{-1}\alpha)$,
we have the following:
\[
  ||\,\psi_m^{-1}(f)\,||_{Z.C}:=
 \int_{\Delta^{\ast}(C)}
 |\psi_m^{-1}(f)|\frac{dr\cdot d\alpha}{r\cdot (-\log r)}
=\frac{1}{m}
 \int_{\Delta^{\ast}(C)}
 \psi_m^{-1}\left(
 |f|\frac{dr\cdot d\alpha}{r\cdot(-\log r)}
 \right)
\]
Since the degree of the map
$\psi_m:\Delta(C)\lrarr \Delta(C^m)$ is $m$,
the right hand side is same as the following:
\[
 \int_{\Delta^{\ast}(C^m)}
 |f|\frac{dr\cdot d\alpha}{r\cdot(-\log r)}
=:
 ||\,f\,||_{Z,C^m}.
\]
Thus we are done.
\hfill\qed

Let us return to the proof of Proposition \ref{prop;12.4.80}.
Note the following equality:
\[
 I_{i\,j}^{(m)}
=
 \psi_m^{-1}\bigl(
 I_{i\,j}
 \bigr)\cdot
 m^{k(v_i)-k(v_{0\,j})}.
\]
Thus we obtain the following equality:
\[
 \zbar\cdot \delbar_z I_{i\,j}^{(m)}
 (-\log|z|)^{k(v_i)-k(v_{0\,j})+1}
=\psi_m^{-1}
 \Bigl(
 \zbar\cdot \delbar_z I_{i\,j}
 (-\log|z|)^{k(v_i)-k(v_{0\,j})+1}
 \Bigr)
\]
Due to the result of Simpson,
we know the finiteness for any $C<1$:
\[
 \big|\big|
 \zbar\cdot \delbar_z I_{i\,j}
 (-\log|z|)^{k(v_i)-k(v_{0\,j})+1}
 \big|\big|_{Z,C}<\infty.
\]
Thus we obtain the following convergence:
\[
 \lim_{m\to\infty}
  \big|\big|
 \zbar\cdot \delbar_z I_{i\,j}^{(m)}
 (-\log|z|)^{k(v_i)-k(v_{0\,j})+1}
 \big|\big|_{Z,C}
=\lim_{m\to\infty}\big|\big|
 \zbar\cdot \delbar_z I_{i\,j}
 (-\log|z|)^{k(v_i)-k(v_{0\,j})+1}
 \big|\big|_{Z,C^m}
=0.
\]
Thus we obtain the following vanishing:
\[
 \big|\big|
 \zbar\cdot \delbar_z I_{i\,j}^{(\infty)}
 (-\log|z|)^{k(v_i)-k(v_{0\,j})+1}
 \big|\big|_{Z,C}=0
\]
It implies the vanishing $\delbar_zI_{i\,j}\superinfty=0$.
Thus $\nbigi\superinfty$ is holomorphic.

We know that 
$I_{i\,j}\cdot (-\log|z|)^{k(v_i)-k(v_{0\,j})}$ is bounded
over $K\times\Delta^{\ast}$
for any compact subset $K\subset \cnum_{\lambda}$.
Thus $I_{i\,j}\superinfty$ is dominated by
a polynomial of $(-\log |z|)$ on such regions.
It implies that $\nbigi\superinfty$ naturally induces
the morphism of the prolongments.

When $k(v_i)-k(v_{0\,j})\neq 0$,
we have the finiteness:
\[
 || I_{i\,j}\superinfty\cdot (-\log|z|)^{k(v_i)-k(v_{0\,j})}||_{W,C}
<\infty.
\]
It implies that $I\superinfty_{i\,j\,|\,\nbigd}=0$
if $k(v_i)>k(v_{0\,j})>0$.
Namely $\nbigi\superinfty$ preserves the weight filtration
on $\nbigd$.

Thus the proof of Proposition \ref{prop;12.4.80} is completed.
\hfill\qed

We have the conjugate $\theta^{(\infty)\,\dagger}$
of $\theta\superinfty$ with respect to the metric $h\superinfty$.
We also have the conjugate $\theta_0^{(\infty)\,\dagger}$
of $\theta_0\superinfty$ with respect to the metric $h_0\superinfty$.
\begin{cor}
We have $\theta^{(\infty)\dagger}=\theta_0^{(\infty)\dagger}$.
\end{cor}
\pf
The holomorphic structures of $\nbigf^{(\infty)}$
and $\nbigf^{(\infty)}_0$
are given by the following:
\[
 \delbar_{\lambda}+\delbar_F+\lambda\cdot\theta^{(\infty)\dagger},
\quad
\delbar_{\lambda}+\delbar_F+\lambda\cdot\theta_0^{(\infty)\dagger}.
\]
Since the $C^{\infty}$-isomorphism
$\nbigf^{(\infty)}\stackrel{=}{\lrarr}
 p_{\lambda}^{-1}(F)
\stackrel{=}{\lrarr}
 \nbigf^{(\infty)}_0$
is holomorphic,
we obtain the equality
$\theta^{(\infty)\,\dagger}=\theta_0^{(\infty)\dagger}$.
\hfill\qed

\begin{cor}
We have the equality:
\[
 \theta\superinfty\cdot\theta^{(\infty)\,\dagger}
+\theta^{(\infty)\,\dagger}\cdot\theta\superinfty
=\theta_0\superinfty\cdot\theta_0^{(\infty)\,\dagger}
+\theta_0^{(\infty)\,\dagger}\cdot\theta_0\superinfty.
\]
We define the section $C\in C^{\infty}(\Delta^{\ast},End(F))$ by
$C\cdot dz\cdot d\zbar:=
\theta\superinfty\cdot\theta^{(\infty)\,\dagger}
+\theta^{(\infty)\,\dagger}\cdot\theta\superinfty$.
\hfill\qed
\end{cor}

By our construction,
we have $\directsummand_{0,h}^0=\directsummand_h^0$ 
on $\nbigx^{0}=\{0\}\times X$.

\begin{lem}
For any $P\in\Delta^{\ast}$,
the subspace $\directsummand_{0,h\,|\,(0,P)}^0$
is the eigenspace of $C_{|(0,P)}$
with the eigenvalue $h\cdot \phi(P)$.
Here $\phi(P)$ denotes the function as follows:
\[
 -\phi(P)=|z(P)|^{-2}\cdot (-\log|z(P)|)^{-2}\neq 0.
\]
\end{lem}
\pf
Since we have
$C\cdot dz\cdot d\zbar=
\theta_0\superinfty\cdot\theta_0^{(\infty)\,\dagger}
+\theta_0^{(\infty)\,\dagger}\cdot\theta_0\superinfty$,
we only have to check the equality in the case of model bundles
$Mod(l+1,1,1)$.

On $Mod(l+1,1,1)$, we have the canonical frame
$\bigl(e_1^p\cdot e_{-1}^q\,|\,p+q=l\bigr)$.
Here $(e_1,e_{-1})$ is the orthogonal frame for $Mod(2,1,1)$,
introduced in the subsubsection \ref{subsubsection;12.4.90}.
We have the following:
\[
 \theta_0(e_1^p\cdot e_{-1}^q)=p\cdot e_1^{p-1}\cdot e_{-1}^{q+1}
\cdot \frac{dz}{z},
\quad\quad
 \theta_0^{\dagger}(e_1^{p}\cdot e_{-1}^q)=
 q\cdot e_1^{p+1}\cdot e_{-1}^{q-1}
\cdot \frac{d\zbar}{\zbar\cdot (-\log|z|)^2}.
\]
Thus we obtain the following:
\[
 C(e_1^p\cdot e_{-1}^q)=-(p-q)\cdot e_1^p\cdot e_{-1}^q\cdot
 \phi(P).
\]
Hence we are done.
\hfill\qed

\begin{lem}
The endomorphism $C_{|P}$ of $E_{|P}$
is anti-self adjoint with respect to
the metric $h_{|P}$.
\end{lem}
\pf
We put $\theta=f\cdot dz$.
Let $f^{\dagger}$
denote the adjoint of $f$ with respect to the metric $h$.
Then $C_{|P}$ is $f\circ f^{\dagger}-f^{\dagger}\circ f$.
Thus it is anti-self adjoint.
\hfill\qed

\vspace{.1in}

In general, the eigenspaces of the anti-self adjoint operator
corresponding to the eigenvalues $\alpha_1$ and $\alpha_2$
are orthogonal, if $\alpha_1\neq \alpha_2$.
Thus we can conclude that
$\directsummand^0_{h|\,P}$ and $\directsummand^0_{h'|P}$ are orthogonal
if $h$ and $h'$ are different.
Hence the proof of Theorem \ref{thm;12.4.100} is completed.
\hfill\qed

\subsubsection{Limiting CVHS in one dimensional case}

We have the following immediate corollary.
\begin{cor} \label{cor;12.7.21}
Let $\harmonicbundle$ be a tame nilpotent harmonic bundle
with trivial parabolic structure over $\Delta^{\ast}$.
Let $(F,\theta\superinfty,h\superinfty)$ be a limiting harmonic bundle,
via the pull backs $\psi_{m,\itibar}$.
Then $(F,\theta^{(\infty)},h^{(\infty)})$
gives a complex variation of polarized Hodge structure,
up to grading.
\end{cor}
\pf
We have the decomposition $F=\bigoplus_{h}\directsummand_h^0$.
We have already known that
$\directsummand_h$ and $\directsummand_{h'}$ are orthogonal
if $h\neq h'$,
and $\theta\superinfty\bigl(\directsummand_{h}\bigr)\subset
 \directsummand_{h-2}\otimes\Omega_{\Delta}(\log O)$.

Let $a$ be an element of $S^1=\{z\in\cnum\,|\,|z|=1\}$.
We have the morphism $\rho_a:F\lrarr F$ given by
$\rho_a:=\bigoplus_{h}a^{-h}\cdot id_{\directsummand_h^0}$.
It gives the isomorphism of the Higgs bundle
$(F,\theta\superinfty)$ and $(F,a\cdot\theta\superinfty)$.
It also gives the isomorphism of the hermitian metrics.
Hence we obtain the $S^{1}$-action on
the harmonic bundle on the harmonic bundle
$(F,\theta\superinfty,h\superinfty)$.
Thus it gives a complex variation of the polarized Hodge structures.
(See \cite{s1} and \cite{s2}. See also Appendix.)
\hfill\qed


\vspace{.1in}

Let use the real coordinate $z=r\cdot \exp(\sqrt{-1}\alpha)$.
\begin{prop} \label{prop;12.5.60}
The $\nbigh(r)$-valued function
$H(h\superinfty,\vecv\superinfty)$
is independent of $\alpha$.
\end{prop}
\pf
Let consider $H^{(m)}=H(\psi_m^{\ast}(h),\vecv^{(m)})$,
and the restriction of $H^{(m)}$
to $S_a^1:=\{\,z\in\cnum\,|\,|z|=a\,\}$ for $0<a<1$.
The components $H_{i\,j}^{(m)}$
is contained in the image of the linear morphism:
$F_m: L^2(S^1_{a^n})\lrarr L^2(S^1_a)$.
The image $J_m$ of $F_m$ is generated by the following:
\[
 \Bigl\{
 \exp\bigl(\sqrt{-1}h\cdot m\cdot\alpha\bigr)
 \,\Big|\, h\in\seisuu
 \Bigr\}.
\]
Since the intersection $\bigcap_{m_i}J_{m_i}$ is $\{0\}$,
we obtain the result.
\hfill\qed

\subsection{The Chern class of the vector bundle
$\graded^{W(N^{\sankaku(1)}(\nibar)) }_{h_1,h_2}$}

\label{subsection;12.15.31}

We use the notation in the subsubsection \ref{subsubsection;12.10.40}.
Let consider the case $Q\in D_{\nibar}$.
We have the vector bundle
$\graded^{W(N^{\sankaku(1)}(\nibar)) }_{h_1,h_2}$ over $\proj^1$.
\begin{lem} \label{lem;1.31.12}
Let $b$ denote the bottom number of $W(\itibar)$.
We have the following:
\begin{equation} \label{eq;12.4.155}
 c_1(\graded^{W(N^{\sankaku(1)}(\nibar))}_{b,h_2})=
 (h_1+h_2)\cdot
 \rank \graded^{W(N^{\sankaku(1)}(\nibar))}_{b,h_2}
\end{equation}
\end{lem}
\pf
We only have to consider the harmonic bundles
as in Assumption \ref{assumption;12.4.120},
because the property on $\graded^{\sankaku(1)}$ of $S(Q,P)$
is not changed if we take a pull back via $\phi$
given in (\ref{eq;12.10.20}).
Thus we assume Assumption \ref{assumption;12.4.120}
in the rest of the proof.

Since the first Chern class is topological invariant,
it does not depend on a choice of points $P\in X-D$
and $Q\in D_{\nibar}$.
We will use it in the following without mention.
One immediate consequence is that
we only have to consider the case
$Q\in D_{\nibar}-
 \bigcup_{i\neq 1,2}D_{\nibar}\cap D_i$.

Let $(E,\theta,h)$ be a harmonic bundle
satisfying Assumption {\rm \ref{assumption;12.4.120}}.
We denote the vector bundle
$\graded^{W(N^{\sankaku(1)}(\nibar))}_{h_1,h_2}$ for $(E,\theta,h)$
by
 $\graded^{W(N^{\sankaku(1)}(\nibar))}_{h_1,h_2}(E,\theta,h)$
to distinguish the dependence on
$(E,\theta,h)$.

Let take a limit $(F,\theta\superinfty,h\superinfty)$
of $(E,\theta,h)$ as in the subsection {\rm\ref{subsection;12.4.70}}.
We obtain the vector bundle
$\graded^{W(N^{\sankaku(1)})(\nibar)}_{h_1,h_2}
     (F,\theta\superinfty,h\superinfty)$.

We need the following lemma.
\begin{lem} \label{lem;2.1.1}
We have the equality of the Chern classes:
\begin{equation} \label{eq;12.4.152}
 c_1\Bigl(
 \graded^{W(N^{\sankaku(1)})(\nibar)}_{h_1,h_2}
     (F,\theta\superinfty,h\superinfty)
\Bigr)=
 c_1\Bigl(
  \graded^{W(N^{\sankaku(1)})(\nibar)}_{h_1,h_2}
     (E,\theta,h)
\Bigr).
\end{equation}
\end{lem} 
\pf
We denote the vector bundle $S(Q,P)$ for $(E,\theta,h)$ by
$S(Q,P,(E,\theta,h))$.
Let $\vecw$ be a normalizing frame of
the deformed holomorphic bundle $(\nbige,\DD)$
of $(E,\theta,h)$
over $\nbigx^{\shikaku}-\nbigd^{\shikaku}$.
Assume that $\vecw$ is compatible with
the filtration $\weightfilt(\itibar)$.
We put, as usual,
$w^{(m)}_i:=\psi_m^{\ast}(w_i)\cdot m^{-k(w_i)}$.
Here we put $2\cdot k(w_i)=\deg^{\weightfilt(\itibar)}(w_i)$.
Since it gives the normalizing frame
of $(\psi_m^{\ast}\nbige,\DD^{(m)})$,
we obtain the following natural isomorphism:
\[
 S(Q,P,\psi_m^{\ast}(E,\theta,h))
\simeq
 S(Q,\psi_m(P),(E,\theta,h)).
\]
The isomorphism preserves the nilpotent maps
$N_1^{\sankaku}$ and $N_2^{\sankaku}$.
Thus we obtain the following:
\[
 c_1\Bigl(
 \graded^{W(N^{\sankaku(1)})(\nibar)}_{h_1,h_2}
 \Bigl(
     \psi_m^{\ast}(E,\theta,h)
 \Bigr)
\Bigr)=
 c_1\Bigl(
  \graded^{W(N^{\sankaku(1)})(\nibar)}_{h_1,h_2}
     (E,\theta,h)
\Bigr).
\]

Let pick the subsequence $\{m_i\}$ of $\{m\}$
for the limit $(F,\theta\superinfty,h\superinfty)$.
We can assume that 
the sequence of the frames $\bigl\{\Phi_{m_i}(\vecw^{(m_i)})\bigr\}$
converges to $\vecw^{(\infty)}$.
Then the sequence of the gluings of the bundle
$S(Q,P,\psi_m^{\ast}(E,\theta,h))$
converging to the gluings of
$S(Q,P,(F,\theta\superinfty,h\superinfty))$.
In Lemma \ref{lem;12.4.51} and Corollary \ref{cor;12.4.50},
we have seen the following:
\begin{itemize}
\item
The conjugacy class of $N(\itibar)$ is not changed
in the limit.
\item
The conjugacy class of $N^{(1)}(\nibar)$ is not changed
in the limit.
\end{itemize}
Then the gluings of the vector bundle
$\graded^{W(N^{\sankaku(1)})(\nibar)}_{h_1,h_2}
 \Bigl(
 \psi_{m_i}^{\ast}(E,\theta,h)\Bigr)$
converges to the gluing of the vector bundle
$\graded^{W(N^{\sankaku(1)})(\nibar)}_{h_1,h_2}
     (F,\theta\superinfty,h\superinfty)$.
Since the Chern class is topological invariant,
we obtain the following for sufficiently large $m_i$:
\[
 c_1\Bigl(
 \graded^{W(N^{\sankaku(1)})(\nibar)}_{h_1,h_2}
 \Bigl(
     \psi_{m_i}^{\ast}(E,\theta,h)
 \Bigr)
 \Bigr)
=c_1\Bigl(
 \graded^{W(N^{\sankaku(1)})(\nibar)}_{h_1,h_2}
     (F,\theta\superinfty,h\superinfty)
\Bigr).
\]
Thus the proof of Lemma \ref{lem;2.1.1} is completed.
\hfill\qed

Let us return to the proof of Lemma \ref{lem;1.31.12}.
We can assume that the harmonic bundle
considered is a limiting harmonic bundle in one direction.
Let $Q$ be a point contained in $D_{\nibar}$.
Let consider the tuple of vector space
$\nbigv:=\prolong{\nbige}_{|(\lambda,Q)}$
and the nilpotent maps
$\nbign_1:=N_{1\,|\,(\lambda,Q)}$
and $\nbign_2:=N_{2\,|\,(\lambda,Q)}$.
The tuple $(\nbigv,\nbign_1,\nbign_2)$
is decomposed as follows:

\begin{lem}\mbox{{}}\label{lem;12.4.140}
There exists the number $M$ and the following data:
\begin{itemize}
\item
Vector spaces $\nbigv_{i\,a}$ for $i=1,2$ and for $a=1,\ldots,M$.
\item
Nilpotent maps $\nbign_{i\,a}\in End(\nbigv_{i\,a})$
for $i=1,2$ and for $a=1,\ldots,M$.
\end{itemize}
 They satisfy the following:
\begin{itemize}
\item
 $\nbigv$ is isomorphic to a direct sum
 $\bigoplus_{a=1}^M\nbigv_{1\,a}\otimes\nbigv_{2\,a}$.
\item
 Under the isomorphism above,
 the nilpotent map $\nbign_1$ is same as
 $\sum_{a=1}^M \nbign_{1\,a}\otimes id_{V_{2\,a}}$.
\item
 Under the isomorphism above,
 the nilpotent map $\nbign_2$ is same as
 $\sum_{a=1}^M id_{V_{1\,a}}\otimes\nbign_{2\,a}$.
\end{itemize}
\end{lem}
\pf
We use the notation in subsection \ref{subsection;12.4.130}.
We have the decomposition
$\prolong{\nbige}=\bigoplus_h\directsummand_h$
satisfying the following:
\[
 f_{A_{1}\superinfty}(\directsummand_h)\subset \directsummand_{h-2},
\quad
 f_{A_2\superinfty}(\directsummand_h)\subset\directsummand_h.
\]
The space $\bigoplus_h\directsummand_{h\,|\,(\lambda,Q)}$
is naturally isomorphic to the graded associated vector space
of the weight filtration of 
$\nbign_1$.
Thus we can decompose 
$\directsummand_{h\,|\,(\lambda,Q)} $
into the primitive parts.
Since $\nbign_1$ and $\nbign_2$
are commutative, the primitive parts are preserved
by $\nbign_2$.
Hence we obtain Lemma \ref{lem;12.4.140}.
\hfill\qed

\vspace{.1in}

Let compare the filtrations
$W(\itibar)_{|(\lambda,Q)}$, $W(\nibar)_{|(\lambda,Q)}$
of $\prolong{\nbige}_{|(\lambda,Q)}$,
and $W(N^{(1)}(\nibar))_{|(\lambda,Q)}$
of $\graded^{(1)}=\graded^{W(\itibar)}$.
For simplicity of the notation,
we omit to denote the notation of the restriction $`|(\lambda,Q)'$.
Let $b$ be the bottom number of the filtration $W(\itibar)$.
Let $v$ be a non-zero element of $W(\itibar)_b$.
We have the degree of $v$ with respect
to the filtration $W(\nibar)$,
which we denote by $\deg^{W(\nibar)}(v)$.
Since we have the natural isomorphism
$W(\itibar)_b\simeq \graded^{(1)}_b$,
we have the degree of $v\in Gr^{(1)}_b$ with respect to
the filtration $W(N^{(1)}(\nibar))$.
The degree is denoted by $\deg^{W(N^{(1)}(\nibar))}(v)$.

\begin{lem}
We have the following equality:
\[
 \deg^{W(N^{(1)}(\nibar))}(v)+b= \deg^{W(\nibar)}(v).
\]
\end{lem}
\pf
Due to the decomposition as in Lemma \ref{lem;12.4.140},
we only have to consider the case:
$\nbigv=\nbigv_1\otimes\nbigv_2$,
$\nbign_{1}=\nbign_1'\otimes id_{\nbigv_2}$,
and
$\nbign_2=id_{\nbigv_1}\otimes \nbign_2'$.
In this case, we can check the equality
by using the decomposition as in 
the subsubsection \ref{subsubsection;12.4.145}.
\hfill\qed

Thus we obtain the following implication
of the vector subbundle over $S(O,P)$ for any point $P$:
\[
 \graded^{W(N^{\sankaku(1)}(\nibar))}_{b,h_2}
\subset
 \graded^{W(N^{\sankaku}(\nibar))}_{b+h_2}.
\]
Let pick an appropriate point $P$
from the curve $\{(t,t,c_3,\ldots,c_n)\in X-D\}$.
Then $S(O,P)$ with the filtration $W^{\sankaku}(\nibar)$
is a mixed twistor.
Then $\graded^{W(N(\nibar))}_{b+h_2}$
is isomorphic to a direct sum
of $\nbigo_{\proj^1}(b+h_2)$.
Thus we obtain the following inequality:
\begin{equation} \label{eq;12.4.150}
 c_1\Bigl(
 \graded^{W(N^{\sankaku(1)}(\nibar))}_{b,h_2}
 \Bigr)
\leq
 (b+h_2)\cdot
 \rank \Bigl(
 \graded^{W(N^{\sankaku(1)}(\nibar))}_{b,h_2}
 \Bigr).
\end{equation}
On the other hand,
we have the following equality:
\[
 \sum_{h_2}
 c_1\Bigl(
 \graded^{W(N^{\sankaku(1)}(\nibar))}_{b,h_2}
 \Bigr)
=c_1\Bigl(
 \graded^{(1)}_{b}
 \Bigr)
=b\cdot\rank \graded^{(1)}_{b}.
\]
Note the following equalities:
\[
 \sum_{h_2} \rank\Bigl(
 \graded^{W(N^{\sankaku(1)}(\nibar))}_{b,h_2}
 \Bigr)=\rank\Bigl(\graded^{(1)}_{b}
 \Bigr),
\quad\quad
 \sum_{h_2} h_2\cdot \rank\Bigl(
 \graded^{W(N^{\sankaku(1)}(\nibar))}_{b,h_2}
 \Bigr)=0.
\]
Thus we have the following equality:
\begin{equation}\label{eq;12.4.151}
\sum_{h_2}c_1\Bigl(
 \graded^{W(N^{\sankaku(1)}(\nibar))}_{b,h_2}
 \Bigr)=
\sum_{h_1} (b+h_2)\cdot
 \rank \Bigl(
 \graded^{W(N^{\sankaku(1)}(\nibar))}_{b,h_2}
 \Bigr).
\end{equation}
The inequality (\ref{eq;12.4.150})
and the equality (\ref{eq;12.4.151})
imply the equality (\ref{eq;12.4.155}).
Hence the proof of Lemma \ref{lem;1.31.12}
is completed.
\hfill\qed


\section{The constantness of the filtrations on the positive cones}

\label{section;12.15.32}

\subsection{Preliminary norm estimate}
\label{subsection;12.5.100}

We put $X=\Delta^n=\{(\zeta_1,\ldots,\zeta_n)\in\Delta^n\}$,
$D_i:=\{\zeta_i=0\}$, and $D=\bigcup_{i=1}^nD_i$.
We also put $D_{\mbar}=\bigcap_{j\leq m}D_j$.
Let $\harmonicbundle$ be a tame nilpotent harmonic bundle
with trivial parabolic structure over $X-D$.
We have the deformed holomorphic bundle
$(\nbige,\DD,h)$ on $\nbigx-\nbigd$,
and the prolongment $\prolong{\nbige}$.
We have the residue $N_i:=\Res_{\nbigd_i}(\DD)$.
For an element $\veca=(a_i)\in \cnum^{n}$,
we put
$N(\veca)_{|(\lambda,O)}
:=\sum_{i=1}^n a_i\cdot N_{i\,|\,(\lambda,O)}$.
We denote the weight filtration of $N(\veca)_{|(\lambda,O)}$
by $\weightfilt(\veca)_{|(\lambda,O)}$.

We have already known the following
(Proposition \ref{prop;12.13.10} and Corollary \ref{cor;12.11.31}):
\begin{lem}
If $N(\veca_1)_{|(\lambda,O)}$ and $N(\veca_2)_{|(\lambda,O)}$
be general in the sense of Definition {\rm\ref{df;12.3.40}},
then $\weightfilt(\veca_1)_{|(\lambda,O)}
=\weightfilt(\veca_2)_{|(\lambda,O)}$.
We also know that 
if $N(\veca_1)_{|(\lambda,O)}$ is general,
then $N(\veca_1)_{|(\lambda',O)}$ is general for any $\lambda'$.
\hfill\qed
\end{lem}
Thus we say $\veca$ is general if $N(\veca)_{|(\lambda,O)}$
is general for all $\lambda$.
In the following of this subsection, we fix a general element $\veca$.

Let $d$ be a positive integer.
Let $\vecb_j=(b_{j\,1},b_{j\,2},\ldots,b_{j\,n})$
be a general element of
$\seisuu_{>0}^n$,
for $j=1,\ldots,d$.
We put $\blowup{X}=\{(z_1,\ldots,z_d)\in\Delta^d\}$,
$\blowup{D}_i=\{z_i=0\}$,
and $\blowup{D}=\bigcup_{i=1}^d\blowup{D}_i$.
We also put $\blowup{D}_{\mbar}=\bigcap_{j\leq m}\blowup{D}_j$.
We have the morphism $f:\blowup{X}\lrarr X$
defined as follows:
\[
 f^{\ast}(\zeta_i)=\prod_{j=1}^d z_j^{b_{j\,i}}.
\]
Note that we have $f(\blowup{D}_i)\subset D_{\nbar}$.
We put $(\blowup{E},\blowup{\theta},\blowup{h}):=f^{\ast}(E,\theta,h)$.
We also put $(\blowup{\nbige},\blowup{\DD})=f^{\ast}(\nbige,\DD)$.

Let $\vecv$ be a holomorphic frame of $\prolong{\nbige}$
over $\nbigx$, which is compatible with
the filtration $W(\veca)$ on $\nbigd_{\dbar}$.
We put $2\cdot k(v_i)=\deg^{W(\veca)}(v_i)$.
We obtain the frame $\blowup{\vecv}:=f^{\ast}\vecv$,
which gives the frame of $\prolong{\blowup{\nbige}}$.
We obtain the $C^{\infty}$-frame
$\blowup{\vecv}'$ defined as follows:
\[
 \blowup{v}_i':=
 \Bigl(-\sum_{j=1}^d\log|z_j|\Bigr)^{-k(v_i)}\cdot \blowup{v}_i.
\]

\begin{prop}\label{prop;12.10.210}
 The frame $\blowup{\vecv}'$ is adapted over
$\Delta_{\lambda}(R)\times (\blowup{X}-\blowup{D})$
for any $R>0$.
\end{prop}
\pf
We use an induction on $d$.
We assume that the claim holds for $d-1$,
and we will prove the claim for $d$,
in the following.
The assumption will be used in Lemma \ref{lem;12.13.30}.

First of all, we note the following:
Let $\vecv_1$ be another holomorphic frame of $\prolong{\nbige}$
over $\nbigd$, which is compatible with
the filtration $W(\veca)$.
Then we obtain the $C^{\infty}$-frame $\blowup{\vecv}'_1$
by the same procedure.
\begin{lem} \label{lem;12.10.230}
If $\blowup{\vecv}'_1$ is adapted over
$\Delta_{\lambda}(R)\times (\blowup{X}-\blowup{D})$,
then $\blowup{\vecv}'$ is adapted over the same region.
\end{lem}
\pf
We have the relation $v_i=\sum_{j}c_{j\,i}\cdot v_{1\,j}$
over $\nbigx$.
Since $\vecv$ and $\vecv_1$ are compatible with the filtration
$W(\veca)$,
we have the following:
\[
  \deg^{W(\veca)}(v_i)<\deg^{W(\veca)}(v_{1\,j})
\Longrightarrow
  c_{j\,i}(\lambda,O)=0.
\]
We have the following relation:
\[
 \blowup{v}_i'=
 \sum_j f^{\ast}(c_{j\,i})\cdot
 \Bigl(-\sum_{m=1}^d\log|z_m|\Bigr)^{-k(v_i)+k(v_{1\,j})}\cdot
 \blowup{v}_{1\,j}'.
\]
If $-k(v_i)+k(v_{1\,j})>0$,
then $f^{\ast}(c_{j\,i})$ is of the form $(\prod_{m=1}^d z_m)\cdot g$
for some holomorphic function $g$ over $\blowup{\nbigx}$.
Hence the transformation matrices of
$\blowup{\vecv}'$ and $\blowup{\vecv}_1'$
are bounded.
Thus we obtain Lemma \ref{lem;12.10.230}.
\hfill\qed

Let us return to the proof of Proposition \ref{prop;12.10.210}.
Let take a holomorphic frame $\vecv$
of $\prolong{\nbige}$ over $\nbigx$,
satisfying the following:
\begin{itemize}
\item $\vecv_{|(\lambda,O)}$ is compatible with the filtration 
$\weightfilt(\veca)_{|(\lambda,O)}$ for the general $\veca$.
\item
For $\vecb_1$,
the representing matrix of the endomorphism
$\sum_i b_{1\,i}\cdot N_{i|\,(\lambda,O)}$
with respect to the frame $\vecv_{|(\lambda,O)}$
is constant,
in other words, independent of $\lambda$.
\end{itemize}
We only have to check the claim of Proposition \ref{prop;12.10.210}
for the frame $\vecv$.

We have the $\lambda$-connection form $\nbiga$
of $\DD$ with respect to $\vecv$.
We decompose $\nbiga$ as follows:
\[
 \nbiga=\sum_{i=1}^n A_i(\lambda,\zeta)\cdot \frac{d\zeta_i}{\zeta_i}.
\]

We denote the $\lambda$-connection form of $\blowup{\DD}$
with respect to the frame $\blowup{\vecv}$
by $\blowup{\nbiga}=\sum\blowup{\nbiga}_j$.
We have the following:
\[
\blowup{\nbiga}= f^{\ast}\nbiga
=\sum_{i=1}^n
  f^{\ast}(A_i)\cdot f^{\ast}\Bigl(\frac{d\zeta_i}{\zeta_i}\Bigr)
=\sum_{i=1}^n
  f^{\ast}(A_i)\cdot\Bigl(\sum_{j=1}^db_{j\,i}\cdot
 \frac{dz_j}{z_j}\Bigr)
=\sum_{j=1}^d\Bigl(
 \sum_{i=1}^n
     b_{j\,i}\cdot f^{\ast}A_i
 \Bigr)\frac{dz_j}{z_j}.
\]
We put $\blowup{A}_j=\sum_{i=1}^n b_{j\,i}\cdot f^{\ast}A_i$.
We also have $f(\blowup{D}_j)\subset \{O\}$.
Then we obtain the following relation.
\begin{equation}
 \blowup{N}_{j}:=
 \Res_{\blowup{\nbigd}_j}(\blowup{\DD})
=\sum_{i=1}^n
  b_{j\,i}\cdot 
  f^{\ast}\Bigl(
  \Res_{\nbigd_i}(\DD)_{|(\lambda,O)}
 \Bigr),
\quad
(j=1,\ldots,d).
\end{equation}
Thus the weight filtration $W(\blowup{N}_j)$ on $\blowup{\nbigd}_j$
is naturally isomorphic to the pull back of $W(\vecb_j)$
for each $j$.

The residue $\blowup{N}_1$
is represented by the following $M(r)$-valued functions
with respect to the frame $\blowup{\vecv}$:
\[
 \blowup{A}_{1\,|\,\blowup{\nbigd}_1}
=\sum_{i=1}^n b_{1\,i}\cdot A_i(\lambda,O).
\]
Note the following lemma.
\begin{lem}
$\blowup{A}_{1\,|\,\blowup{\nbigd}_j}$ is constant, say $A$,
for any $j=1,\ldots,d$.
\end{lem}
\pf
Recall that
we have assumed that $\sum_{i=1}^n b_{1\,i}\cdot A_i(\lambda,O)$
is independent of $\lambda$.
\hfill\qed

We put $V=\prolong{\blowup{\nbige}}_{|(0,\blowup{O})}$
and $N=\Res_{\blowup{\nbigd}_1}(\blowup{\DD})_{|(0,\blowup{O})}$.
We have a model bundle $E(V,N)=(E_0,h_0,\theta_0)$.
We denote the deformed holomorphic bundle
by $(\nbige_0,\DD_0,h_0)$.
We have the canonical frame $\vecv_0$
such that $\DD_0\vecv_0=\vecv_0\cdot A\cdot dz/z$.

Let $q_1$ denote the projection $\Delta^n\lrarr \Delta$
onto the first component.
We put $\blowup{\nbige}_0:=q_1^{\ast}\nbige_0$.
We have the $\lambda$-connection $\blowup{\DD}_0:=q_1^{\ast}\DD_0$
along the $z_1$-direction.
We also put $\blowup{\vecv}_0:=q_1^{\ast}\vecv_0$.

Due to the frames $\blowup{\vecv}$
and $\blowup{\vecv}_0$,
we obtain the holomorphic isomorphism
$\Phi:\prolong{\blowup{\nbige}}_0\lrarr \prolong{\blowup{\nbige}}$.

\begin{lem} \label{lem;12.5.1}
Note the following.
\begin{itemize}
\item
We have
$\Phi\circ \blowup{\DD}_0-\projection_1(\blowup{\DD})\circ\Phi=0$
on $\blowup{\nbigd}_j$
for $j=2,\ldots,d$.
Similarly we have
$\Phi^{-1}\circ\projection_1(\blowup{\DD})-\blowup{\DD}_0\circ\Phi^{-1}=0$
on $\blowup{\nbigd}_j$
for $j=2,\ldots,d$.
\item
We have $\Res(\Phi\circ\blowup{\DD}_0
  -\projection_1(\blowup{\DD})\circ\Phi)=0$
on $\blowup{\nbigd}_1$.
Similarly we have
$\Res(\Phi^{-1}\circ\projection_1(\blowup{\DD})
 -\blowup{\DD}_0\circ\Phi^{-1})=0$
on $\blowup{\nbigd}_1$.
\end{itemize}
\end{lem}
\pf
Clear from our construction.
(See the subsection \ref{subsection;12.14.101}
for $\projection_1$.)
\hfill\qed

We put as follows:
\[
 \tilde{h}_0(\lambda,z_1,\ldots,z_d):=
 h_0\Bigl(\lambda,
 \prod_{i=1}^dz_i
 \Bigr).
\]
Proposition \ref{prop;12.10.210} is a consequence of the following
lemma.
\begin{lem} \label{lem;12.13.30}
The morphisms $\Phi$ and $\Phi^{-1}$ are bounded
with respect to the metrics
$\blowup{h}$ and $\blowup{h}_0$,
over $\Delta_{\lambda}(R)\times (\blowup{X}-\blowup{D})$.
\end{lem}
\pf
Let $a$ be an element of $\Delta^{\ast}$.
We put $\blowup{X}_a=\{(a,z_2,\ldots,z_n)\in\blowup{X}\}$,
and $\blowup{D}_{a\,i}:=\blowup{X}_a\cap \blowup{D}_i$.
Since $\blowup{X}_a$ is $(d-1)$-dimensional,
we can apply the assumption of the induction of
the proof of Proposition \ref{prop;12.10.210}
to $\blowup{X}_a\lrarr X$.
Then we know the following:
{\em
We have the $C^{\infty}$-frame $\blowup{\vecv}^{\clubsuit}$
defined as follows:
\[
 \blowup{v}_i^{\clubsuit}:=
 \Bigl(-\sum_{j=2}^d\log|z_j|\Bigr)^{-k(\blowup{v}_i)}
 \cdot \blowup{v}_i.
\]
Here we put $2\cdot k(\blowup{v}_i)=\deg^{W(\veca)}(v_i)$.
Then $\blowup{\vecv}^{\clubsuit}$ is adapted
over $\Delta_{\lambda}(R)\times(\blowup{X}_a-\blowup{D}_a)$.
}


On the other hand,
we have the $C^{\infty}$-frame $\blowup{\vecv}^{\clubsuit}_0$
over $\blowup{X}_a$ defined as follows:
\[
 \blowup{v}_{0\,i}^{\clubsuit}
=\blowup{v}_{0\,i}\cdot
 \Bigl(-\sum_{j=2}^d\log|z_j|
 \Bigr)^{-k(v_{0\,i})}\blowup{v}_{0\,i}.
\]
Then it is easy to see that
$\blowup{v}_0^{\clubsuit}$ is adapted over
$\Delta_{\lambda}(R)\times (\blowup{X}_a-\blowup{D}_a)$,
due to our construction of $\blowup{v}_0$
and the metric $\blowup{h}_0$.

Thus we obtain the boundedness of the morphisms $\Phi$ and $\Phi^{-1}$
over the boundary
$\Delta_{\lambda}(R)\times
 \bigl\{(z_1,\ldots,z_n)\in\blowup{X}-\blowup{D}\,\big|\,
 |z_1|=C\bigr\}$.
Then we obtain the boundedness of $\Phi$ and $\Phi^{-1}$
on the region
$\Delta_{\lambda}(R)\times
 \bigl\{(z_1,\ldots,z_n)\in\blowup{X}-\blowup{D}\,\big|\,
 |z_1|\leq C\bigr\}$
due to Lemma \ref{lem;12.3.71}.
\hfill\qed

\vspace{.1in}
Thus the induction of the proof of Proposition \ref{prop;12.10.210}
can proceed,
namely, the proof of Proposition \ref{prop;12.10.210}
is completed.
\hfill\qed

\label{subsection;12.5.3}

\subsection{Taking limit}

\label{subsection;12.15.33}

\subsubsection{Replacement of Notation}

We essentially use the setting in the subsection
\ref{subsection;12.5.3},
by putting $d=n$.
For simplicity of the notation,
we replace $\blowup{X}$ with $X$, and make similar replacements for others.
More precisely, we consider as follows:

We put $X=\Delta^n$, $D_i=\{z_i=0\}$ and $D=\bigcup_{i=1}^nD_i$.
Let $\harmonicbundle$ be tame nilpotent harmonic bundle
with trivial parabolic structure over $X-D$.
As usual,
$(\nbige,\DD)$ denotes the deformed holomorphic bundle
with the $\lambda$-connection.
We assume the following:
\begin{condition}\mbox{{}}
\begin{itemize}
\item
There is a tame nilpotent harmonic bundle
$(E_1,\theta_1,h_1)$
with trivial parabolic structure over $X-D$.
We denote the deformed holomorphic bundle by
$(\nbige_1,\DD_1)$.
\item
There are general elements $\vecb_j\in\seisuu_{>0}^n$ $(j=1,\ldots,n)$,
and $f:X\lrarr X$ as in the subsection {\rm\ref{subsection;12.5.3}}.
\item
We have $f^{\ast}(E_1,\theta_1,h_1)=(E,\theta,h)$.
\item
We have a frame $\vecv_1$ of $\prolong{\nbige}_1$
compatible with the filtration
$W(\veca)$ on $\cnum_{\lambda}\times \{O\}$
for general $\veca\in\cnum^n$.
And we have $\vecv=f^{\ast}(\vecv_1)$.
\hfill\qed
\end{itemize}
\end{condition}

By the frame $\vecv$, we decompose $\prolong{\nbige}$ as follows:
\[
 \prolong{\nbige}=
 \bigoplus_h U_h,
\quad\quad
 U_h:=
 \big\langle
  v_i\,\big|\,\deg^{\weightfilt(N_1)}(v_i)=h
 \big\rangle.
\]
We put $2\cdot k(v_i)=\deg^{\weightfilt(N_1)}(v_i)$.
Recall that the degrees of $v_i$
with respect to the weight filtration
$\weightfilt(N_j)$ are independent of $j$,
i.e.,
$\deg^{W(N_j)}(v_i)=\deg^{W(N_1)}(v_i)$.

\subsubsection{Pull backs}

Let $m$ be a non-negative integers.
We have the morphism $\psi_{m,\nbar}:X\lrarr X$ 
or $X-D\lrarr X-D$, defined as follows:
\[
 \psi_{m,\nbar}(z_1,z_2,\ldots,z_n):=(z_1^m,z_2^m,\ldots,z_n^m).
\]
Then we obtain the harmonic bundles
$\psi_{m,\nbar}^{\ast}\harmonicbundle$ on $X-D$.
We also obtain the deformed holomorphic bundles
and the $\lambda$-connections
$(\psi_{m,\nbar}^{\ast}\prolong{\nbige},\psi_{m,\nbar}^{\ast}\DD)$.
We obtain the holomorphic frame $\vecv^{(m)}$
of $\psi_{m,\nbar}^{\ast}\prolong{\nbige}$,
defined as follows:
\[
 v^{(m)}_i:=
 \psi_{m,\nbar}^{-1}(v_i)\cdot m^{-k(v_i)}.
\]
We put $H^{(m)}:=H(\psi_m^{\ast}(h),\vecv^{(m)})$,
which is an $\nbigh(r)$-valued function.
The following lemma can be shown by an argument
similar to the proof of Lemma \ref{lem;12.5.10},
by using Proposition \ref{prop;12.10.210}.
\begin{lem}
On any compact subset $K\subset X-D$,
the $\nbigh(r)$-valued functions
$\{H^{(m)}\}$ and $\{H^{(m)\,-1}\}$
are bounded independently of $m$.
\hfill\qed
\end{lem}

We have the $\lambda$-connection form
$\nbiga\in
    \Gamma(\nbigx,
    M(r)\otimes p_{\lambda}^{\ast}\Omega_{X}^{1,0}(\log X))$
of $\DD$ with respect to the frame $\vecv$.
We decompose $\nbiga$ as follows:
\[
 \nbiga=\sum_j \nbiga_j,
\quad\quad
 \nbiga_j
 \in \Gamma(\nbigx,
   M(r)\otimes
    p_{\lambda}^{\ast}q_j^{\ast}
       \Omega_{\Delta}(\log O)).
\]
We obtain $A_j\in M(r)\otimes\nbigo_{\nbigx}$ satisfying the following:
\[
 \nbiga_j= A_j\cdot \frac{dz_j}{z_j}.
\]
We decompose $A_j$ into $\sum_{h,k}A_{j,(h,k)}$
as in the subsubsection \ref{subsubsection;12.5.10}.
We have the following:
\[
 A_{j,(h,k)}(\lambda,O)=0,
\quad\quad
 \mbox{ if } h>k-2.
\]

We have the $\lambda$-connection form $\nbiga^{(m)}$
of $\psi_{m,\nbar}^{\ast}\DD$
with respect to the frame $\vecv^{(m)}$.
We decompose $\nbiga^{(m)}$ into
$\sum_j\nbiga^{(m)}_j$ as in the case of $\nbiga$.
\begin{lem} \label{lem;12.5.51}
We have the following equalities:
\[
 \nbiga^{(m)}_j=
 \sum_{h,k}\psi_{m,\nbar}^{\ast}A_{j\,(h,k)}\cdot m^{(h-k+2)/2}\cdot
 \frac{dz_j}{z_j}.
\]
In particular, the sequences $\{\nbiga^{(m)}\}$ converges
to $\nbiga^{(\infty)}=\sum_j\nbiga_j^{(\infty)}$ given as follows:
\[
 \nbiga_j^{(\infty)}=
 \sum_h A_{j\,(h,h+2)}(\lambda,O)\cdot 
 \displaystyle  \frac{dz_j}{z_j}.
\]
\end{lem}
\pf
Similar to Lemma \ref{lem;12.4.51}.
\hfill\qed

We decompose $\theta$ into
$\sum_{j=1}^n f_j\cdot dz_j/z_j$
Then we have
$\psi_{m,\nbar}^{\ast}\theta=
 \sum_j
 m\cdot\psi_m^{\ast}(f_j)\cdot dz_j/z_j$.
\begin{lem}
We have the following inequality independent of $m$:
\[
 |m\cdot\psi_m^{\ast}(f_j)|_{\psi_m^{\ast}(h)}
 \leq C\cdot(-\log|z_j|)^{-1}.
\]
\end{lem}
\pf
Similar to Lemma \ref{lem;12.13.31}.
\hfill\qed

In all we obtain the following.
\begin{lem}
The sequence $\{\vecv^{(m)}\}$
satisfies Condition {\rm \ref{condition;11.27.16}}.
\hfill\qed
\end{lem}

\subsubsection{Limit}

We can apply the result in the subsection \ref{subsection;12.4.10}
as in the subsubsection \ref{subsubsection;12.14.2}.
Thus we obtain a limiting harmonic bundle
$(F,\theta\superinfty,h\superinfty)$
and the frame $\vecv\superinfty$
of the prolongment $\prolong{\nbigf\superinfty}$.
We have the following equality:
\[
 \DD^{(\infty)}\vecv^{(\infty)}=
 \vecv^{(\infty)}\cdot\nbiga^{(\infty)},
\quad
 \nbiga^{(\infty)}
=\sum_j\sum_h A_{j\,(h,h+2)}(\lambda,O)\frac{dz_j}{z_j}.
\]

We have the $C^{\infty}$-frame $\vecv^{\prime\,(\infty)}$
defined by 
$v_i^{\prime\,(\infty)}=
 \Bigl(-\sum_{j=1}^n\log|z_j|\Bigr)^{-k(v_i)}\cdot v_i^{(\infty)}$.
Here $\vecv$ denotes our original frame.
\begin{lem}
The frame $\vecv^{\prime\,(\infty)}$ is adapted.
\end{lem}
\pf
Let consider the $C^{\infty}$-frame $\vecv^{\prime\,(m)}$
defined as follows:
\[
 v^{\prime\,(m)}_i:=v^{(m)}_i\cdot
 \Bigl(-\sum_{j=1}^n\log|z_j|\Bigr)^{-k(v_i)}.
\]
Then $\vecv^{\prime\,(m)}$ is adapted,
independently of $m$.
Then the adaptedness of $\vecv^{\prime\,(\infty)}$ follows immediately.
\hfill\qed

We have the real coordinate
$z_i=r_i\cdot\exp(2\pi\sqrt{-1}\alpha_i)$
for $i=1,\ldots,n$.
\begin{lem}
The $\nbigh(r)$-valued function $H(h^{(\infty)},\vecv^{(\infty)})$
is independent of $\alpha_i$
for any $i$.
\end{lem}
\pf
Similar to Proposition \ref{prop;12.5.60}.
\hfill\qed

We put $\vece^{(\infty)}=\vecv^{(\infty)}_{|\nbigx^0}$.
We have the decomposition of $F$ into $\bigoplus_h U_h$,
where $U_h$ denote the vector subbundles
of $F$ generated by
$\{e_i\superinfty\,|\,\deg^{W(N_1)}(e_i)=h\}$.
Note that $A_{j\,(h,h+2)}$ and $\vece\superinfty$
determines the morphism $U_{h+2}\lrarr U_h$.

Let pick an element $\vecc=(c_1,\ldots,c_n)\in \real_{>0}^n$.
Then we have the morphism
$\Xi_{\vecc}:\hyperh\lrarr X-D$ defined as follows:
\[
 \Xi_{\vecc}(\zeta)=
 (\exp(2\pi\sqrt{-1}c_1\zeta),
\ldots,
 \exp(2\pi\sqrt{-1}c_n\zeta)
 )
\]
We put as follows:
\[
 (F_{\vecc},h_{\vecc},\theta_{\vecc})
 :=\Xi_{\vecc}^{\ast}
 (F,h\superinfty,\theta\superinfty),
\quad
 \vecw_{\vecc}:=
 \Xi_{\vecc}^{\ast}\vece^{(\infty)},
\quad
 U_{h,\vecc}:=
 \Xi_{\vecc}^{\ast}U_h.
\]
We use the real coordinate $\zeta=x+\sqrt{-1}y$.
We have the $C^{\infty}$-frame $\vecw_{\vecc}'$
defined by
$w_{\vecc\,i}':=y^{-k(v_i)}\cdot w_{\vecc\,i}$.
Here $\vecv$ is our original frame.

\begin{lem}\mbox{{}}
\begin{itemize}
\item
$H(h_{\vecc},\vecw_{\vecc})$ is independent of $x$.
\item
We have the following equality:
\[
 \theta_{\vecc}(\vecw_{\vecc})=
 \vecw_{\vecc}\cdot
 \Bigl(
 \sum_{j}c_j\cdot A_j^{(\infty)}(0,O)
 \Bigr)\cdot(2\pi\sqrt{-1}d\zeta).
\]
\item
We have $A_j^{(\infty)}(0,O)=\sum_h A_{j,(h,h+2)}(0,O)$.
\item
The frame $\vecw_{\vecc}'$ is adapted.
\hfill\qed
\end{itemize}
\end{lem}

We can descend $(F_{\vecc},h_{\vecc},\theta_{\vecc})$
to the harmonic bundle
over $\Delta^{\ast}$ by using the framing $\vecw_{\vecc}$.
Since $H(h_{\vecc},\vecw_{\vecc})$
does not depend on $x$,
it gives the metric $h_{\vecc}$ defined over $\Delta^{\ast}$.
Similarly we obtain the Higgs field $\theta_{\vecc}$
and the vector subbundles $U_{h,\vech}$
defined over $\Delta^{\ast}$.

On $\Delta^{\ast}$,
we have the following:
\begin{lem}\mbox{{}}
\begin{itemize}
\item
We have the following equality:
\[
 \theta_{\vecc}(\vecw_{\vecc})=
 \vecw_{\vecc}\cdot
 \Bigl(
 \sum_{j}c_j\cdot A_j^{(\infty)}(0,O)
 \Bigr)\cdot\frac{dz}{z}.
\]
In particular, $\theta_{\vecc}$ is tame and nilpotent.
\item
We have the $C^{\infty}$-frame $\vecw'_{\vecc}$
defined by
$w'_{\vecc\,i}:=w_{\vecc\,i}\cdot(-\log|z|)^{-k(v_i)}$,
where $\vecv=(v_i)$ is our original frame.
Then the frame $\vecw_{\vecc}'$ is adapted.
\hfill\qed
\end{itemize}
\end{lem}

Thus $\vecw_{\vecc}$ naturally gives a frame of
$\prolong{F_{\vecv}}$.

\begin{lem}
The weight filtration of 
$\sum_{j}c_j\cdot A_j^{(\infty)}(0,O)$
is same as the filtration 
$\bigl\{W_l:=\bigoplus_{h\leq l}\prolong{U}_{h\,|\,O}\,|\,l\in\seisuu
 \bigr\}$.
\end{lem}
\pf
We know that $|w_{\vecc\,i}|_{h_{\vecc}}\sim (-\log|z|)^{k(v_i)}$,
where $\vecv$ is the original frame.
Thus the norm estimate in one dimensional case
implies our result.
\hfill\qed

We put $A_{h,h+2}(\vecc):=\sum_j c_j\cdot A_{j\,h\,h+2}(0,O)$.
The endomorphism $N(\vecc)_{|(0,O)}\in End(\prolong{\nbige}_{|(0,O)})$
induces the morphism
$g_h(\vecc):Gr_{h\,|\,(0,O)}^{W(N_1)}\lrarr Gr^{W(N_1)}_{h-2\,|\,(0,O)}$.
We also have the frame $\vecv^{(1)}$ of $Gr^{W(N_1)}_{|(0,O)}$
induced by $\vecv$.
The matrix $A_{h,h+2}(\vecc)$ represents
$g_{h+2}(\vecc)$ with respect to the frame $\vecv^{(1)}$.

\begin{cor} \label{cor;12.11.20}
For any $k\geq 0$,
the following isomorphism is isomorphic:
\[
 g_{-k,-k+2}(\vecc)\circ\cdots\circ g_{k-2,k}(\vecc):
 Gr_{k\,|\,(0,O)}^{W(N_1)}\lrarr Gr_{-k\,|\,(0,O)}^{W(N_1)}.
\]
\end{cor}
\pf
We only have to note that $\prolong{U}_{k\,|\,O}$
is naturally isomorphic to $Gr^{W(N_1)}_k$
via the frames $\vecv^{(1)}$ and $\vecw_{\vecc}$.
\hfill\qed

\begin{cor}\label{cor;12.11.21}
The weight filtration of $N(\vecc)_{|(\lambda,O)}$ is same
as $W(N_1)_{|(\lambda,O)}$ for any $\lambda$.
\end{cor}
\pf
In the case $\lambda=0$,
it is a consequence of Corollary \ref{cor;12.11.20}
and the fact that $N(\vecc)_{|(0,O)}$ drops the degree
with respect to $W(N_1)_{|(0,O)}$
by $2$.
In the other case, we only have to use
the fact that the conjugacy classes of $N(\vecc)_{|(\lambda,O)}$
are independent of $\lambda$.
\hfill\qed

\subsection{Constantness of the filtrations on the positive cones}

We put $X=\Delta^n$, $D_i=\{z_i=0\}$
and $D=\bigcup_{i=1}^l D_i$ for $l\leq n$.
Let $(E,\theta,h)$ be a tame nilpotent harmonic bundle over $X-D$.

Let $Q$ be a point of $D_{\mbar}$.
We have the nilpotent maps
$N_{j\,|\,(\lambda,Q)}$ $(j=1,\ldots,m)$.
For any element $\veca\in\real_{\geq 0}^m$,
we put $N(\veca)_{|(\lambda,Q)}:=\sum_{j=1}^ma_j\cdot N_{j\,|\,(\lambda,Q)}$.
Let $W(\veca)_{|(\lambda,Q)}$ 
denote the weight filtration of $N(\veca)_{|(\lambda,Q)}$.

\label{subsection;12.15.34}
\begin{thm} \label{thm;12.7.200}
The constantness of the filtration on the positive cones holds.
Namely we have $W(\veca_1)_{|(\lambda,Q)}=W(\veca_2)_{|(\lambda,Q)}$
for $\veca_1,\veca_2\in\real_{>0}^I$,
where $I$ denotes any subset of $\mbar$.
\end{thm}
\pf
We only have to check the following:
Assume that $a_{1\,j}$ and $a_{2\,j}$ are positive for $j=1,\ldots,k$.
Assume that $a_{1\,j}=a_{2\,j}=0$ for $j=k+1,\ldots,l$.
Then $W(\veca_1)_{|(\lambda,Q)}=W(\veca_2)_{|(\lambda,Q)}$.

Since such property is not depend on $\lambda$,
we only have to check the claim in the case $\lambda=1$.
Then, by using the normalizing frame,
we know that we have to check such claim when $m=k$,
i.e.,
$Q\in D_{\kbar}-\bigcup_{j>k}(D_{\kbar}\cap D_j)$.

We only have to check the claim when $m=k=n$.
For, we take an $m$-dimensional hyperplane,
which intersects with $D_{\mbar}$ transversally at $Q$.

Thus we have reduced the theorem to the following claim:
\begin{quote}
We put $X=\Delta^n$ and $D=\bigcup_{i=1}^nD_i$.
Let $\harmonicbundle$ be a tame nilpotent harmonic bundle
over $X-D$.
Let $O$ be the origin.
We have the residues
$N_{i\,|\,(\lambda,O)}$ $i=1,\ldots,n$.
We have $W(\veca_1)_{|(\lambda,Q)}=W(\veca_2)_{|(\lambda,Q)}$
if $\veca_i\in\real_{>0}^n$.
\end{quote}

The following lemma can be shown by an elementary argument.
\begin{lem}
Let $\veca$ be an element of $\real^n_{>0}$.
Then there exists general elements
$\vecb_1,\ldots,\vecb_n$ contained in $\rnum_{>0}^n$
and
positive numbers $\alpha_1,\ldots,\alpha_n$
satisfying the following:
\[
 \veca=\sum_{j=1}^n\alpha_j\cdot \vecb_j.
\]
\end{lem}
\pf
For example we can argue as follows:
Let take general elements $\vecc_j\in\real_{>0}^n$
sufficiently close to
$(\overbrace{0,\ldots,0}^{j-1},1,\overbrace{0,\ldots,0}^{n-j})$.
Then $\veca$ is contained in the positive cone
generated by $\vecc_i$ $(i=1,\ldots,n)$.

We take $\vecd_i\in\rnum_{>0}^n$ sufficiently close to $\vecc_i$.
Then $\vecd_i$ are general, and $\veca$ is contained
in the positive cone generated by $\vecd_i$.

We take $\vecb_i\in\seisuu_{>0}^n$ such that
$\vecb_i=C_i\cdot\vecd_i$ for some $C_i\in\real_{>0}$.
Then the tuple $(\vecb_1,\ldots,\vecb_n)$ have the desired property.
\hfill\qed

Then we take the morphism $f:\blowup{X}\lrarr X$
for $(\vecb_1,\ldots,\vecb_n)$
as in the subsection \ref{subsection;12.5.100}.
For the harmonic bundle
$f^{\ast}(E,\theta,h)$,
we have already known the constantness of the filtrations
(Corollary \ref{cor;12.11.21}).
Hence we know that $\veca$ is general.
Thus the proof of Theorem \ref{thm;12.7.200} is completed.
\hfill\qed

\vspace{.1in}
We restate a limiting mixed twistor theorem.
Let $Q$ be a point of $D=\bigcup_{i=1}^lD_i$.
We put $I=\{i\,|\,Q\in D_i\}$.
Then we put $N^{\sankaku}(I):=\sum_{i\in I}N^{\sankaku}_i$.
We denote the weight filtration of $N^{\sankaku}(I)$
by $W^{\sankaku}(I)$.

\begin{thm} \label{thm;12.15.121}
\begin{itemize}
\item
Let $U$ be any neighborhood of $Q$.
For an appropriate point $P\in U\cap (X-D)$,
the filtered vector bundle
$(S(Q,P),W^{\sankaku}(I))$ is a mixed twistor.
\item
For any $i\in I$,
the morphism $N_i^{\sankaku}:S(Q,P)\lrarr S(Q,P)\otimes\nbigo_{\proj^1}(2)$
is a morphism of mixed twistor.
\hfill\qed
\end{itemize}
\end{thm}

\section{Strong sequential compatibility}

\label{section;12.15.35}

\subsection{The comparison in the bottom part}
\label{subsection;12.15.36}

\subsubsection{Two dimensional case}
\label{subsubsection;12.15.37}

Let $X$ be $\{(\zeta_1,\zeta_2)\in\Delta^2\}$,
$D_i=\{\zeta_i=0\}$,
and $D=D_1\cup D_2$.
We put $\{O\}=D_1\cap D_2=D_{\nibar}$, and $D_1=D_{\itibar}$.
Let $\harmonicbundle$ be a harmonic bundle over $X-D$.
We have the deformed holomorphic bundle
$(\nbige,\DD)$ and the prolongment $\prolong{\nbige}$.
We put $N_i=\Res_{\nbigd_i}(\DD)$.
On $\nbigd_{\jbar}$,
we put $N(\jbar)=\sum_{i\leq j}N_{i\,|\,\nbigd_{\jbar}}$ as usual.
Let $W(\jbar)$ denote the weight filtration of 
$N(\jbar)$.
From the filtration $W(\itibar)$,
we obtain the graded vector bundle $\graded^{(1)}$
over $\nbigd_{\itibar}$.
On $\graded^{(1)}_{|\nbigd_{\nibar}}$,
we have the induced operator $N^{(1)}(\nibar)$.
Thus we obtain the weight filtration
$W(N^{(1)}(\nibar))$.

Hence we have the filtrations $W(\itibar)$ and $W(\nibar)$
of $\prolong{\nbige}_{|\nbigd_{\nibar}}$.
Let $b$ be the bottom number of $W(\itibar)$.
We also have the filtrations
$W(N^{(1)}(\nibar))_h\cap \graded^{(1)}_{b|\nbigd_{\nibar}}$
and $W^{(1)}(\nibar)_h\cap\graded^{(1)}_{b|\nbigd_{\nibar}}$
of $\graded^{(1)}_{b|\nbigd_{\nibar}}$.
We would like to see the relation of them.

Then we have the natural isomorphism
$\graded_b^{(1)}\simeq W(\itibar)_b$.
Thus we obtain the following inclusion:
\[
 W(N^{(1)}(\nibar))_h\cap \graded_{b|\nbigd_{\nibar}}^{(1)}
\subset
 W(\itibar)_{b\,|\,\nbigd_{\nibar}}\subset
 \prolong{\nbige}_{|\nbigd_{\nibar}}.
\]
We also have the following:
\[
 W^{(1)}(\nibar)_h
\cap
 \gradedrestrictedtoni
=W(\nibar)_h\cap W(\itibar)_{b\,|\,\nbigd_{\nibar}}.
\]

\begin{lem} \label{lem;12.7.130}
We have the following implication:
\[
 W(N^{(1)}(\nibar))_h\cap \gradedrestrictedtoni
\subset
W(\nibar)_{h+b}\cap W(\itibar)_{b\,|\,\nbigd_{\nibar}}.
\]
\end{lem}
\pf
We have already known that both of them are vector bundle
over $\nbigd_{\nibar}$
(subsubsection \ref{subsubsection;12.10.40}).
Thus we only have to prove the claim on $\cnum_{\lambda}^{\ast}$.
Let pick $\lambda\neq 0$.
Let take a normalizing frame $\vecv$ of $\prolong{\nbige}^{\lambda}$
compatible with the sequence of the filtrations
$(W(\itibar),W(N^{(1)}(\nibar)))$.
Namely it satisfies the following:
\begin{itemize}
\item
It is compatible with the filtration $W(\itibar)$
on $\nbigd_{\itibar}$.
\item
The induced frame $\vecv^{(1)}$ is compatible with
$W(N^{(1)}(\nibar))$
on $\nbigd(\nibar)$.
\end{itemize}

Let $\nbiga=A_1\cdot d\zeta_1/\zeta_1+A_2\cdot d\zeta_2/\zeta_2$
denote the $\lambda$-connection form of $\DD^{\lambda}$
with respect to $\vecv$.
Here $A_i$ are elements of $M(r)$.
We put $U_b:=\langle v_i\,|\,\deg^{W(\itibar)}(v_i)=b\rangle$,
which is a subbundle of $\prolong{\nbige^{\lambda}}$.
\begin{lem}
The vector subbundle $U_b$ is stable under the action of
 $\DD^{\lambda}$.
Namely $\DD^{\lambda}(f)$ is a section of $U_b$,
if $f$ is a section of $U_b$.
\end{lem}
\pf
Let $f_{A_i}$ be the endomorphism determined by $\vecv$
and $A_i$.
Then we have $f_{A_1}(U_b)\subset U_{b-2}=0$ and $f_{A_2}(U_b)\subset U_b$.
It implies our claim.
\hfill\qed

We put $\blowup{X}=\Delta^2$, $\blowup{D}_i=\{z_i=0\}$,
and $\blowup{D}=\blowup{D}_1\cup \blowup{D}_2$.
Let consider the morphism $f:\blowup{X}\lrarr X$
defined by $(z_1,z_2)\longmapsto (z_1\cdot z_2,z_2)$.
We take the pull back 
$(\blowup{\nbige},\blowup{\DD}^{\lambda},\blowup{h},\blowup{\vecv})
:=f^{\ast}(\nbige^{\lambda},\DD^{\lambda},h,\vecv)$.
We also put $\blowup{U}_b:=f^{\ast}\blowup{U}_b$.
Let $\blowup{\nbiga}$ denote
the $\lambda$-connection form of $\blowup{\DD}^{\lambda}$
with respect to $\blowup{\vecv}$.
Then we have the following:
\[
\blowup{\nbiga}=
 A_1\cdot\frac{dz_1}{z_1}
+(A_1+A_2)\cdot\frac{dz_2}{z_2}.
\]

We take the projection
$\gminiq_2:\Omega^{1,0}_{X-D}\lrarr
 q_2^{\ast}\Omega^{1,0}_{\Delta^{\ast}}$.
Note that we consider the projection onto the second component,
different from the other cases.
Then we obtain the $\lambda$-connection along the $z_2$-direction:
$\gminiq_2(\blowup{\DD}^{\lambda})$.
By restriction,
we have the holomorphic bundle
$\blowup{U}_b$ and the $\lambda$-connection
$\gminiq_2(\blowup{\DD}^{\lambda})$.
We put as follows:
\[
 \blowup{\vecv}_b
:=(\blowup{v}_i\,|\,\deg^{W(\itibar)}(\blowup{v}_i)=b).
\]
Then  we have $\gminiq_2(\blowup{\DD}^{\lambda})(\blowup{\vecv}_b)
=\blowup{\vecv}_b\cdot A_3\cdot dz_2/z_2 $
for some constant matrix $A_3$.

We put $V=W(\itibar)_{b\,|\,(\lambda,O)}$.
We have the induced operator $N^{(1)}(\nibar)$
on $V$. We denote it by $N$.
Note that it is represented by the matrix $A_3$
for some appropriate base.
Then we have a model bundle $E(V,N)=(E_0,\theta_0,h_0)$.
We have the deformed holomorphic bundle
$(\nbige_0^{\lambda},\DD^{\lambda}_0)$
on $\Delta^{\ast}$.
We also have the canonical frame $\vecv_0$
such that $\DD^{\lambda}_0\vecv_0=\vecv_0\cdot A_3\cdot dz/z$.

Let $q_2$ denote the projection
$\Delta^{\ast}_1\times \Delta^{\ast}_2
\lrarr \Delta_2^{\ast}$
onto the second component.
We put
$(\blowup{\nbige}_0^{\lambda},\DD_0^{\lambda},\blowup{\vecv}_0)
:=q_2^{\ast}(\nbige_0,\DD_0,\vecv_0)$.
We also put as follows:
\[
 \blowup{h}_0(\lambda,z_1,z_2):=
 (-\log|z_1|)^b\cdot
 h_0(\lambda,z_2).
\]

Due to the frames $\blowup{\vecv}_b$
and $\blowup{\vecv}_0$,
we obtain the holomorphic isomorphism
$\Phi:\blowup{\nbige}_0^{\lambda}\lrarr \blowup{U}_b$,
satisfying the following:
\begin{equation} \label{eq;12.7.120}
 \Phi\circ\blowup{\DD}_0
-\gminiq_2(\blowup{\DD})\circ\Phi=0.
\end{equation}
The morphism $\Phi$ induces the morphism
$\blowup{\nbige}_0^{\lambda}\lrarr\blowup{\nbige}^{\lambda}$.
We denote it by the same notation $\Phi$.
It satisfies (\ref{eq;12.7.120}).

\begin{lem} 
The morphism $\Phi$ is bounded
with respect to the metrics
$\blowup{h}_0$ and $\blowup{h}$.
\end{lem}
\pf
We only have to check the boundedness on the boundary
$\{(z_1,z_2)\,|\,|z_2|=C\}$.
Let $\blowup{v}_i$ be an element of $\blowup{\vecv}_b$.
Then we have $\deg^{W(\itibar)}(\blowup{v}_i)=b$.
Thus we obtain the estimate on the boundary,
due to the norm estimate in one dimensional case.
\hfill\qed

Let $\blowup{v}_i$ be an element of $\blowup{\vecv}_b$.
Due to the lemma,
we obtain the following inequality:
\[
 |\blowup{v}_i|_{\blowup{h}}
\leq C\cdot (-\log|z_1|)^{b/2}\cdot (-\log|z_2|)^{k(\blowup{v}_i)}.
\]
Here we put
$2\cdot k(\blowup{v}_i)
:=\deg^{W^{(1)}(\nibar)}(\blowup{v}_i)$.
Let consider the curve $C(1,1)=\{(z,z)\in \blowup{X}\}$.
On the curve $C(1,1)$,
we have the following estimate:
\[
 |\blowup{v}_{i\,|\,C(1,1)}|_{\blowup{h}}
\leq C\cdot (-\log|z|)^{b/2+k(\blowup{v}_i)}.
\]
It implies $\blowup{v}_{i}(O)$ is contained
in $W(\nibar)_{b+2k(\blowup{v}_i)}$.
Thus we obtain our claim.
\hfill\qed

\begin{lem} \label{lem;12.7.151}
We have
$W(N^{(1)}(\nibar))_h\cap \graded^{(1)}_b=
 W^{(1)}(\nibar)_{b+h}\cap \graded^{(1)}_b$.
\end{lem}
\pf
We will use the mixed twistor structure.
Let take an appropriate point $P\in \blowup{X}-\blowup{D}$
such that the filtration $W(\nibar)$ gives the mixed twistor structure
to the filtered vector bundle $S(O,P)$.
Then the filtration
$\{W^{(1)}(\nibar)_h\cap \graded^{(1)}_b\,|\,h\in\seisuu\}$
gives a mixed twistor structure to $\graded^{(1)}_b$.
We denote the filtration by $W^{\star}$ for simplicity.
The graded vector bundle is denoted by $\graded^{\star}$.

On the other hand,
we have the nilpotent map
$N^{(1)}(\nibar):
 \graded^{(1)}_b\lrarr\graded^{(1)}_b\otimes\nbigo(2)$.
The weight filtration is denoted by $W^{\heartsuit}$,
and the associated graded bundle is denoted by $\graded^{\heartsuit}$.
Our purpose is to show that
$W^{\heartsuit}_{h}=W^{\star}_{h+b}$.

Due to Lemma \ref{lem;12.7.130},
we have already known that $W^{\heartsuit}_{h}\subset W^{\star}_{h+b}$.
Due to Lemma \ref{lem;1.31.12},
we have
$c_1(\graded^{\heartsuit}_h)
=(h+b)\cdot\rank(\graded^{\heartsuit}_h)$.
Then, due to Lemma \ref{lem;1.8.10},
the vector bundle $\graded^{\heartsuit}_h$
naturally gives a vector subbundle of $\graded^{\star}_{h+b}$.

We obtain the inequalities
$\rank\graded_h^{\heartsuit}\leq \rank\graded^{\star}_{h+b}$.
We also have the equalities:
\[
 \sum_h\rank\graded_h^{\heartsuit}
=\rank\graded^{(1)}_b
=\sum_h\rank\graded^{\star}_{h+b}.
\]
Thus we obtain the equalities
$\rank\graded_h^{\heartsuit}=\rank\graded^{\star}_{h+b}$.
It implies that $W_{h}^{\heartsuit}=W_{h+b}^{\star}$.
\hfill\qed

\subsubsection{Higher dimensional case}

\label{subsubsection;12.15.40}

We put $X:=\Delta^n$, $D_i=\{z_i=0\}$
and $D=\bigcup_{i=1}^l D_i$ for some $l\leq n$.
Let $\harmonicbundle$ be a tame nilpotent harmonic bundle
with trivial parabolic structure over $X-D$.
We have the deformed holomorphic bundle with $\lambda$-connection
$(\nbige,\DD)$ and the prolongment $\prolong{\nbige}$.
We put $N_i=\Res_{\nbigd_i}(\DD)$.

Let $Q$ be a point of $D_{\mbar}=\bigcap_{j\leq m}D_j$
for $m\leq l$.
We put $N(\jbar)_{|(\lambda,Q)}
=\sum_{i\leq j}N_{i\,|\,(\lambda,Q)}$.
We denote the weight filtration of $N(\jbar)_{|(\lambda,Q)}$
by $W(\jbar)_{|(\lambda,Q)}$.
Let $b$ be the bottom number of
$W(\itibar)$.
On $W(\itibar)_{b\,|\,(\lambda,Q)}
\simeq \graded^{(1)}_{b\,|\,(\lambda,Q)}$,
we have filtrations
$W^{(1)}(\jbar)_{|(\lambda,Q)}\cap \graded^{(1)}_{b\,|\,(\lambda,Q)}$
and
$W(N^{(1)}(\jbar))_{|(\lambda,Q)}\cap\graded^{(1)}_{b\,|\,(\lambda,Q)}$.
\begin{prop} \label{prop;12.7.150}
We have 
$W^{(1)}(\jbar)_{h+b|(\lambda,Q)}\cap \graded^{(1)}_{b\,|\,(\lambda,Q)}
=
W(N^{(1)}(\jbar))_{h|(\lambda,Q)}\cap\graded^{(1)}_{b\,|\,(\lambda,Q)}$.
\end{prop}
\pf
We have already known that both of them 
are vector bundles over $\nbigd_{\mbar}$.
Thus we only have to see the equality
in the case $\lambda\neq 0$.

Let take a morphism $\blowup{X}\lrarr X$
defined  as in the subsection \ref{subsection;12.4.1}.
Let $\blowup{Q}\in \blowup{\nbigd}_1\cap\blowup{\nbigd}_j$
be a point such that $\phi(\blowup{Q})=Q$.
Then we have 
$\blowup{N}_{1\,|\,(\lambda,\blowup{Q})}:=
 \Res_{\blowup{\nbigd}_1}(\blowup{\DD})_{|(\lambda,\blowup{Q})}
=\phi^{\ast}N(\itibar)_{|(\lambda,Q)}$,
and
$\blowup{N}_{j\,|\,(\lambda,\blowup{Q})}:=
 \Res_{\blowup{\nbigd}_j}(\blowup{\DD})_{|(\lambda,\blowup{Q})}
=\phi^{\ast}N(\jbar)_{|(\lambda,Q)}$.
Thus we only have to compare
$\blowup{N}_{1\,|\,(\lambda,Q)} $
and $\blowup{N}_{j\,|\,(\lambda,Q)}$
at $(\lambda,\blowup{Q})$.

Since we have assumed $\lambda\neq 0$,
we can take a normalizing frame.
Then we only have to compare
in the case
$(\lambda,\blowup{Q})\in 
\blowup{\nbigd}_1^{\lambda}
 \cap
\blowup{\nbigd}_j^{\lambda}-
\bigcup_{i\neq j,1}
 \blowup{\nbigd}_1^{\lambda}
 \cap
\blowup{\nbigd}_j^{\lambda}
\cap\blowup{\nbigd}_i^{\lambda}$.

Let take a two dimensional hyperplane of $\blowup{X}$,
which intersects $\blowup{D}_1\cap \blowup{D}_j$
at $\blowup{Q}$ transversally.
Then we can reduce Proposition \ref{prop;12.7.150}
to Lemma \ref{lem;12.7.151}.
\hfill\qed

\subsection{Theorem}

\label{subsection;12.15.41}

We put $X=\Delta^n$, $D_i=\{z_i=0\}$, and $D=\bigcup_{i=1}^l D$
for $l\leq n$.
We put $D_{\mbar}=\bigcap_{j\leq m}D_j$.
Let $\harmonicbundle$ be a tame nilpotent harmonic bundle
with trivial parabolic structure,
over $X-D$.
We have the deformed holomorphic bundle with $\lambda$-connection
$(\nbige,\DD)$,
and the prolongment $\prolong{\nbige}$.
We put $N_i:=\Res_{\nbigd_i}(\DD)$.

\begin{thm}\label{thm;12.13.50}
The tuple $(N_1,\ldots,N_l)$ is sequentially compatible.
\end{thm}
\pf
We have already known the following:
\begin{itemize}
\item
 The conjugacy classes of $N(\jbar)_{|\,Q}$ is independent of
a choice of $Q\in D_{\jbar}$
(Corollary \ref{cor;12.11.31}).
\item
 Let $\vech=(h_1,\ldots,h_m)$ be an $m$-tuple of integers.
The intersections $\bigcap_{j=1}^m W(\jbar)_{h_j}$
form a vector subbundle of $\prolong{\nbige}_{\nbigd_{\mbar}}$
(Lemma \ref{lem;12.11.30}).
\end{itemize}

Let $Q$ be a point contained in $D_{\mbar}$.
We only have to check that
$(N_{1\,|\,(\lambda,Q)},\ldots,N_{m\,|\,(\lambda,Q)})$
is sequentially compatible.
We have already known that the constantness of the filtration
on the positive cone holds
(Theorem \ref{thm;12.7.200}).

We use an induction on $l$.
We put $C(1,1):=\{(t,t)\in\Delta^2\}\subset \Delta^2$,
and $C(1,1)^{\ast}=C(1,1)-\{O\}$.
By considering the restriction of $\harmonicbundle$
to $C(1,1)^{\ast}\times \Delta^{\ast\,l-2}\times \Delta^{n-l}$,
we obtain that the tuple
$((N_1+N_2)_{|\,(\lambda,Q)},
 N_{3\,|\,(\lambda,Q)},\ldots, N_{l\,|\,(\lambda,Q)})$
is sequentially compatible,
due to our assumption of induction.

\begin{lem}
The tuple
$(N^{(1)}_2,\ldots,N^{(1)}_l)$ is sequentially compatible.
\end{lem}
\pf
Let take a limiting harmonic bundle in one dimensional direction
$(F,\theta\superinfty,h\superinfty)$
of $\harmonicbundle$.
We have the deformed holomorphic bundle
$(\nbigf^{(\infty)},\DD\superinfty)$.
Let pick $a\in\Delta^{\ast}$.
We put
$X_a:=\{a\}\times \Delta^{n-1}\subset \Delta^{n}=X$,
and $D_{a\,i}:=D_i\cap X_a$
for $i=2,\ldots,l$.
We have
$\nbigx_a\subset \nbigx$
and $\nbigd_{a\,i}=\nbigd_i\cap \nbigx_a$.
We have the residues
$N\superinfty_{a\,i}:=\Res_{\nbigd_{a\,i}}(\DD\superinfty_{|\nbigx_a})$
for $i=2,\ldots,l$.
By our assumption of the induction on $l$,
we have already known that
$(N\superinfty_{a\,2},\ldots,N\superinfty_{a\,l})$
is sequentially compatible.

On the other hand, we have already seen that
the tuple of the vector space 
$\prolong{\nbigf^{(\infty)\,\lambda}}_{|(\lambda,Q)}$
the residues
$\Res_{\nbigd_i}(\DD^{(\infty)\,\lambda}_{|\nbigx_a^{\lambda}})$
are isomorphic to 
the tuple of the graded vector space
$\graded^{(1)}$
and the residues $N^{(1)}_i$,
due to Lemma  \ref{lem;12.4.51} and Corollary \ref{cor;12.4.50}.
Thus we obtain our claim.
\hfill\qed

We know that 
$(N_1,\ldots,N_l)$ is 
sequentially compatible in the bottom part,
due to Proposition \ref{prop;12.7.150}.
The proposition \ref{prop;12.7.150}
holds for any tame nilpotent harmonic bundle
with trivial parabolic structure.
Thus we obtain the universal sequentially compatibility
in the bottom part.
Then we obtain the sequentially compatibility
due to Proposition \ref{prop;12.10.80}.
\hfill\qed

\begin{thm} \label{thm;12.10.250}
The tuple $(N_1,\ldots,N_l)$ is strongly sequentially compatible.
\end{thm}
\pf
Let $Q$ be a point contained in $D_{\mbar}$.
We only have to check that
$(N_{1\,|\,(\lambda,Q)},\ldots,N_{m\,|\,(\lambda,Q)})$
is strongly sequentially compatible.
We use an induction on $m$.
We assume that the claim for $m-1$ holds,
and then we will prove that the claim for $m$ holds.

Assume that $Q$ is contained in 
$D_{\mbar}-\bigcup_{h>m}(D_{\mbar}\cap D_h)$.
Then we can pick a point $P\in X-D$ 
such that the filtration $W^{\sankaku}(\mbar)$
of $S(Q,P)$ is a mixed twistor.
By our assumption of induction,
the sequence
$(N_{1\,|\,(\lambda,Q)},
 \ldots, N_{m-1\,|\,(\lambda,Q)})$ is
strongly sequentially compatible.
We have already proved that
the sequence
$(N_{1\,|\,(\lambda,Q)},
 \ldots,N_{m-1\,|\,(\lambda,Q)}, N_{m\,|\,(\lambda,Q)})$
is sequentially compatible
(Theorem \ref{thm;12.13.50}).
Thus we obtain that 
$(N_{1\,|\,(\lambda,Q)},
 \ldots,N_{m\,|\,(\lambda,Q)})$
is strongly sequentially compatible in this case,
due to Proposition \ref{prop;12.10.70}.

Let consider the case that $Q$ is contained in 
$\bigcup_{h>m}(D_{\mbar}\cap D_h)$.
Let pick $\lambda\neq 0$.
By using a normalizing frame,
we obtain that 
$(N_{1\,|\,(\lambda,Q)},
 \ldots,N_{m-1\,|\,(\lambda,Q)}, N_{m\,|\,(\lambda,Q)})$
is strongly sequentially compatible,
if $\lambda\neq 0$.
We put $I=\{i\,|\,D_i\ni Q\}$.
Then we can pick an appropriate $P$
such that $(S(Q,P),W^{\sankaku}(I))$ is mixed twistor.
By using the mixed twistor structure,
we can conclude that
$(N_{1\,|\,(0,Q)},
 \ldots, N_{m\,|\,(0,Q)})$
is also strongly sequentially compatible.

Thus the induction on $m$ can proceed,
namely the proof of Theorem \ref{thm;12.10.250} is completed.
\hfill\qed

\begin{thm}\label{thm;12.15.110}
The tuple $(N_1,\ldots,N_l)$ is of Hodge.
\end{thm}
\pf
For any $\sigma\in\gbigl$, 
the tuple $(N_{\sigma(1)},\ldots,N_{\sigma(l)})$
is strongly sequentially compatible.
Thus we are done.
\hfill\qed

\section{Norm estimate and limiting CVHS in higher dimensional case}
\label{section;12.15.80}

\subsection{Norm estimate}
\label{subsection;12.15.45}

\subsubsection{Preliminary norm estimate} \label{subsubsection;12.7.100}
We put $X=\Delta^n=\{(\zeta_1,\ldots,\zeta_n)\}$,
$D_i=\{\zeta_i=0\}$ and $D=\bigcup_{i=1}^l D_i$ for any $l\leq n$.
Let $\harmonicbundle$ be a tame nilpotent harmonic bundle
with trivial parabolic structure over $X-D$.
We have the deformed holomorphic bundle with $\lambda$-connection
$(\nbige,\DD)$ and the prolongment $\prolong{\nbige}$.
We put $N_i=\Res_{\nbigd_i}(\DD)$.
Then $(N_1,\ldots,N_n)$ gives a strongly sequentially compatible tuple.

Let $\vecv$ be a frame of $\prolong{\nbige}$ over $\nbigx$
which is compatible with $(N_1,\ldots,N_n)$.
We put as follows:
\[
 2\cdot k_j(v_i)=\deg^{W(\jbar)}(v_i)-\deg^{W(\jminusitibar)}(v_i)
\]

We put $\blowup{X}:=\{(z_1,\ldots,z_n)\in\Delta^n\}$,
$\blowup{D}_i:=\{z_i=0\}$,
and $\blowup{D}=\bigcup_{i=1}^l\blowup{D}_i$.
Let $\bigl(c(1),\ldots,c(n)\bigr)$ be a sequence of integers
such that
$0\leq c(1)<c(2)<\ldots<c(n)$.
Let consider the morphism $\phi_N:\blowup{X}\lrarr X$
defined as follows:
\begin{equation}
 \phi_N^{\ast}(\zeta_i)=
 \left\{
\begin{array}{ll}
 {\displaystyle
 \left(
 \prod_{j=i}^l z_j
 \right)^{N^{c(i)}}
 }
 & (i\leq l)\\
\mbox{{}}\\
 z_i & (i\geq l+1).
\end{array}
 \right.
\end{equation}
We put $(\blowup{\nbige},\blowup{\DD},\blowup{h})=
 \phi_N^{\ast}(\nbige,\DD,h)$.
We also put $\blowup{\vecv}=\phi_N^{\ast}\vecv$.
We have the $C^{\infty}$-frame $\blowup{\vecv}'$ defined as follows:
\[
 \blowup{v}_i':=
 \blowup{v}_i\cdot
 \prod_{m=1}^l
 \Bigl(-\sum_{t=m}^l\log|z_t|\Bigr)^{-k_m(v_i)}.
\]
\begin{prop} \label{prop;12.10.301}
The frame $\blowup{\vecv}'$ is adapted
over $\Delta_{\lambda}(R)\times (\blowup{X}-\blowup{D})$
for any $R>0$.
\end{prop}
\pf
We use an induction on $l$.
We assume that the claim for $l-1$ holds
and we will prove the claim for $l$.
The assumption of the induction will be used in Lemma \ref{lem;12.13.60}.

First we note the following:
Let $\vecv_1$ be another frame compatible
with the sequence $(N_1,\ldots,N_l)$.
By the same procedure,
we obtain the $C^{\infty}$-frame
$\blowup{\vecv}_1'$.
\begin{lem} \label{lem;12.10.300}
If $\blowup{\vecv}_1'$ is adapted over
$\Delta_{\lambda}(R)\times (\blowup{X}-\blowup{D})$,
then $\blowup{\vecv}'$ is also adapted over the same region.
\end{lem}
\pf
We have the relation $v_i=\sum_j b_{j\,i}\cdot v_{1\,j}$.
Since $\vecv$ and $\vecv_1$ are compatible with
the sequence $(N_1,\ldots,N_l)$,
we have the following:
\[
 \deg^{W(\mbar)}(v_i)<\deg^{W(\mbar)}(v_{1\,j})
 \Longrightarrow
 b_{j\,i\,|\,\nbigd_{\mbar}}=0.
\]
If we put 
$m_0(i,j)=\min\{m\,|\,k_m(v_i)<k_{m}(v_{1\,j})\}$,
then we obtain that 
$b_{j\,i\,|\,D_{\mzeroijbar}}=0$.
Since we have $\phi_N(\blowup{\nbigd}_t)\subset \nbigd_{\mzeroijbar}$
for $t\geq m_0(i,j)$,
the holomorphic function
$\phi_N^{\ast}b_{j\,i}$ is of the form
$g\cdot \prod_{t=m_0(i,j)}^lz_t$ for some holomorphic function $g$.

We have the following relation:
\[
 \blowup{v}_i'=
 \sum_j \blowup{B}'_{j\,i}\cdot\blowup{v}_{1\,j}',
\quad\quad
 \blowup{B}'_{j\,i}=
 \phi_N^{\ast}(b_{j\,i})
\cdot \prod_{m=1}^l
 \Bigl(-\sum_{t=m}^l \log|z_t|\Bigr)^{-k_m(v_i)+k_m(v_{1\,j})}.
\]
Thus $\blowup{B}'_{j\,i}$ is bounded over
$\Delta_{\lambda}(R)\times (\blowup{X}-\blowup{D})$
for any $R>0$.
Thus $\blowup{B}'$ is bounded over the same region.
Similarly, the boundedness of $\blowup{B}^{\prime\,-1}$
over the same region can be shown.
Thus the proof of Lemma \ref{lem;12.10.300} is completed.
\hfill\qed

Let
$\vecv$ be a frame of $\prolong{\nbige}$ over $\nbigx$,
which is strongly compatible with
the tuple $(N_1,\ldots,N_l)$
(See the subsubsection \ref{subsubsection;12.7.50}
 and Definition \ref{df;12.10.90}).
Then $\vecv$ satisfies the following:
\begin{itemize}
\item
 $\vecv=
   (v_{k,\vech,\eta}\,|\,
  k\geq 0,\,\vech\in \seisuu^{l},\,\,\eta=1,\ldots, d(k,\vech)
 )$.
\item
 $2\cdot k_j(v_{k,\vech,\eta})= h_{j}$.
\item
 We put $\vecdelta=(1,0,\ldots,0)$.
 For an element $\vech$, the first component is denoted by $h_1$.
 Then we have the following:
 \[
  N_1(v_{k,\vech,\eta})=
 \left\{
 \begin{array}{ll}
 v_{k,\vech-2\cdot\vecdelta_1,\eta} &(h_1>-k)\\
 0 &(h_1=-k).
 \end{array}
 \right.
 \]
\end{itemize}
We only have to prove the claim of Proposition \ref{prop;12.10.301}
for such a frame $\vecv$.

As a preliminary consideration,
we see the frame $\blowup{\vecv}^{\circ}$
defined as follows:
\[
  \blowup{v}^{\circ}_{k,\vech,\eta}
:=\blowup{v}_{k,\vech,\eta}\cdot 
 \Bigl(-\sum_{j=2}^l\log|z_j|\Bigr)^{-(h_1+h_2)/2}
\cdot 
\prod_{i=3}^l\Bigl(
 -\sum_{j=i}^l\log|z_j|
\Bigr)^{-h_i/2}.
\]

\begin{lem} \label{lem;12.13.60}
Let $C$ be a real number such that $0<C<1$.
On $\Delta_{\lambda}(R)\times
   \{(z_1,\ldots,z_n)\in \blowup{X}-\blowup{D}\,|\,|z_1|=C\}$,
the frame $\blowup{\vecv}^{\circ}$ is adapted
for any $R>0$.
\end{lem}
\pf
Let take $a\in\cnum$ such that $|a|=C$.
We put $\blowup{X}_a:=\{(a,z_2,\ldots,z_n)\in\blowup{X}\}$
and $\blowup{D}_{a\,i}:=\blowup{D}_i\cap\blowup{X}_a$.
Then $\blowup{D}_a=\blowup{D}\cap X=\bigcup_{i=2}^l\blowup{D}_{a,i}$.
We also put
$X_a:=\{(\zeta_1,\ldots,\zeta_n)\in X\,|\, 
  \zeta_1^{N^{c(2)-c(1)}}=a^{N^{c(2)-c(1)}}\zeta_2\}$
and $D_{a\,i}=D_i\cap X_a$.
Then $D_a=D\cap X=\bigcup_{i=2}^l D_{a,i}$.
We have the isomorphisms:
\[
 \blowup{X}_a\simeq \{(z_2,\ldots,z_n)\in \Delta^{n-1}\},
\quad
 X_a\simeq \{(\zeta_1,\zeta_3,\ldots,\zeta_n)\in \Delta^{n-1}\}.
\]
The restriction of $\phi_N$
gives the morphism $\blowup{X}_a\lrarr X_a$.
Applying the assumption of the induction,
we obtain the result.
\hfill\qed

We have the $\lambda$-connection along the $z_1$-direction
$\gminiq_1(\blowup{\DD})$
(the subsection \ref{subsection;12.14.101}).
We have the $\lambda$-connection form $A\cdot dz_1/z_1$
of $\gminiq_1(\blowup{\DD})$ with respect to the frame $\blowup{\vecv}$.
We have the decomposition of the frame as follows:
\[
 \blowup{\vecv}=
\bigcup_{k\geq 0,\,\,\vecl\in\seisuu^{l-1}}
\Bigl(
 v_{k,\vech,\eta}\,|\,\pi_1(\vech)=\vecl,\,
 \eta=1,\ldots,d(k,\vech).
\Bigr)
\]
Here we put
$\pi(\vech)=(h_2,\ldots,h_n)$ for $\vech=(h_1,h_2,\ldots,h_n)$.
Corresponding to the decomposition of the frame,
we obtain the decomposition of
$\prolong{\nbige}$:
\[
 \prolong{\nbige}
=\bigoplus_{k\geq \,0}
 \bigoplus_{\vecl\in\seisuu^{l-1}}  
 \nbigl_{k,\vecl}.
\]

\begin{lem}\mbox{{}}
\begin{itemize}
\item
The restrictions $A_{|\blowup{\nbigd}_i}$ $(i=1,\ldots,l)$
are constant.
\item
Corresponding to the decomposition of
$\prolong{\nbige}$,
it is decomposed into $A_{k,\vecl}$ for $k\geq 0$,
and $\vecl\in\seisuu^{l-1}$.
\end{itemize}
\end{lem}
\pf
The first claim is clear from our construction.
The second claim corresponds to the fact
that $N_1$ preserves the bundles
$\nbigl_{k,\vecl}$ for $k\geq 0$ and $\vecl\in\seisuu^{l-1}$.
\hfill\qed

\vspace{.1in}

For $k\geq 0$ and $\vecl\in\seisuu^{l-1}$,
we put $V_{k,\vecl}
:=\nbigl_{k,\vecl\,|\,(0,O)}$.
The morphism $N_1$ induces $N_{k,\vecl}$
on $\nbigv_{k,\vecl}$.
We have a harmonic bundle
$E(V_{k,\vecl},N_{k,\vecl})=
(E_{0,k,\vecl},\theta_{0,k,\vecl},h_{0,k,\vecl})$.
We have the deformed holomorphic bundle
$(\nbige_{0,k,\vecl},\DD_{0,k\,\vecl})$.
We have the canonical frame
$\vecv_{0,k,\vecl}=
 (v_{0,k,\vech,\eta}\,|\,\pi_1(\vech)=\vecl,\,\,
 \eta=1,\ldots,d(k,\vech))$
satisfying
$\DD_{0,k,\vecl}\vecv_{0,k,\vecl}
=\vecv_{0,k,\vecl}\cdot A_{k,\vecl}dz/z$.

Let $q_1$ denote the projection $\Delta^{n}\lrarr \Delta$
onto the first component.
We put
$(\blowup{\nbige}_{0,k,\vecl},
  \blowup{\DD}_{0,k,\vecl},
  \blowup{\vecv}_{0,k,\vecl})
=q_1^{\ast}(\nbige_{0,k\,\vecl},\DD_{0,k,\vecl},\vecv_{0,k,\vecl})$.
We also put as follows:
\[
 \blowup{h}_{0,k,\vecl}(z_1,\ldots,z_n)=
 \left[
\prod_{i=2}^l
 \Bigl(-\sum_{j=i}^l\log|z_j|\Bigr)^{l_{i-1}/2}
 \right] 
\cdot
 h_{0,k,\vecl}
\Bigl( \prod_{j=1}^lz_j\Bigr).
\]

Let consider the frame $\blowup{\vecv}^{\star}_{0,k,\vecl}$
defined as follows:
\[
 \blowup{v}^{\star}_{0,k,\vech,\eta}:=
 \blowup{v}_{0,k,\vech,\eta}\cdot
 \Bigl(-\sum_{j=2}^l\log|z_j|
 \Bigr)^{-(h_1+h_2)/2}
 \cdot
 \prod_{i=3}^l
 \Bigl(-\sum_{j=i}^l\log|z_j|
 \Bigr)^{-h_i/2}.
\]
Note that we have $h_i=l_{i-1}$ by definition.
\begin{lem} \label{lem;12.13.70}
Let $C$ be a real number such that $0<C<1$.
The $C^{\infty}$-frame $\vecv^{\star}_{0,k,\vecl}$
is adapted over
$\Delta_{\lambda}(R)\times
 \bigl\{(z_1,\ldots,z_n)\in \blowup{X}-\blowup{D}\,\big|\,|z_1|=C\bigr\}$
for any $R>0$.
\end{lem}
\pf
Clear from our construction.
\hfill\qed

We put as follows:
\[
 (\blowup{\nbige}_0,\blowup{\DD}_0,\blowup{h}_0,\blowup{\vecv}_0)
:=
 \bigoplus_{k\geq 0}\bigoplus_{\vecl\in\seisuu^{n-1}}
 (
 \blowup{\nbige}_{0,k,\vecl},
 \blowup{\DD}_{0,k,\vecl},\blowup{h}_{0,k,\vecl},
 \blowup{\vecv}_{0,k,\vecl}
 ).
\]
By using the frames
$\blowup{\vecv}_0$ and $\blowup{\vecv}$,
we obtain the holomorphic isomorphism
$\Phi:\blowup{\nbige}_0\lrarr\blowup{\nbige}$.
The proof of Proposition \ref{prop;12.10.301} 
is reduced to the following lemma.
\begin{lem} \label{lem;12.7.102}
The morphisms $\Phi$ and $\Phi^{-1}$
are bounded with respect to the metrics $\blowup{h}$
and $\blowup{h}_0$,
over $\Delta_{\lambda}(R)\times (\blowup{X}-\blowup{D})$
for any $R>0$
\end{lem}
\pf
We only have to check the boundedness
on the boundary
$\Delta_{\lambda}(R)\times
\{(z_1,\ldots,z_n)\in \blowup{X}-\blowup{D} \,|\,|z_1|=C\}$,
which we have already seen
(Lemma \ref{lem;12.13.60} and Lemma \ref{lem;12.13.70}).
\hfill\qed

Thus we obtain the adaptedness of the frame $\blowup{\vecv}$.
The induction on $l$ can proceed,
namely the proof of Proposition \ref{prop;12.10.301}
is completed.
\hfill\qed

\subsubsection{Replacement of Notation}

We put $X=\{(z_1,\ldots,z_n)\in\Delta^n\}$,
$D_i=\{z_i=0\}$, and $D=\bigcup_{i=1}^lD_i$ for $l\leq n$.
Let $\harmonicbundle$ be a tame nilpotent harmonic bundle
with trivial parabolic structure over $X-D$.
We have the deformed holomorphic bundle with $\lambda$-connection
$(\nbige,\DD)$ and the prolongment $\prolong{\nbige}$.
We put $N_i=\Res_{\nbigd_i}(\DD)$.
We take a holomorphic frame $\vecv$ of
$\prolong{\nbige}$,
which is compatible with $(N_1,\ldots,N_l)$.
We take the $C^{\infty}$-frame $\vecv'$
defined as follows:
\[
 v'_{k,\vech,\eta}=v_{k,\vech,\eta}\cdot
 \prod_{i=1}^l(-\log|z_i|)^{-h_i/2}.
\]
Then the following theorem is the reformulation of 
Proposition \ref{prop;12.10.301}
when we put $c(i)=i-1$.
\begin{thm} \label{thm;12.7.110}
Let consider the following region:
\[
 Z(id,l,N):=
 \bigl\{(z_1,\ldots,z_n)\,\big|\,|z_{i-1}|^N<|z_{i}|,\,\,i=2,\ldots,l
 \bigr\}.
\]
On $Z(id,l,N)$, the frame $\vecv'$ is adapted.
\end{thm}
\pf
By the morphism $\phi_N$,
we have the isomorphism of the following regions:
\[
\begin{CD}
 \Delta^{\ast\,l}\times \Delta^{n-l}
 @>{\phi_N\,\simeq}>>
 Z(id,l,N).
\end{CD}
\]
Moreover, we have the following replacement:
\[
 -N^{i-1}\sum_{j=i}^l\log|z_j|
 \stackrel{\phi_N}{\longleftrightarrow}
 -\log|z_i|,
\quad\quad
 (i\leq l).
\]
Here we put $c(i)=i-1$.
Thus we are done.
\hfill\qed

\subsubsection{Theorem}
\label{subsubsection;12.15.46}

Let $\sigma$ be an element of $\gbigs_l$,
and we put $I_j=\{\sigma(i)\,|\,i\leq j\}$
for $j=1,\ldots,l$.
We put $\nbigd_{I_j}=\bigcap_{i\in I}\nbigd_i$.
On $\nbigd_{I_j}$,
we put $N(I_j):=\sum_{i\in I_j}N_{i\,|\,\nbigd_{I_j}}$.
We denote the weight filtration of $N(I_j)$
by $W(I_j)$.
Then we obtain the strongly sequentially compatible 
commuting tuple $(N_{\sigma(1)},\ldots,N_{\sigma(l)})$
and the compatible sequence of the filtrations
$(W(I_1),W(I_2),\ldots,W(I_l))$.

Let $\vecv$ be a holomorphic frame of $\prolong{\nbige}$
which is compatible with the sequentially compatible tuple
$(N_{\sigma(1)},\ldots,N_{\sigma(l)})$.
Recall that we have the number
$\deg^{W(I(j))}(v_i)$ for any $v_i$
and any $j=1,\ldots,l$.
We put as follows:
\[
 2\cdot k_j(v_i):=\deg^{W(I(j))}(v_i)-\deg^{W(I(j-1))}(v_i).
\]
We obtain the $C^{\infty}$-frame $\vecv'$
defined as follows:
\[
 v_i':=
 v_i\cdot
 \prod_{j=1}^l(-\log|z_{\sigma(j)}|)^{-k_j(v_i)}.
\]
Recall that $Z(\sigma,l,C)$ denotes the following region:
\[
 Z(\sigma,l,C):=
 \bigl\{(z_1,\ldots,z_n)\in\Delta^{\ast\,n}\,\big|\,
 |z_{\sigma(i-1)}|^C<|z_{\sigma(i)}|,\,\,
 i=2,\ldots,l
 \bigr\}.
\]
Then we obtain the following theorem, which is the norm estimate
in higher dimensional case.
\begin{thm}\label{thm;12.15.111}
Let $C$ be a real number such that $C>1$.
The frame $\vecv'$ is adapted
on $Z(\sigma,l,C)$.
\end{thm}
\pf
Clearly we only have to consider the case $\sigma=id$.
We can choose a natural number $N$ larger than $C$.
And then, we only have to apply Theorem \ref{thm;12.7.110}.
\hfill\qed

\begin{cor}\label{cor;12.11.50}
Let $f$ be a holomorphic section
compatible with the sequence of the filtrations:
\[
 (W(I_1),W(I_2),\ldots,W(I_l)).
\]
We put $2\cdot k_{j}(f):=\deg^{W(I_j)}(f)-\deg^{W(I_{j-1})}(f)$.
Then we obtain the following estimate
on $Z(\sigma,l,C)$:
\[
 0<C_1<
 |f|_{h}\cdot\prod_{j=1}^l (-\log |z_{\sigma(j)}|)^{-k_j(f)}
 <C_2.
\]
Here $C_i$ $(i=1,2)$ denote positive constants.
\hfill\qed
\end{cor}

\subsubsection{Norm estimate for flat sections}

We give a brief argument to obtain a norm estimate
for flat sections.
For simplicity, we consider a tame nilpotent harmonic bundle
$\harmonicbundle$
with trivial parabolic structure
over $\Delta^{\ast\,n}$.
Let $(\nbige^{1},\DD^{1})$
denote the corresponding flat bundle,
that is the deformed holomorphic bundle over
$\{1\}\times \Delta^{\ast\,n}$.

Let $P$ be a point of $\Delta^{\ast\,n}$.
Let $\gamma_i$ denote the loop
defined as follows:
\[
 z_j(\gamma_i(t))=
 \left\{
 \begin{array}{ll}
 z_i(P)\cdot exp(2\pi\sqrt{-1}t) & (j=i)\\
 z_j(P) &(j\neq i).
 \end{array}
 \right.
\]
Then we obtain the monodromy $M(\gamma_i)\in End(\nbige^1_{|P})$.
Since it is unipotent, we obtain the logarithm
$N(\gamma_i)\in End(\nbige^{1}_{|P})$.

Let $\sigma$ be an element of $\gbigs_n$.
Then we obtain the sets
$I_j=\{\sigma(i)\,|\,i\leq j\}$.
We put $N(I_j):=\sum_{i\in I_j}N(\gamma_i)$.
It induces the weight filtration $W(I_j)$.

\begin{lem} \label{lem;12.8.50}
We have the following implication:
\[
 \left\{
 \begin{array}{ll}
N(I_j)\Bigl(
 W(I_k)_{h}\Bigr)\subset W(I_k)_h & (k<j)\\
N(I_j)\Bigl(
 W(I_k)_{h}\Bigr)
 \subset W(I_k)_{h-2}&(k\geq j).
 \end{array}
 \right.
\]
\end{lem}
\pf
This is a consequence of Theorem \ref{thm;12.3.50}
and Theorem \ref{thm;12.7.200}.
\hfill\qed

\vspace{.1in}
We put $\hyperh:=\{\zeta=x+\sqrt{-1}y\in\cnum\,|\,y>0\}$.
Then we have the covering
$\pi:\hyperh^n\lrarr \Delta^{\ast\,n}$
given by $z_i=\exp(2\pi\sqrt{-1}\zeta_i)$.
We put as follows:
\[
 \blowup{Z}(\sigma,n,C,A):=
 \Bigl\{
 (x_1+\sqrt{-1}y_1,\ldots,x_n+\sqrt{-1}y_n)
 \in\hyperh^n
\,\Big|\,
 \begin{array}{ll}
 |x_i|<A, &i=1,\ldots,n\\
 C\cdot y_{i-1}>y_i, &i=2,\ldots,n
\end{array}
 \Bigr\}.
\]
Let take the pull back
$\pi^{\ast}(\nbige^1,\DD^1,h)$.
Let $\blowup{P}$ be a point of $\blowup{Z}(\sigma,n,C,A)$
such that $\pi(\blowup{P})=P$.

Let take a non-zero element $v$ of $\nbige^1_{|P}$.
We have the numbers $h_j:=\deg^{W(I_j)}(v)$.
We take a flat section $f$ such that 
$\pi_{\ast}(f_{|\blowup{P}})=v$.
We have the following estimate of the norm of $f$.
\begin{thm}
There exist positive numbers $C_1$ and $C_2$
such that the following equality holds
on $\blowup{Z}(\sigma,l,C,A)$:
\[
0< C_1\leq 
 |f|_h^2\cdot y_1^{h_1}\times\prod_{i=2}^ny_i^{h_i-h_{i-1}}
 \leq C_2.
\]
\end{thm}
\pf
We give an only indication.
We put as follows:
\[
 g=\exp\Bigl(
 \sum 
 \bigl(\zeta_i-\zeta_i(\blowup{P})\bigr)
 \cdot N(\gamma_i)
 \Bigr)\cdot f.
\]
Then $g$ is a holomorphic section of $\nbige^1$
defined over $\Delta^{\ast}$.
In fact, it is a section of $\prolong{\nbige^1}$
compatible with the filtrations on the divisors.
Thus we obtain the norm estimate for $g$.

By using Lemma \ref{lem;12.8.50},
we can show that the logarithmic order of 
the norm of $g-f$ is lower than $f$.
Thus we obtain the norm estimate of $f$,
by using Corollary \ref{cor;12.11.50}.
\hfill\qed

\subsection{Limiting CVHS}

\label{subsection;12.15.47}

\subsubsection{Taking a limit}

\label{subsubsection;12.14.1}

We put $X=\Delta^n$, $D_i=\{z_i=0\}$, and $D=\bigcup_{i=1}^nD_i$.
For any non-negative integer $m$,
we have the morphism $\psi_{m,\nbar}:X\lrarr X$
defined as follows:
\[
 (z_1,\ldots,z_n)\longmapsto
 (z_1^m,\ldots,z_n^m).
\]
We will consider a limit of harmonic bundle
via the pull backs of $\{\psi_{m\,\nbar}\,|\,m\in\seisuu\}$.

Let $\harmonicbundle$ be a tame nilpotent harmonic bundle
over $X-D$.
We have the deformed holomorphic bundle with the $\lambda$-connection
by $(\nbige,\DD)$
and the prolongment $\prolong{\nbige}$.
Let $\vecv$ be a holomorphic frame of $\prolong{\nbige}$,
compatible with the filtration
$W(\nbar)$ on $\nbigd_{\nbar}$.
We put $2\cdot k(v_i)=\deg^{W(\nbar)}(v_i)$.

Then we obtain the harmonic bundles
$\psi_{m,\nbar}^{\ast}\harmonicbundle$ on $X-D$.
We also obtain the deformed holomorphic bundles
and the $\lambda$-connections
$(\psi_{m,\nbar}^{\ast}\prolong{\nbige},\psi_{m,\nbar}^{\ast}\DD)$.
We obtain the holomorphic frame $\vecv^{(m)}$
of $\psi_{m,\nbar}^{\ast}\prolong{\nbige}$,
defined as follows:
\[
 v^{(m)}_i:=
 \psi_{m,\nbar}^{-1}(v_i)\cdot m^{-k(v_i)}.
\]

We put $H^{(m)}:=H(\psi_{m,\nbar}^{\ast}(h),\vecv^{(m)})$,
which is an $\nbigh(r)$-valued function.
\begin{lem}
Let $K$ be a compact subset of $\nbigx-\nbigd$.
On $K$, the functions $H^{(m)}$ and $H^{(m)\,-1}$
are bounded independently of $m$.
Moreover, we have the following estimate for
sufficiently large $M$ and for a positive constant $C$,
independently of $m$:
\begin{equation}\label{eq;12.7.1}
 |v_i^{(m)}|_{\psi_m^{\ast}(h)}
\leq C\cdot \prod_{i=1}^n(-\log|z_i|)^M
\end{equation}
\end{lem}
\pf
We only have to consider the boundedness
on a compact subset $K\subset \cnum_{\lambda}\times Z(\sigma,n,N)$
when $\sigma=id\in\gbigs_n$.
Note that $Z(id,n,N)$ is stable under $\psi_{m,\nbar}$.

Let $\vecw$ be a holomorphic frame of $\prolong{\nbige}$
over $\nbigx$,
which is compatible with the sequence of
the filtrations $(N(\itibar),N(\nibar),\ldots,N(\nbar))$.
We have the $C^{\infty}$-frame $\vecw'$ over $\nbigx-\nbigd$,
defined as follows:
\[
 w_i':=w_i\cdot
 \prod_{j=1}^n(-\log|z_j|)^{-k_j(w_i)}.
\]
Here we put
$2\cdot k_j(w_i):=
   \deg^{W(\jbar)}(w_i)-\deg^{W(\jminusitibar)}(w_i)$.
As we have already seen, the frame $\vecw'$ is adapted
over $\Delta_{\lambda}(R)\times Z(id,n,N)$, namely,
$H(h,\vecw)$ and $H(h,\vecw)^{-1}$ are bounded
over $\Delta_{\lambda}(R)\times Z(id,n,N)$ for any $R>0$.
It implies the boundedness of the functions
$\bigl\{H(\psi_{m,\nbar}^{\ast}h,\psi_{m,\nbar}\vecw')\bigr\}$
and
$\bigl\{H(\psi_{m,\nbar}^{\ast}h,\psi_{m,\nbar}\vecw')^{-1}\bigr\}$.

We put $2\cdot k(w_i)=\deg^{W(\nbar)}(w_i)$.
Then we have the following equality:
\[
 \sum_{j=1}^n k_j(w_i)=k(w_i).
\]
We have the following relation:
\[
 \psi_{m,\nbar}^{-1}w_i'=
 \psi_{m,\nbar}^{-1}(w_i)\cdot
 m^{-k(w_i)}\cdot 
 \prod_{j=1}^n(-\log|z_j|)^{-k_j(w_i)}
=w_i^{(m)}\cdot \prod_{j=1}^n(-\log|z_j|)^{-k_j(w_i)}
\]
Here we put $w_i^{(m)}:=\psi_{m,\nbar}^{-1}(w_i)\cdot  n^{-k(w_i)}$.
Let $L$ denote the diagonal matrix whose $(i,i)$-component
is $\prod_{j=1}^n(-\log|z_j|)^{k_j(w_i)} $.
Then we obtain the following relation:
\[
 H(\psi_{m,\nbar}^{\ast}h,\vecw^{(m)})=
 L\cdot
 H(\psi_{m,\nbar}^{\ast}h,\psi_{m,\nbar}\vecw')
\cdot L.
\]
Thus we obtain the boundedness
of $\{H(\psi_{m,\nbar}^{\ast}h,\vecw^{(m)}) \}$
over a compact subset $K\subset \cnum_{\lambda}\times Z(id,n,N)$,
independently of $m$.
Similarly we obtain the boundedness
of $\{H(\psi_{m,\nbar}^{\ast}h,\vecw^{(m)})^{-1} \}$
over a compact subset $K\subset \cnum_{\lambda}\times Z(id,n,N)$,
independently of $m$.

We have the relation $\vecv=\vecw\cdot B$, namely
we have the following equalities:
\[
 v_j=\sum_{i} b_{i\,j}\cdot w_i.
\]
Note that $b_{i\,j}(\lambda,O)=0$ if $k(v_j)<k(w_i)$,
for both of $\vecv$ and $\vecw$ are compatible
with the filtration $W(\nbar)$.
We have the relation $\vecv^{(m)}=\vecw^{(m)}\cdot B^{(m)}$.
The component of $B^{(m)}$ is determined by the following:
\[
 v_j^{(m)}
=\sum_{i} \psi_{m,\nbar}^{\ast}(b_{i\,j})\cdot
 m^{-k(v_j)+k(w_i)}\cdot w_i^{(m)}.
\]
Namely we have
$B^{(m)}_{i\,j}=
\psi_{m,\nbar}^{\ast}(b_{i\,j})\cdot
 m^{-k(v_j)+k(w_i)}$.
It implies the boundedness of $B^{(m)}$
independently of $m$.
Similarly we obtain the boundedness of $B^{(m)\,-1}$.
Thus we obtain the boundedness of $H^{(m)}$
and $H^{(m)\,-1}$
over any compact subset $K$,
independently of $m$.
The estimates (\ref{eq;12.7.1}) are also obtained.
\hfill\qed

Let consider $\Omega(\vecv)=v_1\wedge\cdots\wedge v_r$
and $\Omega(\vecv^{(m)})=v_1^{(m)}\wedge\cdots\wedge v_r^{(m)}$.
\begin{lem}
Let $R$ be any positive number.
There exist positive constants $C_1$ and $C_2$
satisfying the following inequality over
$\{(\lambda,P)\in\nbigx-\nbigd\,|\,|\lambda|<R\}$,
for any $m$:
\[
 0<C_1<|\Omega(\vecv^{(m)})|_{\psi_{m,\nbar}^{\ast}(h)}<C_2.
\]
\end{lem}
\pf
Since $\Omega(\vecv)$ gives a holomorphic frame
of $\det(\prolong{\nbige})$,
we have the following inequality for some $0<C_1<C_2$:
\[
 0<C_1<|\Omega(\vecv)|_{h}<C_2.
\]
Then it is clear we obtain the result.
\hfill\qed

We put $U_h:=\langle v_i\,|\, \deg^{W(\nbar)}(v_i)=h\rangle$.
We have the $\lambda$-connection form 
$\nbiga$ of $\DD$ with respect to the frame $\vecv$.
We decompose $\nbiga$ into $\sum_j \nbiga_j=\sum A_j\cdot dz_j/z_j$.
We also decompose $A_j$ into $A_{j\,(h,k)}$
corresponding to the decomposition $\prolong{\nbige}=\bigoplus U_h$
(see the subsubsection \ref{subsubsection;12.5.10}).
By our choice of $U_h$,
we have the vanishing
$A_{j\,(h,k)}(\lambda,O)=0$ when $h>k-1$.

Let $\nbiga^{(m)}=\sum \nbiga^{(m)}_j$ denote the $\lambda$-connection form
of $\psi_{m,\nbar}^{-1}\DD$
with respect to the frame $\vecv^{(m)}$.
\begin{lem} \label{lem;12.7.51}
We have the following equalities:
\[
 \nbiga^{(m)}_j=
 \sum_{h,k}\psi_{m,\nbar}^{\ast}A_{j\,(h,k)}\cdot m^{(h-k+2)/2}\cdot
 \frac{dz_j}{z_j}.
\]
In particular, the sequences $\{\nbiga^{(m)}\}$ converges
to $\nbiga^{(\infty)}=\sum_j\nbiga_j^{(\infty)}$ given as follows:
\[
 \nbiga_j^{(\infty)}=
 \sum_h A_{j\,(h,h+2)}(\lambda,O)\cdot \frac{dz_j}{z_j}
\]
\end{lem}
\pf
Similar to Lemma \ref{lem;12.4.51}.
\hfill\qed

We decompose $\theta$ into $\sum f_j\cdot dz_j/z_j$.
Then $\psi_m^{\ast}\theta
=\sum_j m\cdot \psi_{m,\nbar}^{-1}(f_j)\cdot dz_j/z_j$.
Since we have the estimate
$|f_j|_h\leq C\cdot (-\log|z_j|)^{-1}$,
we obtain the following estimate independently of $m$:
\[
 |m\cdot \psi_{m,\nbar}^{\ast}(f_j)|_{\psi_m^{\ast}(h)}
\leq C\cdot(-\log|z_j|)^{-1}.
\]

In all we obtain the following.
\begin{lem}
The sequence $\{\vecv^{(m)}\}$
satisfies Condition {\rm \ref{condition;11.27.16}}.
\hfill\qed
\end{lem}

We can apply the result in subsection \ref{subsection;12.4.10},
as in the subsubsection \ref{subsubsection;12.14.2}.
Then we can pick a subsequence $\{m_i\}$ of $\{m\}$
and the limiting harmonic bundle $(F,\theta\superinfty,h\superinfty)$
for $\{m_i\}$.
We also obtain the deformed holomorphic bundle
with $\lambda$-connection
$(\nbigf^{(\infty)},\DD^{(\infty)})$
and the holomorphic frame $\vecv\superinfty$.
Due to the estimates (\ref{eq;12.7.1}),
we obtain the estimates
$|v_i\superinfty|_{h\superinfty}<C\cdot \prod_{j=1}^n(-\log|z_j|)^M$
for sufficiently large $M$.
We also obtain the estimate
$0<C_1<|\Omega(\vecv^{(\infty)})|_{h\superinfty}$.
Thus $\vecv\superinfty$ naturally gives
a holomorphic frame of the prolongment $\prolong{\nbigf^{(\infty)}}$.

We have the decomposition
$\nbigf\superinfty=\bigoplus_h\directsummand_h$,
where $\directsummand_h$ denotes the vector subbundle
generated by $\{v_i\superinfty\,|\,\deg^{W(\nbar)}v_i=h\}$.

\begin{lem}
The subbundles $\directsummand_h$ are independent of 
a choice of the original frame compatible with $W(\nbar)$.
\end{lem}
\pf
Similar to Lemma \ref{lem;12.7.10}.
\hfill\qed

We have the $\lambda$-connection form
$\nbiga^{(\infty)}=\sum A_j^{(\infty)}dz_j/z_j$
of $\DD^{(\infty)}$ with respect to the frame $\vecv\superinfty$.
Let $f_{A_j\superinfty}$ denote the endomorphism
determined by $A_j\superinfty$ and $\vecv\superinfty$.
\begin{lem}
We have $f_{A_j\superinfty}(\directsummand_h)\subset \directsummand_{h-2}$.
In particular,
we obtain
$\theta(\directsummand_h)
\subset\directsummand_{h-2}\otimes\Omega^{1,0}_{X-D}$.
\end{lem}
\pf
The claim immediately follows from Lemma \ref{lem;12.7.51}.
\hfill\qed

We have the $\nbigh(r)$-valued function $H(h\superinfty,v\superinfty)$.
We use the real coordinate $z_i=r_i\cdot \exp(2\pi\sqrt{-1}\alpha_i)$
for $i=1,\ldots,n$.
\begin{lem}
The function $H(h\superinfty,v\superinfty)$
is independent of $\alpha_i$ for any $i$.
\end{lem}
\pf
Similar to Proposition \ref{prop;12.5.60}.
\hfill\qed

We put $\directsummand^{\lambda}
:=\directsummand_{|\nbigx^{\lambda}-\nbigd^{\lambda}}$.
\begin{thm} \label{thm;12.7.20}
If $h\neq h'$,
then $\directsummand^0_{h}$ and $\directsummand^0_{h'}$
are orthogonal.
\end{thm}
\pf
We put
$\Delta_{\real,+}^{\ast,n}=\real_{>0}^n\cap \Delta^{\ast\,n}$.
We have the following:
\[
 H(h\superinfty,\vecv\superinfty)
 (r_1\cdot \exp(2\pi\sqrt{-1}\alpha_1),\ldots
 r_n\cdot \exp(2\pi\sqrt{-1}\alpha_n) )
=H(h\superinfty,\vecv\superinfty)
 (r_1
 \ldots
 r_n).
\]
Thus we only have to check the orthogonality
over $\Delta_{\real,+}^{\ast\,n}$.

Let $(r_1,\ldots,r_n)$ be an element of $\Delta_{\real,+}^{\ast\,n}$.
Then we have the real numbers $\alpha_j>0$ satisfying
$r_j=r_1^{\alpha_j}$.
The real numbers $\alpha_j$ can be approximated by rational numbers.
For any $\vech\in\seisuu_{>0}^n$,
we put $C_{\vech,\real,+}:=
 \{(t^{h_1},\ldots,t^{h_n})\,|\,0<t<1\}
\subset \Delta_{\real,+}^{\ast\,n}$.
Then we know the following set is dense in $\Delta_{\real,+}^{\ast\,n}$:
\[
 \bigcup_{\vech\in\seisuu_{>0}^n}
 C_{\vech,\real,+}.
\]
Thus we only have to check the orthogonality
on $C_{\vech,\real,+}$.

We put $C_{\vech}:=\{(t_1^{h_1},\ldots,t_n^{h_n})\in\Delta^{\ast}\}$.
Since we have $C_{\vech,\real,+}\subset C_{\vech}$,
we only have to check the orthogonality
on $C_{\vech}$.
Note that $C_{\vech}$ is the image of the morphism
$g_{\vech}:\Delta^{\ast}\lrarr X-D$,
and that we have the following commutative diagramm:
\[
 \begin{CD}
 \Delta^{\ast} @ >{\psi_{m,\itibar}}>> \Delta^{\ast}\\
 @V{g_{\vech}}VV @V{g_{\vech}}VV \\
 X-D @>{\psi_{m,\nbar}}>> X-D.
 \end{CD}
\]
Thus the theorem \ref{thm;12.7.20} is reduced to the theorem
\ref{thm;12.4.100}.
\hfill\qed

As a direct corollary, we obtain the following theorem.
\begin{thm}\label{thm;12.15.120}
The tuple $(F,\theta^{(\infty)},h^{(\infty)})$
gives a complex variation of polarized Hodge structure,
up to grading.
\end{thm}
\pf
Similar to Corollary \ref{cor;12.7.21}.
\hfill\qed

\begin{df}
The tuple $(F,\theta\superinfty,h\superinfty)$
is called a limiting CVHS of $(E,\theta,h)$.
\hfill\qed
\end{df}

\label{subsubsection;12.11.60}

\subsubsection{Real structure}

Let consider the real structure of the harmonic bundles.
Let $\harmonicbundle$ is a harmonic bundle
over a complex manifold $X$.
\begin{df} \label{df;12.14.10}
Let $\iota:E\lrarr E$ be an anti-linear isomorphism.
We say that $\iota$ is a real structure of
$\harmonicbundle$, if the following holds:
\begin{itemize}
\item $\iota^2$ is the identity map.
\item
 $\iota$ preserves the metric $h$.
\item
 $\iota$ replaces $\del_E$ and $\delbar_E$.
 Namely we have $\iota(\del_Ef)=\delbar_E\iota(f)$
 and $\iota(\delbar_Ef)=\del_E\iota(f)$.
\item
 $\iota$ replaces $\theta$ and $\theta^{\dagger}$.
\hfill\qed
\end{itemize}
\end{df}

By the natural anti-linear isomorphism $E\simeq E^{\lor}$
induced by the hermitian metric $h$,
we can also regard $\iota$ as the holomorphic isomorphism
$E\lrarr E^{\lor}$.

\begin{rem}
Recall that the harmonic bundle is regarded as a variation
of pure polarized twistor structure, according to Simpson
({\rm\cite{s3}}).
Although the real structure of harmonic bundle should be seen
from such view point,
we do not discuss such issue.
\hfill\qed
\end{rem}

We put $X=\Delta^n$, $D_i=\{z_i=1\}$ and $D=\bigcup_{i=1}^n D_i$.
Let $\harmonicbundle$ be a tame nilpotent harmonic bundle
with trivial parabolic structure over $X-D$.
Assume that $\harmonicbundle$ has a real structure $\iota$.
Let consider the sequence
$\bigl\{\psi_{m,\nbar}^{\ast}\harmonicbundle\bigr\}$
in the subsubsection \ref{subsubsection;12.14.1}.
As in \ref{subsubsection;12.14.2},
Let $F=\bigoplus_{i=1}^r\nbigo_{X-D}\cdot u_i$ be a holomorphic bundle
with the frame $\vecu=(u_i)$,
over $X-D$.
We put $\vece^{(m)}:=\vecv^{(m)}_{|\nbigx^0}$.
It is the frame of $\prolong{\psi_{m,\itibar}^{\ast}\nbige^{0}}$
over $\nbigx^0=\{0\}\times X$.
The frames $\vece^{(m)}$ and $\vecu$ give the holomorphic isomorphism
$\Phi_m:\nbige^0\lrarr F$
over $X-D$.
The morphism $\Phi_m$ induces the structure of harmonic bundle
on $F$.
Moreover we have the sequence of real structure
$\iota^{(m)}$,
which is the image of $\psi_{m,\nbar}^{\ast}(\iota)$
via the morphism $\Phi_m$.

Note that $|\psi_{m,\nbar}(\iota)|_{\psi_{m,\nbar}^{\ast}(h)}=1$,
and they are holomorphic as morphisms
$\psi_{m,\nbar}^{\ast}E\lrarr \psi_{m,\nbar}^{\ast}E^{\lor}$.
Thus we can pick a subsequence $\{m_i\}$ for a limiting CVHS
$(F,\theta\superinfty,h\superinfty)$
such that the sequence $\{\iota^{(m_i)}\}$ also converges.
We denote the limit by $\iota\superinfty$.

\begin{lem}
$\iota\superinfty$ gives a real structure
to the limiting CVHS $(F,\theta\superinfty,h\superinfty)$.
\end{lem}
\pf
It is clear from our construction that
the conditions in \ref{df;12.14.10} are satisfied.
Thus $\iota\superinfty$ gives a real structure of a harmonic bundle.

Thus we only have to check
$\iota\superinfty(\directsummand^0_h)=\directsummand^0_{-h}$.
Let consider the filtration
$\nbigw$ and $\nbigw^{\dagger}$ given as follows:
\[
 \nbigw_h:=\bigoplus_{l\leq h}\directsummand^0_l,
\quad
 \nbigw_h^{\dagger}:=\bigoplus_{l\geq -h}\directsummand^0_l.
\]
From $\iota(\theta\superinfty)=\theta^{(\infty)\,\dagger}$,
it is easy to see 
$\iota(\nbigw_h)=\nbigw_{-h}^{\dagger}$.
Then we obtain $\iota(\directsummand_h^0)=\directsummand_{-h}^0$,
by using the orthogonality of them and
an easy ascending induction on $h$.
\hfill\qed

\begin{df}
The tuple $(F,\theta\superinfty,h\superinfty,\iota\superinfty)$
is called a limiting RVHS.
\hfill\qed
\end{df}

\subsubsection{Purity theorem}
\label{subsubsection;12.15.50}

Recall the setting of the purity theorem, following Kashiwara and Kawai
(See \cite{k4} and \cite{k3}. But our filtration
is slightly different from theirs. See remark \ref{rem;12.14.15}.)
Let $V$ be a finite dimensional vector space
and $(N_1,\ldots,N_n)$ be a commuting tuple of nilpotent maps.
We put $N(\jbar)=\sum_{i\leq j}N_i$.
We denote the weight filtration of $N(\nbar)$
by $W(\nbar)$.
Let $(e_1,\ldots,e_n)$ be the standard base of $\seisuu^n$.
We have the partial Koszul complex $\partialkoszul$.
The $k$-th part $\partialkoszul^k$ is defined as follows:
\[
 \partialkoszul^k=
 \bigoplus_{\substack{|J|=k,\\J\subset\nbar}}
 \Image N_J\otimes (\cnum \cdot e_J).
\]
Here we put $N_J=\prod_{j\in J}N_j$
and $e_J=\bigwedge_{j\in J}e_j$.
The differential $d$ is given by $\sum_{j=1}^n N_j\wedge e_j$.
Namely we put as follows:
\[
 d(v\otimes e_J):=
 \sum_{j=1}^n N_j(v)\otimes (e_j\wedge e_J).
\]
The filtration $W$ of $\Image(N_j)$ is given as follows:
\[
 W_h(\Image(N_J))=N_J(W_h(\nbar))\subset \Image(N_J).
\]
Then we obtain the filtration $W$
of $\partialkoszul^k$ as follows:
\[
 W_h(\partialkoszul^k):=
 \bigoplus_{\substack{|J|=k,\\J\subset\nbar}}
 W_h(\Image(N_J))\otimes (\cnum\cdot e_J).
\]

We put $X=\Delta^n$, $D_i=\{z_i=0\}$
and $D=\bigcup_{i=1}^n D_i$.
Let $\harmonicbundle$ be a tame nilpotent harmonic bundle
with trivial parabolic structure
over $X-D$.
We put $V=\prolong{\nbige}_{|(\lambda,O)}$.
Then we have the nilpotent maps
$N_i:=\Res_{\nbigd_i}(\DD)_{|(\lambda,O)}$
for $i=1,\ldots,n$.

\begin{lem} \label{lem;12.7.35}
In this case, we have the following:
\[
 W_k(\Image(N_J))=\Image(N_J)\cap W_{k-2|J|}.
\]
\end{lem}
\pf
Let take a point $P\in X-D$, and then 
the filtration $W^{\sankaku}(\nbar)$ gives a mixed twistor structure
to the vector bundle $S(O,P)$.
Then Lemma  \ref{lem;12.7.35} is a consequence of 
Lemma \ref{lem;12.7.30}.
\hfill\qed

We have the induced filtration on the cohomology group
$H^{\ast}(\partialkoszul)$.
\begin{thm} \label{thm;12.7.40}
Assume that $\harmonicbundle$ has a real structure.
Then the purity theorem holds for the tuple $(V,N_1,\ldots,N_n)$.
Namely we have the following:
\[
 H^k(\partialkoszul)=
 W_k(H^k(\partialkoszul)).
\]
In other words, the following morphism is surjective.
\[
 W_k(\partialkoszul)^k\cap\ker(d)\lrarr H^k(\partialkoszul).
\]
\end{thm}
\pf
Let $Gr_{\cdot}(\partialkoszul)$
denote the associated graded complex of the filtered complex
$\partialkoszul$.
Then we obtain the following natural isomorphism:
\[
 Gr_h(\partialkoszul)^a
:=\bigoplus_{\substack{|J|=a,\\J\subset\nbar}}
 \frac{N_J(W_h)}{N_J(W_{h-1})}
 \otimes (\cnum\cdot e_J)
 \simeq
 \bigoplus_{\substack{|J|=a,\\J\subset\nbar}}
 \frac{N_J(W_h)}{N_J(W_h)\cap W_{h-1-2a}}
\otimes(\cnum\cdot e_J).
\]
\begin{lem} \label{lem;12.13.102}
To show Theorem {\rm\ref{thm;12.7.40}},
we only have to show the following:
\[
 H^a(Gr_k(\partialkoszul))=0,
\quad
\mbox{if } a<k.
\]
\end{lem}
\pf
It is shown by an elementary homological algebraic argument.
\hfill\qed

We put $V\superinfty_h=Gr^{W(\nbar)}_h$,
and $V\superinfty=\bigoplus_h V\superinfty_h$.
The morphism $N_i$ induces
$N\superinfty_{i\,h}:V_h\superinfty\lrarr V_{h-2}\superinfty$
for any $h$.
We denote the direct sum $\bigoplus_h N\superinfty_{i\,h}$
by $N\superinfty_i$.
Then we obtain the commuting tuple
$(N\superinfty_1,\ldots,N\superinfty_n)$.

Let consider a limiting RVHS $(F,\theta\superinfty,h\superinfty)$
of $(E,\theta,h)$
obtained in the subsubsection \ref{subsubsection;12.11.60}.
Let $\nbigf\superinfty$ denote the deformed holomorphic bundle.
Due to Lemma \ref{lem;12.7.51}
and the fact that $\vecv\superinfty$ gives a frame of
$\prolong{\nbigf\superinfty}$,
the tuple $(V\superinfty,N_1\superinfty,\ldots,N_n\superinfty)$
is isomorphic to the following tuple:
\[
 (\prolong{\nbigf\superinfty}_{|(\lambda,O)},
 \Res_{\nbigd_1}(\DD\superinfty)_{|(\lambda,O)},
 \ldots
 \Res_{\nbigd_n}(\DD\superinfty)_{|(\lambda,O)},
).
\]
For a RVHS, the purity theorem holds,
due to the result  of Cattani-Kaplan-Schmid
\cite{cks2} or Kashiwara-Kawai \cite{k3}.

Note that the complex $\partialkoszulinfty$ is graded.
The grading is induced by the grading of $V\superinfty$.
Namely we have the following decomposition:
\begin{equation} \label{eq;12.10.100}
 \partialkoszulinfty^a
=\bigoplus_h
\Bigl(
 \bigoplus_{|J|=a,J\subset\nbar}
 N_J\superinfty(V\superinfty_h)
\Bigr).
\end{equation}
It is easy to see that the differential preserves the grading.
Thus we obtain the following complex:
\[
   \partialkoszulinfty_h^a=
 \bigoplus_{|J|=a,J\subset\nbar}
 N_J\superinfty(V\superinfty_h),
\quad
 d=\sum_{j=1}^n N_j\superinfty\wedge e_j.
\]
Since the decomposition of $V\superinfty$
gives the splitting of the filtration of $V\superinfty$,
the decomposition (\ref{eq;12.10.100})
gives the splitting of the filtration of $\partialkoszulinfty$.
Namely we have the following:
\[
 W_k(\partialkoszulinfty)^a=
\bigoplus_{h\leq k}
\partialkoszulinfty_h^a.
\]
Then the claim of the purity theorem
for the tuple $(V\superinfty,N_1\superinfty,\ldots,N_n\superinfty)$
implies the following vanishing:
\begin{equation} \label{eq;12.13.100}
 H^a(\partialkoszulinfty_k)=0,
\quad
 \mbox{if }a<k.
\end{equation}
By our construction,
we have the following:
\[
 \partialkoszulinfty_k^a
\simeq
 \bigoplus_{|J|=a,J\subset\nbar}
 \frac{N_J(W_k)}{N_J(W_k)\cap W_{k-1-2a}}.
\]
Namely we have the isomorphism of the complexes:
\begin{equation}\label{eq;12.13.101}
 Gr_k(\partialkoszul)
\simeq
 \partialkoszulinfty_k.
\end{equation}
Due to (\ref{eq;12.13.100}), (\ref{eq;12.13.101})
and Lemma \ref{lem;12.13.102},
we obtain Theorem \ref{thm;12.7.40}.
\hfill\qed

\begin{rem}\label{rem;12.14.15}
We should remark the relations between
the weight filtration $W$ of $\partialkoszul$
considered here,
$W^{cks}$ given by Cattani-Kaplan-Schmid,
and $W^{kk}$ given by Kashiwara-Kawai.
The weight filtration $W^{cks}$ of Cattani-Kaplan-Schmid 
is as follows:
\[
 W^{cks}_h(\partialkoszul^a)
=W_{a+h}(\partialkoszul^a).
\]
The weight filtration $W^{kk}$ of $\partialkoszul^a$
for the RVHS of weight $w$ is given as follows:
\[
 W^{kk}_{h}(\partialkoszul^a)
=W_{h-w}(\partialkoszul^a).
\]
The results can be stated as
$W^{cks}_{0}H^h=W^{kk}_{h+w}H^h=W_{h}H^h=H^h$.
\hfill\qed
\end{rem}

\begin{rem}
We should remark that we do not obtain another proof
of the purity theorem for RVHS.
We just reduced the purity theorem for tame nilpotent harmonic bundles
with trivial parabolic structure 
to the purity theorem for RVHS.
\hfill\qed
\end{rem}

\section{Appendix}

We recall the definition of complex variation of polarized Hodge
structure (CVHS),
and the relation of CVHS with a harmonic bundle.
We also recall the real variation of polarized Hodge structure
(RVHS).

Let $X$ be a complex manifold.
Recall the definition of complex variation of polarized Hodge
structures (See section 8 of \cite{s1}, for example).
\begin{df}
Let $V$ be a $C^{\infty}$-vector bundle over $X$
with a decomposition $V=\bigoplus_{p+q=w}V^{p,q}$.
Let $\DD^1$ be a flat connection of $V$,
and $\langle\cdot,\cdot\rangle$ be a sesqui-linear form
on $V$.
The tuple
$\Bigl(\bigoplus_{p+q=w} V^{p,q},\DD^1,\langle\cdot,\cdot\rangle\Bigr)$
is called a complex variation of Hodge structure of weight $w$,
if the following conditions are satisfied:
\begin{itemize}
\item
 We have the following implication:
 \begin{multline}
 \DD^1
 \Bigl(
  C^{\infty}(X,V^{p,q})
 \Bigr)
\subset
  C^{\infty}\Bigl(X,
 V^{p+1,q-1}\otimes\Omega^{0,1}\Bigr)
 \oplus
 C^{\infty}\Bigl(X,
 V^{p,q}\otimes\Omega^{0,1}
 \Bigr)\\
 \oplus
 C^{\infty}\Bigl(X,
 V^{p,q}\otimes\Omega^{1,0}
 \Bigr)
 \oplus
 C^{\infty}\Bigl(X,
 V^{p-1,q+1}\otimes\Omega^{1,0}
 \Bigr)
 \end{multline}
\item
 $\langle u,v\rangle=(-1)^w\overline{\langle v,u\rangle}$.
\item
 $V^{p,q}$ and $V^{p',q'}$ are orthogonal
 if $(p,q)$ and $(p',q')$ are different.
\item
 We denote the restriction of $\langle\cdot,\cdot\rangle$
by $\langle\cdot,\cdot\rangle_{(p,q)}$.
 Then $(\sqrt{-1})^{p-q}\langle\cdot,\cdot\rangle_{(p,q)}$
 is positive definite on $V^{p,q}$.
\end{itemize}
A complex variation of Hodge structure is called
CVHS in this paper, for simplicity.
\hfill\qed
\end{df}

Let
$\Bigl(\bigoplus_{p+q=w} V^{p,q},\DD^1,\langle\cdot,\cdot\rangle\Bigr)$
be a CVHS on $X$.
We obtain the metric $h$ given as follows:
\[
 h=\bigoplus_{p+q=w}
 (\sqrt{-1})^{p-q}\langle\cdot,\cdot\rangle_{(p,q)}.
\]
We have the decomposition of $\DD^1$ corresponding
to the decomposition
$V\otimes\Omega^{1}=
 \bigoplus V^{p,q}\otimes (\Omega^{0,1}\oplus \Omega^{1,0})$:
\[
 \DD^1=\theta^{\dagger}+\delbar_V+\del_V+\theta.
\]

The following proposition can be shown by a direct calculation.
\begin{prop}
A tuple $(V,\delbar_V,\theta,h)$ is a harmonic bundle.
The $(1,0)$-part of the unitary metric
associated with $\delbar_V$ and $h$ is given by $\del_V$,
and the adjoint of $\theta$ is $\theta^{\dagger}$.
\hfill\qed
\end{prop}

On the other hand,
we obtain a CVHS from a harmonic bundle
with a nice grading.
Let consider the harmonic bundle
$\harmonicbundle$, with a holomorphic decomposition
$E=\bigoplus_l E_l$.
Assume that the following:
\begin{itemize}
\item $\theta(E_l)\subset E_{l-1}\otimes \Omega^{1,0}$.
\item If $l\neq l'$, then $E_{l}$ and $E_{l'}$ are orthogonal
 with respect to the metric $h$.
\end{itemize}
Let $w$ be an integer.
For example,
we can put $V^{p,q}=E_{p}$ for a pair $(p,q)$ such that $p+q=w$.
We have the flat connection $\DD^1$ of
$C^{\infty}$-bundle $V=\bigoplus V^{p,q}=E$.
By reversing the construction above,
we obtain the flat sesqui linear form
$\langle\cdot,\cdot\rangle$ on $V$.
Then the tuple
$\Bigl(\bigoplus V^{p,q},\DD^1,\langle\cdot,\cdot\rangle\Bigr)$
gives a CVHS.
We do not have to care
a choice of the weight $w$ and the way of superscripts $(p,q)$.

Let consider the harmonic bundle
$\harmonicbundle$, with a holomorphic decomposition
$E=\bigoplus_l E_l$.
Assume that the following:
\begin{itemize}
\item $\theta(E_l)\subset E_{l-2}\otimes \Omega^{1,0}$.
\item If $l\neq l'$, then $E_{l}$ and $E_{l'}$ are orthogonal
 with respect to the metric $h$.
\end{itemize}
In this case, we can put, for example, as follows:
For a pair $(p,q)\in\seisuu^2$ such that $p+q=0$,
we put $V_{0}^{p,q}=E_{p-q}$.
For a pair $(p,q)\in\seisuu^2$ such that $p+q=1$,
we put $V_{1}^{p,q}=E_{p-q}$.
We put $V_0=\bigoplus_{p+q=0}V^{p,q}$
and $V_1=\bigoplus_{p+q=1}V^{p,q}$.
Then $E$ is decomposed into $V_0\oplus V_1$.
We have the induced structures of harmonic bundles
on $V_0$ and $V_1$.
They are complex variations of polarized Hodge structures
weight $0$ and $1$ respectively.

\begin{df}
Let
 $\Bigl(\bigoplus_{p+q=w}V^{p,q},\DD^1,\langle\cdot,\cdot\rangle\Bigr)$
be a CVHS on $X$.
A real structure $\iota$ of
$\Bigl(\bigoplus_{p+q=w}V^{p,q},\DD^1,\langle\cdot,\cdot\rangle\Bigr)$
is an anti-linear $C^{\infty}$-isomorphism
$V\lrarr V$ satisfying the next conditions:
\begin{itemize}
 \item $\iota^2=id_V$ and  $\iota(V^{p,q})=V^{q,p}$.
 \item
 $\iota$ preserves the metric
 $\bigoplus_{p+q=w}(-\sqrt{-1})^{p-q}\langle\cdot,\cdot\rangle$.
 \item
 $\iota$ is flat with respect to $D$.
 Equivalently,
 $\iota$ replaces $(\delbar_V,\theta^{\dagger})$ and $(\del_V,\theta)$.
\end{itemize}
Such a tuple
$\Bigl(\bigoplus_{p+q=w}V^{p,q},
 \DD^1,\langle\cdot,\cdot\rangle,
 \iota\Bigr)$
is called a real variation of polarized Hodge structures.
For simplicity, it is called RVHS in this paper.
\hfill\qed
\end{df}

\noindent
{\it Address\\
Department of Mathematics,
Osaka City University,
Sugimoto, Sumiyoshi-ku,
Osaka 558-8585, Japan.
takuro@sci.osaka-cu.ac.jp\\
School of Mathematics,
Institute for Advanced Study,
Einstein Drive, Princeton,
NJ, 08540, USA. \\
takuro@math.ias.edu
}

\end{document}